\documentclass[12pt]{article}
\usepackage{graphicx}
\usepackage{amsmath}
\usepackage{amsthm}
\usepackage{graphics}
\usepackage{epsfig}
\usepackage{amssymb}
\usepackage{mathtools}
\usepackage{enumitem}
\usepackage{url}
\usepackage{hyperref}
\usepackage{xcolor}
\usepackage{multirow}
 \usepackage{algorithm}
\usepackage{algpseudocode}
\usepackage{float}
\usepackage{titletoc}

\usepackage[toc,page]{appendix}
\usepackage{xcolor}
 \usepackage{setspace}
\usepackage{extarrows}

\usepackage[authoryear]{natbib}
\usepackage{subfig}
\usepackage{lmodern}

\usepackage{pdflscape}
\usepackage{subfig}
\usepackage{enumitem}
\allowdisplaybreaks

\makeatletter
\newcommand{\myitem}[1]{%
\item[#1]\protected@edef\@currentlabel{#1}%
}
\makeatother

\newcommand{\R}{\mathbb{R}}

\newcommand{\E}{\mathbb{E}}

\newcommand{\tr}{\operatorname{tr}}

\usepackage{xcolor}

\newtheorem{theorem}{Theorem}[section]
\newtheorem{lemma}{Lemma}[section]
\newtheorem{proposition}{Proposition}[section]
\newtheorem{assertion}{Assertion}
\newtheorem{corollary}{Corollary}[section]

\newtheorem{remark}{Remark}[section]
\newtheorem{assumption}{Assumption}[section]

\numberwithin{equation}{section}


\begin{document}
\title{Quantitative Resolvent Calculus for Sample Covariance Matrices}
\date{\today}

\author{%
    Thomas Stark\\
    Department of Mathematics\\
    Aarhus University\\
    \texttt{thomas.stark@math.au.dk}
    \and
    Nina Dörnemann\\
    Department of Mathematics\\
    Aarhus University\\
    \texttt{ndoernemann@math.au.dk}
}

\maketitle

\begin{abstract}
In this work, we develop a quantitative resolvent calculus for sample
covariance matrices tailored to applications in high-dimensional
statistics.
We treat single resolvents, their derivatives, and products
of two resolvents with different spectral parameters
in the proportional-dimensional regime.
Our bounds hold uniformly over compact subsets of
$\mathbb C\setminus\mathbb R_{\geq0}$ with approximation errors
of order $O_{\mathbb P}(\|\Theta_n\|_F/\sqrt n)$ for deterministic
test matrices $\Theta_n$. These bounds require only uniformly bounded fourth moments
and a uniformly bounded population covariance operator norm.
Neither regular spectral edges nor convergence of the
population spectral distribution is required.
As a key  technical ingredient, we develop an averaging argument for centered
quadratic forms involving one or two leave-one-out resolvents. Applications of our results include precision-matrix shrinkage, prediction
risk in ridge regression, and generalized cross-validation.
For these problems, we obtain quantitative guarantees
for estimation and data-driven parameter selection.
\end{abstract}

MSC2020 subject classifications: 
Primary 60B20; Secondary 62H12, 62J07.

\noindent 
Keywords and phrases: Deterministic equivalents; sample covariance matrices;
resolvents; high-dimensional statistics; ridge regression;
precision-matrix estimation; generalized cross-validation.

\allowdisplaybreaks

\section{Introduction}
Deterministic equivalents for linear functionals of sample
covariance resolvents have become a central tool in modern
high-dimensional statistics. They underpin advances across a wide range of problems,
including high-dimensional prediction, distributed estimation,
regularized linear discriminant analysis, and principal component
regression \citep{HastieMontanariRossetTibshirani2022,
DobribanSheng2021, WangJiang2018,
Lecoiu:Mukherjee:Sur:2026, Green2025HighDimPCR}.
To fix notation, let $x_1,\ldots,x_n$ be centered independent random vectors
in $\mathbb R^p$ with common covariance matrix $\Sigma_n$,
and write
\[
    \hat\Sigma_n=\frac1n\sum_{i=1}^n x_ix_i^\top.
\]
We say that a possibly random matrix
$A_n\in\mathbb C^{p\times p}$ has a \textit{deterministic equivalent}
$B_n\in\mathbb C^{p\times p}$ if, for a suitable class of
deterministic test matrices $\Theta_n$,
\begin{align}\label{det_eq}
    \left|
    \operatorname{tr}\!\left[\Theta_n(A_n-B_n)\right]
    \right|
\end{align}
converges to zero in probability or almost surely as
$n\to\infty$.
A well-known deterministic equivalent is
available for the sample covariance resolvent
$A_n=R_n(\zeta)$, where $R_n(\zeta)=(\hat\Sigma_n-\zeta I_p)^{-1}$ and $\zeta\in\mathbb C\setminus\mathbb R_{\geq0}$.
Indeed, for suitable test matrices, a deterministic equivalent is
obtained by choosing $B_n=R_n^{\operatorname{det}}(\zeta)$, where
\begin{align} \label{eq:def_R}
    R_n^{\operatorname{det}}(\zeta)
    \coloneqq\bigl(-\zeta v_n(\zeta)\Sigma_n-\zeta I_p\bigr)^{-1};
\end{align}
see, for example, \citet{Hachem:Loubaton:Najim:2007},
\citet{RubioMestre2011}, and \citet{LedoitPeche2011}.
Here, $v_n$ denotes the Stieltjes transform of the companion
Marchenko--Pastur distribution associated with $p/n$ and
the empirical spectral distribution of $\Sigma_n$.
Its definition is given in Section~\ref{sec:RM-repair}.

Quantitative resolvent approximations are also central to
the theory of local laws
\citep{
Lin:Pan:2024:MultiResolvent, knowles2017anisotropic};
see Section~\ref{sec:related_work}.
These results control resolvents at spectral parameters
approaching the spectrum as the sample size $n$ increases.
Such estimates are often established under strong moment
assumptions and regularity conditions on the spectrum,
including its edges.
Many statistical applications, however, require resolvents
only at spectral parameters that remain a fixed positive
distance from $\mathbb R_{\geq0}$.
At the same time, weaker distributional assumptions
are desirable in these applications.
Uniform control over the regularization parameter is
particularly useful for analyzing data-driven parameter choices.

\noindent\textbf{Main contributions.}
Motivated by these considerations, we develop a quantitative
resolvent calculus for sample covariance matrices tailored
to statistical applications.
\begin{itemize}
    \item
    We establish deterministic approximations for single
    resolvents, their derivatives, and products of two resolvents
    with the population covariance matrix between them,
    allowing different spectral parameters $\zeta, \xi \in \mathbb C \setminus \mathbb R_{\geq 0}$.
    The bounds in Sections \ref{sec:RM-repair} and \ref{sec:two-resolvent} hold uniformly over compact subsets of
    $\mathbb C\setminus\mathbb R_{\geq0}$ and depend explicitly
    on the Frobenius norm of the deterministic test matrix.
    Under mild moment assumptions and a uniform bound on
    the operator norm of $\Sigma_n$, our quantitative results
    require neither regular spectral edges nor convergence
    of the empirical spectral distribution of $\Sigma_n$.

    \item
    In Section~\ref{sec:applications}, we study
    precision-matrix shrinkage, prediction risk in ridge
    regression, and generalized cross-validation.
    We obtain quantitative guarantees for estimation
    and for the excess loss of data-driven parameter choices
    relative to their oracle counterparts.

       \item
    We develop an averaging argument for centered quadratic
    forms involving one or two leave-one-out resolvents
    of sample covariance matrices.
    Inspired by the approach of \citet{RubioMestre2011},
    we use the structure of these resolvents to control
    the dependence between the summands.
    The resulting bounds, presented in
    Section~\ref{sec:Rubmes_erroranalysis}, are a key
    ingredient in the proofs of our main results.
\end{itemize}
Table~\ref{table_overview} summarizes our main quantitative
bounds in probability.
Almost-sure versions are given in
Proposition~\ref{prop:RM-recovered},
Corollary~\ref{cor:GreenRom-derivatives}, and
Theorem~\ref{thm:double-resolvent-limit2}.

\begin{table}[H]
    \centering
    \small
    \renewcommand{\arraystretch}{1.5}
    \setlength{\tabcolsep}{4pt}
    \begin{tabular}{c|c|c|c}
        $A_n$
        & Approximation to $\tr(\Theta_n A_n)$
        & Uniform error bound
        & Result \\
        \hline
        $R_n(\zeta)$
        & $\tr\{\Theta_nR_n^{\operatorname{det}}(\zeta)\}$
        & $O_{\mathbb P}(\|\Theta_n\|_F/\sqrt n)$
        & Thm.~\ref{Lemma:GreenRom}
        \\[2mm]
        $\frac{d^q}{d\zeta^q} R_n(\zeta)$
        & $\tr\left\{\Theta_n\frac{d^q}{d\zeta^q}
          R_n^{\operatorname{det}}(\zeta)\right\}$
        & $O_{\mathbb P}(\|\Theta_n\|_F/\sqrt n)$
        & Cor.~\ref{cor:GreenRom-derivatives:rates}
        \\[2mm]
        $R_n(\zeta)\Sigma_nR_n(\xi)$
        & $K_n^\Theta(\zeta,\xi)$
        & $O_{\mathbb P}(\|\Theta_n\|_F/\sqrt n)$
        & Thm.~\ref{thm:double-resolvent-limit}
        \\
    \end{tabular}
    \caption{Overview of our main quantitative
    results for the approximation in \eqref{det_eq}.
    The deterministic matrix $R_n^{\operatorname{det}}(\zeta)$
    and the scalar $K_n^\Theta(\zeta,\xi)$ are defined in
    \eqref{eq:def_R} and \eqref{def_K_theta}, respectively.
    The integer $q\geq1$ is fixed.
    The error bounds hold uniformly over compact sets of
    spectral parameters in $\mathbb C\setminus\mathbb R_{\geq0}$.}
    \label{table_overview}
\end{table}

\subsection{Related work} \label{sec:related_work}
Our results are related to several lines of research on
deterministic approximations in random matrix theory and high-dimensional statistics.
Given the extensive literature on this topic, we focus
on a selection of works most closely related to our results. We first review these developments and then discuss
how our results relate to the closest contributions.

Foundational contributions include
\citet{Hachem:Loubaton:Najim:2007} on deterministic equivalents
for normalized resolvent traces,
\citet{Hachem:Loubaton:Najim:Vallet:2013} on bilinear
resolvent forms, and
\citet{LedoitPeche2011} on weighted resolvent traces
and their connection to sample and population eigenvectors. More recently, \citet{Chouard:2022} obtains quantitative
deterministic equivalents for sample covariance resolvents
with independent identically distributed columns and dependence
within columns, with explicit bounds in terms of the dimension
and the spectral parameter. \citet{Louart:2026} establishes quantitative $L^q$-bounds
for single- and double-resolvent expressions under moment
assumptions on centered quadratic forms, allowing different
covariance matrices across observations.

Another closely related line of research concerns local laws,
which provide high-probability approximations for resolvents
on fine spectral scales, including spectral parameters $\zeta$ approaching
the spectrum; see, for example,
\citet{Bloemendal:Erdos:Knowles:Yau:Yin:2016,
knowles2017anisotropic}.
They are fundamental tools for establishing eigenvalue rigidity,
eigenvector delocalization, and universality of local eigenvalue
statistics, including Tracy-Widom fluctuations at regular
spectral edges; see, e.g. \cite{Lee:Schnelli:2016, Bao:Pan:Zhou:2015, Ding:Yang:2021}.
Related developments concern products of resolvents with
deterministic matrices inserted between them. \citet{Cipolloni:Erdos:Schroeder:2022} establish optimal
multi-resolvent local laws for Wigner matrices.
For sample covariance matrices,
\citet{Lin:Pan:Zhao:Zhou:2024} establish a two-resolvent
isotropic law as part of their analysis of spiked eigenvalue
fluctuations in signal-plus-noise models.
\citet{Lin:Pan:2024:MultiResolvent} develop local laws for
two- and three-resolvent expressions and use them to obtain
quantitative bounds for eigenvector overlaps.

Related random matrix methods are widely used to study
estimation and prediction in high-dimensional statistics.
\citet{dobriban:2018} derive limiting prediction risks
for ridge regression and regularized discriminant analysis
under a random-effects model.
\citet{HastieMontanariRossetTibshirani2022} study the
asymptotic prediction risk of ridge and minimum-norm
least-squares estimators, while
\citet{DobribanSheng2020,DobribanSheng2021} analyze
distributed linear and ridge regression.
More recently, \citet{Green2025HighDimPCR} use
multi-resolvent traces and eigenvector-overlap measures
to characterize the limiting prediction risk of
principal component regression.
Random matrix methods also provide a basis for
data-driven selection of regularization parameters.
\citet{WangPanTongZhu2015} derive asymptotically optimal
shrinkage parameters for precision-matrix estimation,
while \citet{PatilWeiRinaldoTibshirani2021} establish
almost-sure uniform consistency of generalized
cross-validation for ridge regression.

\vspace{0.25cm}
We conclude this discussion by positioning our contribution
within this literature.
Our results complement the work of \citet{Chouard:2022} on one-resolvent expressions by providing
uniform approximations under a uniform finite fourth-moment assumption and by additionally treating
two-resolvent expressions.
We also revisit the averaging argument of Lemma 4 in
\citet{RubioMestre2011}.
In Section~\ref{sec:Rubmes_erroranalysis}, we give a
counterexample to the underlying assertion and establish
averaging bounds for centered quadratic forms involving
leave-one-out resolvents in the sample covariance setting. Our results also provide quantitative guarantees for estimation
and data-driven parameter selection in precision-matrix
shrinkage and ridge regression, including generalized
cross-validation, as discussed in Section~\ref{sec:applications}. These guarantees complement the asymptotic
risk formulas and consistency results in the statistical
literature discussed above.

While completing this manuscript, we became aware of the
latest version of \citet{Louart:2026}, which now also includes
deterministic approximations for double-resolvent expressions. There are several differences in scope and emphasis
between the two works.
Our double-resolvent approximations allow two different
spectral parameters $\zeta,\xi$, whereas
\citet[Corollary~1.7]{Louart:2026} treats a common
spectral parameter $\zeta = \xi$.
Different spectral parameters are needed for the statistical
application in \cite{StarkKoestenbergerDoernemann2026};
see also Section~\ref{sec:app-ridge}.
The two works also differ in their model assumptions:
we consider a common population covariance matrix,
while \citet{Louart:2026} allows different covariance
matrices across observations. Moreover, the proofs of the two works take different routes.
For double-resolvent approximations, \citet{Louart:2026}
differentiates resolvents and their deterministic equivalents
with respect to an additive matrix perturbation.
This naturally yields two resolvents at the same
spectral parameter $\zeta = \xi$.
Our proof builds on averaging bounds for centered quadratic
forms involving leave-one-out resolvents and directly
treats two different spectral parameters $\zeta, \xi$ (see Section \ref{sec:Rubmes_erroranalysis}).

\vspace{0.5cm}

The remainder of this paper is structured as follows. Section~\ref{sec:RM-repair} presents deterministic equivalents
and quantitative bounds for single resolvents and their derivatives.
These results are extended to two-resolvent expressions in
Section~\ref{sec:two-resolvent}, followed by statistical
applications in Section~\ref{sec:applications}.
Finally, Section~\ref{sec:Rubmes_erroranalysis} revisits the
averaging argument of \citet{RubioMestre2011} and establishes
averaging bounds for leave-one-out resolvent quadratic forms,
providing insight into our proof strategy.
Technical lemmas and most proofs are deferred to the appendix.

\vspace{0.5cm}

\noindent \textbf{Notation and terminology.}
We write
\[
    \mathbb C^\pm=\{z\in\mathbb C:\pm\operatorname{Im}(z)>0\},
    \qquad \mathbb R_{\geq0}=[0,\infty).
\]
Let $I_q$ denote the identity matrix of size $q$.
For a real symmetric matrix $A$ of size $q$, write its eigenvalues
as $\lambda_1(A),\ldots,\lambda_q(A)$. Its empirical spectral measure is
\[
    F_A\coloneqq\frac1q\sum_{i=1}^q\delta_{\lambda_i(A)},
\]
where $\delta_x$ is the Dirac measure at $x$.
For a real or complex matrix $A$, we write $A^\top$ for its transpose
and $A^*$ for its conjugate transpose. We denote its operator norm by
$\|A\|_2$, its Frobenius norm by
$\|A\|_F=\{\tr(A^*A)\}^{1/2}$, and its trace norm by
$\|A\|_1=\tr\{(A^*A)^{1/2}\}$.
For vectors, $\|\cdot\|_2$ denotes the Euclidean norm.
For a nonnegative definite matrix $A$, $A^{1/2}$ denotes its unique
nonnegative definite square root.
We use $\Rightarrow$ to denote weak convergence of measures.
All limits are taken as $n\to\infty$, unless stated otherwise.
All random variables are defined on a common probability space
$(\Omega,\mathcal F,\mathbb P)$.

\section{Deterministic Equivalents for Sample Covariance Resolvents}\label{sec:RM-repair}
In this section, we establish convergence results for resolvent
expressions in forms commonly encountered in the literature. After introducing the notation and
assumptions, we present two technical lemmas that allow us to consider
spectral parameters in $\mathbb{C}\setminus\mathbb{R}_{\geq 0}$ rather than
only in $\mathbb{C}^+$, followed by our main results.

We denote the $n \times p$ data matrix by $X_n = (x_1, \ldots, x_n)^\top.$ 
Before proceeding to the theoretical results, we collect the required assumptions and notations. We emphasize that not all of these assumptions are required for each result in this section.
\begin{assumption}\label{ass:main}
\begin{enumerate}[label=(A-\arabic*)]
    \item\label{A1}
    The dimension $p=p_n$ satisfies
    \[
            \gamma_n\coloneqq\frac pn \longrightarrow \gamma\in(0,\infty), \quad n\to\infty.
    \]

\item\label{A3}
The data matrix $X_n$ has the form $X_n = Z_n \Sigma_n^{1/2}$ and 
 the entries of $Z_n =(z_{ij})_{1\leq i\leq n, 1\leq j\leq p}$ are, for each $n\in \mathbb{N}$, real-valued independently distributed random variables satisfying $\mathbb{E}(z_{ij})=0$, $\mathbb{E}(z_{ij}^2)=1$ and
\begin{equation}\label{X:unif4}
            \sup_{n\ge1}
            \max_{1\le i\le n,\ 1\le j\le p}
            \mathbb E (|z_{ij}|^{4})
            <\infty .
    \end{equation}

    \item\label{A2}
    The matrices $\Sigma_n$ are deterministic and nonnegative
    definite, and satisfy
    \[
            \sup_{n\ge1}\|\Sigma_n\|_2<\infty.
    \]

    \item\label{A4}
    The empirical spectral measure $H_n\coloneqq F_{\Sigma_n}$ converges weakly to a measure $H$ supported on $[0,\infty)$, that is,
    \[
            H_n\Rightarrow H , \quad n\to\infty.
    \]
\end{enumerate}
\end{assumption}

Denote the $n\times n$ companion matrix of $\hat\Sigma_n$ by
\[
        \underline{\hat\Sigma}_n
        \coloneqq
        \frac1n X_nX_n^\top .
\]
For $j=1,\ldots,n$, let $z_j^\top$ denote the $j$th row of $Z_n$;
then $x_j=\Sigma_n^{1/2}z_j$.

We first recall a version of the Mar\v{c}enko--Pastur theorem \citep{MarPas1967}, which characterizes the limiting spectral measure of $\underline{\hat \Sigma}_n.$  
The formulation given here is a special case of
\citet[Theorem~1]{Pan2010}.
\begin{lemma}
\label{thm:MP}
Suppose Assumption~\ref{ass:main} holds. Then, almost surely,
\[
        F_{\underline{\hat\Sigma}_n}
        \Rightarrow
        \bar F_{\gamma,H},
\]
where $\bar F_{\gamma,H}$ is a nonrandom probability measure supported on
$[0,\infty)$. For $z\in\mathbb C^+$, its Stieltjes transform $v$ is the unique solution in $\mathbb C^+$ of
\begin{equation}\label{stielt1}
        v(z)
        =
        -\left(
        z-\gamma\int_0^\infty
        \frac{t}{1+t  v(z)}\,dH(t)
        \right)^{-1}.
\end{equation}
\end{lemma}
Deterministic equivalents for traces involving the resolvent of the sample covariance matrix $\hat\Sigma_n$ are naturally expressed in terms of the
Stieltjes transform of the companion Mar\v{c}enko--Pastur distribution. The following lemma extends the characterization in \eqref{stielt1} to
$\mathbb C\setminus\mathbb R_{\geq0}$, and its proof can be found in Section~\ref{sec:proof:lem}.

\begin{lemma}
\label{lem:companion-uniq-D}
Let $\gamma_0\in(0,\infty)$, and let $H_0$ be a probability measure supported on a compact interval
$[0,C_0]\subseteq[0,\infty)$. Let $\overline F_{\gamma_0,H_0}$ be the nonrandom limiting measure associated with the pair
$(\gamma_0,H_0)$ in Lemma~\ref{thm:MP}, and define
\[
        v_{\gamma_0,H_0}(\zeta)
        :=
        \int_0^\infty
        \frac{1}{s-\zeta}\,d\bar F_{\gamma_0,H_0}(s),
        \qquad
        \zeta\in\mathbb C\setminus\mathbb R_{\ge0}.
\]
Then, for every $\zeta\in\mathbb C\setminus\mathbb R_{\ge0}$,
$v_{\gamma_0,H_0}(\zeta)$ is the unique admissible $w\in\mathbb C$ satisfying
\begin{align}\label{stielt}
        w
        =
        -\left(
        \zeta
        -
        \gamma_0\int_0^\infty
        \frac{t}{1+tw}\,dH_0(t)
        \right)^{-1},
\end{align}
where admissible means that $\Im(w)\Im(\zeta)>0$, if
$\zeta\in\mathbb C\setminus\mathbb R$, and $w>0$, if
$\zeta\in\mathbb R_{<0}$.
\end{lemma}
 For each $n\in\mathbb N$, we define
\begin{align}\label{eq:fixed_point}
        v_n(\zeta)
        \coloneqq
        v_{\gamma_n,H_n}(\zeta),
        \qquad
        \zeta\in\mathbb C\setminus\mathbb R_{\geq0}.
\end{align}
By Lemma~\ref{lem:companion-uniq-D}, $v_n(\zeta)$ is equivalently the unique admissible solution of \eqref{stielt} with $(\gamma_0,H_0)=(\gamma_n,H_n)$.
Next, we turn to our main results on deterministic equivalents for resolvent matrices
\begin{align} \label{eq:def_resolvent}
    R_n (\zeta) := (\hat\Sigma_n - \zeta I_p )^{-1}, \quad \zeta \in \mathbb{C} \setminus \R_{\geq 0}.
\end{align}
The corresponding deterministic matrix $R_n^{\operatorname{det}}(\zeta)$ is defined in \eqref{eq:def_R}.
The following theorem gives a deterministic equivalent uniformly on
compact subsets of
\(\mathbb C\setminus\mathbb R_{\geq0}\), with an error bound that
depends explicitly on the Frobenius norm of the deterministic test matrix.

\begin{theorem}
\label{Lemma:GreenRom}
Suppose Assumptions~\ref{A3} and~\ref{A2} hold and
\begin{align}\label{gamma_n:cond}
0<\liminf_{n\to\infty}\gamma_n
\le
\limsup_{n\to\infty}\gamma_n
<\infty .
\end{align}
Let
\(\{\Theta_n\}_{n\ge1}\) be a sequence of deterministic \(p\times p\) matrices.
Then, for any compact set \(K\subseteq \mathbb{C}\setminus\mathbb{R}_{\geq 0}\),
\begin{align}\label{eq:error_bound_rates}
        \sup_{\zeta\in K}\biggl|\tr\bigl[\Theta_n\bigl\{R_n(\zeta)
        -
        R_n^{\operatorname{det}}(\zeta)\bigr\}\bigr]\biggr|=O_{\mathbb{P}}\left(\frac{\|\Theta_n\|_{F}}{n^{1/2}}\right).
\end{align}
Here, the \(O_{\mathbb P}(\cdot)\)-term may depend on the compactum $K.$
\end{theorem}
We provide a proof outline later in this section; the
complete proof is given in Section~\ref{sec:RM-repair_proofs}.

Next, we illustrate how the rate in
Theorem~\ref{Lemma:GreenRom} depends on the rank of the
test matrix. Consider $\Theta_n=r_n^{-1}P_n$,
where $P_n$ is an orthogonal projection of rank $r_n$.
Then $\|\Theta_n\|_F=r_n^{-1/2}$, and
\eqref{eq:error_bound_rates} yields the rate $O_{\mathbb P}(1/\sqrt{nr_n})$.
For $r_n=1$, the trace functional reduces to a quadratic
form, and the bound is of order $n^{-1/2}$.
At the other extreme, $r_n=p$ gives
$\Theta_n=I_p/p$, so the trace functional is the
normalized trace. In this case, \eqref{gamma_n:cond}
implies a bound of order $n^{-1}$.
Thus, the bound interpolates between the $n^{-1/2}$
rate for quadratic forms and the $n^{-1}$ rate for
normalized traces.
Applications are discussed in
Section~\ref{sec:applications}.

The deterministic equivalent of Theorem \ref{Lemma:GreenRom} is stable under differentiation with
respect to the spectral parameter, as stated in the following corollary. Its proof is postponed to Section~\ref{sec:RM-repair_proofs}.
\begin{corollary}
\label{cor:GreenRom-derivatives:rates}
Consider the assumptions of Theorem~\ref{Lemma:GreenRom}, and let
$q\in\mathbb N_0$ be fixed. Then, 
\[
\sup_{\zeta\in K}\left|
\frac{d^q}{d\zeta^q}
\tr\!\left[
\Theta_n
\left\{
R_n(\zeta)
-
R_n^{\operatorname{det}}(\zeta)
\right\}
\right]
\right|=O_{\mathbb{P}}\left(\frac{\|\Theta_n\|_{F}}{n^{1/2}} \right).
\]
\end{corollary}
Next, we outline the proof of Theorem \ref{Lemma:GreenRom}, before turning to the almost-sure convergence results. To this end, denote the empirical companion Stieltjes transform by
\begin{equation}\label{eq:emp-companion-definition}
    \widehat v_n(z)\coloneqq
    \frac1n\tr\bigl[(\underline{\widehat\Sigma}_n-zI_n)^{-1}\bigr],
    \qquad z\in\mathbb C\setminus\mathbb R_{\geq0}.
\end{equation}

\noindent \textbf{Proof outline for Theorem~\ref{Lemma:GreenRom}.}
The proof consists of two main ingredients. First, we provide a uniform approximation of the companion
Stieltjes transform $\widehat v_n$ from \eqref{eq:emp-companion-definition} by its deterministic counterpart
$v_n$ defined in \eqref{eq:fixed_point}. We establish this approximation in two steps. 
First, under the additional assumptions
$\gamma_n\to\gamma\in(0,\infty)$ and $H_n\Rightarrow H$,
we obtain uniform bounds for the fluctuations
$\hat v_n(\zeta)-\E(\hat v_n(\zeta))$ and the
deterministic bias $\E(\hat v_n(\zeta))-v_n(\zeta)$
in Lemmas~\ref{stielt:var} and~\ref{stielt:bias},
respectively, in Appendix~\ref{sec:RM-repair_proofs}.
Lemma~\ref{lem:uniform-companion-rate} then combines these two estimates
and extends the resulting bound to the weaker assumption
\eqref{gamma_n:cond}, by a subsequence argument. This yields
\begin{align} \label{n11}
        \sup_{\zeta\in K}
        |\hat v_n(\zeta)-v_n(\zeta)|
        =
        O_{\mathbb P}(n^{-1}).
\end{align}
We note that the fluctuation--bias decomposition is akin to that used by
\cite{Silversteinbook2010} in their proof of Lemma 9.11. Establishing uniform bounds on arbitrary compact
subsets $K$ of $\mathbb C\setminus\mathbb R_{\geq 0}$ requires a different
technical approach from their treatment of a specially constructed
compact set $\mathcal{C}^+ \subset \mathbb{C}^+$.

Second, we proceed with the decomposition
\begin{equation} \label{n22}
\begin{split} 
\tr\!\left[
\Theta_n\{R_n(\zeta)-R_n^{\operatorname{det}}(\zeta)\}
\right]
&=
-\zeta\{v_n(\zeta)-\hat v_n(\zeta)\}
\tr\!\left[
\Theta_nR_n(\zeta)\Sigma_n
R_n^{\operatorname{det}}(\zeta)
\right]
\\
&\qquad+
\tr\!\left[
\Theta_nA_n(\zeta)
R_n^{\operatorname{det}}(\zeta)
\right],
\end{split}
\end{equation}
where 
\[
A_n(\zeta)=-\zeta\hat v_n(\zeta)R_n(\zeta)\Sigma_n-R_n(\zeta)\hat\Sigma_n.
\]

Using \eqref{n11}, the first summand in \eqref{n22} is seen to be of order 
     $   O_{\mathbb P}\left(n^{-1}
        \|\Theta_n\|_1 
        \right).$
Repeated use of the Sherman--Morrison identity decomposes the second
summand in \eqref{n22} into a leading averaged term and remainder terms.
Lemma~\ref{RubMes_Lem4_fix} controls the leading term at the rate
$O_{\mathbb P}(n^{-1/2}\|\Theta_n\|_F)$. Together with uniform bounds on the remainder terms, this yields
the asserted rate. This concludes the proof outline. The complete proof is given in
Section~\ref{sec:RM-repair_proofs}.

\medskip

We conclude this section by complementing the preceding quantitative
bounds in Theorem~\ref{Lemma:GreenRom} and
Corollary~\ref{cor:GreenRom-derivatives:rates} with almost-sure convergence
results. For these results, we strengthen the moment condition in
Assumption~\ref{A3}. The aspect-ratio condition remains \eqref{gamma_n:cond}.

\begin{proposition}
\label{prop:RM-recovered}
Suppose Assumptions~\ref{A3} and~\ref{A2} and
\eqref{gamma_n:cond} hold, and that
\eqref{X:unif4} is strengthened to
\[
        \sup_{n\ge1}
        \max_{1\le i\le n,\ 1\le j\le p}
        \mathbb E\bigl(|z_{ij}|^{8}\bigr)
        <\infty .
\]
Let $\{\Theta_n\}_{n\ge1}$ be a sequence of deterministic
$p\times p$ matrices satisfying $\sup_{n\ge1}\|\Theta_n\|_{1}<\infty$.
Then, almost surely, for all
$\zeta\in\mathbb C\setminus\mathbb R_{\geq0}$,
\[
\operatorname{tr}\!\left[
\Theta_n
\left\{
R_n(\zeta)
-
R_n^{\operatorname{det}}(\zeta)
\right\}
\right]
\longrightarrow 0.
\]
\end{proposition}

\begin{remark} {\rm
\rmfamily
Proposition~\ref{prop:RM-recovered} provides a version of
\citet[Theorem~1]{RubioMestre2011} for the sample-covariance setting
considered here. The two results differ in several respects.
First, expressed in our notation and normalization,
the more general matrix model considered by
\citet{RubioMestre2011} takes the form
\[
    B_n
    =
    A_n+\frac{1}{n}
    \Sigma_n^{1/2}Z_n^\top T_nZ_n\Sigma_n^{1/2},
\]
where $A_n$ is a $p\times p$ matrix and $T_n$ is an
$n\times n$ matrix, subject to the assumptions in
\citet{RubioMestre2011}.
Our sample-covariance model corresponds to $A_n=0$ and $T_n=I_n$,
a setting used in many statistical applications; see, for example,
\citet{HastieMontanariRossetTibshirani2022,Green2025HighDimPCR}.
Second, our assumptions on the entries of $Z_n$ are weaker in two
respects: the entries need not be identically distributed, and
Proposition~\ref{prop:RM-recovered} requires only uniformly bounded
eighth moments, rather than $(8+\delta)$ moments for some $\delta>0$.

Our main motivation for including Proposition~\ref{prop:RM-recovered}
is that the proof of \citet[Theorem~1]{RubioMestre2011} relies on
their Lemma~4, whose averaging argument is not valid in the stated
generality; see Section~\ref{sec:Rubmes_erroranalysis} for a counterexample
and a detailed discussion.
We therefore provide an independent proof of the corresponding
deterministic equivalent in the sample-covariance setting considered here. }
\end{remark}

\begin{remark}[Deterministic limits for resolvent functionals]
\label{rem:single-resolvent-limits}
{\rm
The approximation in Proposition~\ref{prop:RM-recovered} yields
convergence to a deterministic limit whenever the corresponding
deterministic equivalent converges. We illustrate this for normalized
traces and quadratic forms.

Consider the assumptions of Proposition~\ref{prop:RM-recovered}, together with Assumptions~\ref{A1} and~\ref{A4}. Using Lemma~\ref{thm:MP} together with dominated convergence, we have $\E(\widehat v_n(\zeta))\to v(\zeta)$ for every $\zeta\in \mathbb C\setminus\mathbb R_{\geq 0}$. 
Combining this with Lemma~\ref{stielt:bias} in Appendix~\ref{sec:RM-repair_proofs}
and applying the triangle inequality, we conclude that
\begin{align}\label{eq:pointwise_vn_v}
        v_n(\zeta)\longrightarrow v(\zeta),
        \qquad
        \zeta\in\mathbb C\setminus\mathbb R_{\geq0}.
\end{align}
For $\Theta_n=p^{-1}I_p$, the definition of
$R_n^{\operatorname{det}}(\zeta)$ and the fixed-point equation
\eqref{stielt} give
\[
\begin{split}
        \frac1p\tr\{R_n^{\operatorname{det}}(\zeta)\}
        =
        -\frac1\zeta
        \int_0^\infty
        \frac{1}{1+t v_n(\zeta)}\,dH_n(t)
        =
        -\left(1-\frac1{\gamma_n}\right)\frac1\zeta
        +
        \frac1{\gamma_n}v_n(\zeta).
\end{split}
\]
Consequently, Proposition~\ref{prop:RM-recovered} and
the preceding two displays yield, almost surely, for all
$\zeta\in\mathbb C\setminus\mathbb R_{\geq0}$,
\[
        \frac1p\tr\{R_n(\zeta)\}
        \longrightarrow
        -\left(1-\frac1\gamma\right)\frac1\zeta
        +
        \frac1\gamma v(\zeta)
        \eqqcolon        m(\zeta),
\]
where $m$ is the Stieltjes transform of $F_{\gamma,H}$;
see \eqref{limitStil} in the appendix.

Next, let $\beta_n\in\mathbb R^p\setminus\{0\}$ be deterministic and consider $\Theta_n = \beta_n \beta_n^\top$. Then, the trace functional in \eqref{eq:error_bound_rates} reduces to a quadratic form. Its deterministic
equivalent can be expressed as an integral with respect to
the spectral measure of $\Sigma_n$ associated with $\beta_n$.To define this measure, consider an orthonormal eigendecomposition
$\Sigma_n=\sum_{\ell=1}^p t_\ell q_\ell q_\ell^\top$ and write
\[
        G_n(t)
        \coloneqq
        \frac1{\|\beta_n\|_2^2}
        \sum_{\ell=1}^p
        \langle\beta_n,q_\ell\rangle^2
        \mathbf 1_{\{t_\ell\leq t\}}.
\]
Suppose additionally that $\|\beta_n\|_2^2\longrightarrow\tau^2\in(0,\infty)$, and $G_n\Rightarrow G$.
For $\Theta_n=\beta_n\beta_n^\top$, we have
\[
    \beta_n^\top R_n^{\operatorname{det}}(\zeta)\beta_n
    =
    -\frac{\|\beta_n\|_2^2}{\zeta}
    \int_0^\infty
    \frac{1}{1+t v_n(\zeta)}\,dG_n(t).
\]
By Assumption~\ref{A2}, the measures $G_n$ are supported
on a common compact interval $[0,M]$.
Since $v_n(\zeta)\to v(\zeta)$ and
$1+t v(\zeta)\neq0$ for all $t\in[0,M]$,
the integrands in the preceding display converge uniformly on $[0,M]$.
Together with $G_n\Rightarrow G$,
$\|\beta_n\|_2^2\to\tau^2$, and
Proposition~\ref{prop:RM-recovered}, this yields,
almost surely, for every
$\zeta\in\mathbb C\setminus\mathbb R_{\geq0}$,
\[
    \beta_n^\top R_n(\zeta)\beta_n
    \longrightarrow
    -\frac{\tau^2}{\zeta}
    \int_0^\infty
    \frac{1}{1+t v(\zeta)}\,dG(t).
\]
}
\end{remark}

Again, the approximation is stable under differentiation, as stated in the following result.
\begin{corollary}
\label{cor:GreenRom-derivatives}
Suppose the assumptions of Proposition~\ref{prop:RM-recovered} hold,
and let \(q\in\mathbb N_0\) be fixed. Then, almost surely, for all
\(\zeta\in\mathbb C\setminus\mathbb R_{\geq0}\),
\[
\frac{d^q}{d\zeta^q}
\operatorname{tr}\!\left[
\Theta_n
\left\{
R_n(\zeta)
-
R_n^{\operatorname{det}}(\zeta)
\right\}
\right]
\longrightarrow 0.
\]
\end{corollary}
In the next section, we extend the preceding deterministic-equivalent results
to trace functionals involving products of two resolvents.

\section{Deterministic Equivalents for Two-Resolvent Expressions}\label{sec:two-resolvent}

We now turn to trace functionals containing two resolvents,
\begin{align}\label{eq:double_res_def}
        \tr\!\left\{
        \Theta_nR_n(\zeta)\Sigma_nR_n(\xi)
        \right\},
        \qquad
        \zeta,\xi\in\mathbb C\setminus\mathbb R_{\geq0}.
\end{align}
Such expressions arise naturally in covariance-weighted quadratic forms,
prediction risks, regularized discriminant analysis, and eigenvector-overlap
problems; see, for example, \cite{Green2025HighDimPCR}, \cite{HastieMontanariRossetTibshirani2022}, \cite{WangJiang2018}, \cite{Lin:Pan:2024:MultiResolvent}.
In contrast to products
of resolvents alone, the presence of the population covariance matrix
$\Sigma_n$ in \eqref{eq:double_res_def} prevents a direct reduction to the single-resolvent
deterministic equivalent, since $\Sigma_n$ does not generally commute with
$\widehat\Sigma_n$. The following lemma justifies the replacement of
$\Sigma_n$ by $\widehat\Sigma_n$, up to the deterministic factors
 $\{-\zeta v_n(\zeta)\}$ and $\{-\xi v_n(\xi)\}$.

\begin{lemma}
\label{lem:ridge-signal-comparison}
Suppose that Assumptions~\ref{A3}, \ref{A2}, and
\eqref{gamma_n:cond} are satisfied.
Let \(K\subset\mathbb C\setminus\mathbb R_{\geq0}\) be compact,
and let \(\Theta_n\) be a sequence of real, deterministic
\(p\times p\) matrices.
\begin{enumerate}[label=(\alph*),ref=\thelemma(\alph*)]
    \item\label{lem:ridge-signal-equivalence:rates}
    Then,
    \begin{align*}
   \sup_{\zeta,\xi\in K}
\Bigg|
\{-\zeta v_n(\zeta)\}\{-\xi v_n(\xi)\}
\tr\!\left\{\Theta_nR_n(\zeta)\Sigma_nR_n(\xi)\right\}
-
\tr\!\left\{\Theta_nR_n(\zeta)\hat\Sigma_nR_n(\xi)\right\}
\Bigg|=
O_{\mathbb P}\!\left(\frac{\|\Theta_n\|_F}{\sqrt n}\right).
    \end{align*}

    \item\label{lem:ridge-signal-equivalence}
    Suppose additionally that
    \(\sup_{n\in\mathbb N}\|\Theta_n\|_F<\infty\) and
    \[
        \sup_{n\ge1}
        \max_{1\le i\le n,\ 1\le j\le p}
        \mathbb E\bigl(|z_{ij}|^8\bigr)
        <\infty.
    \]
    Then, for every
    \(\zeta,\xi\in\mathbb C\setminus\mathbb R_{\geq0}\),
    \begin{align}\label{eq:ridge-signal-equivalence}
    \begin{split}
    \Bigg|
    \{-\zeta v_n(\zeta)\}\{-\xi v_n(\xi)\}
    \tr\!\left\{
        \Theta_nR_n(\zeta)\Sigma_nR_n(\xi)
    \right\}-
    \tr\!\left\{
        \Theta_nR_n(\zeta)\hat\Sigma_nR_n(\xi)
    \right\}
    \Bigg|
    \overset{\mathrm{a.s.}}{\longrightarrow}0.
    \end{split}
    \end{align}
\end{enumerate}
\end{lemma}

The preceding lemma enables us to reduce
construction of a deterministic equivalent for \eqref{eq:double_res_def}
to single-resolvent deterministic equivalents, for which the results of
Section~\ref{sec:RM-repair} apply. Indeed, for $\zeta\neq\xi$, we have
\[
       \tr\{\Theta_nR_n(\zeta)\widehat\Sigma_nR_n(\xi)\}
        =
        \frac{
        \zeta\tr\{\Theta_nR_n(\zeta)\}
        -
        \xi\tr\{\Theta_nR_n(\xi)\}
        }{\zeta-\xi}.
\]

This motivates the following definition of the finite-$n$
deterministic equivalent.
For \(\zeta,\xi\in\mathbb C\setminus\mathbb R_{\geq0}\), let
\begin{align}\label{def_K_theta}
K_n^\Theta(\zeta,\xi)
\coloneqq
\begin{cases}
\displaystyle
\frac{
        \zeta\tr\{\Theta_n R_n^{\operatorname{det}}(\zeta)\}
        -
        \xi\tr\{\Theta_n R_n^{\operatorname{det}}(\xi)\}
}{
        (\zeta-\xi)
        \{-\zeta v_n(\zeta)\}
        \{-\xi v_n(\xi)\}
},
& \zeta\neq \xi, \\[4mm]
\displaystyle
\frac{
        \frac{d}{d\zeta}
        \left[
        \zeta\tr\{\Theta_n R_n^{\operatorname{det}}(\zeta)\}
        \right]
}{
        \{-\zeta v_n(\zeta)\}^2
},
& \zeta=\xi .
\end{cases}
\end{align}
The definition for $\zeta=\xi$ is obtained by continuous extension from the case $\zeta\neq\xi$. Notably, $K_n^\Theta$ is determined entirely by the same
finite-$n$ deterministic resolvent $R_n^{\operatorname{det}}$ that appears in the one-resolvent
theory.

The following theorem shows that the deterministic equivalent $K_n^\Theta(\zeta,\xi)$ defined above
approximates \eqref{eq:double_res_def} with the same rate as the single-resolvent
deterministic equivalent in Section~\ref{sec:RM-repair}.

\begin{theorem}
\label{thm:double-resolvent-limit}
Suppose that Assumptions~\ref{A3} and~\ref{A2} are satisfied, and suppose that \eqref{gamma_n:cond} holds. Let
\(\Theta_n\) be a sequence of real, deterministic \(p\times p\) matrices and $K$ a compact set satisfying $K\subset\mathbb{C}\setminus\mathbb{R}_{\geq 0}$.
Then,
\[
\sup_{\zeta,\xi\in K}\left|
\tr\{\Theta_nR_n(\zeta)\Sigma_nR_n(\xi)\}
-
K_n^\Theta(\zeta,\xi)
\right|=O_{\mathbb{P}}\left(\frac{\|\Theta_n\|_F}{n^{1/2}}\right).
\]
\end{theorem}

Under stronger moment assumptions, the corresponding almost-sure
approximation holds:

\begin{theorem}
\label{thm:double-resolvent-limit2}
Suppose that Assumptions~\ref{A3} and~\ref{A2} are satisfied, and that
\eqref{X:unif4} is strengthened to $\sup_{n\ge1}
        \max_{1\le i\le n,\ 1\le j\le p}
        \mathbb E\bigl(|z_{ij}|^{8}\bigr)
        <\infty$, and suppose that \eqref{gamma_n:cond} holds. Let
\(\Theta_n\) be a sequence of real, deterministic \(p\times p\) matrices
with $\sup_{n\in\mathbb N}\|\Theta_n\|_{1}<\infty$.
Then, almost surely, for every
\(\zeta,\xi\in\mathbb C\setminus\mathbb R_{\geq0}\),
\[
\left|
\tr\{\Theta_nR_n(\zeta)\Sigma_nR_n(\xi)\}
-
K_n^\Theta(\zeta,\xi)
\right|
\longrightarrow0 .
\]
\end{theorem}

The proofs of Theorems~\ref{thm:double-resolvent-limit} and
\ref{thm:double-resolvent-limit2} are given in
Appendix~\ref{sec:two-resolvent-proofs}.

\begin{remark}[Deterministic limits for two-resolvent functionals]
\label{rem:double-resolvent-limits}
{\rm
Similarly to Remark~\ref{rem:single-resolvent-limits} for single
resolvents, we can obtain deterministic limits for two-resolvent
expressions under additional assumptions. As before, we consider the cases $\Theta_n=p^{-1} I_p,$ corresponding to the normalized trace, and $\Theta_n = \beta_n \beta_n^\top$, corresponding to a quadratic form. 
Suppose the assumptions of
Theorem~\ref{thm:double-resolvent-limit2} hold, and assume
additionally that Assumptions~\ref{A1} and~\ref{A4} are
satisfied. Define
\[
        \Delta(\zeta,\xi)
        \coloneqq
        \begin{cases}
        \displaystyle
        \frac{v(\zeta)-v(\xi)}{\zeta-\xi},
        & \zeta\neq\xi,\\[3mm]
        \displaystyle
        \frac{d}{d\zeta}v(\zeta),
        & \zeta=\xi.
        \end{cases}
\]

For $\Theta_n=p^{-1}I_p$, the definition in
\eqref{def_K_theta} gives, when $\zeta\neq\xi$,
\[
        K_n^\Theta(\zeta,\xi)
        =
        \frac{1}{\zeta\xi}
        \frac{v_n(\zeta)-v_n(\xi)}
        {(\zeta-\xi)v_n(\zeta)v_n(\xi)}
        \int_0^\infty
        \frac{t}
        {\{1+t v_n(\zeta)\}\{1+t v_n(\xi)\}}
        \,dH_n(t),
\]
with the divided difference replaced by the derivative when
$\zeta=\xi$. Hence, almost surely, for every
$\zeta,\xi\in\mathbb C\setminus\mathbb R_{\geq0}$,
\begin{align} \label{eq:limit}
        \frac1p
        \tr\!\left[
        R_n(\zeta)\Sigma_nR_n(\xi)
        \right]
        \longrightarrow
        \frac{\Delta(\zeta,\xi)}
        {\zeta\xi v(\zeta)v(\xi)}
        \int_0^\infty
        \frac{t}
        {\{1+t v(\zeta)\}\{1+t v(\xi)\}}
        \,dH(t).
\end{align}
Note that the fixed-point equation \eqref{stielt} implies
\[
        \int_0^\infty
        \frac{1}{1+t v(\zeta)}\,dH(t)
        =
        1-\frac{1+\zeta v(\zeta)}{\gamma}.
\]
Then, the limit in \eqref{eq:limit} can
equivalently be written as
\[
        \begin{cases}
        \displaystyle
        \frac{\zeta v(\zeta)-\xi v(\xi)}
        {\gamma(\zeta-\xi)\zeta\xi v(\zeta)v(\xi)},
        & \zeta\neq\xi,\\[4mm]
        \displaystyle
        \frac{
        v(\zeta)+\zeta\frac{d}{d\zeta}v(\zeta)
        }{
        \gamma\zeta^2v(\zeta)^2
        },
        & \zeta=\xi.
        \end{cases}
\]

For $\Theta_n=\beta_n\beta_n^\top$, let $\beta_n$ and $G_n$
satisfy the additional assumptions in
Remark~\ref{rem:single-resolvent-limits}.
The same calculation, with $\|\beta_n\|_2^2\,dG_n$ in place
of $dH_n$, yields, almost surely, for every
$\zeta,\xi\in\mathbb C\setminus\mathbb R_{\geq0}$,
\begin{align}\label{eq:double_rank_one}
        \beta_n^\top
        R_n(\zeta)\Sigma_nR_n(\xi)\beta_n
        \longrightarrow
        \frac{\tau^2\Delta(\zeta,\xi)}
        {\zeta\xi v(\zeta)v(\xi)}
        \int_0^\infty
        \frac{t}
        {\{1+t v(\zeta)\}\{1+t v(\xi)\}}
        \,dG(t).
\end{align}
In particular, for $\zeta = \xi$,
\[
        \beta_n^\top
        R_n(\zeta)\Sigma_nR_n(\zeta)\beta_n
        \longrightarrow
        \frac{\tau^2}{\zeta^2}
        \frac{\frac{d}{d\zeta}v(\zeta)}{v(\zeta)^2}
        \int_0^\infty
        \frac{t}{\{1+t v(\zeta)\}^2}\,dG(t).
\]

These conclusions follow by passing to the limit in
$K_n^\Theta$ and applying
Theorem~\ref{thm:double-resolvent-limit2}.
For distinct spectral parameters $\zeta \neq \xi $, the required convergence
follows from $v_n\to v$ and the common compact support of
$H_n$ and $G_n$, as in
Remark~\ref{rem:single-resolvent-limits}.
For $\zeta=\xi$, choose a circle $\Gamma$ centered at $\zeta$
whose closed interior is contained in
$\mathbb C\setminus\mathbb R_{\geq0}$.
The Stieltjes bound
$|v_n(w)|\leq\operatorname{dist}(w,\mathbb R_{\geq0})^{-1}$
provides a common integrable bound for
$v_n(w)/(w-\zeta)^2$ on $\Gamma$.
Hence \eqref{eq:pointwise_vn_v}, Cauchy's integral formula,
and dominated convergence give
\[
    \frac{d}{d\zeta}v_n(\zeta)
    \longrightarrow
    \frac{d}{d\zeta}v(\zeta).
\]
Thus, the same conclusion holds for $\zeta=\xi$.
}
\end{remark}

\section{Applications}\label{sec:applications}
The single- and double-resolvent expressions studied in Sections \ref{sec:RM-repair} and \ref{sec:two-resolvent} appear naturally in high-dimensional statistical learning.
The uniform bounds over compact sets are particularly helpful for common tasks such as varying a regularization parameter or selecting it by minimizing a risk criterion.
We illustrate their usefulness through three applications.
Section~\ref{sec:app-precision} considers precision-matrix shrinkage.
In Section~\ref{sec:app-ridge}, we derive a deterministic equivalent
for the prediction risk of high-dimensional ridge regression,
and Section~\ref{sec:app-tuning} studies generalized cross-validation
for ridge regression.

\subsection{Precision-matrix shrinkage}
\label{sec:app-precision}

For a population covariance matrix $\Sigma_n$, consider
the regularized inverse covariance estimator
\[
    \widehat\Omega_n(\alpha,\lambda)
    \coloneqq
    \alpha R_n(-\lambda),
    \qquad
    \alpha>0,\quad
    \lambda\in\Lambda=[\lambda_-,\lambda_+]\subset(0,\infty).
\]
Following \citet[Section~3]{WangPanTongZhu2015}, we consider
minimization over $(\alpha,\lambda)$ of the normalized Frobenius loss
\begin{equation}\label{eq:app-precision-loss}
\begin{split}
    \mathcal Q_n(\alpha,\lambda)
    &\coloneqq
    \frac{1}{p}
    \tr\!\left[
        \left(
            \Sigma_n^{1/2}
            \widehat\Omega_n(\alpha,\lambda)
            \Sigma_n^{1/2}-I_p
        \right)^2
    \right]
    =
    \alpha^2T_{2,n}(\lambda)
    -
    2\alpha T_{1,n}(\lambda)
    +
    1,
\end{split}
\end{equation}
where
$T_{1,n}(\lambda)
\coloneqq p^{-1}\tr[\Sigma_nR_n(-\lambda)]$
and
$T_{2,n}(\lambda)
\coloneqq p^{-1}\tr[
\Sigma_nR_n(-\lambda)\Sigma_nR_n(-\lambda)]$. 
The loss remains well-defined when $\Sigma_n$ is singular.
Indeed, writing $\Sigma_n^\dagger$ for the Moore--Penrose
pseudoinverse, we have
\[
    \mathcal Q_n(\alpha,\lambda)
    =
    \frac{1}{p}
    \left\|
        \Sigma_n^{1/2}
        (\widehat\Omega_n(\alpha,\lambda)-\Sigma_n^\dagger)
        \Sigma_n^{1/2}
    \right\|_F^2+
    \frac{p-\operatorname{rank}(\Sigma_n)}{p}.
\]
The additional term is independent of $\alpha$ and $\lambda$
and therefore does not affect their optimal choice.

Since the loss $\mathcal 
Q_n (\alpha, \lambda)$ depends on the unknown covariance matrix
$\Sigma_n$, it cannot be minimized directly. Our results in Sections~\ref{sec:RM-repair}
and~\ref{sec:two-resolvent} provide deterministic
approximations of these functionals.
We then use these approximations to construct an observable
criterion whose minimizer is asymptotically optimal
over the prescribed parameter set.

Define
\[
    T_{1,n}^{\mathrm{det}}(\lambda)
    \coloneqq
    \frac{1}{p}
    \tr\!\left[
        \Sigma_nR_n^{\mathrm{det}}(-\lambda)
    \right],
    \qquad
    T_{2,n}^{\mathrm{det}}(\lambda)
    \coloneqq
    K_n^{\Sigma_n/p}(-\lambda,-\lambda),\qquad\text{for }\lambda>0.
\]
By the definitions of $R_n^{\mathrm{det}}$ in \eqref{eq:def_R} and $K_n^\Theta$ in \eqref{def_K_theta},
together with \eqref{stielt} for $(\gamma_0,H_0)=(\gamma_n,H_n)$, it follows that

\begin{equation}\label{eq:app-precision-det-traces}
    T_{1,n}^{\mathrm{det}}(\lambda)
    =
    \frac{1}{\gamma_n}
    \left\{
        \frac{1}{\lambda v_n(-\lambda)}-1
    \right\}\qquad\text{and}\qquad
    T_{2,n}^{\mathrm{det}}(\lambda)
    =
    \frac{
        -\frac{d}{d\lambda}v_n(-\lambda)
        -v_n(-\lambda)^2
    }{
        \gamma_n\lambda^2v_n(-\lambda)^4
    }.
\end{equation}
Moreover, define, in analogy to \eqref{eq:app-precision-loss},
\[
    \mathcal Q_n^{\mathrm{det}}(\alpha,\lambda)
    \coloneqq
    \alpha^2T_{2,n}^{\mathrm{det}}(\lambda)
    -
    2\alpha T_{1,n}^{\mathrm{det}}(\lambda)
    +
    1.
\]
The representations in \eqref{eq:app-precision-det-traces}
suggest replacing $v_n$ by $\widehat v_n$ from
\eqref{eq:emp-companion-definition}. The latter also satisfies
\begin{align}\label{eq:emp_fixed_point}
    \widehat v_n(z)
    =-\frac{1-\gamma_n}{z}
      +\frac{\gamma_n}{p}\tr[R_n(z)],
      \qquad z\in\mathbb C\setminus\mathbb R_{\geq0}.
\end{align}
Hence, we define estimators for $T_{1,n}^{\mathrm{det}}$ and $T_{2,n}^{\mathrm{det}}$ in \eqref{eq:app-precision-det-traces} by
\begin{equation}\label{eq:app-precision-empirical-traces-def}
    \widehat T_{1,n}(\lambda)
    \coloneqq
    \frac{1}{\gamma_n}
    \left\{
        \frac{1}{\lambda\widehat v_n(-\lambda)}-1
    \right\}\qquad\text{and}\qquad
    \widehat T_{2,n}(\lambda)
    \coloneqq
    \frac{
        -\frac{d}{d\lambda}\widehat v_n(-\lambda)
        -\widehat v_n(-\lambda)^2
    }{
        \gamma_n\lambda^2
        \widehat v_n(-\lambda)^4
    },\qquad \lambda>0.
\end{equation}
These estimators yield the empirical loss criterion
\begin{equation}\label{eq:app-precision-empirical-loss-def}
    \widehat{\mathcal Q}_n(\alpha,\lambda)
    \coloneqq
    \alpha^2\widehat T_{2,n}(\lambda)
    -
    2\alpha\widehat T_{1,n}(\lambda)
    +
    1.
\end{equation}
The following proposition establishes its uniform approximation
to the deterministic loss and controls the excess loss of
the resulting data-driven parameter choice.

\begin{proposition}\label{prop:app-precision-loss}
Suppose Assumptions~\ref{A3}--\ref{A2} and
\eqref{gamma_n:cond} hold.
Then, for every pair of fixed compact intervals
$\mathcal A,\Lambda\subset(0,\infty)$,
\begin{equation}\label{eq:app-precision-empirical-de}
    \sup_{(\alpha,\lambda)\in\mathcal A\times\Lambda}
    \left\{\left|
        \widehat{\mathcal Q}_n(\alpha,\lambda)
        -
        \mathcal Q_n^{\mathrm{det}}(\alpha,\lambda)
    \right|+\left|
        \mathcal Q_n^{\mathrm{det}}(\alpha,\lambda)- \mathcal Q_n(\alpha,\lambda)
    \right|\right\}
    =
    O_{\mathbb P}(n^{-1}).
\end{equation}
Moreover, any measurable minimizer
$(\widehat\alpha_n,\widehat\lambda_n)
\in\operatorname*{arg\,min}_{(\alpha,\lambda)
\in\mathcal A\times\Lambda}
\widehat{\mathcal Q}_n(\alpha,\lambda)$ satisfies
\begin{equation}\label{eq:app-precision-data-driven}
    \mathcal Q_n(\widehat\alpha_n,\widehat\lambda_n)
    -
    \inf_{(\alpha,\lambda)\in\mathcal A\times\Lambda}
    \mathcal Q_n(\alpha,\lambda)
    =
    O_{\mathbb P}(n^{-1}).
\end{equation}
\end{proposition}

Proposition~\ref{prop:app-precision-loss} shows that the
deterministic equivalent of the loss in
\eqref{eq:app-precision-loss}, and hence the loss itself, can be approximated uniformly
by an observable criterion. Minimizing this criterion yields
data-driven parameters whose excess loss relative to the optimal oracle parameter
choice in $\mathcal A\times\Lambda$ is of order
$O_{\mathbb P}(n^{-1})$. Since $\widehat{\mathcal Q}_n$ is measurable in the data
and continuous in $(\alpha,\lambda)$ on the compact set
$\mathcal A\times\Lambda$, a measurable minimizer exists; see \citet[Exercise~7.2.3]{GineNickl2015}. Neither knowledge of $\Sigma_n$
nor convergence of $\gamma_n$ or $H_n$ is required.

For implementation, let
$\mathcal A=[\alpha_-,\alpha_+]\subset(0,\infty)$.
Since $\widehat v_n$ is the Stieltjes transform of a
probability measure on $[0,\infty)$, we have
$\lambda\widehat v_n(-\lambda)\leq1$ and
$-\frac{d}{d\lambda}\widehat v_n(-\lambda)
\geq\widehat v_n(-\lambda)^2$.
Thus, both $\hat T_{1,n}(\lambda)$ and $\hat T_{2,n}(\lambda)$ in
\eqref{eq:app-precision-empirical-traces-def} are nonnegative and for each $\lambda\in\Lambda$, a minimizer of
$\widehat{\mathcal Q}_n(\alpha,\lambda)$ in \eqref{eq:app-precision-empirical-loss-def} over
$\alpha\in\mathcal A$ is given by
\[
    \widehat\alpha_n(\lambda)
    \coloneqq
    \begin{cases}
        \displaystyle
        \min\!\left\{
            \alpha_+,
            \max\!\left\{
                \alpha_-,
                \frac{\widehat T_{1,n}(\lambda)}
                     {\widehat T_{2,n}(\lambda)}
            \right\}
        \right\},
        & \widehat T_{2,n}(\lambda)>0,\\[3mm]
        \alpha_+,
        & \widehat T_{2,n}(\lambda)=0.
    \end{cases}
\]
Hence, the joint minimization over
$\mathcal A\times\Lambda$ reduces to
\[
    \widehat\lambda_n
    \in
    \operatorname*{arg\,min}_{\lambda\in\Lambda}
    \widehat{\mathcal Q}_n(
        \widehat\alpha_n(\lambda),\lambda),
    \qquad
    \widehat\alpha_n
    =
    \widehat\alpha_n(\widehat\lambda_n).
\]
The resulting procedure requires only the empirical companion
Stieltjes transform $\widehat v_n$ and its derivative.
The remaining one-dimensional optimization over $\Lambda$
can be carried out numerically, for instance by grid search.

\subsection{Ridge regression and continuous spectral filters}
\label{sec:app-ridge}

Consider the linear model
\begin{align}\label{eq:lin_mod}
    y_n=X_n\beta_n+u_n,
\end{align}
where $X_n$ is as in Assumption~\ref{A3},
$\beta_n\in\mathbb R^p\setminus\{0\}$ is deterministic, and
$\mathbb E[u_n\mid X_n]=0$ and
$\mathbb E[u_nu_n^\top\mid X_n]=\sigma_n^2I_n$,
with $\sigma_n\geq0$ for every $n$.
The ridge estimator is given by
\[
    \widehat\beta_n(\lambda)
    =R_n(-\lambda)\frac{X_n^\top y_n}{n},\qquad \lambda>0.
\]
Let $x_0$ be a test vector independent of $(X_n,u_n)$, with
mean zero and covariance $\Sigma_n$. We define the prediction
risk conditional on $X_n$ by
\[
    \mathcal R_n(\lambda)
    \coloneqq
    \mathbb E\!\left[
        \bigl\{x_0^\top
        (\widehat\beta_n(\lambda)-\beta_n)\bigr\}^2
        \,\middle|\,X_n
    \right],
\]
which admits the
following decomposition
\begin{equation}\label{eq:app-ridge-risk}
    \mathcal R_n(\lambda)
    =
    \lambda^2
    \beta_n^\top R_n(-\lambda)\Sigma_nR_n(-\lambda)\beta_n
    +
    \frac{\sigma_n^2 \gamma_n}{p}
    \operatorname{tr}\!\left[
        \Sigma_n\widehat\Sigma_n R_n(-\lambda)^2
    \right];
\end{equation}
see for example \cite{dobriban:2018, HastieMontanariRossetTibshirani2022}. 
The first term is the expression on the left-hand side of \eqref{eq:double_rank_one}, with
$\zeta=\xi=-\lambda$.
For the second term in \eqref{eq:app-ridge-risk}, the identity
\[
    \widehat\Sigma_nR_n(-\lambda)^2
    =R_n(-\lambda)-\lambda R_n(-\lambda)^2 = R_n(-\lambda)
      +\lambda\frac{d}{d\lambda}R_n(-\lambda)
\]
reduces the calculation to a single-resolvent trace and its
derivative. Define $h_n(z)\coloneqq p^{-1}\tr[\Sigma_nR_n^{\mathrm{det}}(z)]$
and
\begin{equation}\label{eq:app-ridge-risk-de}
    \mathcal R_n^{\mathrm{det}}(\lambda)
    \coloneqq
    \lambda^2
    K_n^{\beta_n\beta_n^\top}(-\lambda,-\lambda)
    +
    \sigma_n^2\gamma_n
    \left\{
        h_n(-\lambda)
        +\lambda\frac{d}{d\lambda}h_n(-\lambda)
    \right\}.
\end{equation}

\begin{proposition}\label{prop:app-ridge-risk}
Suppose Assumptions~\ref{A3}--\ref{A2} and \eqref{gamma_n:cond} hold.
Then, for every compact interval
$\Lambda=[\lambda_-,\lambda_+]\subset(0,\infty)$,
\begin{equation}\label{eq:app-ridge-risk-rate}
    \sup_{\lambda\in\Lambda}
    \left|
        \mathcal R_n(\lambda)
        -\mathcal R_n^{\mathrm{det}}(\lambda)
    \right|
    =
    O_{\mathbb P}\!\left(
        \frac{\|\beta_n\|_2^2}{\sqrt n}+\frac{\sigma_n^2}{n}
    \right).
\end{equation}
\end{proposition}

\begin{proof}
Using \eqref{eq:app-ridge-risk} and \eqref{eq:app-ridge-risk-de} the claim follows directly by applying Theorem~\ref{Lemma:GreenRom},
Corollary~\ref{cor:GreenRom-derivatives:rates},
and Theorem~\ref{thm:double-resolvent-limit},
with $\|\beta_n\beta_n^\top\|_F=\|\beta_n\|_2^2$
and $\|\Sigma_n/p\|_F=O(n^{-1/2})$.
\end{proof}
Quantitative deterministic approximations of ridge prediction
risk are also established in
\citet[Theorem~5]{HastieMontanariRossetTibshirani2022},
albeit under different assumptions and with different
probabilistic error bounds. Their analysis relies on the
anisotropic local law for sample covariance matrices of
\citet[Theorem~3.16(i) and Remark~3.17]{knowles2017anisotropic}.
Their result assumes i.i.d. observations with independent
standardized entries having uniformly bounded moments of all
orders, a bounded inverse spectral moment
$p^{-1}\tr(\Sigma_n^{-1})$, and
$M^{-1}\leq\gamma_n\leq M$ with
$|1-\gamma_n|\geq M^{-1}$ for some constant $M>1$.
In contrast, our results require only uniformly bounded fourth moments, allow non-identically distributed observations and
singular population covariance matrices, and impose no
separation of $\gamma_n$ from one.
\cite{HastieMontanariRossetTibshirani2022} provide high-probability error bounds, whereas
our results yield the stochastic-order approximation
\eqref{eq:app-ridge-risk-rate}.

The preceding ridge-risk analysis can be extended to a broader
class of estimators defined by continuous functions of the
sample covariance matrix. Let $f:[0,\infty)\to\mathbb R$ be continuous
and consider
\[
    \widehat\beta_{n,f}
    \coloneqq
    f(\widehat\Sigma_n)\frac{X_n^\top y_n}{n}
    .
\]
Here, $f(\widehat\Sigma_n)$ is defined by applying $f$ to
the eigenvalues of $\widehat\Sigma_n$ while keeping its
eigenvectors unchanged.
In particular, ridge regression corresponds to the continuous
spectral filter $f_\lambda(s)=(s+\lambda)^{-1}$.
Writing $r_f(s)\coloneqq1-sf(s)$, the
 prediction risk of $\hat\beta_{n,f}$ can be written as
\begin{equation}\label{eq:app-filter-risk}
\begin{split}
    \mathcal R_{n,f}
    =
    \beta_n^\top
    r_f(\widehat\Sigma_n)\Sigma_n
    r_f(\widehat\Sigma_n)\beta_n+
    \frac{\sigma_n^2\gamma_n}{p}
    \operatorname{tr}\!\left[
        \Sigma_n\widehat\Sigma_n
        f(\widehat\Sigma_n)^2
    \right].
\end{split}
\end{equation}

In a companion paper \cite{StarkKoestenbergerDoernemann2026},
we use the present results to study the prediction risk
in \eqref{eq:app-filter-risk} for arbitrary continuous spectral
filters. A central difficulty is that convergence of Stieltjes
transforms does not, in general, imply weak convergence of signed
measures. Building on the deterministic equivalents established
here, we overcome this difficulty and characterize the limiting
prediction risk without requiring weak convergence of the
underlying signed eigenvector-overlap measures.

\subsection{Generalized cross-validation for ridge regression}
\label{sec:app-tuning}

The deterministic risk approximation in
Proposition~\ref{prop:app-ridge-risk} can be combined with
an observable criterion to select the regularization
parameter from the data. We consider generalized
cross-validation (GCV) in the linear model \eqref{eq:lin_mod}
and restrict attention to a fixed compact interval
$\Lambda=[\lambda_-,\lambda_+]\subset(0,\infty)$.

Unlike the preceding ridge-risk analysis, where
$\mathcal R_n(\lambda)$ conditions only on $X_n$ and
averages over the training noise, regularization tuning
concerns the prediction risk conditional on the entire
training sample $(X_n,y_n)$.
Let $(x_0,y_0)$ be an independent test observation,
where $x_0$ has mean zero and covariance $\Sigma_n$, and
$y_0=x_0^\top\beta_n+u_0$, with $u_0$ independent of
$x_0$, having mean zero and variance $\sigma_n^2$.
We define
\begin{equation}\label{eq:app-cv-risk}
\begin{split}
    \mathcal L_n(\lambda)
    &\coloneqq
    \mathbb E\!\left[
        \bigl(x_0^\top\widehat\beta_n(\lambda)-y_0\bigr)^2
        \,\middle|\,X_n,y_n
    \right]
\end{split}
\end{equation}
Thus, unlike $\mathcal R_n(\lambda)$, the conditional
prediction risk $\mathcal L_n(\lambda)$
also accounts for the realized training noise. The GCV criterion is defined by 
\begin{equation}\label{eq:app-gcv-def}
    \operatorname{GCV}_n(\lambda)
    \coloneqq
    \frac{
        n^{-1}\|y_n-X_n\widehat\beta_n(\lambda)\|_2^2
    }{
        \left\{
            1-\frac{1}{n}
            \tr[\widehat\Sigma_nR_n(-\lambda)]
        \right\}^2
    },
\end{equation}
(see, e.g., \citet{HastieMontanariRossetTibshirani2022}, \cite{PatilWeiRinaldoTibshirani2021}).
Write
\[
    \underline R_n(z)
    \coloneqq
    \left(\frac{X_nX_n^\top}{n}-zI_n\right)^{-1},
    \qquad
    \widehat v_n(z)
    =
    \frac{1}{n}\tr[\underline R_n(z)], \quad z\in \mathbb C \setminus \mathbb R_{\geq0}.
\]
Since
$y_n-X_n\widehat\beta_n(\lambda)
=\lambda\underline R_n(-\lambda)y_n$
and
$1-n^{-1}\tr[\widehat\Sigma_nR_n(-\lambda)]
=\lambda\widehat v_n(-\lambda)$,
we obtain
\begin{equation}\label{eq:app-gcv-companion-form}
    \operatorname{GCV}_n(\lambda)
    =
    \frac{
        n^{-1}y_n^\top\underline R_n(-\lambda)^2y_n
    }{
        \widehat v_n(-\lambda)^2
    }.
\end{equation}
The following proposition establishes a quantitative
guarantee for data-driven regularization tuning by
controlling the excess conditional prediction risk
incurred by minimizing GCV.

\begin{proposition}\label{prop:app-ridge-tuning}
Suppose Assumptions~\ref{A3}--\ref{A2} and
\eqref{gamma_n:cond} hold, and
$\sup_n\|\beta_n\|_2<\infty$.
Assume additionally that $u_n$ is independent of $X_n$
and has independent entries satisfying
$\mathbb E[u_{n,j}]=0$,
$\mathbb E[u_{n,j}^2]=\sigma_n^2$, and
$\sup_{n,j}\mathbb E|u_{n,j}|^4<\infty$.
Then, for every fixed compact interval
$\Lambda=[\lambda_-,\lambda_+]\subset(0,\infty)$, 
\begin{equation}\label{eq:app-gcv-empirical-de}
    \sup_{\lambda\in\Lambda}
    \left|
        \operatorname{GCV}_n(\lambda)
        -
       \mathcal L_n(\lambda)
    \right|
    =
    O_{\mathbb P}(n^{-1/2}).
\end{equation}
Moreover, any measurable minimizer
$\widehat\lambda_n
\in\operatorname*{arg\,min}_{\lambda\in\Lambda}
\operatorname{GCV}_n(\lambda)$ satisfies
\begin{equation}\label{eq:app-gcv-risk-optimality}
    \mathcal L_n(\widehat\lambda_n)
    -
    \inf_{\lambda\in\Lambda}\mathcal L_n(\lambda)
    =
    O_{\mathbb P}(n^{-1/2}).
\end{equation}
\end{proposition}

Proposition~\ref{prop:app-ridge-tuning} shows that GCV
uniformly approximates the conditional prediction risk
$\mathcal L_n$ with an error of order $O_{\mathbb P}(n^{-1/2})$.
Consequently, minimizing GCV yields a data-driven
regularization parameter whose excess prediction risk
is of the same order. As in Section \ref{sec:app-precision}, a measurable minimizer
exists because $\operatorname{GCV}_n(\lambda)$ is measurable in the data and continuous on
$\Lambda$; see \citet[Exercise~7.2.3]{GineNickl2015}. 

A related result is established by
\citet[Theorem~4.1]{PatilWeiRinaldoTibshirani2021},
who prove almost-sure uniform consistency of GCV,
including for zero and certain negative regularization
parameters. Note that the proof of \citet{PatilWeiRinaldoTibshirani2021} uses an adaptation
of Assertion~\ref{assrt:RM-lemma4}
(Lemma~4 of \citet{RubioMestre2011}).
However, this assertion does not hold under its stated assumptions; see Proposition~\ref{prop1}.
Its use in their proof therefore requires additional justification.
In contrast, our quantitative bounds hold
under uniformly bounded fourth moments, allow
non-identically distributed entries of $Z_n$ and singular
population covariance matrices, and do not require
convergence of $\gamma_n$ or $H_n$.

\section{Averaging of centered leave-one-out quadratic forms}
\label{sec:Rubmes_erroranalysis}

In this section, we study conditions under which averages
of centered quadratic forms converge to zero.
This problem arises in the proofs of our main results in
Sections~\ref{sec:RM-repair} and~\ref{sec:two-resolvent},
but is also of independent interest.
We first show that such averages need not converge to zero,
even when the random matrix in each quadratic form is
independent of the corresponding observation.
We then establish convergence in our sample covariance
setting by exploiting the additional structure of
leave-one-out resolvent matrices.

For $D\subseteq\{1,\ldots,n\}$, define the sample covariance
matrix obtained by excluding the observations indexed by
$D$, together with its corresponding resolvent, by
\begin{align}\label{eq:D_del_res_cov}
    \widehat\Sigma_{n,-D}
    \coloneqq
    \frac1n\sum_{i\notin D}x_ix_i^\top,
    \qquad
    R_{n,-D}(\zeta)
    \coloneqq
    (\widehat\Sigma_{n,-D}-\zeta I_p)^{-1},
    \quad
    \zeta\in\mathbb C\setminus\mathbb R_{\geq0}.
\end{align}
We also write $\mathcal F_{n,-D}
    \coloneqq
    \sigma(x_i:i\notin D)$, for the $\sigma$-algebra generated by all observations
except those indexed by $D$.
For convenience, for pairwise distinct indices
$j_1,\ldots,j_m\in\{1,\ldots,n\}$, we use the notation
\[
    R_{n,-j_1\cdots j_m}(\zeta)
    \coloneqq
    R_{n,-\{j_1,\ldots,j_m\}}(\zeta),
    \qquad
    \mathcal F_{n,-j_1\cdots j_m}
    \coloneqq
    \mathcal F_{n,-\{j_1,\ldots,j_m\}}.
\]
In particular,
$R_{n,-\varnothing}(\zeta)=R_n(\zeta)$.
For the martingale arguments, we write
\[
    \mathcal F_j=\sigma(x_1,\ldots,x_j),\quad j=1,\ldots,n,
    \qquad \mathcal F_0=\{\varnothing,\Omega\}.
\]

Our starting point is the following assertion from
\citet[Lemma~4]{RubioMestre2011}, stated in the notation
of the present paper.

\begin{assertion}
\label{assrt:RM-lemma4}
For each $n\in\mathbb N$, let $x_1,\ldots,x_n$ be i.i.d. random
vectors taking values in $\mathbb R^p$, where $p=p_n$ may depend on $n$.
Suppose that their entries are i.i.d. with mean zero, variance
one, and uniformly bounded $(4+\delta)$-moments for some
$\delta>0$. For $j=1,\ldots,n$, let $C_{n,-j}$ be a random
$p\times p$ matrix measurable with respect to $\mathcal F_{n,-j}$.
Suppose that
\[
        \sup_{n\geq1}\max_{1\leq j\leq n}
        \|C_{n,-j}\|_{1}<\infty.
\]
Then, as $n\to\infty,$
\[
        \frac1n\sum_{j=1}^n
        \left\{
        x_j^\top C_{n,-j}x_j-\tr(C_{n,-j})
        \right\}
        \longrightarrow0
        \qquad\text{almost surely}.
\]
\end{assertion}
The proof of Assertion~\ref{assrt:RM-lemma4} in
\citet{RubioMestre2011} treats the centered quadratic forms $$\eta_j
        \coloneqq
        x_j^\top C_{n,-j}x_j-\tr(C_{n,-j})$$ as
martingale differences and applies a Burkholder inequality;
see \citet[Lemma~2.12]{Silversteinbook2010}.
However, $C_{n,-j}$ may depend on
$x_{j+1},\ldots,x_n$, so $\eta_j$ need not be
$\mathcal F_j$-measurable. Consequently, the sequence
$\{\eta_j\}_{j=1}^n$ need not be a martingale-difference sequence
with respect to $\{\mathcal F_j\}_{j=0}^n$. Requiring $C_{n,-j}$ to be $\mathcal F_{j-1}$-measurable
would justify this martingale-difference argument. Such a requirement would, however, exclude the quadratic
forms involving the leave-one-out resolvent matrices
$R_{n,-j}$ used in the proofs of
Sections~\ref{sec:RM-repair} and~\ref{sec:two-resolvent}. 

The lack of adaptedness is not merely a technical issue in the proof: Noting that Assertion~\ref{assrt:RM-lemma4} imposes no relation between $p$ and $n$, the following result shows that its conclusion already fails for $p=2$. We give a short constructive proof below. However, the failure of convergence to zero is not specific to the fixed-dimensional setting. In Remark~\ref{rem:proportional-counterexample}, we explain the extension to arbitrary dimension sequences satisfying $p_n\geq2$.
\begin{proposition}
\label{prop1}
There exist i.i.d. random vectors $\{x_j\}_{j\geq1}$ in
$\mathbb R^2$ and random matrices
$C_{n,-j}$, $n\in\mathbb N$, $1\leq j \leq n$, satisfying
the assumptions of Assertion~\ref{assrt:RM-lemma4}, such that
\[
        \left|
        \frac1n\sum_{j=1}^n
        \left\{
        x_j^\top C_{n,-j}x_j-\tr(C_{n,-j})
        \right\}
        \right|
        =1
        \qquad\text{almost surely, for every }n\in\mathbb N.
\]
\end{proposition}

\begin{proof}
Let $\{\varepsilon_j,\eta_j:j\geq1\}$ be independent Rademacher
random variables, and set
$
        x_j=(\varepsilon_j,\eta_j)^\top.
$
The vectors $x_j$ are i.i.d., and their entries have mean zero,
variance one, and bounded moments of every order.
For $n\in\mathbb N$ and $j=1,\ldots,n$, define
\[
        a_{n,-j}
        \coloneqq
        \prod_{\substack{1\leq k\leq n\\k\neq j}}
        \varepsilon_k\eta_k,
        \qquad
        C_{n,-j}
        \coloneqq
        \begin{pmatrix}
        0 & a_{n,-j}/2\\
        a_{n,-j}/2 & 0
        \end{pmatrix},
\]
where an empty product equals one.
The matrix $C_{n,-j}$ is measurable with respect to
$\mathcal F_{n,-j}$.
Moreover, since $a_{n,-j}^2=1$,
\[
        C_{n,-j}^\top C_{n,-j}=\frac14 I_2.
\]
Its two singular values are therefore both equal to $1/2$,
and hence $\|C_{n,-j}\|_{1}=1$.

Since $\tr(C_{n,-j})=0$, we have
\[
\begin{split}
        x_j^\top C_{n,-j}x_j-\tr(C_{n,-j})
        &=
        a_{n,-j}\varepsilon_j\eta_j
        =
        \prod_{k=1}^n\varepsilon_k\eta_k.
\end{split}
\]
Thus all $n$ centered quadratic forms are equal to the same
random variable, and
\[
        \frac1n\sum_{j=1}^n
        \left\{
        x_j^\top C_{n,-j}x_j-\tr(C_{n,-j})
        \right\}
        =
        \prod_{k=1}^n\varepsilon_k\eta_k.
\]
The right-hand side takes values in $\{-1,1\}$, which proves
the assertion.
\end{proof}

 \begin{remark}[Extension to arbitrary dimensions]
 \label{rem:proportional-counterexample}
{\rm For any dimension sequence with $p_n\geq2$, extend each $x_j$
by independent Rademacher coordinates and each $C_{n,-j}$
by zero rows and columns. All assumptions remain satisfied,
and the quadratic forms are unchanged. Thus, the counterexample
also applies in growing dimensions, including the proportional regime \eqref{gamma_n:cond}.}
\end{remark}

Proposition \ref{prop1} is relevant beyond the original argument of
\citet{RubioMestre2011}, since their deterministic-equivalent framework
has been used extensively in high-dimensional probability and
statistics. Direct applications or adaptations of
\citet[Lemma~4]{RubioMestre2011} appear, for example, in
\citet{BodnarGuptaParolya2014},
\citet{ParolyaHeinyKurowicka2024}, \citet{DobribanSheng2021}, and
\citet{PatilWeiRinaldoTibshirani2021}. More broadly, the
deterministic-equivalent \citet[Theorem~1]{RubioMestre2011} and related
resolvent approximations have been used in high-dimensional regression,
distributed linear and ridge regression, covariance shrinkage,
portfolio theory, and overparameterized learning; see, among others,
\citet{HastieMontanariRossetTibshirani2022},
\citet{DobribanSheng2021}, \citet{DobribanSheng2020},
\citet{Green2025HighDimPCR}, \citet{BodnarParolyaSchmid2018},
\citet{BodnarOkhrinParolya2022}, \citet{BodnarParolyaThorsen2024},
and \citet{PatilDuKuchibhotla2023}.
However, Proposition~\ref{prop1} does not imply that the conclusions of these works relying on \cite{RubioMestre2011} are false. It shows that the averaging step in Assertion~\ref{assrt:RM-lemma4} requires additional information about the dependence between the coefficient matrices.

In our sample covariance setting, the leave-one-out
resolvent matrices $R_{n,-j}$ have additional structure
that allows us to establish the desired convergence.
Specifically, the Sherman--Morrison formula allows us
to control how $R_{n,-j}(\zeta)$ changes when another
observation is removed.
In the next subsection, we establish averaging bounds
for centered quadratic forms involving one or two
leave-one-out resolvents.

\subsection{Averaging for resolvent quadratic forms} \label{sec:multi_resolvent}

Recall the sample covariance setting of Section~\ref{sec:RM-repair}. The first lemma gives a quantitative averaging bound for
coefficients of the form
\[
    E_n(\zeta)R_{n,-j}(\zeta)F_n(\zeta),
\]
where $E_n(\zeta)$ and $F_n(\zeta)$ are deterministic.
\begin{lemma}
\label{RubMes_Lem4_fix}
Suppose that Assumptions~\ref{A3} and~\ref{A2} hold.
Let $K\subset\mathbb C\setminus\mathbb R_{\geq0}$ be compact,
and let $K_1\subset\mathbb C\setminus\mathbb R_{\geq0}$ be
compact with $K\subset\operatorname{int}(K_1)$.
Let
\[
        E_n,F_n:
        \mathbb C\setminus\mathbb R_{\geq0}
        \longrightarrow\mathbb C^{p\times p}
\]
be deterministic matrix-valued functions whose entries
are holomorphic.
Suppose that, for a deterministic sequence $r_n>0$,
either
\begin{align}\label{cond:lem4}
        \sup_{n\geq1}\sup_{\zeta\in K_1}
        \|E_n(\zeta)\|_2<\infty,
        \qquad
        \sup_{\zeta\in K_1}
        \|F_n(\zeta)\|_F\leq C_{K_1}r_n,
\end{align}
or
\begin{align}\label{cond:lem4:2}
        \sup_{n\geq1}\sup_{\zeta\in K_1}
        \|F_n(\zeta)\|_2<\infty,
        \qquad
        \sup_{\zeta\in K_1}
        \|E_n(\zeta)\|_F\leq C_{K_1}r_n.
\end{align}
Then
\[
        \sup_{\zeta\in K}
        \left|
        \frac1n\sum_{j=1}^n
        \left\{
        x_j^\top E_n(\zeta)R_{n,-j}(\zeta)F_n(\zeta)x_j
        -
        \tr\!\left[
        \Sigma_nE_n(\zeta)R_{n,-j}(\zeta)F_n(\zeta)
        \right]
        \right\}
        \right|
        =
        O_{\mathbb P}\!\left(
        r_n\left\{
        \frac1{\sqrt n}\vee\frac{\sqrt p}{n}
        \right\}
        \right).
\]
\end{lemma}
Lemma~\ref{RubMes_Lem4_fix} does not require a restriction
on the aspect ratio. Under \eqref{gamma_n:cond}, its rate
reduces to $O_{\mathbb P}(r_n/\sqrt n)$.
The next lemma gives pointwise almost-sure convergence
under a uniform eighth-moment bound.

\begin{lemma}
\label{RubMes_Lem4_fix_2}
Suppose that Assumptions~\ref{A3} and~\ref{A2} hold,
that $p/n=O(1)$, and that \eqref{X:unif4} is strengthened to
\[
    \sup_{n\geq1}\max_{\substack{1\leq i\leq n\\1\leq j\leq p}}
    \mathbb E|z_{ij}|^8<\infty.
\]
Let $E_n,F_n:\mathbb C\setminus\mathbb R_{\geq0}
\to\mathbb C^{p\times p}$ be deterministic matrix-valued
functions.
Suppose that, for every
$\zeta\in\mathbb C\setminus\mathbb R_{\geq0}$, either
\[
        \sup_{n\geq1}\|E_n(\zeta)\|_2<\infty,
        \qquad
        \sup_{n\geq1}\|F_n(\zeta)\|_F<\infty,
\]
or
\[
        \sup_{n\geq1}\|F_n(\zeta)\|_2<\infty,
        \qquad
        \sup_{n\geq1}\|E_n(\zeta)\|_F<\infty.
\]
Then, for every fixed
$\zeta\in\mathbb C\setminus\mathbb R_{\geq0}$,
\[
        \frac1n\sum_{j=1}^n
        \left\{
        x_j^\top E_n(\zeta)R_{n,-j}(\zeta)F_n(\zeta)x_j
        -
        \tr\!\left[
        \Sigma_nE_n(\zeta)R_{n,-j}(\zeta)F_n(\zeta)
        \right]
        \right\}
        \longrightarrow0
        \qquad\text{almost surely}.
\]
\end{lemma}
The proofs of Lemmas~\ref{RubMes_Lem4_fix}
and~\ref{RubMes_Lem4_fix_2} are given in
Appendix~\ref{sec:averaging-proofs}.

We next establish a quantitative averaging bound for
quadratic forms involving matrices of the form
\[
    R_{n,-j}(\zeta)\Theta_nR_{n,-j}(\xi),
\]
where $\Theta_n$ is deterministic.
The next two lemmas extend the preceding averaging results
to this setting. In particular, the quantitative estimate
is used in the proof of
Lemma~\ref{lem:ridge-signal-equivalence:rates} to compare
trace functionals containing $\Sigma_n$ with those
containing $\hat\Sigma_n$.

\begin{lemma}
\label{lem:ridge_sandwich}
Suppose that Assumptions~\ref{A3} and~\ref{A2} hold.
Let $K\subset\mathbb C\setminus\mathbb R_{\geq0}$ be compact,
and let $\Theta_n$ be a sequence of real, deterministic
$p\times p$ matrices.
Then
\[
        \sup_{\zeta,\xi\in K}
        \left|
        \frac1n\sum_{j=1}^n
        \left\{
        x_j^\top
        R_{n,-j}(\zeta)\Theta_nR_{n,-j}(\xi)x_j
        -
        \tr\!\left[
        \Sigma_nR_{n,-j}(\zeta)\Theta_nR_{n,-j}(\xi)
        \right]
        \right\}
        \right|
        =
        O_{\mathbb P}\!\left(
        \|\Theta_n\|_F
        \left\{
        \frac1{\sqrt n}\vee\frac{\sqrt p}{n}
        \right\}
        \right).
\]
\end{lemma}

\begin{lemma}
\label{lem:ridge_sandwich2}
Suppose that Assumptions~\ref{A3} and~\ref{A2} hold,
that $p/n=O(1)$, and that \eqref{X:unif4} is strengthened to
\[
    \sup_{n\geq1}\max_{\substack{1\leq i\leq n\\1\leq j\leq p}}
    \mathbb E|z_{ij}|^8<\infty.
\]
Let $\Theta_n$ be a sequence of real, deterministic
$p\times p$ matrices satisfying $\sup_{n\geq1}\|\Theta_n\|_F<\infty$.
Then, for every fixed
$\zeta,\xi\in\mathbb C\setminus\mathbb R_{\geq0}$,
\[
        \frac1n\sum_{j=1}^n
        \left\{
        x_j^\top
        R_{n,-j}(\zeta)\Theta_nR_{n,-j}(\xi)x_j
        -
        \tr\!\left[
        \Sigma_nR_{n,-j}(\zeta)\Theta_nR_{n,-j}(\xi)
        \right]
        \right\}
        \longrightarrow0
        \qquad\text{almost surely}.
\]
\end{lemma}

The proofs of Lemmas~\ref{lem:ridge_sandwich}
and~\ref{lem:ridge_sandwich2} are given in
Appendix~\ref{sec:averaging-proofs}.
As in the single-resolvent case, the quantitative bound
retains the dependence on both $n$ and $p$.
Under \eqref{gamma_n:cond}, it becomes
$O_{\mathbb P}(\|\Theta_n\|_F/\sqrt n)$.

  \section*{Acknowledgments}
The work of TS and ND is partially supported by  the Aarhus University Research
Foundation (AUFF),  project numbers 47221 and 47388.
ND gratefully acknowledges the hospitality of the Department of
Statistics at Harvard University during the preparation of this
manuscript. ND's stay is supported by a RECRUIT Sabbatical
Grant from the Novo Nordisk Foundation, project number 52629. We thank Georg Köstenberger for a helpful proof suggestion.
  
   \setlength{\bibsep}{1pt}
\begin{small}
\bibliography{references}
\end{small}

\newpage
\appendix
\section*{Appendix}
\startcontents[appendix]

\printcontents[appendix]{}{1}{
    \subsection*{Contents}
}

\newpage
\section{Preliminaries and notations}
In the following proofs, $C=C_K$ denotes a generic positive constant
whose value may change from line to line. It may depend on $K$, on
an explicitly chosen compact neighborhood $K_1$, and on the constants
in the assumptions, but not on $n$, $p=p_n$, or the observation indices.
The leave-out matrices and sigma-algebras are those defined in
Section~\ref{sec:Rubmes_erroranalysis}.
We also use the Stieltjes transform of $F_{\widehat\Sigma_n}$,
\[
    \widehat m_n(z)\coloneqq\frac1p\tr[R_n(z)],
    \qquad z\in\mathbb C\setminus\mathbb R_{\geq0}.
\]
The matrices $X_n^\top X_n/n$ and $X_nX_n^\top/n$ have the same
nonzero eigenvalues. Consequently, their empirical spectral
measures satisfy
\begin{equation}\label{disrel}
    F_{\widehat\Sigma_n}
    =\left(1-\frac1{\gamma_n}\right)\delta_0
      +\frac1{\gamma_n}F_{\underline{\widehat\Sigma}_n}.
\end{equation}
Taking Stieltjes transforms gives
\begin{equation}\label{Stielcon}
    \widehat m_n(z)
    =-\left(1-\frac1{\gamma_n}\right)\frac1z
      +\frac1{\gamma_n}\widehat v_n(z),
    \qquad z\in\mathbb C\setminus\mathbb R_{\geq0}.
\end{equation}
Under Assumption~\ref{ass:main}, Lemma~\ref{thm:MP}
and \eqref{disrel} imply
\[
    F_{\widehat\Sigma_n}\Rightarrow F_{\gamma,H}
    \quad\text{almost surely},
    \qquad
    F_{\gamma,H}
    =\left(1-\frac1\gamma\right)\delta_0
      +\frac1\gamma\bar F_{\gamma,H}.
\]
Writing $m$ for the Stieltjes transform of $F_{\gamma,H}$,
we obtain
\begin{equation}\label{limitStil}
    m(z)=-\left(1-\frac1\gamma\right)\frac1z
         +\frac1\gamma v(z),
    \qquad z\in\mathbb C\setminus\mathbb R_{\geq0}.
\end{equation}

\section{Proofs of the results in Section~\ref{sec:RM-repair}}\label{sec:RM-repair_proofs}

\subsection{Lemma~\ref{stielt:var} and its proof}

\begin{lemma}
\label{stielt:var}
Suppose that Assumptions~\ref{A1}, \ref{A3}, and~\ref{A2} are satisfied.
Let \(K\subset\mathbb C\setminus\mathbb R_{\geq0}\) be compact. Then
\[
        \sup_{\zeta\in K}
        \left|
        \hat{v}_n(\zeta)-\mathbb E\bigl(\hat v_n(\zeta)\bigr)
        \right|
        =
        O_{\mathbb{P}}(n^{-1}).
\]
\end{lemma}

\begin{proof}[Proof of Lemma~\ref{stielt:var}]
For $0<\rho<\operatorname{dist}(K,\mathbb{R}_{\geq 0})/2$, define
\[
K_1\coloneqq \{w\in \mathbb C:\operatorname{dist}(w,K)\leq\rho\}.
\]

First, we are going to show that,
\[
\sup_{\zeta\in K_1}\mathbb{E}\bigl(\bigl|
        n\{\hat v_n(\zeta)-\mathbb E(\hat v_n(\zeta))\}
    \bigr|^2\bigr)\leq C.
\]
To this end, we derive a decomposition for $n\{\hat v_n(\zeta)-\mathbb E(\hat v_n(\zeta))\}$. To this end, we reproduce some arguments given in Section 9.10 of \cite{Silversteinbook2010} in our setting. Readers familiar with their work may skip directly to \eqref{eq:centered-UV-decomp}. From that point onward, our approach departs from \cite{Silversteinbook2010}.
Let
\[
        \mathcal F_j\coloneqq\sigma(x_1,\ldots,x_j),
        \qquad
        \mathbb E_j(\,\cdot\,)\coloneqq\mathbb E(\,\cdot\,\mid\mathcal F_j),
        \qquad j=0,1,\ldots,n,
\]
with \(\mathcal F_0\) the trivial sigma-field.  
Using \eqref{Stielcon}, we obtain
\begin{align}
\label{lem:b1_1}
        n\{\hat v_n(\zeta)-\mathbb E(\hat v_n(\zeta))\}
        =
        p\{\hat m_n(\zeta)-\mathbb E(\hat m_n(\zeta))\}
        =
        \sum_{j=1}^n
        (\mathbb E_j-\mathbb E_{j-1})
        \{\tr(R_n(\zeta))\}.
\end{align}
Here,
\[
(\E_j-\E_{j-1})\{\tr(R_n(\zeta))\}
=
\E_j\{\tr(R_n(\zeta))\}
-
\E_{j-1}\{\tr(R_n(\zeta))\}.
\]
Using the independence between \(R_{n,-j}(\zeta)\) and \(x_j\), we obtain
\begin{align}
\label{lem:b1_2}
(\E_j-\E_{j-1})
\bigl\{\tr(R_n(\zeta))\bigr\}
=
(\E_j-\E_{j-1})
\bigl\{
\tr(R_n(\zeta))
-
\tr(R_{n,-j}(\zeta))
\bigr\}.
\end{align}
By the Sherman--Morrison formula,
\begin{align}
\label{lem:b1_3}
\tr(R_n(\zeta))
-
\tr(R_{n,-j}(\zeta))
=
-
\frac{1}{
1+\frac1n x_j^\top R_{n,-j}(\zeta)x_j
}
\frac{
x_j^\top R_{n,-j}(\zeta)^2x_j
}{n}
=
-\beta_j(\zeta)q_j(\zeta),
\end{align}
where
\begin{align}
\label{eq:beta_jandq_j}
q_j(\zeta)
&\coloneqq
\frac1n
x_j^\top R_{n,-j}(\zeta)^2x_j,
\qquad
\beta_j(\zeta)
\coloneqq
\left(
1+\frac1n
x_j^\top R_{n,-j}(\zeta)x_j
\right)^{-1}.
\end{align}

Next, define
\begin{align}
\label{eq:a_jandb_j}
a_j(\zeta)
&\coloneqq
\frac1n
\tr\bigl\{
\Sigma_nR_{n,-j}(\zeta)
\bigr\},
\qquad
b_j(\zeta)
\coloneqq
\{1+a_j(\zeta)\}^{-1},
\end{align}
and
\begin{align}
\label{eq:gamma_def}
\gamma_j(\zeta)
\coloneqq
\frac1n
x_j^\top R_{n,-j}(\zeta)x_j
-
\frac1n
\tr\bigl\{
\Sigma_nR_{n,-j}(\zeta)
\bigr\}.
\end{align}
Note that
\begin{align}\label{eq:b_j_representation}
    \beta_j(\zeta)^{-1}
    =
    1+a_j(\zeta)+\gamma_j(\zeta)
    =
    b_j(\zeta)^{-1}+\gamma_j(\zeta).
\end{align}
Multiplying the preceding display by $\beta_j(\zeta)b_j(\zeta)$ and rearranging gives
\begin{align}\label{eq:beta_b_gamma}
    \beta_j(\zeta)
    =
    b_j(\zeta)-b_j(\zeta)\beta_j(\zeta)\gamma_j(\zeta).
\end{align}
Using \eqref{lem:b1_1}, \eqref{lem:b1_2}, \eqref{lem:b1_3} and \eqref{eq:beta_b_gamma}, we have for each $\zeta\in K_1$
\begin{align}
\begin{split}\label{dec1_Thm_main}
    n\{\hat{v}_n(\zeta)-\E(\hat{v}_n(\zeta))\}
    &=-\sum_{j=1}^n (\E_j-\E_{j-1})\{q_j(\zeta)b_j(\zeta)\}\\
    &\qquad+ \sum_{j=1}^n (\E_j-\E_{j-1})\{q_j(\zeta)b_j(\zeta)\beta_j(\zeta)\gamma_j(\zeta)\}.
    \end{split}
\end{align}
Furthermore, define, for each $\zeta\in K_1$,
\begin{align}\label{eq:tildeq_and_eta}
    \widetilde q_j(\zeta)
    \coloneqq
    \frac1n\operatorname{tr}\{\Sigma_nR_{n,-j}(\zeta)^2\}\qquad\text{and}\qquad
    \eta_j(\zeta)\coloneqq q_j(\zeta)-\widetilde q_j(\zeta).
\end{align}
Since \(b_j(\zeta)\widetilde q_j(\zeta)\) depends on
\(\{x_i:i\ne j\}\) but not on \(x_j\), and \(x_j\) is independent of
\(\{x_i:i\ne j\}\), we have
\[
    (\E_j-\E_{j-1})\{b_j(\zeta)\widetilde q_j(\zeta)\}=0.
\]
Therefore, using the decomposition in \eqref{dec1_Thm_main}, we get
\begin{align}\label{eq:centered-UV-decomp}
    n\{\hat v_n(\zeta)-\mathbb E(\hat v_n(\zeta))\}
    =
    -U_n(\zeta)+V_n(\zeta),
\end{align}
where
\[
    U_n(\zeta)
    \coloneqq
    \sum_{j=1}^n (\E_j-\E_{j-1})\{b_j(\zeta)\eta_j(\zeta)\},
\]
and 
\[
    V_n(\zeta)
    \coloneqq
    \sum_{j=1}^n
    (\E_j-\E_{j-1})\{b_j(\zeta)\beta_j(\zeta)q_j(\zeta)\gamma_j(\zeta)\}.
\]
Using \eqref{eq:general-trace-denominator-bound} with $c=n^{-1}$, $A=\widehat\Sigma_{n,-j}$ and $B=\Sigma_n$, we get
\begin{align}\label{eq:b-bound}
\sup_{\zeta\in K_1}\max_{1\leq j\leq n}|b_j(\zeta)|\leq C.
\end{align}

Next, we show that,
\begin{align}\label{eq:beta-q-bound}
    \sup_{\zeta\in K_1}\max_{1\le j\le n}
    |\beta_j(\zeta)q_j(\zeta)|
    \le C.
\end{align}
Recall \eqref{lem:b1_3}.
Since, $\operatorname{rank}\{R_n(\zeta)-R_{n,-j}(\zeta)\}\leq 1$, for each $n\geq 1$, $j\in \{1,\dots,n\}$ and $\zeta\in K_1$, we conclude that
\[
 | \beta_j(\zeta) q_j(\zeta)| 
 \leq \operatorname{rank}\bigl\{R_n(\zeta)-R_{n,-j}(\zeta)\bigr\}\|R_n(\zeta)-R_{n,-j}(\zeta)\|_2\leq \frac{2}{\operatorname{dist}(\zeta,\mathbb{R}_{\geq 0})}\leq  \frac{2}{\operatorname{dist}(K_1,\mathbb{R}_{\geq 0})},
\]
and therefore the statement in \eqref{eq:beta-q-bound} follows.

Fix \(j\in\{1,\ldots,n\}\) and \(\zeta\in K_1\). Applying
Lemma \ref{lem_b26_bai_silverstein} conditionally on
\(R_{n,-j}(\zeta)\), and using \(x_j=\Sigma_n^{1/2}z_j\), we obtain
\begin{align*}
    \mathbb E\!\left(
        |\gamma_j(\zeta)|^2
        \,\middle|\, R_{n,-j}(\zeta)
    \right)
    &\le
    \frac{C}{n^2}
    \tr\!\left[
        \Sigma_n^{1/2}R_{n,-j}(\zeta)\Sigma_n
        R_{n,-j}(\zeta)^*\Sigma_n^{1/2}
    \right]                                                        \\
    &\le
    \frac{C}{n^2}
    \|\Sigma_n\|_2^2
    \|R_{n,-j}(\zeta)\|_2^2\,p
    \le
    \frac{C}{n},
\end{align*}
where the last inequality follows by Assumptions~\ref{A1} and~\ref{A2} and \eqref{eq:general-resolvent-bound}. Hence,
\begin{align}\label{eq:gamma-bound}
    \sup_{\zeta\in K_1}\max_{1\le j\le n}
    \mathbb E\bigl(|\gamma_j(\zeta)|^2\bigr)
    \le
    \frac{C}{n}.
\end{align}
Similarly, almost surely, 
\begin{align*}
    \mathbb E\left(
        |\eta_j(\zeta)|^2
        \mid R_{n,-j}(\zeta)
    \right)
    \le
    \frac{C}{n^2}
    \operatorname{tr}\left[
        \Sigma_n^{1/2}R_{n,-j}(\zeta)^2\Sigma_n
        (R_{n,-j}(\zeta)^{2})^\ast\Sigma_n^{1/2}
    \right]                      \le
    \frac{C}{n}.
\end{align*}
Therefore,
\begin{align}\label{eq:eta-bound}
    \sup_{\zeta\in K_1}\max_{1\le j\le n}
    \mathbb E\bigl(|\eta_j(\zeta)|^2\bigr)
    \le
    \frac{C}{n}.
\end{align}
We can now bound \(U_n\). Since the summands in \(U_n\) are martingale
differences, it is straightforward to verify that
\[
    \mathbb E \bigl(|U_n(\zeta)|^2\bigr)
    =
    \sum_{j=1}^n
    \mathbb E \bigl(\left|
        (\E_{j}-\E_{j-1})\{b_j(\zeta)\eta_j(\zeta)\}
    \right|^2\bigr).
\]
Moreover, for any square-integrable random variable \(Y\), using $(a+b)^2\leq 2a^2+2b^2$, for $a,b\geq 0$ and Jensen's inequality, we obtain
\begin{align*}
    \mathbb E\bigl\{\bigl|(\E_{j}-\E_{j-1})\{Y\}\bigr|^2\bigr\}
    \leq C\E\{|Y|^2\}.
\end{align*}
Thus, combining the previous two displays together with \eqref{eq:b-bound} and \eqref{eq:eta-bound}, we have, for any $\zeta\in K_1$,
\[
\begin{aligned}
    \mathbb E\bigl(|U_n(\zeta)|^2\bigr)
    &\le
    C\sum_{j=1}^n
    \mathbb E\bigl(|b_j(\zeta)\eta_j(\zeta)|^2\bigr)                  \le
    C\sum_{j=1}^n
    \mathbb E\bigl(|\eta_j(\zeta)|^2\bigr)
    \le C.
\end{aligned}
\]
Consequently,
\begin{align}\label{eq:U-bound}
   \sup_{\zeta\in K_1}\mathbb E\bigl(|U_n(\zeta)|^2\bigr)
    \le C.
\end{align}

The term \(V_n\) is handled similarly, which gives
\begin{align*}
    \mathbb E\bigl(|V_n(\zeta)|^2\bigr)
    &\leq\sum_{j=1}^n
    \mathbb E\left[\bigl|(\E_j-\E_{j-1})\{
        b_j(\zeta)\beta_j(\zeta)q_j(\zeta)\gamma_j(\zeta)
    \}\bigr|^2\right]\\[2 mm]
    &\leq
    4\sum_{j=1}^n
    \mathbb E\left(\bigl|
        b_j(\zeta)\beta_j(\zeta)q_j(\zeta)\gamma_j(\zeta)
    \bigr|^2\right).
\end{align*}
Combining \eqref{eq:b-bound} and \eqref{eq:beta-q-bound}, we obtain
\[
    \sup_{\zeta\in K_1}\max_{1\leq j\leq n}|b_j(\zeta)\beta_j(\zeta)q_j(\zeta)|
    \le C.
\]
Therefore, by
\eqref{eq:gamma-bound},
\[
\begin{aligned}
    \mathbb E\bigl(|V_n(\zeta)|^2\bigr)
    &\le
    C
    \sum_{j=1}^n
    \mathbb E\bigl(|\gamma_j(\zeta)|^2\bigr)
    \le C.
\end{aligned}
\]
Hence,
\begin{align}\label{eq:V-bound}
    \sup_{\zeta\in K_1}\mathbb E\bigl(|V_n(\zeta)|^2\bigr)
    \le C.
\end{align}
Combining \eqref{eq:centered-UV-decomp}, \eqref{eq:U-bound}, and
\eqref{eq:V-bound}, we obtain
\[
    \sup_{\zeta\in K_1}
    \mathbb E\bigl(|
        n\{\hat v_n(\zeta)-\mathbb E(\hat v_n(\zeta))\}
    |^2\bigr)
    \le C.
\]

Fix an \(\omega\in\Omega\) and note that for each \(n\), the function
\[
        D_n(\zeta,\omega)
        \coloneqq
        n\{\hat v_n(\zeta,\omega)-\E(\hat v_n(\zeta))\}
\]
is analytic on \(\mathbb C\setminus\mathbb R_{\ge0}\). Indeed,
\(\zeta\mapsto \hat v_n(\zeta,\omega)\) is analytic on this domain, and
\(\zeta\mapsto \E(\hat v_n(\zeta))\) is analytic there as well. Therefore, we can use
Lemma~\ref{lem_sup_measurable} and conclude that,
\(\sup_{\zeta\in K}|D_n(\zeta)|\)
is \(\mathcal F\)-measurable for every \(n\ge1\).
Applying Lemma~\ref{lem:cauchy_reduction_random} to $D_n$,
with the compact neighborhood $K_1$ fixed at the beginning
of the proof, gives
\begin{align}\label{eq:var:square}
        \E\left[
        \sup_{\zeta\in K}
        \left|
        n\{\hat v_n(\zeta)-\E(\hat v_n(\zeta))\}
        \right|^2
        \right]
        &\le
         C
        \sup_{\zeta\in K_1}
        \E
        \bigl(|
        n\{\hat v_n(\zeta)-\E(\hat v_n(\zeta))\}
        |^2\bigr) \leq C.
\end{align}
Using Markov's inequality we conclude that 
\[
\sup_{\zeta\in K}\bigl|\{\hat{v}_n(\zeta)-\E(\hat{v}_n(\zeta))\}\bigr|=O_{\mathbb{P}}(n^{-1}).
\]
\end{proof}

\subsection{Lemma~\ref{stielt:bias} and its proof}

\begin{lemma}
\label{stielt:bias}
Suppose that Assumptions~\ref{A1}--\ref{A4} are satisfied.
Let $K\subset\mathbb C\setminus\mathbb R_{\geq0}$ be compact. Then
\[
        \sup_{\zeta\in K}
        \left|
        \mathbb E\bigl(\widehat v_n(\zeta)\bigr)
        -v_n(\zeta)
        \right|
        =
        O(n^{-1}).
\]
\end{lemma}

\begin{proof}[Proof of Lemma~\ref{stielt:bias}]

\textit{Proof outline for Lemma~\ref{stielt:bias}.}
The proof is divided into five steps.

\textbf{Step 1:}\label{step1:lemb2} First, we will show that it suffices to show that
\begin{align} \label{stielt:bias:goal}
        \sup_{\zeta\in\Gamma}
        \left|
        n\{\E(\hat v_n(\zeta))-v_n(\zeta)\}
        \right|
        \le C ,
\end{align}
where \(\Gamma\subseteq\mathbb C\setminus\mathbb R_{\ge0}\) denotes a finite union
of contours  surrounding
\(K\).

\textbf{Step 2:}\label{step2:lemb2} Second, we establish that there exists a constant $c>0$, such that
\begin{align}\label{eq:lower}
        \inf_{\zeta\in\Gamma}|v_n(\zeta)|\ge c,
        \qquad
        \inf_{\zeta\in\Gamma}|\E(\hat v_n(\zeta))|\ge c
\end{align}
for all sufficiently large \(n\). These bounds allow us to divide by
\(v_n(\zeta)\) and \(\E(\hat v_n(\zeta))\) uniformly on \(\Gamma\) for all sufficiently large $n$.

\textbf{Step 3:}\label{step3:lemb2} Third, we derive the perturbed equation
\[
        -\frac1{\E(\hat v_n(\zeta))}
        =
        \zeta
        -
        \frac1n
        \tr\!\left[
        \Sigma_n(I_p+\E(\hat v_n(\zeta))\Sigma_n)^{-1}
        \right]
        +
        \varepsilon_n(\zeta),
\]
where
\[
        \sup_{\zeta\in\Gamma}|n\varepsilon_n(\zeta)|\le C .
\]

\textbf{Step 4:}\label{step4:lemb2} Fourth, we obtain a deterministic function $L_n$, defined below, such that
\begin{align}
\label{stielt:bias:step4}
\{\E(\hat v_n(\zeta))-v_n(\zeta)\}L_n(\zeta)
=
\varepsilon_n(\zeta).
\end{align}

\textbf{Step 5:}\label{step5:lemb2} Finally, we verify that
\begin{align}
\label{eq:Ln-lower-bound}
        \inf_{\zeta\in\Gamma}|L_n(\zeta)|\ge c>0
\end{align}
for all sufficiently large $n$. Consequently,
\[
        \sup_{\zeta\in\Gamma}
        \left|
        n\{\E(\hat v_n(\zeta))-v_n(\zeta)\}
        \right|
        \le C.
\]
This implies the desired bound \eqref{stielt:bias:goal}.
Here and below the dependence on \(\zeta\) is often suppressed when no confusion is
possible.
\bigskip

\textbf{\hyperref[step1:lemb2]{Proof of Step~1}}. For some $0<\rho<\operatorname{dist}(K,\mathbb{R}_{\geq 0})$, 
we may choose finitely many points \(a_1,\ldots,a_m\in K\) such that
\[
        K\subseteq \bigcup_{\ell=1}^m B(a_\ell,\rho/4).
\]
For each \(\ell\), let
\[
        \Gamma_\ell\coloneqq\{w\in \mathbb{C}:|w-a_\ell|=\rho/2\},
        \qquad
        \Gamma\coloneqq\bigcup_{\ell=1}^m\Gamma_\ell .
\]
By the choice of $\rho$, all $\Gamma_l$ and their interiors are contained in $\mathbb{C}\setminus \R_{\geq 0}.$
Define the function
\[
        B_n(\zeta)\coloneqq
        n\{\E(\hat{v}_n(\zeta))-v_n(\zeta)\}
\]
and note that $B_n$ is holomorphic on $\mathbb{C}\setminus\mathbb{R}_{\geq 0}$. 
Let $\zeta \in K$, which implies \(\zeta\in B(a_\ell,\rho/4)\) for some $\ell$. 
Then, Cauchy's integral formula gives
\[
        B_n(\zeta)
        =
        \frac1{2\pi i}
        \int_{\Gamma_\ell}
        \frac{B_n(w)}{w-\zeta}\,dw .
\]
Since \(|w-\zeta|\ge \rho/4\) for \(w\in\Gamma_\ell\), we have
\[
        |B_n(\zeta)|
        \le
        \frac{2}{\pi\rho}
        \int_{\Gamma_\ell}|B_n(w)|\,|dw|
        \le
        C\sup_{w\in\Gamma}|B_n(w)|.
\]
Taking the supremum over \( \zeta\in K\) in the previous display, it is enough to prove \eqref{stielt:bias:goal}, which completes the first step of the proof.

\bigskip
\textbf{\hyperref[step2:lemb2]{Proof of Step~2}} We next establish \eqref{eq:lower}. In order to do so, we record two uniform
bounds that will be used repeatedly. Since $\Gamma$ is bounded away
from $\mathbb R_{\ge0}$, it is straightforward to verify that, for
every $n\in\mathbb N$,
\begin{align}\label{eq:v_nandEv_n}
        \sup_{\zeta\in\Gamma}|v_n(\zeta)|\le C,
        \qquad
        \sup_{\zeta\in\Gamma}|\E(\hat v_n(\zeta))|\le C .
\end{align}        
We start by verifying the lower bound for \(v_n\) in \eqref{eq:lower}. Suppose, to the contrary, that
the bound in \eqref{eq:lower} fails. Then, for every \(c>0\) and every
\(N\in\mathbb N\), there exists \(n\ge N\) such that
\(\inf_{\zeta\in\Gamma}|v_n(\zeta)|<c\). Taking \(c=1/m\) and \(N=m\), we can choose \(n_m\ge m\) and a sequence
\(\{\zeta_m\}_{m\ge1}\subseteq\Gamma\) such that
\[
        |v_{n_m}(\zeta_m)|<\frac1m
        \quad\text{for all }m\ge1 .
\]
Thus, along this subsequence,
\[
        |v_{n_m}(\zeta_m)|\to0 .
\]
By the implicit representation in \eqref{stielt}, applied to the
finite-sample pair \((\gamma_{n_m},H_{n_m})\), \(v_{n_m}(\zeta)\) satisfies
\eqref{stielt} for every
\(\zeta\in\mathbb C\setminus\mathbb R_{\ge0}\). %

Hence, after rearranging the implicit representation of
\(v_{n_m}(\zeta_m)\), we obtain
\begin{align}\label{eq:vn-identity-lower}
        1
        =
        -\zeta_m v_{n_m}(\zeta_m)
        +
        \gamma_{n_m}
        \int_0^\infty
        \frac{t v_{n_m}(\zeta_m)}
             {1+t v_{n_m}(\zeta_m)}
        \,dH_{n_m}(t).
\end{align}
Since \(\Gamma\) is compact, \(\sup_{\zeta\in\Gamma}|\zeta|<\infty\), and hence
\[
        |\zeta_m v_{n_m}(\zeta_m)|
        \le
        \sup_{\zeta\in\Gamma}|\zeta|\,
        |v_{n_m}(\zeta_m)|
        \longrightarrow 0,
        \qquad \text{as }m\to\infty .
\]
Consequently, $-\zeta_m v_{n_m}(\zeta_m)\longrightarrow 0$.
Moreover, by Assumption~\ref{A2}, $\operatorname{supp}(H_n)\subseteq[0,C]$ for all $n\in \mathbb{N}$.
Hence, for all sufficiently large \(m\),
$C |v_{n_m}(\zeta_m)|\le \frac12.$
By the reverse triangle inequality,
\[
        |1+t v_{n_m}(\zeta_m)|
        \ge
        \left|1-t|v_{n_m}(\zeta_m)|\right|,\quad t\in [0,C].
\]
Since \(t|v_{n_m}(\zeta_m)|\le C|v_{n_m}(\zeta_m)|\le 1/2\) for all sufficiently large \(m\), the
quantity \(1-t|v_{n_m}(\zeta_m)|\) is nonnegative. Hence,
\[
        \left|1-t|v_{n_m}(\zeta_m)|\right|
        =
        1-t|v_{n_m}(\zeta_m)|
        \ge \frac12 .
\]
Consequently, for all sufficiently large \(m\) and all \(t\in[0,C]\),
\[
        \left|
        \frac{t v_{n_m}(\zeta_m)}
             {1+t v_{n_m}(\zeta_m)}
        \right|
        \le
        2C |v_{n_m}(\zeta_m)|.
\]
Since \(\gamma_{n_m}=O(1)\) by Assumption~\ref{A1}, it follows that
\[
        \lim_{m\to\infty}
        \gamma_{n_m}
        \int_0^\infty
        \frac{t v_{n_m}(\zeta_m)}
             {1+t v_{n_m}(\zeta_m)}
        \,dH_{n_m}(t)
        =0 .
\]
Taking limits in \eqref{eq:vn-identity-lower} gives a contradiction.
Therefore, the first lower bound in \eqref{eq:lower} holds.

It remains to prove the second lower bound in \eqref{eq:lower}. Recall that
\[
        \hat v_n(\zeta)
        =
        \frac1n\tr\{(\hat{\underline\Sigma}_n-\zeta I_n)^{-1}\}=\int_0^\infty \frac{1}{x-\zeta}dF_{\hat{\underline{\Sigma}}_n}(x).
\]
By \citet[Theorem~2]{Latala2005} and Assumption~\ref{A3},
\[
    \mathbb E(\|Z_n\|_2)
    \leq
    C\left(\sqrt n+\sqrt p+(np)^{1/4}\right)
    \leq C\sqrt n,
\]
where the last inequality follows from
Assumption~\ref{A1}.
Consequently, Markov's inequality and Assumption~\ref{A2}
give
\[
    \|\underline{\widehat\Sigma}_n\|_2
    \leq
    \|\Sigma_n\|_2\frac{\|Z_n\|_2^2}{n}
    =
    O_{\mathbb P}(1).
\]

Hence, there exists an \(M<\infty\) such that, for all sufficiently large \(n\),
\begin{align}\label{eq:lmax-tight-companion}
        \mathbb P\{\|\hat{\underline\Sigma}_n\|_2\le M\}
        \ge \frac12 .
\end{align}
Let \(\zeta\in\Gamma\) and set $M_\Gamma\coloneqq\sup_{\xi\in\Gamma}|\xi|<\infty$.
We distinguish two cases.
First suppose that
\begin{align} \label{b1}
        |\Im(\zeta)|
        \ge
        \frac12\operatorname{dist}(\Gamma,\mathbb R_{\ge0}).
\end{align}
Denote the eigenvalues of $\hat{\underline{\Sigma}}_n$ by $\lambda_{1,n},\ldots,\lambda_{n,n}$ and define $D_n\coloneqq\{\|\hat{\underline\Sigma}_n\|_2\leq M\}$. 
Then, for each \(k=1,\ldots,n\), we have
\[
        \Im\left(
        \frac1{\lambda_{k,n}-\zeta}
        \right)
        =
        \frac{\Im(\zeta)}
             {|\lambda_{k,n}-\zeta|^2}.
\]
This implies
\begin{align*}
        |\E(\hat v_n(\zeta))|
        &\ge
        \left|
        \Im\{\E(\hat v_n(\zeta))\}
        \right|                   
        =
        |\Im(\zeta)|
        \frac1n\sum_{k=1}^n
        \mathbb E\left[
        \frac1{|\lambda_{k,n}-\zeta|^2}
        \right]                      \ge
        |\Im(\zeta)|
        \frac1n\sum_{k=1}^n
        \E\left[
        \frac{1}
             {|\lambda_{k,n}-\zeta|^2}\mathbf 1_{D_n}
        \right].
\end{align*}
On \(D_n\), we have \(\lambda_{k,n}\in[0,M]\) for every
\(k=1,\ldots,n\). Hence, for \(\zeta\in\Gamma\),
\begin{align} \label{b3}
        |\lambda_{k,n}-\zeta|
        \le
        |\lambda_{k,n}|+|\zeta|
        \le
        M+\sup_{\xi\in\Gamma}|\xi|
        =
        M+M_\Gamma .
\end{align}
Thus, on \(D_n\),
\[
        \frac1{|\lambda_{k,n}-\zeta|^2}
        \ge
        \frac1{(M+M_\Gamma)^2}.
\]
Consequently,
\[
\begin{aligned}
        |\mathbb E(\hat v_n(\zeta))|\geq
        \frac{|\Im(\zeta)|}{(M+M_\Gamma)^2}
        \mathbb P(D_n).
\end{aligned}
\]
Since we assumed \eqref{b1}
and since \(\mathbb P(D_n)\ge1/2\) for all sufficiently large \(n\), it follows
that
\[
        |\mathbb E(\hat v_n(\zeta))|
        \ge
        \frac{
        \operatorname{dist}(\Gamma,\mathbb R_{\ge0})
        }{
        4(M+M_\Gamma)^2
        }.
\]

Now suppose that
\begin{align}\label{bias:case2}
        |\Im(\zeta)|
        <
        \frac12\operatorname{dist}(\Gamma,\mathbb R_{\ge0}).
\end{align}
If \(\Re(\zeta)\ge0\), then 
\[
        \operatorname{dist}(\zeta,\mathbb R_{\ge0})
        =
        |\Im(\zeta)|\geq  \operatorname{dist}(\Gamma,\mathbb R_{\ge0})>\frac{1}{2}\operatorname{dist}(\Gamma,\mathbb R_{\ge0}),
\]
which contradicts \eqref{bias:case2}. Therefore \(\Re(\zeta)<0\). %

Write \(\zeta=a+ib\) for an arbitrary $\zeta\in \Gamma$. Since \(a<0\), we have
\[
        \operatorname{dist}(\zeta,\mathbb R_{\ge0})^2= \inf_{x\ge0}\{(x-a)^2\}+b^2
        =a^2+b^2=|\zeta|^2.
\]
Using
$
        \operatorname{dist}(\zeta,\mathbb R_{\ge0})
        \ge
        \operatorname{dist}(\Gamma,\mathbb R_{\ge0})$, \eqref{bias:case2} and the previous display, we obtain
\[
        |a|^2
        =
        |\zeta|^2-|b|^2
        \ge
        \operatorname{dist}(\Gamma,\mathbb R_{\ge0})^2
        -
        \frac14\operatorname{dist}(\Gamma,\mathbb R_{\ge0})^2
        =
        \frac34\operatorname{dist}(\Gamma,\mathbb R_{\ge0})^2.
\]
In particular,
\[
        |a|
        \ge
        \frac12\operatorname{dist}(\Gamma,\mathbb R_{\ge0}).
\]

Furthermore, for \(x\ge0\), note that
\[
        \Re\left(\frac1{x-\zeta}\right)
        =
        \frac{x-a}{(x-a)^2+b^2}.
\]
Since \(a<0\), we have \(x-a
>0\), and hence
$\operatorname{Re}\left(\frac1{x-\zeta}\right)\ge0$.
Therefore,
\begin{align}
        |\E(\hat v_n(\zeta))|
        &\ge
        \Re\{\E(\hat v_n(\zeta))\}                  =
        \frac1n\sum_{k=1}^n
        \mathbb E\left[
        \frac{\lambda_{k,n}-a}
             {(\lambda_{k,n}-a)^2+b^2}
        \right]                      
        \ge
        \frac1n\sum_{k=1}^n
        \mathbb E\left[
        \mathbf 1_{D_n}
        \frac{\lambda_{k,n}-a}
             {(\lambda_{k,n}-a)^2+b^2}
        \right]. \label{b2}
\end{align}
Since on \(D_n\), we have \(\lambda_{k,n}\in[0,M]\) for each $k\in \{1,\dots,n\}$, we conclude that
\[
        \lambda_{k,n}-a
        \ge
        -a
        =
        |a|
        \ge
        \frac12\operatorname{dist}(\Gamma,\mathbb R_{\ge0}).
\]
Thus, on \(D_n\), using also \eqref{b3},
\[
        \frac{\lambda_{k,n}-a}
             {(\lambda_{k,n}-a)^2+b^2}
        \ge
        \frac{
        \operatorname{dist}(\Gamma,\mathbb R_{\ge0})
        }{
        2(M+M_\Gamma)^2
        }.
\]
Combining this with \eqref{b2},
\[
\begin{aligned}
        |\mathbb E(\hat v_n(\zeta))|
        &\ge
        \frac{
        \operatorname{dist}(\Gamma,\mathbb R_{\ge0})
        }{
        4(M+M_\Gamma)^2
        }.
\end{aligned}
\]
Combining the two cases yields
\[
        \inf_{\zeta\in\Gamma}
        |\mathbb E(\hat v_n(\zeta))|
        \ge
        \frac{
        \operatorname{dist}(\Gamma,\mathbb R_{\ge0})
        }{
        4\{M+ M_\Gamma\}^2
        }
\]
for all sufficiently large \(n\). This proves the second lower bound in
\eqref{eq:lower} and completes the second step.

\bigskip
\textbf{\hyperref[step3:lemb2]{Proof of Step~3}}
In the third step of the proof, we establish that, for all sufficiently large $n$,
\begin{align}\label{eq:un-perturbed-equation}
-\frac1{\E(\hat{v}_n(\zeta))}
        =
        \zeta
        -
        \frac1n
        \tr\!\left[
        \Sigma_n(I_p+\E(\hat{v}_n(\zeta))\Sigma_n)^{-1}
        \right]
        +
        \varepsilon_n(\zeta),
\end{align}
where we will show that
\begin{align}\label{eq:epsilon-rate}
\sup_{\zeta\in\Gamma}|n\varepsilon_n(\zeta)|
        \le C_\Gamma .
\end{align}
Recall the definition of 
$\beta_j$ in \eqref{eq:beta_jandq_j} and note that 
\begin{align}\label{eq:beta_j_stielt}
        -\zeta \hat v_n(\zeta)
        =
        \frac1n\sum_{j=1}^n\beta_j(\zeta).
\end{align}        
Therefore,
\begin{align*}%
        -\zeta \E(\hat{v}_n(\zeta))
        =
        \frac1n\sum_{j=1}^n\mathbb E(\beta_j(\zeta)).
\end{align*}

Note that
\[
        1-\beta_j(\zeta)
        =
        \frac{1}{n}\beta_j(\zeta) x_j^\top R_{n,-j}(\zeta)x_j .
\]
Using the previous two displays, we obtain
\begin{align}\label{eq:basic-difference}
        1+\zeta \E(\hat{v}_n(\zeta))
        =
        1-\frac1n\sum_{j=1}^n\E(\beta_j(\zeta))=
        \frac1n\sum_{j=1}^n
        \mathbb E\!\left[
        \frac{\beta_j(\zeta)}{n}x_j^\top R_{n,-j}(\zeta)x_j
        \right].
\end{align}
For brevity, we write
\[
       A_n(\zeta)
        \coloneqq
        (I_p+\E(\hat{v}_n(\zeta))\Sigma_n)^{-1}.
\]
Since $\mathbb E(\widehat v_n)$ is the Stieltjes transform
of a probability measure on $[0,\infty)$,
\eqref{eq:deterministic-resolvent-bound} gives
\begin{align}\label{eq:bias-An-uniform-bound}
\sup_{\zeta\in\Gamma}\|A_n(\zeta)\|_2\leq C_\Gamma.
\end{align}
Define, for $\zeta\in \Gamma$,
\[
S_n(\zeta)\coloneqq \frac1n\sum_{j=1}^n
        \mathbb E\!\left[
        \frac{\beta_j}{n} x_j^\top R_{n,-j}(\zeta) x_j
        \right]
        -\E(\hat{v}_n(\zeta))\frac1n\tr\{\Sigma_n A_n(\zeta)\}.
\]
We shall prove that
\begin{align}\label{eq:S-rate}
        \sup_{\zeta\in\Gamma}\left|S_n(\zeta)\right|
        \le \frac{C_\Gamma}{n}.
\end{align}

Since $\E(\hat v_n(\zeta))$, is
bounded away from zero uniformly for \(\zeta\in\Gamma\) by \eqref{eq:lower}, we may divide \eqref{eq:basic-difference} by this term for sufficiently large $n$. By setting \[
-\E(\hat{v}_n(\zeta))^{-1}S_n(\zeta)\eqqcolon\varepsilon_n(\zeta),
\]
the identity
\eqref{eq:un-perturbed-equation} follows together with the bound
\eqref{eq:epsilon-rate}. 
Using \eqref{eq:a_jandb_j}, \eqref{eq:gamma_def}, \eqref{eq:b_j_representation}, we have

\[
        \beta_j^{-1}
        =
        1+a_j+\gamma_j
        =
        b_j^{-1}+\gamma_j,
\]
and therefore
\begin{align}\label{eq:beta-expansion}
        \beta_j
        =
        b_j
        -
        b_j^2\gamma_j
        +
        b_j^2\beta_j\gamma_j^2 .
\end{align}
Recalling the arguments from the proof of Lemma~\ref{stielt:var} we have that
\begin{align}\label{eq:bias:unif_bounds}
\begin{split}
        \sup_{\zeta\in\Gamma}\max_{1\leq j\leq n}
        |b_j(\zeta)|&\le C_\Gamma,\\
       \sup_{\zeta\in\Gamma}\max_{1\leq j\leq n}
        |\beta_j(\zeta)|&\le C_\Gamma,\qquad\text{and}\\
          \sup_{\zeta\in\Gamma}\max_{1\leq j\leq n}
        \mathbb E\bigl(|\gamma_j(\zeta)|^2\bigr)
        &\le \frac{C_\Gamma}{n}.
        \end{split}
\end{align}
We shall also need the following centered quadratic form. Put
\[
        \delta_j(\zeta)
        \coloneqq
        \frac{1}{n}x_j^\top R_{n,-j}(\zeta) A_n(\zeta) x_j
        -
        \frac1n
        \tr\{\Sigma_nR_{n,-j}(\zeta)A_n(\zeta)\}.
\]
Since \(R_{n,-j}\) and \(A_n\) are independent of \(x_j\), using Lemma \ref{lem_b26_bai_silverstein} gives
\begin{align}\label{eq:delta-L2}
        \sup_{\zeta\in\Gamma}\max_{1\leq j\leq n}
        \mathbb E\bigl(|\delta_j(\zeta)|^2\bigr)
        \le \frac{C_\Gamma}{n}.
\end{align}
Moreover, using Cauchy-Schwarz's inequality, \eqref{eq:bias:unif_bounds} and \eqref{eq:delta-L2}, we obtain
\begin{align}\label{eq:gamma-delta-cov}
\sup_{\zeta\in\Gamma}\max_{1\leq j\leq n}
        \mathbb E\bigl(\left|\gamma_j(\zeta)\delta_j(\zeta)\right|\bigr)
        \le \frac{C_\Gamma}{n}.
\end{align}
We now turn to the decomposition of $S_n(\zeta)$. The argument
uses the following elementary observations. Since
\[
        I_p=(I_p+\E(\hat{v}_n(\zeta))\Sigma_n)A_n
        =
        A_n+\E(\hat{v}_n(\zeta))\Sigma_n A_n,
\]
we have
\[
        \frac{1}{n}x_j^\top R_{n,-j}x_j
        =
        \frac{1}{n}x_j^\top R_{n,-j}A_n x_j
        +
        \frac{1}{n}\E(\hat{v}_n(\zeta)) x_j^\top R_{n,-j}\Sigma_n A_n x_j.
\]
Using the Sherman--Morrison identity, we get
\[
        R_nx_j=\left(I_p - \frac{\beta_j}{n}x_j^\top R_{n,-j} x_j I_p\right)R_{n,-j}x_j=\beta_j R_{n,-j}x_j,
\]
and using \(R_n\hat\Sigma_n=I_p+\zeta R_n\), we obtain
\begin{align*}
        \frac1n\sum_{j=1}^n
        \frac{\beta_j}{n} x_j^\top R_{n,-j}\Sigma_n A_n x_j
        &=
        \frac1n
        \tr\!\left[
        R_n\Sigma_nA_n\sum_{j=1}^n \frac{x_jx_j^\top}{n}
        \right]                                                   \\
        &=
        \frac1n
        \tr\!\left[
        R_n\Sigma_nA_n\hat\Sigma_n
        \right]                                                   \\
        &=
        \frac1n
        \tr(\Sigma_nA_n)
        +
        \zeta
        \frac1n
        \tr(R_n\Sigma_nA_n).
\end{align*}
Putting these observations together, yields
\begin{align*}
       S_n(\zeta)
        &=
        \frac1n\sum_{j=1}^n
        \frac{1}{n} \E\{\beta_jx_j^\top R_{n,-j} x_j\}
        -
        \E(\hat{v}_n(\zeta))
        \frac1n\tr\{\Sigma_n A_n\}\\[2 mm]
        &=\frac{1}{n}\sum_{j=1}^n \E\left[\frac{\beta_j}{n}x_j^\top R_{n,-j}A_n x_j + \frac{\beta_j}{n}\E(\hat{v}_n(\zeta)) x_j^\top R_{n,-j}\Sigma_n A_n x_j\right] -
        \E(\hat{v}_n(\zeta))
        \frac1n\tr\{\Sigma_n A_n\}\\[2 mm]
        &=\frac{1}{n}\sum_{j=1}^n \E\left[\frac{\beta_j}{n}x_j^\top R_{n,-j}A_n x_j\right]+\frac{\zeta \E(\hat{v}_n(\zeta))}{n}\E\bigl[\tr\{R_n \Sigma_n A_n\}\bigr]
\end{align*}
Next write
\[
        \frac{1}{n}x_j^\top R_{n,-j}A_n x_j
        =
        \frac1n\tr(\Sigma_nR_{n,-j}A_n)+\delta_j.
\]
Thus,
\begin{align}\label{eq:Sn-decomp-pre-cancel}
         S_n(\zeta)
        =
        \frac1n\sum_{j=1}^n
        \mathbb E\left[
        \beta_j
        \frac1n\tr(\Sigma_nR_{n,-j}A_n)
        \right]
        +
        \frac1n\sum_{j=1}^n\mathbb E(\beta_j\delta_j) 
        +
        \frac{\zeta\,\E(\hat{v}_n(\zeta))}{n}
        \mathbb E\{\tr(R_n\Sigma_nA_n)\}.
\end{align}
We next show that the first and third terms in
\eqref{eq:Sn-decomp-pre-cancel} cancel up to an error of order \(n^{-1}\).
By the Sherman--Morrison formula,
the matrix \(R_n-R_{n,-j}\) has rank at most one. Since the matrices
\(\Sigma_n\), \(A_n\), \(R_n\), and \(R_{n,-j}\) have uniformly bounded spectral
norms on \(\Gamma\), uniformly in \(n\), we have
\begin{align*}
\left|
        \frac1n\tr\{\Sigma_n(R_n-R_{n,-j})A_n\}
\right|
&\le
\frac1n
\|\Sigma_n(R_n-R_{n,-j})A_n\|_1                                      \\
&\le
\frac1n
\|\Sigma_n\|_2\|A_n\|_2\|R_n-R_{n,-j}\|_1                              \\
&\le
\frac1n
\|\Sigma_n\|_2\|A_n\|_2
\operatorname{rank}(R_n-R_{n,-j})\|R_n-R_{n,-j}\|_2                        \\
&\le
\frac1n
\|\Sigma_n\|_2\|A_n\|_2
\{\|R_n\|_2+\|R_{n,-j}\|_2\}                                           \\
&\le
\frac{C_\Gamma}{n}.
\end{align*}
Hence, for each \(j\in\{1,\dots,n\}\) and \(\zeta\in\Gamma\), we have
\[
       \frac1n\tr(\Sigma_nR_{n,-j}A_n)
        =
        \frac1n\tr(\Sigma_nR_nA_n)+O(n^{-1}).
\]
Using this bound, the uniform bound for \(\beta_j\) in \eqref{eq:bias:unif_bounds}, and the identity in\eqref{eq:beta_j_stielt}
we conclude that
\begin{align} \label{n1} 
\begin{split}\frac1n\sum_{j=1}^n
        \mathbb E\left[
        \beta_j
        \frac1n\tr(\Sigma_nR_{n,-j}A_n)
        \right]
        &=
        \mathbb E\left[
        \frac1n\tr(\Sigma_nR_nA_n)
        \frac1n\sum_{j=1}^n\beta_j
        \right]
        +
        O(n^{-1})                                         \\
        &=
        -\zeta
        \mathbb E\left[
        \hat v_n(\zeta)
        \frac1n\tr(\Sigma_nR_nA_n)
        \right]
        +
        O(n^{-1}).
        \end{split}
\end{align}
Applying \eqref{eq:var:square} with $\Gamma$ in place of $K$
and then Cauchy--Schwarz gives
\[
    \sup_{\zeta\in\Gamma}
    \mathbb E|\widehat v_n(\zeta)-\mathbb E\widehat v_n(\zeta)|
    \leq
    \left\{\mathbb E\sup_{\zeta\in\Gamma}
    |\widehat v_n(\zeta)-\mathbb E\widehat v_n(\zeta)|^2\right\}^{1/2}
    \leq\frac{C_\Gamma}{n}.
\]
Furthermore, for each \(n\in\mathbb N\) and \(\zeta\in\Gamma\),
\[
        \frac1n|\tr(\Sigma_nR_nA_n)|
        \le
        \frac1n\|\Sigma_nR_nA_n\|_1
        \le
        \gamma_n\|\Sigma_n\|_2\|R_n\|_2\|A_n\|_2
        \le C_\Gamma.
\]
Hence,
\begin{align} \label{n2}
           \mathbb E\left[
        \hat v_n(\zeta)
        \frac1n\tr(\Sigma_nR_nA_n)
        \right]
        =
        \E(\hat{v}_n(\zeta))
        \E\left[
        \frac1n\tr(\Sigma_nR_nA_n)
        \right]
        +
        O(n^{-1}),
\end{align}
uniformly for \(\zeta\in\Gamma\).
Thus, by combining \eqref{eq:Sn-decomp-pre-cancel}, \eqref{n1} and \eqref{n2}
\[
    S_n(\zeta)
        =
        \frac1n\sum_{j=1}^n\mathbb E(\beta_j\delta_j)
        +
        O(n^{-1}).
\]
It remains to control the average of \(\mathbb E(\beta_j\delta_j)\). Since $b_j(\zeta)$ is $\mathcal F_{n,-j}$-measurable
and $\mathbb E[\delta_j(\zeta)\mid\mathcal F_{n,-j}]=0$,
we get
\[
        \mathbb E\{b_j(\zeta)\delta_j(\zeta)\}=0.
\]
Consequently, using the decomposition of $\beta_j$ in 
\eqref{eq:beta-expansion},
\[
        \mathbb E\{\beta_j(\zeta)\delta_j(\zeta)\}
        =
        -\mathbb E\{b_j(\zeta)^2\gamma_j(\zeta)\delta_j(\zeta)\}
        +
        \mathbb E\{b_j(\zeta)^2\beta_j(\zeta)\gamma_j(\zeta)^2
        \delta_j(\zeta)\}.
\]
The first term in the previous display is \(O(n^{-1})\) by \eqref{eq:gamma-delta-cov} and the uniform boundedness of \(b_j\) from \eqref{eq:bias:unif_bounds}. For the second term, we use the identity  $\beta_j^{-1}=b_j^{-1}+\gamma_j$. 
Multiplying both sides by \(b_j\beta_j\) gives $b_j=\beta_j+b_j\beta_j\gamma_j$,
and hence
\[
        b_j\beta_j\gamma_j=b_j-\beta_j .
\]
Since both \(b_j\) and
\(\beta_j\) are uniformly bounded, we conclude with the previous display that \(b_j\beta_j\gamma_j\) is uniformly
bounded. Hence,
\[
        |b_j^2\beta_j\gamma_j^2\delta_j|
        =
        |b_j\beta_j\gamma_j|\,|b_j|\,|\gamma_j\delta_j|
        \le
        C_\Gamma|\gamma_j\delta_j|.
\]
Therefore, by the Cauchy--Schwarz inequality, the bounds in \eqref{eq:bias:unif_bounds} and \eqref{eq:delta-L2}, 
we obtain
\begin{align*}
        \left|
        \E(b_j^2\beta_j\gamma_j^2\delta_j)
        \right|
        \le
        C\,\E(|\gamma_j\delta_j|)\le
        C\{\E(|\gamma_j|^2)\}^{1/2}
          \{\E(|\delta_j|^2)\}^{1/2}\le
        \frac{C}{n}.
\end{align*}
Consequently,
\[
        \sup_{\zeta\in\Gamma}\max_{1\le j\le n}
        |\mathbb E\{\beta_j(\zeta)\delta_j(\zeta)\}|
        \le \frac{C}{n}.
\]
Therefore,
\[
        \sup_{\zeta\in\Gamma}|S_n(\zeta)|
        \le
        \frac1n\sum_{j=1}^n\frac{C}{n}
        +
        \frac{C}{n}
        \le
        \frac{C}{n}.
\]
This proves \eqref{eq:S-rate} and concludes step 3.
\bigskip

\textbf{\hyperref[step4:lemb2]{Proof of Step~4}}
Next, we show \eqref{stielt:bias:step4}. Note that, using the implicit
representation of $v_n$, it follows that, for
$\zeta\in\mathbb C\setminus\mathbb R_{\geq0}$,
\[
-\frac{1}{v_n(\zeta)}
=
\zeta
-
\frac{1}{n}
\tr\!\left[
\Sigma_n(I_p+v_n(\zeta)\Sigma_n)^{-1}
\right]
=
\zeta
+
\frac{\zeta}{n}
\tr\!\left[
\Sigma_nR_n^{\operatorname{det}}(\zeta)
\right].
\]
Subtracting the preceding identity from \eqref{eq:un-perturbed-equation}, we obtain, for all sufficiently large $n$,
\begin{align}\label{eq:resolvent_dec1}
\begin{split}
\frac{\E(\hat v_n(\zeta))-v_n(\zeta)}
     {\E(\hat v_n(\zeta))v_n(\zeta)}
&=
-\frac1{\E(\hat v_n(\zeta))}
+
\frac1{v_n(\zeta)}
\\
&=
-\frac1n
\tr\!\left[
\Sigma_n
\left\{
A_n(\zeta)
+
\zeta R_n^{\operatorname{det}}(\zeta)
\right\}
\right]
+
\varepsilon_n(\zeta).
\end{split}
\end{align}
Furthermore, by the resolvent identity,
\begin{align}
\label{eq:resol-ident}
        A_n(\zeta)
        +
        \zeta R_n^{\operatorname{det}}(\zeta)
        &=
        \zeta
        \{\E(\hat v_n(\zeta))-v_n(\zeta)\}
        A_n(\zeta)\Sigma_n
        R_n^{\operatorname{det}}(\zeta).
\end{align}
Substituting \eqref{eq:resol-ident} into
\eqref{eq:resolvent_dec1}, we obtain \eqref{stielt:bias:step4} with
\[
\begin{split}
        L_n(\zeta)
        \coloneqq
        \frac1{\E(\hat v_n(\zeta))v_n(\zeta)}
        +
        \frac{\zeta}{n}
        \tr\!\left[
        \Sigma_n A_n(\zeta)\Sigma_n
        R_n^{\operatorname{det}}(\zeta)
        \right].
\end{split}
\]

\bigskip
\textbf{\hyperref[step5:lemb2]{Proof of Step~5}}
It remains to prove the fifth step and establish \eqref{eq:Ln-lower-bound}. We derive this bound, using
\begin{align}\label{eq:weak-bias-on-Gamma}
        \sup_{\zeta\in\Gamma}
        |\E(\hat v_n(\zeta))-v_n(\zeta)|
        \longrightarrow 0 .
\end{align}
Suppose, to the contrary, that \eqref{eq:weak-bias-on-Gamma} fails. Then
there exists an \(\eta>0\) such that, for every \(N\in\mathbb N\), there is
\(n\ge N\) with
\[
        \sup_{\zeta\in\Gamma}
        |\E(\hat v_n(\zeta))-v_n(\zeta)|
        \ge \eta .
\]
Choosing the indices recursively, we may find a subsequence
\((n_k)_{k\ge1}\) such that \(n_k\to\infty\) and
\[
        \sup_{\zeta\in\Gamma}
        |\E(\hat v_{n_k}(\zeta))-v_{n_k}(\zeta)|
        \ge \eta
        \qquad\text{for all }k\ge1 .
\]
Since \(\Gamma\) is compact and, for each \(k\), the map
\[
        \zeta\mapsto
        \E(\hat v_{n_k}(\zeta))-v_{n_k}(\zeta)
\]
is continuous on \(\Gamma\), the supremum is attained for each \(k\). Hence, for each
\(k\ge1\), there exists \(\zeta_k\in\Gamma\) such that
\[
        |\E(\hat v_{n_k}(\zeta_k))-v_{n_k}(\zeta_k)|
        =
        \sup_{\zeta\in\Gamma}
        |\E(\hat v_{n_k}(\zeta))-v_{n_k}(\zeta)|
        \ge \eta .
\]

Since $\{\zeta_k\}_{k\ge1}\subseteq\Gamma$ and $\Gamma$ is compact, there
exists a subsequence $\{k_\ell\}_{\ell\ge1}\subseteq\mathbb N$ with
$k_\ell\to\infty$ and a point $\zeta_\ast\in\Gamma$ such that $\zeta_{k_\ell}\to\zeta_\ast$. 
We now pass from the subsequence $\{n_k\}_{k\ge1}$ to the further subsequence
$\{n_{k_\ell}\}_{\ell\ge1}$. To simplify notation, we relabel this further
subsequence as $\{n_k\}_{k\ge1}$ and the corresponding points
$\{\zeta_{k_\ell}\}_{\ell\ge1}$ as $\{\zeta_k\}_{k\ge1}$. Thus, without loss
of generality, we may assume that
\[
        \zeta_k\to\zeta_\ast\in\Gamma .
\]

Recall that the quantities $\E(\hat v_n(\zeta))$ and $v_n(\zeta)$ are
uniformly bounded on $\Gamma$. Therefore, after passing to a further
subsequence, still denoted by $\{n_k\}_{k\ge1}$, if necessary, we may assume
that
\[
        \E(\hat v_{n_k}(\zeta_k))\to w_1,
        \qquad
        v_{n_k}(\zeta_k)\to w_2 .
\]
By \eqref{eq:lower}, with \(c>0\), we have, for all sufficiently large \(k\),
\[
        |\E(\hat v_{n_k}(\zeta_k))|\ge c,
        \qquad
        |v_{n_k}(\zeta_k)|\ge c .
\]
Letting \(k\to\infty\), we obtain
\[
        |w_1|\ge c,
        \qquad
        |w_2|\ge c .
\]
In particular, \(w_1\neq0\) and \(w_2\neq0\).

We now pass to the limit in the perturbed equation. From
\eqref{eq:un-perturbed-equation}, evaluated at $\zeta=\zeta_k$ and
$n=n_k$, we have
\[
        -\frac1{\E(\hat v_{n_k}(\zeta_k))}
        =
        \zeta_k
        -
        \frac1{n_k}
        \tr\!\left[
        \Sigma_{n_k}
        A_{n_k}
        \right]
        +
        \varepsilon_{n_k}(\zeta_k),
\]
which can be equivalently written as
\begin{align}\label{subseq:implicit}
        -\frac1{\E(\hat v_{n_k}(\zeta_k))}
        =
        \zeta_k
        -
        \gamma_{n_k}
        \int_0^\infty
        \frac{t}{1+t\E(\hat v_{n_k}(\zeta_k))}
        \,dH_{n_k}(t)
        +
        \varepsilon_{n_k}(\zeta_k).
\end{align}
Since $\sup_{\zeta\in\Gamma}|\varepsilon_n(\zeta)|\to0$, we have
$\varepsilon_{n_k}(\zeta_k)\to0$. Moreover, by Assumptions~\ref{A1}
and~\ref{A4}, 
\[\gamma_{n_k}\to\gamma,\quad \text{and}\quad H_{n_k}\Rightarrow H.
\]
Next, we establish that 
\[
        \int_0^\infty
        \frac{t}{1+t\E(\hat v_{n_k}(\zeta_k))}
        \,dH_{n_k}(t)
        \longrightarrow
        \int_0^\infty
        \frac{t}{1+tw_1}\,dH(t).
\]
To this end, we first record a few facts needed for the passage to the limit.
By Assumption~\ref{A2}, the measures $H_n$ are supported on a common
compact interval, say $[0,C]$, for every $n\in \mathbb{N}$.

We now show that $1+tw_1\neq0$ for every $t\in[0,C]$. It is enough to show that $w_1\notin\mathbb R_{<0}$. Indeed, if
$1+tw_1=0$ for some $t\in[0,C]$, then necessarily $t>0$ and
$w_1=-1/t\in\mathbb R_{<0}$. We use the estimates obtained in the proof of the second lower bound in
\eqref{eq:lower}. They imply that, for all sufficiently large $n$ and all
$\zeta\in\Gamma$,
\[
        \Im(\zeta)\ge
        \frac12\operatorname{dist}(\Gamma,\mathbb R_{\ge0})
        \quad\Longrightarrow\quad
        \Im\{\E(\hat v_n(\zeta))\}
        \ge
        \frac{\operatorname{dist}(\Gamma,\mathbb R_{\ge0})}
        {4\{M+\sup_{\xi\in\Gamma}|\xi|\}^2}.
\]

The case $\Im(\zeta)\le
        -\frac12\operatorname{dist}(\Gamma,\mathbb R_{\ge0})$
is identical, with the inequalities reversed. Indeed, for $x\geq 0$,
\[
        \Im\left(\frac1{x-\zeta}\right)
        =
        \frac{\Im(\zeta)}{|x-\zeta|^2},
\]
so \(\Im\{\E(\hat v_n(\zeta))\}\) has the same sign as \(\Im(\zeta)\). The same estimate on $D_n$ therefore yields
\[
        \Im\{\E(\hat v_n(\zeta))\}
        \le
        -\frac{\operatorname{dist}(\Gamma,\mathbb R_{\ge0})}
        {4\{M+\sup_{\xi\in\Gamma}|\xi|\}^2}.
\]

Moreover, in the remaining case
\[
        |\Im(\zeta)|
        <
        \frac12\operatorname{dist}(\Gamma,\mathbb R_{\ge0}),
\]
the same proof gives
\[
        \Re\{\E(\hat v_n(\zeta))\}
        \ge
        \frac{\operatorname{dist}(\Gamma,\mathbb R_{\ge0})}
        {4\{M+\sup_{\xi\in\Gamma}|\xi|\}^2}.
\]
Applying these alternatives to \(\zeta=\zeta_k\in\Gamma\) and \(n=n_k\), and
passing to a further subsequence if necessary, at least one of the three
alternatives holds along the whole subsequence. Letting \(k\to\infty\) and
using \(\E(\hat v_{n_k}(\zeta_k))\to w_1\), we obtain one of the following
three conclusions:
\[
        \Im(w_1)>0,
        \qquad
        \Im(w_1)<0,
        \qquad
        \Re(w_1)>0.
\]
In each case, \(w_1\notin\mathbb R_{<0}\). Thus \(1+tw_1\neq0\) for every
\(t\in[0,C]\). Since the map $t\mapsto |1+tw_1|$ is continuous on the compact interval
$[0,C]$, it attains a positive minimum. Therefore,
\[
        \delta\coloneqq \inf_{t\in[0,C]}|1+tw_1|>0 .
\]

Since $\E(\hat v_{n_k}(\zeta_k))\to w_1$, we have
\[
\sup_{t\in[0,C]}
t\left|\E(\hat v_{n_k}(\zeta_k))-w_1\right|
\le
C\left|\E(\hat v_{n_k}(\zeta_k))-w_1\right|
\longrightarrow0 .
\]
Hence, for all sufficiently large $k$,
\[
\sup_{t\in[0,C]}
t\left|\E(\hat v_{n_k}(\zeta_k))-w_1\right|
\le
\frac{\delta}{2}.
\]
Consequently, for every $t\in[0,C]$ and sufficiently large $k$, 
\begin{align*}
|1+t\E(\hat v_{n_k}(\zeta_k))|&=|1+tw_1-tw_1+t\E(\hat v_{n_k}(\zeta_k))|\\
&\ge
|1+tw_1|
-
t|\E(\hat v_{n_k}(\zeta_k))-w_1|\ge
\frac{\delta}{2}.
\end{align*}
Thus,
\[
        \inf_{t\in[0,C]}
        |1+t\E(\hat v_{n_k}(\zeta_k))|
        \ge
        \frac{\delta}{2}>0
\]
for all sufficiently large $k$. Moreover,
\begin{align*}
\sup_{t\in[0,C]}
\left|
        \frac{t}{1+t\E(\hat v_{n_k}(\zeta_k))}
        -
        \frac{t}{1+tw_1}
\right|                                      &=
\sup_{t\in[0,C]}
\left|
        \frac{
        t^2\{w_1-\E(\hat v_{n_k}(\zeta_k))\}
        }{
        \{1+t\E(\hat v_{n_k}(\zeta_k))\}(1+tw_1)
        }
\right| \\[2 mm]                                            &\le
C\left|
        \E(\hat v_{n_k}(\zeta_k))-w_1
\right|
\longrightarrow0 .
\end{align*}
Therefore,
\[
        t\mapsto
        \frac{t}{1+t\E(\hat v_{n_k}(\zeta_k))}
\]
converges uniformly on \([0,C]\) to
\[
        t\mapsto
        \frac{t}{1+tw_1}.
\]
Since, with $H_{n_k}\Rightarrow H$ and $t\mapsto t(1+tw_1)^{-1}$ is continuous and bounded, we get
\begin{align*}
&\left|
\int_0^\infty
        \frac{t}{1+t\E(\hat v_{n_k}(\zeta_k))}
        \,dH_{n_k}(t)
        -
        \int_0^\infty
        \frac{t}{1+tw_1}\,dH(t)
\right|                                                        \\
&\quad\le
\left|
\int_0^\infty
        \left\{
        \frac{t}{1+t\E(\hat v_{n_k}(\zeta_k))}
        -
        \frac{t}{1+tw_1}
        \right\}
        \,dH_{n_k}(t)
\right|                                                        \\
&\qquad+
\left|
\int_0^\infty
        \frac{t}{1+tw_1}\,dH_{n_k}(t)
        -
        \int_0^\infty
        \frac{t}{1+tw_1}\,dH(t)
\right|                                                        \\
&\quad\le
\sup_{t\in[0,C]}
\left|
        \frac{t}{1+t\E(\hat v_{n_k}(\zeta_k))}
        -
        \frac{t}{1+tw_1}
\right|+
\left|
\int_0^\infty
        \frac{t}{1+tw_1}\,dH_{n_k}(t)
        -
        \int_0^\infty
        \frac{t}{1+tw_1}\,dH(t)
\right|
\longrightarrow 0 .
\end{align*}
Letting $k\to\infty$ in \eqref{subseq:implicit} and using the convergence
established in the preceding display, we obtain
\begin{align}\label{eq:w1-limit-equation}
        -\frac1{w_1}
        =
        \zeta_\ast
        -
        \gamma
        \int_0^\infty
        \frac{t}{1+tw_1}\,dH(t).
\end{align}

By \eqref{eq:scalar-deterministic-resolvent-bound},
applied to the $1\times1$ matrix $\Sigma=t$, we have
\[
    \inf_{t\in[0,C]}|1+t v_{n_k}(\zeta_k)|
    \geq c
    \coloneqq
    \frac{\operatorname{dist}(\Gamma,\mathbb R_{\geq0})}
         {\sup_{\zeta\in\Gamma}|\zeta|}
    >0.
\]
Since $v_{n_k}(\zeta_k)\to w_2$, the same lower bound
holds for $|1+tw_2|$ on $[0,C]$. Consequently,
\[
\sup_{t\in[0,C]}
\left|
    \frac{t}{1+t v_{n_k}(\zeta_k)}
    -
    \frac{t}{1+tw_2}
\right|
\leq
\frac{C^2}{c^2}
|v_{n_k}(\zeta_k)-w_2|
\longrightarrow0.
\]
Together with $H_{n_k}\Rightarrow H$ and the common compact
support of these measures, this yields
\begin{align}\label{eq:w2-limit-equation}
    -\frac1{w_2}
    =
    \zeta_\ast
    -
    \gamma
    \int_0^\infty
    \frac{t}{1+tw_2}\,dH(t).
\end{align}
Next we show that, both solutions are admissible in the sense of Lemma~\ref{lem:companion-uniq-D}. Indeed, if $\zeta_\ast$ is
nonreal, the fact that $\E\hat v_n$ and $v_n$ are Stieltjes
transforms of probability measures on $[0,\infty)$ implies that $\Im(w_j)\Im(\zeta_\ast)\geq0$ for $j=1,2$.
Equality is excluded by \eqref{eq:w1-limit-equation}
and~\eqref{eq:w2-limit-equation}, since otherwise
taking imaginary parts would give $\Im(\zeta_\ast)=0$. 
If $\zeta_\ast<0$, the Stieltjes representations of
$\widehat v_n$ in \eqref{eq:emp-companion-definition} and of $v_n$ in
Lemma~\ref{lem:companion-uniq-D} imply that
$w_1,w_2$ are real and nonnegative.
Since neither is zero by \eqref{eq:lower},
both are strictly positive.

Thus, Lemma~\ref{lem:companion-uniq-D} gives $w_1=w_2$.
Consequently,
\[
    |\E(\hat v_{n_k}(\zeta_k))-v_{n_k}(\zeta_k)|
    \longrightarrow |w_1-w_2|=0,
\]
contradicting the choice of the subsequence.
This proves \eqref{eq:weak-bias-on-Gamma}.

We now use \eqref{eq:weak-bias-on-Gamma} to prove
\eqref{eq:Ln-lower-bound}.

By Lemma~\ref{lem:companion-uniq-D}, we may write
\[
v_n(\zeta)=\int_0^\infty\frac1{x-\zeta}\,d\bar F_{\gamma_n,H_n}(x).
\]
By \cite[Theorems~10.6--10.7]{Rudin1987}, \(v_n\) is holomorphic on
\(\mathbb C\setminus\mathbb R_{\geq0}\). Moreover, differentiation
under the integral sign yields
\[
v_n'(\zeta)\coloneqq\frac{d}{d\zeta}v_n(\zeta)
=
\int_0^\infty\frac{1}{(x-\zeta)^2}\,
d\bar F_{\gamma_n,H_n}(x).
\]
Next, we differentiate \eqref{stielt}, written as an equation for $v_n(\zeta)^{-1}$ with $(\gamma_0,H_0)=(\gamma_n,H_n)$.
This yields
\begin{align}
\label{eq:finite-n-stability-derivative}
        \left\{
        \frac1{v_n(\zeta)^2}
        -
        \frac{\zeta^2}{n}
        \tr\!\left[
        \Sigma_n^2
        \{R_n^{\operatorname{det}}(\zeta)\}^2
        \right]
        \right\}
        v_n'(\zeta)
        =
        1 .
\end{align}
Since, for each \(\zeta\in\Gamma\), $|v_n'(\zeta)|
        \leq
        \operatorname{dist}(\Gamma,\mathbb R_{\geq0})^{-2}$,
it follows from \eqref{eq:finite-n-stability-derivative} that
\begin{align}
\label{eq:det-stability}
        \inf_{\zeta\in\Gamma}
        \left|
        \frac1{v_n(\zeta)^2}
        -
        \frac{\zeta^2}{n}
        \tr\!\left[
        \Sigma_n^2
        \{R_n^{\operatorname{det}}(\zeta)\}^2
        \right]
        \right|
        \geq
        \operatorname{dist}(\Gamma,\mathbb R_{\geq0})^2 .
\end{align}
It remains to compare \(L_n(\zeta)\) with the quantity in
\eqref{eq:det-stability}. By \eqref{eq:lower} and
\eqref{eq:weak-bias-on-Gamma},
\begin{align*}
\sup_{\zeta\in\Gamma}
\left|
        \frac1{\E(\hat v_n(\zeta))v_n(\zeta)}
        -
        \frac1{v_n(\zeta)^2}
\right|
&=
\sup_{\zeta\in\Gamma}
\left|
        \frac{v_n(\zeta)-\E(\hat v_n(\zeta))}
        {\E(\hat v_n(\zeta))v_n(\zeta)^2}
\right|
\\
&\leq
C
\sup_{\zeta\in\Gamma}
|\E(\hat v_n(\zeta))-v_n(\zeta)|
\longrightarrow0 .
\end{align*}

Furthermore, using \eqref{eq:resol-ident}, we obtain
\begin{align*}
&\frac{\zeta}{n}
\tr\!\left[
        \Sigma_n A_n(\zeta)\Sigma_n
        R_n^{\operatorname{det}}(\zeta)
\right]
+
\frac{\zeta^2}{n}
\tr\!\left[
        \Sigma_n^2
        \{R_n^{\operatorname{det}}(\zeta)\}^2
\right]
\\
&\quad=
\frac{\zeta}{n}
\tr\!\left[
        \Sigma_n
        \{A_n(\zeta)+\zeta R_n^{\operatorname{det}}(\zeta)\}
        \Sigma_n
        R_n^{\operatorname{det}}(\zeta)
\right]
\\
&\quad=
\frac{\zeta^2}{n}
\{\E(\hat v_n(\zeta))-v_n(\zeta)\}
\tr\!\left[
        \Sigma_n A_n(\zeta)\Sigma_n
        R_n^{\operatorname{det}}(\zeta)
        \Sigma_n
        R_n^{\operatorname{det}}(\zeta)
\right].
\end{align*}
Therefore, using Assumption~\ref{A2}, the uniform bounds on
\(A_n(\zeta)\) and \(R_n^{\operatorname{det}}(\zeta)\) on \(\Gamma\), and
\(p/n=\gamma_n=O(1)\), which follows from Assumption~\ref{A1}, we obtain
\begin{align*}
\sup_{\zeta\in\Gamma}
\left|
        \frac{\zeta}{n}
        \tr\!\left[
        \Sigma_n A_n(\zeta)\Sigma_n
        R_n^{\operatorname{det}}(\zeta)
        \right]
        +
        \frac{\zeta^2}{n}
        \tr\!\left[
        \Sigma_n^2
        \{R_n^{\operatorname{det}}(\zeta)\}^2
        \right]
\right|
\leq
C
\sup_{\zeta\in\Gamma}
|\E(\hat v_n(\zeta))-v_n(\zeta)|
\longrightarrow0 ,
\end{align*}
where the convergence follows from
\eqref{eq:weak-bias-on-Gamma}.

Combining this with the preceding bound, we conclude that
\begin{align*}
\sup_{\zeta\in\Gamma}
\Biggl|
        L_n(\zeta)
        -
        \left\{
        \frac1{v_n(\zeta)^2}
        -
        \frac{\zeta^2}{n}
        \tr\!\left[
        \Sigma_n^2
        \{R_n^{\operatorname{det}}(\zeta)\}^2
        \right]
        \right\}
\Biggr|
\longrightarrow0 .
\end{align*}
Hence, by the reverse triangle inequality,
\begin{align*}
\inf_{\zeta\in\Gamma}|L_n(\zeta)|
&\geq
\inf_{\zeta\in\Gamma}
\left|
        \frac1{v_n(\zeta)^2}
        -
        \frac{\zeta^2}{n}
        \tr\!\left[
        \Sigma_n^2
        \{R_n^{\operatorname{det}}(\zeta)\}^2
        \right]
\right|
\\
&\quad-
\sup_{\zeta\in\Gamma}
\Biggl|
        L_n(\zeta)
        -
        \left\{
        \frac1{v_n(\zeta)^2}
        -
        \frac{\zeta^2}{n}
        \tr\!\left[
        \Sigma_n^2
        \{R_n^{\operatorname{det}}(\zeta)\}^2
        \right]
        \right\}
\Biggr|.
\end{align*}
Combining the previous two displays with
\eqref{eq:det-stability}, we obtain
\[
        \inf_{\zeta\in\Gamma}|L_n(\zeta)|
        \geq
        \operatorname{dist}(\Gamma,\mathbb R_{\geq0})^2-o(1),
\]
which proves \eqref{eq:Ln-lower-bound}.
\end{proof}

\subsection{Lemma~\ref{lem:uniform-companion-rate} and its proof}

\begin{lemma}
\label{lem:uniform-companion-rate}
Suppose Assumptions~\ref{A3} and~\ref{A2} hold, and suppose that
\begin{align}
\label{eq:uniform-companion-aspect-ratio}
        0
        <
        \liminf_{n\to\infty}\gamma_n
        \leq
        \limsup_{n\to\infty}\gamma_n
        <
        \infty.
\end{align}
Let \(K\subseteq\mathbb C\setminus\mathbb R_{\geq0}\) be compact. Then
\begin{align}
\label{clt:stat:1}
        \sup_{\zeta\in K}
        \left|
        \hat v_n(\zeta)-v_n(\zeta)
        \right|
        =
        O_{\mathbb P}(n^{-1}).
\end{align}
\end{lemma}

\begin{proof}
Suppose, to the contrary, that \eqref{clt:stat:1} does not hold.
Equivalently, the sequence
\[
        n\sup_{\zeta\in K}
        \left|
        \hat v_n(\zeta)-v_n(\zeta)
        \right|
\]
is not bounded in probability. Hence, there exists
\(\varepsilon_0>0\) such that, for every \(M>0\) and every
\(N\in\mathbb N\), there is an \(n\geq N\) satisfying
\[
        \mathbb P\left(
        n\sup_{\zeta\in K}
        \left|
        \hat v_n(\zeta)-v_n(\zeta)
        \right|
        >M
        \right)
        >
        \varepsilon_0.
\]
Choose a sequence \((M_j)_{j\geq1}\) such that \(M_j\to\infty\).
We may then choose recursively a strictly increasing sequence
\((n_j)_{j\geq1}\) such that
\begin{align}
\label{eq:uniform-companion-contradiction}
        \mathbb P\left(
        n_j\sup_{\zeta\in K}
        \left|
        \hat v_{n_j}(\zeta)-v_{n_j}(\zeta)
        \right|
        >M_j
        \right)
        >
        \varepsilon_0
\end{align}
for every \(j\geq1\).

By \eqref{eq:uniform-companion-aspect-ratio}, the sequence
\((\gamma_{n_j})_{j\geq1}\) is bounded and bounded away from zero.
Therefore, after passing to a subsequence, there exists
\(\gamma_\ast\in(0,\infty)\) such that
\[
        \gamma_{n_j}\longrightarrow\gamma_\ast.
\]

By Assumption~\ref{A2}, there exists \(C<\infty\) such that
\[
        \operatorname{supp}(H_n)\subseteq[0,C]
\]
for every \(n\in\mathbb N\). In particular, the sequence
\((H_{n_j})_{j\geq1}\) is tight. Hence, by Prohorov's theorem, after
passing to a further subsequence, there exists a probability
distribution \(H_\ast\) such that
\[
        H_{n_j}\Rightarrow H_\ast.
\]
Since every \(H_{n_j}\) is supported on the closed interval \([0,C]\),
the limiting probability distribution \(H_\ast\) is also supported on
\([0,C]\).

After passing to these subsequences and relabeling them together with
the corresponding thresholds \(M_j\), we may therefore assume that
\begin{equation}
\label{eq:uniform-companion-subsequence-limits}
        \gamma_{n_j}\longrightarrow\gamma_\ast\in(0,\infty),
        \qquad
        H_{n_j}\Rightarrow H_\ast.
\end{equation}
Consequently, the relabeled subsequence satisfies
Assumptions~\ref{A1}--\ref{A4}, with \((\gamma_\ast,H_\ast)\) in place
of \((\gamma,H)\).

Applying Lemma~\ref{stielt:var} to this subsequence gives
\[
        \sup_{\zeta\in K}
        \left|
        \hat v_{n_j}(\zeta)
        -
        \mathbb E\bigl(\hat v_{n_j}(\zeta)\bigr)
        \right|
        =
        O_{\mathbb P}(n_j^{-1}),
\]
whereas Lemma~\ref{stielt:bias} gives
\[
        \sup_{\zeta\in K}
        \left|
        \mathbb E\bigl(\hat v_{n_j}(\zeta)\bigr)
        -
        v_{n_j}(\zeta)
        \right|
        =
        O(n_j^{-1}).
\]
Thus, by the triangle inequality,
\[
\begin{split}
        \sup_{\zeta\in K}
        \left|
        \hat v_{n_j}(\zeta)-v_{n_j}(\zeta)
        \right|
        \leq
        \sup_{\zeta\in K}
        \left|
        \hat v_{n_j}(\zeta)
        -
        \mathbb E\bigl(\hat v_{n_j}(\zeta)\bigr)
        \right|
        +
        \sup_{\zeta\in K}
        \left|
        \mathbb E\bigl(\hat v_{n_j}(\zeta)\bigr)
        -
        v_{n_j}(\zeta)
        \right|
        =
        O_{\mathbb P}(n_j^{-1}).
\end{split}
\]
It follows that
\[
        n_j\sup_{\zeta\in K}
        \left|
        \hat v_{n_j}(\zeta)-v_{n_j}(\zeta)
        \right|
        =
        O_{\mathbb P}(1).
\]
Therefore, there exist \(M<\infty\) and \(j_0\in\mathbb N\) such that
\[
        \mathbb P\left(
        n_j\sup_{\zeta\in K}
        \left|
        \hat v_{n_j}(\zeta)-v_{n_j}(\zeta)
        \right|
        >M
        \right)
        <
        \frac{\varepsilon_0}{2}
\]
for every \(j\geq j_0\).

Since \(M_j\to\infty\), we have \(M_j\geq M\) for all sufficiently
large \(j\). Hence, for all such \(j\),
\[
\mathbb P\left(
        n_j\sup_{\zeta\in K}
        \left|
        \hat v_{n_j}(\zeta)-v_{n_j}(\zeta)
        \right|
        >M_j
        \right)
\leq
\mathbb P\left(
        n_j\sup_{\zeta\in K}
        \left|
        \hat v_{n_j}(\zeta)-v_{n_j}(\zeta)
        \right|
        >M
        \right)
<
\frac{\varepsilon_0}{2}.
\]
This contradicts \eqref{eq:uniform-companion-contradiction}. Therefore,
\eqref{clt:stat:1} holds.
\end{proof}

\subsection{Proof of Theorem~\ref{Lemma:GreenRom}}

\begin{proof}[Proof of Theorem~\ref{Lemma:GreenRom}]
Fix a compact set
$K\subseteq\mathbb C\setminus\mathbb R_{\geq0}$.
Choose $K_1$ as in Lemma~\ref{lem:cauchy_reduction_random}.

Using Lemma~\ref{lem:resolvent-denominator-bounds},
there exists a constant $C_{K_1}<\infty$ such that,
for every $\zeta\in K_1$,
\begin{align}
\label{RM:bounds}
        \|R_n(\zeta)\|_2\leq C_{K_1},
        \qquad
        \|R_{n,-j}(\zeta)\|_2\leq C_{K_1},
        \qquad
        \|R_n^{\operatorname{det}}(\zeta)\|_2\leq C_{K_1}.
\end{align}

By the resolvent identity,
\[
        R_n(\zeta)-R_n^{\operatorname{det}}(\zeta)
        =
        R_n(\zeta)
        \bigl(
        -\zeta v_n(\zeta)\Sigma_n-\hat\Sigma_n
        \bigr)
        R_n^{\operatorname{det}}(\zeta).
\]
Adding and subtracting \(-\zeta\hat v_n(\zeta)\Sigma_n\), we obtain
\begin{align}
\label{dec:RM:1}
\begin{split}
        \tr\!\left[
        \Theta_n\{R_n(\zeta)-R_n^{\operatorname{det}}(\zeta)\}
        \right]
        &=
        -\zeta\{v_n(\zeta)-\hat v_n(\zeta)\}
        \tr\!\left[
        \Theta_nR_n(\zeta)\Sigma_n
        R_n^{\operatorname{det}}(\zeta)
        \right]
        \\[1mm]
        &\qquad+
        \tr\!\left[
        \Theta_nA_n(\zeta)
        R_n^{\operatorname{det}}(\zeta)
        \right],
\end{split}
\end{align}
where
\[
        A_n(\zeta)
        \coloneqq
        -\zeta\hat v_n(\zeta)R_n(\zeta)\Sigma_n
        -
        R_n(\zeta)\hat\Sigma_n .
\]

We shall bound the two terms on the right-hand side of
\eqref{dec:RM:1} separately. First, we show that
\begin{align}
\label{Main:thm:step1}
\sup_{\zeta\in K}
\left|
        \zeta\{v_n(\zeta)-\hat v_n(\zeta)\}
        \tr\!\left[
        \Theta_nR_n(\zeta)\Sigma_n
        R_n^{\operatorname{det}}(\zeta)
        \right]
\right|
=
O_{\mathbb P}\left(
        \frac{\|\Theta_n\|_1}{n}
\right),
\end{align}
and, second, that
\begin{align}
\label{Main:thm:step2}
\sup_{\zeta\in K}
\left|
        \tr\!\left[
        \Theta_nA_n(\zeta)R_n^{\operatorname{det}}(\zeta)
        \right]
\right|
=
O_{\mathbb P}\left(
        \frac{\|\Theta_n\|_1}{n}
        +
        \frac{\|\Theta_n\|_F}{n^{1/2}}
\right).
\end{align}
Together, \eqref{Main:thm:step1} and \eqref{Main:thm:step2}
imply the desired claim, because
\[
    \frac{\|\Theta_n\|_1}{n}
    \leq\sqrt{\frac pn}\frac{\|\Theta_n\|_F}{\sqrt n}
    =O\left(\frac{\|\Theta_n\|_F}{\sqrt n}\right)
\]
under \eqref{gamma_n:cond}. Using Assumption~\ref{A2} and \eqref{RM:bounds}, we obtain, for every
\(\zeta\in K\),
\[
\left|
        \tr\!\left[
        \Theta_nR_n(\zeta)\Sigma_n
        R_n^{\operatorname{det}}(\zeta)
        \right]
\right|
\leq
\|\Theta_n\|_1
\left\|
        R_n(\zeta)\Sigma_n
        R_n^{\operatorname{det}}(\zeta)
\right\|_2
\leq
C_{K_1}\|\Theta_n\|_1.
\]
Therefore,
\[
\sup_{\zeta\in K}
\left|
        -\zeta\{v_n(\zeta)-\hat v_n(\zeta)\}
        \tr\!\left[
        \Theta_nR_n(\zeta)\Sigma_n
        R_n^{\operatorname{det}}(\zeta)
        \right]
\right|
\leq
C_{K_1}\|\Theta_n\|_1
\sup_{\zeta\in K}
\left|
        v_n(\zeta)-\hat v_n(\zeta)
\right|.
\]
Lemma~\ref{lem:uniform-companion-rate} therefore establishes
\eqref{Main:thm:step1}.

We continue with \eqref{Main:thm:step2}. By the Sherman--Morrison formula,
\[
        R_n(\zeta)
        =
        R_{n,-j}(\zeta)
        -
        \frac{
        R_{n,-j}(\zeta)(x_jx_j^\top/n)R_{n,-j}(\zeta)
        }{
        1+\frac1n x_j^\top R_{n,-j}(\zeta)x_j
        },
\]
and using \eqref{eq:beta_jandq_j}, we have
\[
        R_n(\zeta)x_j
        =
        \beta_j(\zeta)R_{n,-j}(\zeta)x_j .
\]
Substituting these identities into \(A_n(\zeta)\), we obtain
\begin{align*}
        A_n(\zeta)
        &=
        \frac1n\sum_{j=1}^n
        \beta_j(\zeta)
        R_{n,-j}(\zeta)
        \bigl(\Sigma_n-x_jx_j^\top\bigr)
        \\
        &\qquad-
        \frac1{n^2}\sum_{j=1}^n
        \beta_j(\zeta)^2
        R_{n,-j}(\zeta)x_jx_j^\top
        R_{n,-j}(\zeta)\Sigma_n .
\end{align*}
Consequently,
\[
        \tr\!\left[
        \Theta_n A_n(\zeta)R_n^{\operatorname{det}}(\zeta)
        \right]
        =
        T_1(\zeta)+T_2(\zeta),
\]
where
\begin{align*}
        T_1(\zeta)
        &\coloneqq
        \frac1n\sum_{j=1}^n
        \beta_j(\zeta)
        \tr\!\left[
        \Theta_n R_{n,-j}(\zeta)
        (\Sigma_n-x_jx_j^\top)
        R_n^{\operatorname{det}}(\zeta)
        \right],
\end{align*}
and
\begin{align*}
        T_2(\zeta)
        &\coloneqq
        -\frac1{n^2}\sum_{j=1}^n
        \beta_j(\zeta)^2
        \tr\!\left[
        \Theta_n R_{n,-j}(\zeta)
        x_jx_j^\top
        R_{n,-j}(\zeta)\Sigma_n
        R_n^{\operatorname{det}}(\zeta)
        \right].
\end{align*}

The remaining part of the proof consists of showing that 
\begin{align}
   \label{eq:T1:OP}
        \sup_{\zeta\in K}|T_1(\zeta)|
        & =
        O_{\mathbb P}\left(\frac{\|\Theta_n\|_F}{n^{1/2}}\right), \\ 
        \label{eq:T2:OP}
        \sup_{\zeta\in K}|T_2(\zeta)|
       & =
         O_{\mathbb P}\left(
        \frac{\|\Theta_n\|_1}{n}
        \right).
\end{align}

Indeed, combining the bounds in \eqref{eq:T1:OP} and \eqref{eq:T2:OP} gives \eqref{Main:thm:step2}, which completes the proof.

\noindent \textbf{Proof of \eqref{eq:T1:OP}.}
We decompose
\begin{align}\label{T1:dec}
        T_1(\zeta)
        =
        T_{1,1}(\zeta)+T_{1,2}(\zeta),
\end{align}
where
\begin{align*}
        T_{1,1}(\zeta)
        &\coloneqq
        \{-\zeta\hat v_n(\zeta)\}
        \frac1n\sum_{j=1}^n Q_j(\zeta),\qquad
        T_{1,2}(\zeta)
        \coloneqq
        \frac1n\sum_{j=1}^n
        \left[
        \beta_j(\zeta)-\{-\zeta\hat v_n(\zeta)\}
        \right]Q_j(\zeta),
\end{align*}
and
\[
        Q_j(\zeta)
        \coloneqq
        \tr\!\left[
        \Theta_nR_{n,-j}(\zeta)
        (\Sigma_n-x_jx_j^\top)
        R_n^{\operatorname{det}}(\zeta)
        \right].
\]
We are now going to bound the terms \(T_{1,1}(\zeta)\) and
\(T_{1,2}(\zeta)\) on the right-hand side of \eqref{T1:dec} separately.
By cyclicity of the trace,
\[
    Q_j(\zeta)=-\bigl\{
    x_j^\top R_n^{\operatorname{det}}(\zeta)\Theta_nR_{n,-j}(\zeta)x_j
    -\tr[\Sigma_nR_n^{\operatorname{det}}(\zeta)
                \Theta_nR_{n,-j}(\zeta)]\bigr\}.
\]
We apply Lemma~\ref{RubMes_Lem4_fix} with
\[
    E_n(\zeta)=R_n^{\operatorname{det}}(\zeta)\Theta_n,
    \qquad F_n(\zeta)=I_p.
\]
These matrices are deterministic and their entries are holomorphic.
By \eqref{RM:bounds},
\[
    \sup_{n\geq1}\sup_{\zeta\in K_1}\|F_n(\zeta)\|_2=1,
    \qquad
    \sup_{\zeta\in K_1}\|E_n(\zeta)\|_F
    \leq C_{K_1}\|\Theta_n\|_F.
\]
When $\Theta_n\neq0$, the second condition of the lemma holds
with $r_n=\|\Theta_n\|_F$. The case $\Theta_n=0$ is immediate.
Since $p/n=O(1)$ by \eqref{gamma_n:cond}, the lemma yields
\[
    \sup_{\zeta\in K}
    \left|\frac1n\sum_{j=1}^nQ_j(\zeta)\right|
    =O_{\mathbb P}\left(\frac{\|\Theta_n\|_F}{\sqrt n}\right).
\]
Also note that, 
\[
        \sup_{\zeta\in K}|\{-\zeta\hat v_n(\zeta)\}|
        \le
        \sup_{\zeta\in K}|\{-\zeta v_n(\zeta)\}|
        +
        \sup_{\zeta\in K}|\zeta|\,|\hat v_n(\zeta)-v_n(\zeta)|
        =
        O_{\mathbb P}(1),
\]
where we used \eqref{clt:stat:1} and that $-\zeta v_n(\zeta)$ is uniformly bounded on $K$.
Hence,
\begin{align}\label{T1:first:avg}
        \sup_{\zeta\in K}
        \left|
        \{-\zeta\hat v_n(\zeta)\}
        \frac1n\sum_{j=1}^n Q_j(\zeta)
        \right|
        =
        O_{\mathbb P}\left(\frac{\|\Theta_n\|_F}{n^{1/2}}\right).
\end{align}
It remains to bound 
$T_{1,2}(\zeta)$
in \eqref{T1:dec}.
By the Cauchy-Schwarz inequality, we have
\begin{align}\label{T1:fluct:CS}
\begin{split}
\left|
        \frac1n\sum_{j=1}^n
        \left[
        \beta_j(\zeta)-\{-\zeta\hat v_n(\zeta)\}
        \right]Q_j(\zeta)
\right|                                      \le
\left\{
        \frac1n\sum_{j=1}^n
        \left|
        \beta_j(\zeta)-\{-\zeta\hat v_n(\zeta)\}
        \right|^2
\right\}^{1/2}
\left\{
        \frac1n\sum_{j=1}^n |Q_j(\zeta)|^2
\right\}^{1/2}.
\end{split}
\end{align}
We first control the first factor on the right-hand side of \eqref{T1:fluct:CS}. Define
\[
        \tilde\beta_n(\zeta)
        \coloneqq
        \left\{1+ \frac1n\tr\{\Sigma_nR_n(\zeta)\}\right\}^{-1}.
\]
Applying \eqref{eq:general-trace-denominator-bound}
with $A=\widehat\Sigma_n$, $B=\Sigma_n$, and $c=n^{-1}$,
and \eqref{eq:general-denominator-bound}
with $A=\widehat\Sigma_{n,-j}$, $x=x_j$, and $c=n^{-1}$,
gives
\[
        |\tilde\beta_n(\zeta)|
        \leq
        \frac{|\zeta|}
        {\operatorname{dist}(\zeta,\mathbb R_{\geq0})},\qquad\text{and}\qquad |\beta_j(\zeta)|
        \leq
        \frac{|\zeta|}
        {\operatorname{dist}(\zeta,\mathbb R_{\geq0})}..
\]
Since \(K_1\subset\mathbb C\setminus\mathbb R_{\geq0}\) is compact,
\begin{align}\label{eq:tilde-beta-bound}
        \sup_{\zeta\in K_1}
        |\tilde\beta_n(\zeta)|
        \leq C_{K_1}\qquad\text{and}\qquad \sup_{\zeta\in K_1}
        |\beta_j(\zeta)|
        \leq C_{K_1}.
\end{align}
Moreover, since
\[
        \{-\zeta\hat v_n(\zeta)\}
        =
        \frac1n\sum_{j=1}^n\beta_j(\zeta),
\]
and the empirical mean minimizes the quadratic loss, we have
\[
        \frac1n\sum_{j=1}^n
        \left|
        \beta_j(\zeta)-\{-\zeta\hat v_n(\zeta)\}
        \right|^2
        \le
        \frac1n\sum_{j=1}^n
        |\beta_j(\zeta)-\tilde\beta_n(\zeta)|^2 .
\]
Moreover,
\begin{align*}
        \beta_j(\zeta)-\tilde\beta_n(\zeta)
        &=\frac{1}{1+\frac{1}{n}x_j^\top R_{n,-j}(\zeta)x_j}-\frac{1}{1+\frac{1}{n}\tr(\Sigma_n R_n(\zeta))}\\&=
        -\beta_j(\zeta)\tilde\beta_n(\zeta)
        \left\{
        \frac1n x_j^\top R_{n,-j}(\zeta)x_j
        -
        \frac1n\tr\{\Sigma_nR_n(\zeta)\}
        \right\}.
\end{align*}
Recall $\gamma_j(\zeta)$ from \eqref{eq:gamma_def} and put
\[
    \rho_j(\zeta)=\frac1n
        \tr\{\Sigma_n(R_{n,-j}(\zeta)-R_n(\zeta))\}.
\]
Using \eqref{eq:tilde-beta-bound}, we obtain
\[
        \left|
        \beta_j(\zeta)-\tilde\beta_n(\zeta)
        \right|
        \le
        C_{K_1}\{|\gamma_j(\zeta)|+|\rho_j(\zeta)|\}.
\]
Hence,
\begin{align}\label{r1}
\sup_{\zeta\in K}
        \frac1n\sum_{j=1}^n
        \left|
        \beta_j(\zeta)-\{-\zeta\hat v_n(\zeta)\}
        \right|^2                            \le
C_{K_1}\sup_{\zeta\in K}
        \frac1n\sum_{j=1}^n|\gamma_j(\zeta)|^2
+
C_{K_1}\sup_{\zeta\in K}
        \frac1n\sum_{j=1}^n|\rho_j(\zeta)|^2.
\end{align}
By Lemma~\ref{lem_b26_bai_silverstein}, applied conditionally on
$\mathcal F_{n,-j}$, together with \eqref{RM:bounds} and
Assumption~\ref{A2},
\[
    \sup_{\zeta\in K_1}\max_{1\leq j\leq n}
       \mathbb E|\gamma_j(\zeta)|^2\leq C_{K_1}\frac p{n^2}.
\]
Each $\gamma_j$ is holomorphic, so
Lemma~\ref{lem:cauchy_reduction_random} gives
\[
    \mathbb E\left[\sup_{\zeta\in K}
       \frac1n\sum_{j=1}^n|\gamma_j(\zeta)|^2\right]
    \leq\frac1n\sum_{j=1}^n
       \mathbb E\sup_{\zeta\in K}|\gamma_j(\zeta)|^2
    \leq C_{K_1}\frac p{n^2}.
\]
Markov's inequality and \eqref{gamma_n:cond} therefore yield
\begin{equation}\label{r2}
    \sup_{\zeta\in K}\frac1n\sum_{j=1}^n|\gamma_j(\zeta)|^2
    =O_{\mathbb P}(pn^{-2})=O_{\mathbb P}(n^{-1}).
\end{equation}
For $\rho_j(\zeta)$, the Sherman--Morrison identity gives
\[
        R_{n,-j}(\zeta)-R_n(\zeta)
        =
        \frac{\beta_j(\zeta)}{n}
        R_{n,-j}(\zeta)x_jx_j^\top R_{n,-j}(\zeta).
\]
Therefore,
\[
           \sup_{\zeta\in K_1}|\rho_j(\zeta)|
        \le
        \frac{C_{K_1}}{n^2}\|x_j\|_{2}^2.
\]
By the uniform moment bound in Assumption~\ref{A3}, we have
\begin{align} \label{r3}
        \sup_{\zeta\in K_1}\frac1n\sum_{j=1}^n|\rho_j(\zeta)|^2\leq \frac1n\sum_{j=1}^n\sup_{\zeta\in K_1}|\rho_j(\zeta)|^2
        \le
        \frac{C_{K_1}}{n^5}\sum_{j=1}^n\|x_j\|_{2}^4
        =O_{\mathbb P}(p^2n^{-4})=
        O_{\mathbb P}(n^{-2}),
\end{align}
where the last equality follows by \eqref{gamma_n:cond}.
Consequently, by combining \eqref{r1}, \eqref{r2} and \eqref{r3},
\begin{align}\label{beta:emp:fluct}
        \sup_{\zeta\in K}
        \left\{
        \frac1n\sum_{j=1}^n
        \left|
        \beta_j(\zeta)-\{-\zeta\hat v_n(\zeta)\}
        \right|^2
        \right\}^{1/2}
        =
        O_{\mathbb P}(n^{-1/2}).
\end{align}
We now control the second factor in the upper bound \eqref{T1:fluct:CS} for $|T_{1,2}(\zeta)|$. For each fixed
\(j\), the process \(Q_j\) satisfies the assumptions of
Lemma~\ref{lem:cauchy_reduction_random}. Indeed, for each fixed
\(\zeta\in\mathbb C\setminus\mathbb R_{\ge0}\), \(Q_j(\zeta)\) is
\(\mathcal F\)-measurable. Moreover, for each \(\omega\in\Omega\), the map
\(\zeta\mapsto Q_j(\zeta,\omega)\) is holomorphic on
$\mathbb C\setminus\mathbb R_{\ge0}$, since both
$\zeta\mapsto R_{n,-j}(\zeta,\omega)$ and
$\zeta\mapsto  R_n^{\operatorname{det}}(\zeta)$ are holomorphic on this domain. Hence,
by Lemma~\ref{lem:cauchy_reduction_random},
\[
        \E\left(\sup_{\zeta\in K}|Q_j(\zeta)|^2\right)
        \le
        C_{K_1}
        \sup_{\zeta\in K_1}\E(|Q_j(\zeta)|^2).
\]
Therefore,
\[
        \E\left[
        \sup_{\zeta\in K}
        \frac1n\sum_{j=1}^n |Q_j(\zeta)|^2
        \right]
        \le
        \frac1n\sum_{j=1}^n
        \E\left(\sup_{\zeta\in K}|Q_j(\zeta)|^2\right)       \le
        C_{K_1}\frac1n\sum_{j=1}^n
        \sup_{\zeta\in K_1}\E(|Q_j(\zeta)|^2).
\]
By Lemma~\ref{lem_b26_bai_silverstein}, applied conditionally on
\(\mathcal F_{n,-j}\), together with Assumption~\ref{A2} and the uniform
resolvent bounds on \(K_1\),
\[
        \sup_{\zeta\in K_1}\E(|Q_j(\zeta)|^2)
        \le
        C_{K_1}\|\Theta_n\|_F^2 .
\]
Consequently,
\[
        \E\left[
        \sup_{\zeta\in K}
        \frac1n\sum_{j=1}^n |Q_j(\zeta)|^2
        \right]
        \le
        C_{K_1}\|\Theta_n\|_F^2,
\]
and by Markov's inequality,
\begin{align}\label{Q:emp:square}
        \sup_{\zeta\in K}
        \left\{
        \frac1n\sum_{j=1}^n |Q_j(\zeta)|^2
        \right\}^{1/2}
        =
        O_{\mathbb P}(\|\Theta_n\|_F).
\end{align}
Combining \eqref{T1:fluct:CS}, \eqref{beta:emp:fluct} and
\eqref{Q:emp:square}, we get
\begin{align}\label{T1:fluct:bound}
        \sup_{\zeta\in K}
        \left| T_{1,2} (\zeta)
        \right|
        =
        O_{\mathbb P}\left(\frac{\|\Theta_n\|_F}{n^{1/2}}\right).
\end{align}
Finally, \eqref{T1:dec}, \eqref{T1:first:avg} and \eqref{T1:fluct:bound} imply \eqref{eq:T1:OP}.

\vspace{0.5cm}
\noindent \textbf{Proof of \eqref{eq:T2:OP}.}
Using $x_j=\Sigma_n^{1/2}z_j$, write for each $\zeta\in K_1$
\[
        T_2(\zeta)
        =
        -\frac1{n^2}\sum_{j=1}^n
        \beta_j(\zeta)^2
        z_j^\top D_{-j}(\zeta)z_j,
\]
where
\[
      D_{-j}(\zeta)
        \coloneqq
        \Sigma_n^{1/2}
        R_{n,-j}(\zeta)\Sigma_nR_n^{\operatorname{det}}(\zeta)\Theta_n
        R_{n,-j}(\zeta)
        \Sigma_n^{1/2}.
\]
We first prove a second-moment bound for \(T_2(\zeta)\). Using the uniform
bound on $\beta_j(\zeta)$ from \eqref{eq:tilde-beta-bound} and Jensen's inequality, we have, for each \(\zeta\in K_1\),
\begin{align} \label{r5}
        |T_2(\zeta)|^2
        \le
        \frac{C_{K_1}}{n^2}
        \left(\frac{1}{n}
        \sum_{j=1}^n
        \left|
        z_j^\top D_{-j}(\zeta)z_j
        \right|
        \right)^2                                 \le
        \frac{C_{K_1}}{n^3}
        \sum_{j=1}^n
        \left|
        z_j^\top D_{-j}(\zeta)z_j
        \right|^2.
\end{align}
By Lemma~\ref{lem_b26_bai_silverstein}, applied conditionally on
\(\mathcal F_{n,-j}\), together with Assumption~\ref{A2} and \eqref{eq:tilde-beta-bound},
\begin{align*}
        \E\left(
        \left|
        z_j^\top D_{-j}(\zeta)z_j
        \right|^2
        \,\middle|\,\mathcal F_{n,-j}
        \right)
        &\le  2\E\left(
        \left|
        z_j^\top D_{-j}(\zeta)z_j-\tr(D_{-j}(\zeta))
        \right|^2
        \,\middle|\,\mathcal F_{n,-j}
        \right)+2|\tr(D_{-j}(\zeta))|^2 \\
        &\leq C_{K_1}\|D_{-j}(\zeta)\|_F^2
        +
        C_{K_1}|\tr(D_{-j}(\zeta))|^2 .
\end{align*}

Moreover, by assumption \ref{A2} and \eqref{RM:bounds}
\[
        \|D_{-j}(\zeta)\|_F
        \le
        C_{K_1}\|\Theta_n\|_F,
\]
and, similarly, 
\begin{align*}
    |\tr(D_{-j}(\zeta))|
        &=
        \left|
        \tr\!\left(
        R_{n,-j}(\zeta)\Sigma_n R_{n,-j}(\zeta)\Sigma_n
        R_n^{\operatorname{det}}(\zeta)\Theta_n
        \right)
        \right|                                      \\
        &\le
        \|R_{n,-j}(\zeta)\Sigma_n R_{n,-j}(\zeta)\Sigma_n
        R_n^{\operatorname{det}}(\zeta)\|_2\|\Theta_n\|_1          \\
        &\le
        C_{K_1}\|\Theta_n\|_1 .
\end{align*}
Consequently,
\begin{align} \label{r6}
        \sup_{\zeta\in K_1}
        \E\bigl(|
        z_j^\top D_{-j}(\zeta)z_j
        |^2\bigr)
        \le
        C_{K_1}\left(\|\Theta_n\|_F^2+\|\Theta_n\|_1^2\right).
\end{align}
Hence, by combining \eqref{r5} and \eqref{r6}
\begin{align*}
        \sup_{\zeta\in K_1}\E\bigl(|T_2(\zeta)|^2\bigr)
        \le
        \frac{C}{n^3}
        \sum_{j=1}^n
        \sup_{\zeta\in K_1}
        \E\bigl(|
        z_j^\top D_{-j}(\zeta)z_j
        |^2\bigr)                                  \le
        \frac{C_{K_1}}{n^2}
        \left(\|\Theta_n\|_F^2+\|\Theta_n\|_1^2\right).
\end{align*}
To verify \eqref{eq:T2:OP}, it remains to transfer this bound to a bound uniform over $K$ via Lemma \ref{lem:cauchy_reduction_random}. For each fixed
$\zeta\in\mathbb C\setminus\mathbb R_{\ge0}$, $T_2(\zeta)$ is
$\mathcal F$-measurable. Moreover, for every $\omega\in\Omega$, the maps
\[
        \zeta\mapsto R_{n,-j}(\zeta,\omega),
        \qquad
       \zeta\mapsto R_n^{\operatorname{det}}(\zeta)
\]
are holomorphic on \(\mathbb C\setminus\mathbb R_{\ge0}\). Since
\[
        \beta_j(\zeta)
        =
        \left(
        1+\frac1n x_j^\top R_{n,-j}(\zeta)x_j
        \right)^{-1}
\]
and the denominator is nonzero on \(\mathbb C\setminus\mathbb R_{\ge0}\), the
map \(\zeta\mapsto \beta_j(\zeta,\omega)\) is holomorphic on this domain as
well. Hence \(\zeta\mapsto D_{-j}(\zeta,\omega)\), and therefore
\(\zeta\mapsto T_2(\zeta,\omega)\), are holomorphic on
\(\mathbb C\setminus\mathbb R_{\ge0}\). Thus \(T_2\) satisfies the assumptions
of Lemma~\ref{lem:cauchy_reduction_random}. Applying this lemma gives
\begin{align*}
        \E\left(\sup_{\zeta\in K}|T_2(\zeta)|^2\right)
        \le
        C_{K_1}
        \sup_{\zeta\in K_1}\E\bigl(|T_2(\zeta)|^2\bigr)\le
        \frac{C_{K_1}}{n^2}
        \left(\|\Theta_n\|_F^2+\|\Theta_n\|_1^2\right).
\end{align*}
By Markov's inequality,
\[
        \sup_{\zeta\in K}|T_2(\zeta)|
        =
        O_{\mathbb P}\left(
        \frac{\|\Theta_n\|_F}{n}
        +
        \frac{\|\Theta_n\|_1}{n}
        \right).
\]
Since \(\|\Theta_n\|_F\leq \|\Theta_n\|_1\), this implies
\eqref{eq:T2:OP}.

\end{proof}

\subsection{Proof of Corollary~\ref{cor:GreenRom-derivatives:rates}}

\begin{proof}[Proof of Corollary~\ref{cor:GreenRom-derivatives:rates}]
For $q=0$, the assertion follows directly from
Theorem~\ref{Lemma:GreenRom}. We therefore suppose that $q\geq1$.

Choose $\rho>0$ and the compact set $K_1\subseteq
\mathbb C\setminus\mathbb R_{\geq0}$ as in
Lemma~\ref{lem:cauchy_reduction_random}, so that
\[
        K_1
        =
        \{w\in\mathbb C:\operatorname{dist}(w,K)\leq\rho\}.
\]
In particular, for every $\zeta\in K$,
\begin{align}\label{eq:inclusion_sets}
        \overline{B(\zeta,\rho)}\eqqcolon
         \{w\in\mathbb C:|w-\zeta|\leq\rho\}
        \subseteq
        K_1
        \subseteq
        \mathbb C\setminus\mathbb R_{\geq0}.
\end{align}

For each $n\in\mathbb N$, define
\[
        f_n(w)
        \coloneqq
        \tr\!\left[
        \Theta_n
        \left\{
        R_n(w)
        -
        R_n^{\operatorname{det}}(w)
        \right\}
        \right],
        \qquad
        w\in\mathbb C\setminus\mathbb R_{\geq0}.
\]
The function $f_n$ is holomorphic on
$\mathbb C\setminus\mathbb R_{\geq0}$. Indeed, $R_n(w)$ is the
resolvent of a nonnegative definite matrix, while $v_n$ is holomorphic
and $R_n^{\operatorname{det}}(w)$ is well defined and holomorphic on
$\mathbb C\setminus\mathbb R_{\geq0}$. Hence, both matrix-valued
functions appearing in the definition of $f_n$ are holomorphic on
this domain.

Let $\zeta\in K$. Using \eqref{eq:inclusion_sets}, 
Cauchy's integral formula for derivatives
\cite[Chapter~10, Exercise~7]{Rudin1987} gives
\[
        \frac{d^q}{d\zeta^q}f_n(\zeta)
        =
        \frac{q!}{2\pi i}
        \int_{|w-\zeta|=\rho}
        \frac{f_n(w)}{(w-\zeta)^{q+1}}\,dw.
\]
Consequently,
\begin{align*}
        \left|
        \frac{d^q}{d\zeta^q}f_n(\zeta)
        \right|
        &\leq
        \frac{q!}{2\pi}
        \int_{|w-\zeta|=\rho}
        \frac{|f_n(w)|}{|w-\zeta|^{q+1}}\,|dw|\\
        &=
        \frac{q!}{2\pi\rho^{q+1}}
        \int_{|w-\zeta|=\rho}
        |f_n(w)|\,|dw|
        \leq
        \frac{q!}{\rho^q}
        \sup_{w\in K_1}|f_n(w)|.
\end{align*}
Taking the supremum over $\zeta\in K$, we obtain
\[
\sup_{\zeta\in K}
\left|
\frac{d^q}{d\zeta^q}
\tr\!\left[
\Theta_n
\left\{
R_n(\zeta)
-
R_n^{\operatorname{det}}(\zeta)
\right\}
\right]
\right|
\leq
\frac{q!}{\rho^q}
\sup_{w\in K_1}
\left|
\tr\!\left[
\Theta_n
\left\{
R_n(w)
-
R_n^{\operatorname{det}}(w)
\right\}
\right]
\right|.
\]
Applying Theorem~\ref{Lemma:GreenRom} to the compact set $K_1$ gives
\[
\sup_{w\in K_1}
\left|
\tr\!\left[
\Theta_n
\left\{
R_n(w)
-
R_n^{\operatorname{det}}(w)
\right\}
\right]
\right|
=
O_{\mathbb P}\left(\frac{\|\Theta_n\|_F}{\sqrt n}\right).
\]
Since $q$ and $\rho$ are fixed, the deterministic factor
$q!\rho^{-q}$ can be absorbed into the $O_{\mathbb P}$-constant.
Therefore,
\[
\sup_{\zeta\in K}
\left|
\frac{d^q}{d\zeta^q}
\tr\!\left[
\Theta_n
\left\{
R_n(\zeta)
-
R_n^{\operatorname{det}}(\zeta)
\right\}
\right]
\right|
=
O_{\mathbb P}\left(\frac{\|\Theta_n\|_F}{\sqrt n}\right),
\]
which proves the assertion.
\end{proof}

\subsection{Proof of Proposition~\ref{prop:RM-recovered}}

\begin{proof}[Proof of Proposition~\ref{prop:RM-recovered}]
Fix $\zeta\in\mathbb C\setminus\mathbb R_{\geq0}$. We first show that
\begin{align}\label{eq:RM-centered-as}
\tr\{\Theta_nR_n(\zeta)\}
-
\E\left(\tr\{\Theta_nR_n(\zeta)\}\right)
\longrightarrow0
\qquad\text{almost surely}.
\end{align}
Let
\[
        \mathcal F_j=\sigma(x_1,\ldots,x_j),
        \qquad
        \E_j(\,\cdot\,)=\E(\,\cdot\,\mid\mathcal F_j),
        \qquad
        j=0,\ldots,n,
\]
where $\mathcal F_0$ is the trivial sigma-field. Then
\[
\tr\{\Theta_nR_n(\zeta)\}
-
\E\tr\{\Theta_nR_n(\zeta)\}
=
\sum_{j=1}^n
(\E_j-\E_{j-1})
\tr\{\Theta_nR_n(\zeta)\}.
\]
Since $R_{n,-j}(\zeta)$ does not depend on $x_j$, independence of the
rows implies
\[
(\E_j-\E_{j-1})
\tr\{\Theta_nR_{n,-j}(\zeta)\}
=
0.
\]
Consequently,
\[
\tr\{\Theta_nR_n(\zeta)\}
-
\E\tr\{\Theta_nR_n(\zeta)\}
=
\sum_{j=1}^n
(\E_j-\E_{j-1})
\tr\!\left[
\Theta_n\{R_n(\zeta)-R_{n,-j}(\zeta)\}
\right].
\]
By the Sherman--Morrison identity,
\[
\tr\!\left[
\Theta_n\{R_n(\zeta)-R_{n,-j}(\zeta)\}
\right]
=
-\frac{\beta_j(\zeta)}{n}
x_j^\top
R_{n,-j}(\zeta)\Theta_nR_{n,-j}(\zeta)x_j.
\]
We next establish a fourth-moment bound for the quadratic form on the
right-hand side. Conditionally on $\mathcal F_{n,-j}$, the matrix $\Sigma_n^{1/2}R_{n,-j}(\zeta)\Theta_n R_{n,-j}(\zeta)\Sigma_n^{1/2}$
is deterministic and independent of $z_j$. Hence, using
$x_j=\Sigma_n^{1/2}z_j$ and the inequality
$|a+b|^4\leq 8(|a|^4+|b|^4)$, we obtain
\begin{align*}
&\E\left[
\left|
x_j^\top
R_{n,-j}(\zeta)\Theta_nR_{n,-j}(\zeta)x_j
\right|^4
\,\middle|\,
\mathcal F_{n,-j}
\right]
\\
&\qquad\leq
8\E\left[
\left|
z_j^\top
\Sigma_n^{1/2}R_{n,-j}(\zeta)\Theta_n
R_{n,-j}(\zeta)\Sigma_n^{1/2}z_j
-
\tr\!\left\{
\Sigma_nR_{n,-j}(\zeta)\Theta_nR_{n,-j}(\zeta)
\right\}
\right|^4
\,\middle|\,
\mathcal F_{n,-j}
\right]
\\
&\qquad\quad+
8\left|
\tr\!\left\{
\Sigma_nR_{n,-j}(\zeta)\Theta_nR_{n,-j}(\zeta)
\right\}
\right|^4.
\end{align*}
By Lemma~\ref{lem_b26_bai_silverstein}, applied conditionally on
$\mathcal F_{n,-j}$ with $q=4$, the first term is bounded by
\[
C
\left\|
\Sigma_n^{1/2}R_{n,-j}(\zeta)\Theta_n
R_{n,-j}(\zeta)\Sigma_n^{1/2}
\right\|_F^4.
\]
Using Assumption~\ref{A2} and the resolvent bounds, we have
\[
\left\|
\Sigma_n^{1/2}R_{n,-j}(\zeta)\Theta_n
R_{n,-j}(\zeta)\Sigma_n^{1/2}
\right\|_F
\leq
C_\zeta\|\Theta_n\|_F
\leq
C_\zeta\|\Theta_n\|_{1}.
\]
Moreover, by cyclicity of the trace and the nuclear--operator norm
inequality,
\begin{align*}
\left|
\tr\!\left\{
\Sigma_nR_{n,-j}(\zeta)\Theta_nR_{n,-j}(\zeta)
\right\}
\right|
&=
\left|
\tr\!\left\{
\Theta_nR_{n,-j}(\zeta)\Sigma_nR_{n,-j}(\zeta)
\right\}
\right|
\\
&\leq
\|\Theta_n\|_{1}
\|R_{n,-j}(\zeta)\Sigma_nR_{n,-j}(\zeta)\|_2
\\
&\leq
C_\zeta\|\Theta_n\|_{1}.
\end{align*}
Since $\sup_{n\geq1}\|\Theta_n\|_{1}<\infty$, it follows that
\[
\E\left[
\left|
x_j^\top
R_{n,-j}(\zeta)\Theta_nR_{n,-j}(\zeta)x_j
\right|^4
\,\middle|\,
\mathcal F_{n,-j}
\right]
\leq
C_\zeta.
\]
Finally, using the uniform bound on $\beta_j(\zeta)$ from
\eqref{eq:tilde-beta-bound}, we conclude that
\[
\E\left(\left|
\tr\!\left[
\Theta_n\{R_n(\zeta)-R_{n,-j}(\zeta)\}
\right]
\right|^4\right)
\leq
\frac{C_\zeta}{n^4}.
\]
Applying the Burkholder inequality
\citep[Lemma~2.12]{Silversteinbook2010} to the preceding martingale
difference decomposition, together with conditional Jensen's
inequality, yields
\[
\E\left(\left|
\tr\{\Theta_nR_n(\zeta)\}
-
\E\tr\{\Theta_nR_n(\zeta)\}
\right|^4\right)
\leq
\frac{C_\zeta}{n^2}.
\]
Therefore, for every $\varepsilon>0$,
\[
\sum_{n=1}^\infty
\mathbb P\left(
\left|
\tr\{\Theta_nR_n(\zeta)\}
-
\E\tr\{\Theta_nR_n(\zeta)\}
\right|>\varepsilon
\right)
\leq
\frac{C_\zeta}{\varepsilon^4}
\sum_{n=1}^\infty\frac1{n^2}
<\infty.
\]
The Borel--Cantelli lemma proves \eqref{eq:RM-centered-as}.

We now use Theorem~\ref{Lemma:GreenRom}. Since
$\sup_n\|\Theta_n\|_{1}<\infty$ and
$\|\Theta_n\|_F\leq\|\Theta_n\|_{1}$, that theorem implies
\[
\tr\!\left[
\Theta_n
\{R_n(\zeta)-R_n^{\operatorname{det}}(\zeta)\}
\right]
\overset{\mathbb{P}}{\longrightarrow}0.
\]
Therefore, together with \eqref{eq:RM-centered-as},
\begin{align*}
&\E\left(\tr\{\Theta_nR_n(\zeta)\}\right)
-
\tr\{\Theta_nR_n^{\operatorname{det}}(\zeta)\}
\\
&\qquad=
\tr\!\left[
\Theta_n
\{R_n(\zeta)-R_n^{\operatorname{det}}(\zeta)\}
\right]
-
\left[
\tr\{\Theta_nR_n(\zeta)\}
-
\E\left(\tr\{\Theta_nR_n(\zeta)\}\right)
\right]
\overset{\mathbb{P}}{\longrightarrow}0.
\end{align*}
The quantity on the left-hand side is deterministic.
Consequently, convergence in probability is equivalent to ordinary
convergence, and hence
\[
\E\left(\tr\{\Theta_nR_n(\zeta)\}\right)
-
\tr\{\Theta_nR_n^{\operatorname{det}}(\zeta)\}
\longrightarrow0.
\]
Finally, writing
\begin{align*}
&\tr\!\left[
\Theta_n
\{R_n(\zeta)-R_n^{\operatorname{det}}(\zeta)\}
\right]
\\
&\qquad=
\left[
\tr\{\Theta_nR_n(\zeta)\}
-
\E\tr\{\Theta_nR_n(\zeta)\}
\right]
+
\left[
\E\tr\{\Theta_nR_n(\zeta)\}
-
\tr\{\Theta_nR_n^{\operatorname{det}}(\zeta)\}
\right],
\end{align*}
the first term on the right-hand side converges to zero almost surely by
\eqref{eq:RM-centered-as}, while the second term converges to zero
deterministically. Thus,
\begin{align}\label{eq:RM-finite-n-as}
\tr\!\left[
\Theta_n
\{R_n(\zeta)-R_n^{\operatorname{det}}(\zeta)\}
\right]
\longrightarrow0
\qquad\text{almost surely}.
\end{align}

It remains to extend \eqref{eq:RM-finite-n-as} to a common event of
probability one for all
$\zeta\in\mathbb C\setminus\mathbb R_{\geq0}$. Set $X\coloneqq\mathbb C\setminus\mathbb R_{\geq0}$
and, for $m\geq1$, define
\[
K_m
\coloneqq
\left\{
\zeta\in\mathbb C:
|\zeta|\leq m,\ 
\operatorname{dist}(\zeta,\mathbb R_{\geq0})\geq m^{-1}
\right\}.
\]
Then each $K_m$ is compact and 
$X=\bigcup_{m\geq1}K_m$.

Fix $m\geq1$ and consider the restrictions to $K_m$ of
\[
f_n(\omega,\zeta)
\coloneqq
\tr\!\left[
\Theta_n
\{R_n(\zeta)-R_n^{\operatorname{det}}(\zeta)\}
\right].
\]
For every fixed $\zeta\in K_m$, the map
$\omega\mapsto f_n(\omega,\zeta)$ is measurable and
\eqref{eq:RM-finite-n-as} yields
\[
f_n(\cdot,\zeta)\longrightarrow0
\qquad\text{almost surely}.
\]
Moreover, the resolvent structure, together with
the bounds \eqref{eq:general-resolvent-bound} and
\eqref{eq:deterministic-resolvent-bound},
Assumption~\ref{A2}, and
$\sup_n\|\Theta_n\|_{1}<\infty$,
implies that, on an event of probability one,
the family $\{f_n(\omega,\cdot)\}_{n\geq1}$
is equicontinuous on $K_m$.
Lemma~\ref{lem:Georg_global_nullset}, applied with
$X=K_m$ and $Y=\mathbb C$, therefore yields an event
$\widetilde C_m$ of probability one such that, for every
$\omega\in\widetilde C_m$, the sequence $f_n(\omega,\zeta)$ converges
for every $\zeta\in K_m$, and the resulting limit is continuous on
$K_m$. Let $D_m\subset K_m$ be countable and dense, and, for each
$\zeta\in D_m$, let $A_\zeta$ be an event of probability one on which
$f_n(\cdot,\zeta)\to0$. Then $C_m
\coloneqq
\widetilde C_m\cap\bigcap_{\zeta\in D_m}A_\zeta$ has probability one. For every $\omega\in C_m$, the limiting function
is continuous on $K_m$ and vanishes on the dense subset $D_m$.
Consequently,
\[
f_n(\omega,\zeta)\longrightarrow0
\qquad\text{for every }\zeta\in K_m.
\]

Finally, the event $C\coloneqq\bigcap_{m\geq1}C_m$ has probability one. Since $X=\bigcup_{m\geq1}K_m$, it follows that,
for every $\omega\in C$,
\[
\tr\!\left[
\Theta_n
\{R_n(\zeta)-R_n^{\operatorname{det}}(\zeta)\}
\right]
\longrightarrow0
\qquad\text{for every }
\zeta\in\mathbb C\setminus\mathbb R_{\geq0}.
\]
\end{proof}

\subsection{Proof of Corollary~\ref{cor:GreenRom-derivatives}}

\begin{proof}[Proof of Corollary~\ref{cor:GreenRom-derivatives}]
For $q=0$, the assertion follows directly from
Proposition~\ref{prop:RM-recovered}. We therefore suppose that $q\geq1$.

For each $n\in\mathbb N$, define
\[
        f_n(w)
        \coloneqq
        \tr\!\left[
        \Theta_n
        \left\{
        R_n(w)
        -
        R_n^{\operatorname{det}}(w)
        \right\}
        \right],
        \qquad
        w\in\mathbb C\setminus\mathbb R_{\geq0}.
\]
As in the proof of Corollary~\ref{cor:GreenRom-derivatives:rates},
$f_n$ is holomorphic on
$\mathbb C\setminus\mathbb R_{\geq0}$. By Proposition~\ref{prop:RM-recovered}, there exists an event
$\Omega_0$ with $\mathbb P(\Omega_0)=1$ such that, for every
$\omega\in\Omega_0$,
\[
        f_n(w,\omega)\longrightarrow0
        \qquad\text{for all }
        w\in\mathbb C\setminus\mathbb R_{\geq0}.
\]

Fix $\omega\in\Omega_0$ and
$\zeta\in\mathbb C\setminus\mathbb R_{\geq0}$. Choose
$0<\rho<\operatorname{dist}(\zeta,\mathbb R_{\geq0})$ and set
\[
        K_1
        \coloneqq
        \{w\in\mathbb C:|w-\zeta|\leq\rho\}.
\]
Then $K_1 \subseteq\mathbb C\setminus\mathbb R_{\geq0}$.
Hence, Cauchy's integral formula for derivatives
\cite[Chapter~10, Exercise~7]{Rudin1987} gives
\[
        \frac{d^q}{d\zeta^q}f_n(\zeta,\omega)
        =
        \frac{q!}{2\pi i}
        \int_{|w-\zeta|=\rho}
        \frac{f_n(w,\omega)}{(w-\zeta)^{q+1}}\,dw.
\]

We next obtain a bound on the integrand which is uniform in $n$. Since
$K_1$ is compact and bounded away from $\mathbb R_{\geq0}$,
\eqref{eq:general-resolvent-bound} and
\eqref{eq:deterministic-resolvent-bound} give
\[
        \sup_{n\geq1}\sup_{w\in K_1}
        \left\{
        \|R_n(w)\|_2
        +
        \|R_n^{\operatorname{det}}(w)\|_2
        \right\}
        \leq C_{K_1}.
\]
Together with $\sup_n\|\Theta_n\|_{1}<\infty$, this yields
\[
        \sup_{n\geq1}\sup_{w\in K_1}
        |f_n(w,\omega)|
        \leq C_{K_1}.
\]
Therefore, for $|w-\zeta|=\rho$,
\[
\left|
        \frac{f_n(w,\omega)}
        {(w-\zeta)^{q+1}}
\right|
\leq
\frac{C_{K_1}}{\rho^{q+1}}.
\]
On the other hand, since $\omega\in\Omega_0$,
$f_n(w,\omega)\longrightarrow0$ for every $w$ on the circle
$|w-\zeta|=\rho$. Hence, by dominated convergence,
\[
\begin{split}
        \frac{d^q}{d\zeta^q}f_n(\zeta,\omega)
        &=
        \frac{q!}{2\pi i}
        \int_{|w-\zeta|=\rho}
        \frac{f_n(w,\omega)}{(w-\zeta)^{q+1}}\,dw
        \longrightarrow0.
\end{split}
\]
Since $\zeta\in\mathbb C\setminus\mathbb R_{\geq0}$ was arbitrary and
$\Omega_0$ does not depend on $\zeta$, we conclude that, almost surely,
for all $\zeta\in\mathbb C\setminus\mathbb R_{\geq0}$,
\[
\frac{d^q}{d\zeta^q}
\tr\!\left[
\Theta_n
\left\{
R_n(\zeta)
-
R_n^{\operatorname{det}}(\zeta)
\right\}
\right]
\longrightarrow0.
\]
This proves the assertion.
\end{proof}

\subsection{Proof of Lemma~\ref{lem:companion-uniq-D}} \label{sec:proof:lem}

\begin{proof}[Proof of Lemma~\ref{lem:companion-uniq-D}]
We first establish existence. %
By Lemma~\ref{thm:MP}, with
$(\gamma,H)=(\gamma_0,H_0)$ and $\zeta\in\mathbb C^+$,
$v_{\gamma_0,H_0}(\zeta)$ satisfies
\[
        v_{\gamma_0,H_0}(\zeta)
        =
        -\left(
        \zeta
        -
        \gamma_0\int_0^\infty
        \frac{t}{1+t v_{\gamma_0,H_0}(\zeta)}\,dH_0(t)
        \right)^{-1}.
\]
We show that the fixed-point equation continues to hold on
$\mathbb C\setminus\mathbb R_{\ge0}$.

Since \(\bar F_{\gamma_0,H_0}\) is supported on \([0,\infty)\),
\(v_{\gamma_0,H_0}\) is the Stieltjes transform of a probability measure
supported on \(\mathbb R_{\ge0}\). Hence
\(v_{\gamma_0,H_0}\) is holomorphic on
\(\mathbb C\setminus\mathbb R_{\ge0}\) by
\citet[Proposition~2.2]{Hachem:Loubaton:Najim:2007}. Next, we claim that

\[
        1+t v_{\gamma_0,H_0}(\zeta)\neq0
\]
for every $\zeta\in\mathbb C\setminus\mathbb R_{\ge0}$ and every
$t\in[0,C_0]$. The assertion is immediate for $t=0$. If $t>0$ and
$1+t v_{\gamma_0,H_0}(\zeta)=0$, then
$v_{\gamma_0,H_0}(\zeta)=-(1/t)\in\mathbb R_{<0}$. However, this is
impossible. Indeed, if $\zeta=a+ib\in\mathbb C\setminus\mathbb R$, then
\[
        \Im\{v_{\gamma_0,H_0}(\zeta)\}
        =
        b\int_0^\infty
        \frac1{(s-a)^2+b^2}\,d\bar F_{\gamma_0,H_0}(s)\neq 0.
\] 
If
$\zeta\in\mathbb R_{<0}$, then $s-\zeta>0$ for every $s\ge0$, and therefore
$v_{\gamma_0,H_0}(\zeta)>0$. Thus
$v_{\gamma_0,H_0}(\zeta)$ cannot be a negative real number, and the claim
follows.

We now verify that
\begin{align}\label{eq:hol:check}
        h(\zeta)
        \coloneqq
        \int_0^\infty
        \frac{t}{1+t v_{\gamma_0,H_0}(\zeta)}\,dH_0(t)
\end{align}
is holomorphic on $\mathbb C\setminus\mathbb R_{\ge0}$. Let
$\zeta_0\in\mathbb C\setminus\mathbb R_{\ge0}$ be arbitrary and choose
$r>0$ such that
$\overline{B(\zeta_0,r)}\subset\mathbb C\setminus\mathbb R_{\ge0}$.
As shown above, the function $(\zeta,t)\mapsto |1+t v_{\gamma_0,H_0}(\zeta)|$
is strictly positive on $\overline{B(\zeta_0,r)}\times[0,C_0]$. Moreover,
since $v_{\gamma_0,H_0}$ is continuous on
$\mathbb C\setminus\mathbb R_{\ge0}$, this map is continuous on the compact
set $\overline{B(\zeta_0,r)}\times[0,C_0]$. Hence, there exists
$c_{\zeta_0,r}>0$ such that
\[
        \inf_{\zeta\in\overline{B(\zeta_0,r)}}\inf_{t\in[0,C_0]}
        |1+t v_{\gamma_0,H_0}(\zeta)|
        \ge c_{\zeta_0,r}.
\]
Consequently, the integrand in \eqref{eq:hol:check} is uniformly bounded on
$\overline{B(\zeta_0,r)}\times[0,C_0]$, and the dominated convergence
theorem implies that $h$ is continuous on $B(\zeta_0,r)$.

Moreover, for each fixed $t\in[0,C_0]$, the map
\begin{align}\label{eq:hol:map}
        \zeta\mapsto
        \frac{t}{1+t v_{\gamma_0,H_0}(\zeta)}
\end{align}
is holomorphic on $B(\zeta_0,r)$, since $v_{\gamma_0,H_0}$ is holomorphic on
$B(\zeta_0,r)$ and $1+t v_{\gamma_0,H_0}(\zeta)\neq0$ there.

Let $\Delta$ be any closed triangle contained in $B(\zeta_0,r)$. We show
that
\[
        \int_{\partial\Delta}h(\zeta)\,d\zeta=0.
\]
Since the integrand in \eqref{eq:hol:check} is uniformly bounded on
$\overline{B(\zeta_0,r)}\times[0,C_0]$, Fubini's theorem gives
\begin{align*}
        \int_{\partial\Delta}
        \left(
        \int_0^\infty
        \frac{t}{1+t v_{\gamma_0,H_0}(\zeta)}\,dH_0(t)
        \right)\,d\zeta
        &=
        \int_0^\infty
        \left(
        \int_{\partial\Delta}
        \frac{t}{1+t v_{\gamma_0,H_0}(\zeta)}\,d\zeta
        \right)\,dH_0(t).
\end{align*}
Since, for each fixed $t\in[0,C_0]$, the map in \eqref{eq:hol:map} is
holomorphic on $B(\zeta_0,r)$, the Cauchy--Goursat theorem
\citep[Theorem~10.13]{Rudin1987} gives
\[
        \int_{\partial\Delta}
        \frac{t}{1+t v_{\gamma_0,H_0}(\zeta)}\,d\zeta
        =
        0,
        \qquad t\in[0,C_0].
\]
Combining the preceding two displays, we obtain
\[
        \int_{\partial\Delta}h(\zeta)\,d\zeta=0
\]
for every closed triangle $\Delta\subset B(\zeta_0,r)$. Since $h$ is
continuous on $B(\zeta_0,r)$, Morera's theorem
\citep[Theorem~10.17]{Rudin1987} implies that $h$ is holomorphic on
$B(\zeta_0,r)$. Since $\zeta_0$ was arbitrary, $h$ is holomorphic on
$\mathbb C\setminus\mathbb R_{\ge0}$.

Define
\begin{align*}
        G(\zeta)
        &\coloneqq
        v_{\gamma_0,H_0}(\zeta)
        \left(
        \zeta
        -
        \gamma_0\int_0^\infty
        \frac{t}{1+t v_{\gamma_0,H_0}(\zeta)}\,dH_0(t)
        \right)
        +1,
        \quad\text{and}\\[2mm]
        Z(G)
        &\coloneqq
        \{\zeta\in\mathbb C\setminus\mathbb R_{\ge0}:G(\zeta)=0\}.
\end{align*}
Combining the preceding arguments, we conclude that $G$ is holomorphic on
$\mathbb C\setminus\mathbb R_{\ge0}$. Moreover,
$\mathbb C\setminus\mathbb R_{\ge0}$ is open and connected. Hence, by
\citet[Theorem~10.18]{Rudin1987}, either
\[
        Z(G)=\mathbb C\setminus\mathbb R_{\ge0},
\]
or $Z(G)$ has no limit point in $\mathbb C\setminus\mathbb R_{\ge0}$.
However, Lemma~\ref{thm:MP} implies that $G(\zeta)=0$ for
$\zeta\in\mathbb C^+$. Therefore, $\mathbb C^+\subset Z(G)$, and since
$\mathbb C^+$ is a nonempty open subset of
$\mathbb C\setminus\mathbb R_{\ge0}$, the set $Z(G)$ has a limit point in
$\mathbb C\setminus\mathbb R_{\ge0}$. Consequently,
\[
        Z(G)=\mathbb C\setminus\mathbb R_{\ge0}.
\]
Equivalently,
\[
        v_{\gamma_0,H_0}(\zeta)
        =
        -\left(
        \zeta
        -
        \gamma_0\int_0^\infty
        \frac{t}{1+t v_{\gamma_0,H_0}(\zeta)}\,dH_0(t)
        \right)^{-1}
\]
for every $\zeta\in\mathbb C\setminus\mathbb R_{\ge0}$. Thus
$v_{\gamma_0,H_0}(\zeta)$ is an admissible solution for every
$\zeta\in\mathbb C\setminus\mathbb R_{\ge0}$.

It remains to prove uniqueness. Let first $\zeta\in\mathbb C^+$. Then
uniqueness of the admissible solution follows from Lemma~\ref{thm:MP} with
the pair $(\gamma_0,H_0)$.

Next, let $\zeta\in\mathbb C^-$. Let $w_1,w_2\in\mathbb C^-$ be two
admissible solutions. Thus, for $j=1,2$,
\[
        w_j
        =
        -\left(
        \zeta
        -
        \gamma_0\int_0^\infty
        \frac{t}{1+tw_j}\,dH_0(t)
        \right)^{-1}.
\]
Taking complex conjugates, we get
\[
        \overline w_j
        =
        -\left(
        \overline\zeta
        -
        \gamma_0\int_0^\infty
        \frac{t}{1+t\overline w_j}\,dH_0(t)
        \right)^{-1},
        \qquad j=1,2.
\]
Since $\overline\zeta\in\mathbb C^+$ and $\overline w_i\in\mathbb C^+$, the
uniqueness result on $\mathbb C^+$ implies
$\overline w_1=\overline w_2$. Hence $w_1=w_2$.

It remains to treat the case $\zeta\in\mathbb R_{<0}$. Write
$\zeta=-\lambda$ with $\lambda>0$, and let $w_1,w_2>0$ be two admissible
solutions. The fixed-point equation is equivalent to
\[
        \frac1w
        =
        \lambda
        +
        \gamma_0\int_0^\infty
        \frac{t}{1+tw}\,dH_0(t).
\]
Therefore, for $j=1,2$,
\begin{align}\label{eq:ridge:type}
        \frac1{w_j}
        =
        \lambda
        +
        \gamma_0\int_0^\infty
        \frac{t}{1+tw_j}\,dH_0(t).
\end{align}
Subtracting the equations for $w_1$ and $w_2$ gives
\[
        \frac1{w_1}-\frac1{w_2}
        =
        \gamma_0\int_0^\infty
        \left\{
        \frac{t}{1+tw_1}
        -
        \frac{t}{1+tw_2}
        \right\}\,dH_0(t).
\]
Since
\[
        \frac{t}{1+tw_1}
        -
        \frac{t}{1+tw_2}
        =
        \frac{t^2(w_2-w_1)}{(1+tw_1)(1+tw_2)}
\]
and
\[
        \frac1{w_1}-\frac1{w_2}
        =
        \frac{w_2-w_1}{w_1w_2},
\]
we obtain
\[
        (w_2-w_1)
        \left\{
        \frac1{w_1w_2}
        -
        \gamma_0\int_0^\infty
        \frac{t^2}{(1+tw_1)(1+tw_2)}\,dH_0(t)
        \right\}
        =
        0.
\]
We show that the term in braces is strictly positive. Multiplying \eqref{eq:ridge:type} for $j=1$ by $w_1$ yields
\[
        1
        =
        \lambda w_1
        +
        \gamma_0\int_0^\infty
        \frac{tw_1}{1+tw_1}\,dH_0(t).
\]
Since $\lambda w_1>0$, it follows that
\[
        \gamma_0\int_0^\infty
        \frac{tw_1}{1+tw_1}\,dH_0(t)<1.
\]
Moreover, for $t\ge0$ and $w_2>0$,
\[
        0\le \frac{tw_2}{1+tw_2}\le1.
\]
Consequently,
\begin{align*}
        \gamma_0 w_1w_2
        \int_0^\infty
        \frac{t^2}{(1+tw_1)(1+tw_2)}\,dH_0(t)
        &=
        \gamma_0\int_0^\infty
        \frac{tw_1}{1+tw_1}
        \frac{tw_2}{1+tw_2}\,dH_0(t)\\
        &\le
        \gamma_0\int_0^\infty
        \frac{tw_1}{1+tw_1}\,dH_0(t)
        <1.
\end{align*}
Hence,
\[
        \frac1{w_1w_2}
        -
        \gamma_0\int_0^\infty
        \frac{t^2}{(1+tw_1)(1+tw_2)}\,dH_0(t)
        >0.
\]
Thus the preceding factorization implies $w_2-w_1=0$, and therefore
$w_1=w_2$.

This proves uniqueness of the admissible solution for every
$\zeta\in\mathbb C\setminus\mathbb R_{\ge0}$. Together with the existence
established above, this completes the proof.
\end{proof}

\begin{remark}
For each $n\in\mathbb N$, we have $\gamma_n\in(0,\infty)$, while $H_n$ is a
nonrandom probability measure supported on $[0,C]$ by Assumption~\ref{A2}. Therefore,
Lemma~\ref{lem:companion-uniq-D}, applied with $(\gamma_0,H_0)=(\gamma_n,H_n)$,
guarantees that $v_n$ is well-defined as the unique admissible solution
of \eqref{stielt} with $(\gamma_0,H_0)=(\gamma_n,H_n)$.
Similarly, choosing $(\gamma,H)$ as in Assumptions~\ref{A1} and~\ref{A4},
we may apply Lemma~\ref{lem:companion-uniq-D} with
\((\gamma_0,H_0)=(\gamma,H)\). Hence, the implicit representation on $\mathbb{C}\setminus\mathbb{R}_{\geq 0}$ also holds
for the limiting pair.
\end{remark}

\section{Proofs of the results in Section~\ref{sec:two-resolvent}}\label{sec:two-resolvent-proofs}

\subsection{Proof of Lemma~\ref{lem:ridge-signal-equivalence:rates}}

\begin{proof}[Proof of Lemma~\ref{lem:ridge-signal-equivalence:rates}]
Fix a compact set
$K\subseteq\mathbb C\setminus\mathbb R_{\geq0}$, and choose a compact set
$K_1\subseteq\mathbb C\setminus\mathbb R_{\geq0}$ as in
Lemma~\ref{lem:cauchy_reduction_random}. 
By \eqref{eq:general-denominator-bound}, \eqref{eq:general-resolvent-bound}, the definition of $\beta_j$ in \eqref{eq:beta_jandq_j}, and the compactness of $K_1$,
\begin{align}\label{eq:signal-uniform-scalar-bounds}
\sup_{w\in K_1}\max_{1\leq j\leq n}|\beta_j(w)|
+
\sup_{w\in K_1}|\{-wv_n(w)\}|
+
\sup_{w\in K_1}|\{-w\hat v_n(w)\}|
\leq C_{K_1}.
\end{align}
For $j=1,\ldots,n$ and $\zeta,\xi\in K_1$, define
\[
\begin{split}
        Q_j(\zeta,\xi)
        &\coloneqq
        x_j^\top
        R_{n,-j}(\xi)\Theta_nR_{n,-j}(\zeta)x_j
        -
        T_j(\zeta,\xi),\\
        T_j(\zeta,\xi)
        &\coloneqq
        \tr\!\left[
        \Theta_nR_{n,-j}(\zeta)\Sigma_nR_{n,-j}(\xi)
        \right],
\end{split}
\]
and
\[
        T(\zeta,\xi)
        \coloneqq
        \tr\!\left[
        \Theta_nR_n(\zeta)\Sigma_nR_n(\xi)
        \right].
\]
Note that
\[
        T_j(\zeta,\xi)
        =
        \tr\!\left[
        \Sigma_nR_{n,-j}(\xi)\Theta_nR_{n,-j}(\zeta)
        \right],
\]
and hence $Q_j(\zeta,\xi)$ is conditionally centered given
$\mathcal F_{n,-j}$.
Since
\[
        R_n(w)x_j
        =
        \beta_j(w)R_{n,-j}(w)x_j,\qquad\text{and}\qquad \hat\Sigma_n=\frac1n\sum_{j=1}^n x_jx_j^\top,
\]
we have
\begin{align}\label{eq:signal-exact-empirical}
\begin{split}
\tr\!\left[
\Theta_nR_n(\zeta)\hat\Sigma_nR_n(\xi)
\right]
&=
\frac1n\sum_{j=1}^n
x_j^\top R_n(\xi)\Theta_nR_n(\zeta)x_j\\
&=
\frac1n\sum_{j=1}^n
\beta_j(\zeta)\beta_j(\xi)
\left\{
Q_j(\zeta,\xi)+T_j(\zeta,\xi)
\right\}.
\end{split}
\end{align}

We first collect the estimates needed to analyze the right-hand side of
\eqref{eq:signal-exact-empirical}. By
\eqref{beta:emp:fluct} and Lemma~\ref{lem:uniform-companion-rate},
applied to $K_1$,
\begin{align}\label{eq:beta_minus_v_n}
\begin{split}
&\sup_{w\in K_1}
\left\{
\frac1n\sum_{j=1}^n
\left|
\beta_j(w)-\{-wv_n(w)\}
\right|^2
\right\}^{1/2}\\
&\qquad\leq
\sup_{w\in K_1}
\left\{
\frac1n\sum_{j=1}^n
\left|
\beta_j(w)-\{-w\hat v_n(w)\}
\right|^2
\right\}^{1/2}
+
\sup_{w\in K_1}
|w|\,|\hat v_n(w)-v_n(w)|\\
&\qquad=
O_{\mathbb P}(n^{-1/2}).
\end{split}
\end{align}
Define $d_j(\zeta,\xi)
        \coloneqq
        \beta_j(\zeta)\beta_j(\xi)
        -
        \{-\zeta v_n(\zeta)\}\{-\xi v_n(\xi)\}$.
Then,
\[
d_j(\zeta,\xi)
=
\beta_j(\zeta)
\left[
\beta_j(\xi)-\{-\xi v_n(\xi)\}
\right]+
\{-\xi v_n(\xi)\}
\left[
\beta_j(\zeta)-\{-\zeta v_n(\zeta)\}
\right].
\]
Hence, by \eqref{eq:signal-uniform-scalar-bounds},
\[
|d_j(\zeta,\xi)|
\leq
C_{K_1}
\left\{
\left|
\beta_j(\xi)-\{-\xi v_n(\xi)\}
\right|
+
\left|
\beta_j(\zeta)-\{-\zeta v_n(\zeta)\}
\right|
\right\}.
\]
Therefore, Minkowski's inequality, yields
\[
\begin{split}
\left\{
\frac1n\sum_{j=1}^n
|d_j(\zeta,\xi)|^2
\right\}^{1/2}
&\leq
C_{K_1}
\left\{
\frac1n\sum_{j=1}^n
\left|
\beta_j(\xi)-\{-\xi v_n(\xi)\}
\right|^2
\right\}^{1/2}\\
&\quad+
C_{K_1}
\left\{
\frac1n\sum_{j=1}^n
\left|
\beta_j(\zeta)-\{-\zeta v_n(\zeta)\}
\right|^2
\right\}^{1/2}.
\end{split}
\]
Taking the supremum over $\zeta,\xi\in K_1$ and using \eqref{eq:beta_minus_v_n} yields
\begin{align}\label{eq:signal-beta-product-square}
\sup_{\zeta,\xi\in K_1}
\left\{
\frac1n\sum_{j=1}^n
|d_j(\zeta,\xi)|^2
\right\}^{1/2}
=
O_{\mathbb P}(n^{-1/2}).
\end{align}

We next prove an upper bound for
$|n^{-1}\sum_{j=1}^nd_j(\zeta,\xi)|$. Using \eqref{eq:beta_j_stielt},
we have
\[
\begin{split}
\frac1n\sum_{j=1}^n d_j(\zeta,\xi)
&=
\frac1n\sum_{j=1}^n
\left[
\beta_j(\zeta)-\{-\zeta\hat v_n(\zeta)\}
\right]
\left[
\beta_j(\xi)-\{-\xi\hat v_n(\xi)\}
\right]\\
&\qquad+
\{-\zeta\hat v_n(\zeta)\}\{-\xi\hat v_n(\xi)\}
-
\{-\zeta v_n(\zeta)\}\{-\xi v_n(\xi)\}.
\end{split}
\]
The first term is $O_{\mathbb P}(n^{-1})$ uniformly on $K_1\times K_1$
by the Cauchy--Schwarz inequality and \eqref{beta:emp:fluct}. The second
term has the same order by Lemma~\ref{lem:uniform-companion-rate} and
\eqref{eq:signal-uniform-scalar-bounds}. Consequently,
\begin{align}\label{eq:signal-beta-product-mean}
\sup_{\zeta,\xi\in K_1}
\left|
\frac1n\sum_{j=1}^n d_j(\zeta,\xi)
\right|
=
O_{\mathbb P}(n^{-1}).
\end{align}

Next, by Lemma~\ref{lem:ridge_sandwich}, we have 
\begin{align}\label{eq:signal-Q-average}
\sup_{\zeta,\xi\in K}
\left|
\frac1n\sum_{j=1}^nQ_j(\zeta,\xi)
\right|
=
O_{\mathbb P}\left(
\|\Theta_n\|_F
\left\{
\frac1{\sqrt n}\vee\frac{\sqrt p}{n}
\right\}
\right).
\end{align}
We shall also use
\begin{align}\label{eq:signal-Q-square}
\sup_{\zeta,\xi\in K}
\left\{
\frac1n\sum_{j=1}^n|Q_j(\zeta,\xi)|^2
\right\}^{1/2}
=
O_{\mathbb P}(\|\Theta_n\|_F).
\end{align}
To see this, conditionally on $\mathcal F_{n,-j}$, write
\[
        Q_j(\zeta,\xi)
        =
        z_j^\top M_{n,j}(\zeta,\xi)z_j
        -
        \tr\{M_{n,j}(\zeta,\xi)\},
\]
where
\[
        M_{n,j}(\zeta,\xi)
        =
        \Sigma_n^{1/2}
        R_{n,-j}(\zeta)\Theta_n^\top
        R_{n,-j}(\xi)\Sigma_n^{1/2}.
\]
By Assumption~\ref{A2} and the uniform resolvent bounds on $K_1$,
\[
        \sup_{\zeta,\xi\in K_1}
        \|M_{n,j}(\zeta,\xi)\|_F
        \leq
        C_{K_1}\|\Theta_n\|_F.
\]
Lemma~\ref{lem_b26_bai_silverstein}, applied conditionally on
$\mathcal F_{n,-j}$, therefore gives
\[
        \sup_{\zeta,\xi\in K_1}
        \E\bigl(|Q_j(\zeta,\xi)|^2\bigr)
        \leq
        C_{K_1}\|\Theta_n\|_F^2.
\]
The two-dimensional part of
Lemma~\ref{lem:cauchy_reduction_random}, followed by Markov's
inequality, proves \eqref{eq:signal-Q-square}.

We finally compare $T_j(\zeta,\xi)$ with $T(\zeta,\xi)$. The
Sherman--Morrison identity gives
\[
        R_{n,-j}(w)-R_n(w)
        =
        \frac{\beta_j(w)}n
        R_{n,-j}(w)x_jx_j^\top R_{n,-j}(w).
\]
Expanding both resolvent factors yields
\begin{align}\label{eq:signal-T-difference}
\begin{split}
T_j(\zeta,\xi)-T(\zeta,\xi)
&=
\frac{\beta_j(\zeta)}n
x_j^\top
R_{n,-j}(\zeta)\Sigma_nR_{n,-j}(\xi)
\Theta_nR_{n,-j}(\zeta)x_j\\
&\quad+
\frac{\beta_j(\xi)}n
x_j^\top
R_{n,-j}(\xi)\Theta_nR_{n,-j}(\zeta)
\Sigma_nR_{n,-j}(\xi)x_j\\
&\quad-
\frac{\beta_j(\zeta)\beta_j(\xi)}{n^2}
\left\{
x_j^\top
R_{n,-j}(\zeta)\Sigma_nR_{n,-j}(\xi)x_j
\right\}
\left\{
x_j^\top
R_{n,-j}(\xi)\Theta_nR_{n,-j}(\zeta)x_j
\right\}.
\end{split}
\end{align}
The first two terms on the right-hand side satisfy by \eqref{eq:signal-uniform-scalar-bounds} and
Lemma~\ref{lem_b26_bai_silverstein},
\[
\sup_{\zeta,\xi\in K_1}
\E\left|
\frac{\beta_j(\zeta)}n
x_j^\top
R_{n,-j}(\zeta)\Sigma_nR_{n,-j}(\xi)
\Theta_nR_{n,-j}(\zeta)x_j
\right|^2
\leq
\frac{C_{K_1}p}{n^2}\|\Theta_n\|_F^2,
\]
and the same bound holds for the second term.

For the last term in \eqref{eq:signal-T-difference}, the Sherman--Morrison formula gives
\[
\frac{|\beta_j(w)|}{n}
\|R_{n,-j}(w)x_j\|_2^2
=
\|R_{n,-j}(w)-R_n(w)\|_2
\leq
\frac{2}{\operatorname{dist}(K_1,\mathbb R_{\geq0})}.
\]
Hence, by the Cauchy--Schwarz inequality,
\[
\begin{split}
&\frac{\sqrt{|\beta_j(\zeta)\beta_j(\xi)|}}n
\left|
x_j^\top
R_{n,-j}(\zeta)\Sigma_nR_{n,-j}(\xi)x_j
\right|\\
&\qquad\leq
\|\Sigma_n\|_2
\left\{
\frac{|\beta_j(\zeta)|}{n}
\|R_{n,-j}(\zeta)x_j\|_2^2
\right\}^{1/2}
\left\{
\frac{|\beta_j(\xi)|}{n}
\|R_{n,-j}(\xi)x_j\|_2^2
\right\}^{1/2}\\
&\qquad\leq
\frac{2\|\Sigma_n\|_2}
{\operatorname{dist}(K_1,\mathbb R_{\geq0})}
\leq
C_{K_1},
\end{split}
\]
where the last inequality follows from Assumption~\ref{A2}.
Consequently, the absolute value of the last term in
\eqref{eq:signal-T-difference} is bounded by,
\[
\frac{C_{K_1}\sqrt{|\beta_j(\zeta)\beta_j(\xi)|}}{n}
\left|
x_j^\top
R_{n,-j}(\xi)\Theta_nR_{n,-j}(\zeta)x_j
\right|.
\]
Using \eqref{eq:signal-uniform-scalar-bounds}, and applying
Lemma~\ref{lem_b26_bai_silverstein} conditionally on
$\mathcal F_{n,-j}$, we therefore obtain
\[
\sup_{\zeta,\xi\in K_1}
\E\left|
\frac{\beta_j(\zeta)\beta_j(\xi)}{n^2}
\left\{
x_j^\top
R_{n,-j}(\zeta)\Sigma_nR_{n,-j}(\xi)x_j
\right\}
\left\{
x_j^\top
R_{n,-j}(\xi)\Theta_nR_{n,-j}(\zeta)x_j
\right\}
\right|^2
\leq
\frac{C_{K_1}p}{n^2}\|\Theta_n\|_F^2.
\]
It follows that
\[
        \sup_{\zeta,\xi\in K_1}
        \E\bigl(
        |T_j(\zeta,\xi)-T(\zeta,\xi)|^2
        \bigr)
        \leq
        \frac{C_{K_1}p}{n^2}\|\Theta_n\|_F^2.
\]
Using the two-dimensional part of
Lemma~\ref{lem:cauchy_reduction_random} and Markov's inequality, we
obtain
\begin{align}\label{eq:signal-T-square}
\sup_{\zeta,\xi\in K}
\left\{
\frac1n\sum_{j=1}^n
|T_j(\zeta,\xi)-T(\zeta,\xi)|^2
\right\}^{1/2}
=
O_{\mathbb P}\left(
\frac{\sqrt p}{n}\|\Theta_n\|_F
\right).
\end{align}
Furthermore,
\begin{align}\label{eq:signal-T-bound}
\sup_{\zeta,\xi\in K}
|T(\zeta,\xi)|
\leq
\|\Theta_n\|_F
\sup_{\zeta,\xi\in K}
\|R_n(\zeta)\Sigma_nR_n(\xi)\|_F\leq
C_K\sqrt p\,\|\Theta_n\|_F.
\end{align}

We now combine the preceding estimates. By
\eqref{eq:signal-exact-empirical},
\begin{align}\label{eq:signal-final-decomposition}
\begin{split}
\tr\!\left[
\Theta_nR_n(\zeta)\hat\Sigma_nR_n(\xi)
\right]
&-
\{-\zeta v_n(\zeta)\}\{-\xi v_n(\xi)\}
T(\zeta,\xi)\\
&\quad=
\{-\zeta v_n(\zeta)\}\{-\xi v_n(\xi)\}
\frac1n\sum_{j=1}^nQ_j(\zeta,\xi)\\
&\qquad+
\frac1n\sum_{j=1}^n
d_j(\zeta,\xi)Q_j(\zeta,\xi)\\
&\qquad+
\{-\zeta v_n(\zeta)\}\{-\xi v_n(\xi)\}
\frac1n\sum_{j=1}^n
\{T_j(\zeta,\xi)-T(\zeta,\xi)\}\\
&\qquad+
\frac1n\sum_{j=1}^n
d_j(\zeta,\xi)
\{T_j(\zeta,\xi)-T(\zeta,\xi)\}\\
&\qquad+
T(\zeta,\xi)
\frac1n\sum_{j=1}^n d_j(\zeta,\xi).
\end{split}
\end{align}

The first term on the right-hand side of
\eqref{eq:signal-final-decomposition} has the asserted order by
\eqref{eq:signal-uniform-scalar-bounds} and
\eqref{eq:signal-Q-average}. By the Cauchy--Schwarz inequality,
\eqref{eq:signal-beta-product-square}, and
\eqref{eq:signal-Q-square}, the second term is $O_{\mathbb P}(n^{-1/2}\|\Theta_n\|_F)$.
The third term is $O_{\mathbb P}(\sqrt{p}n^{-1}\|\Theta_n\|_F)$
by \eqref{eq:signal-T-square}. The fourth term is, again by the
Cauchy--Schwarz inequality,
\[
        O_{\mathbb P}\left(
        \frac{\sqrt p}{n^{3/2}}\|\Theta_n\|_F
        \right),\qquad
\text{and the final term is}\qquad
        O_{\mathbb P}\left(
        \frac{\sqrt p}{n}\|\Theta_n\|_F
        \right)
\]
by \eqref{eq:signal-beta-product-mean} and
\eqref{eq:signal-T-bound}. Combining these estimates proves
\begin{align*}
\begin{split}
&\sup_{\zeta,\xi\in K}
\Bigg|
\{-\zeta v_n(\zeta)\}\{-\xi v_n(\xi)\}
\tr\!\left[
\Theta_nR_n(\zeta)\Sigma_nR_n(\xi)
\right]-
\tr\!\left[
\Theta_nR_n(\zeta)\hat\Sigma_nR_n(\xi)
\right]
\Bigg|\\
&\qquad=
O_{\mathbb P}\left(
\|\Theta_n\|_F
\left\{
\frac1{\sqrt n}\vee\frac{\sqrt p}{n}
\right\}
\right),
\end{split}
\end{align*}
as asserted.
\end{proof}

\subsection{Proof of Lemma~\ref{lem:ridge-signal-equivalence}}

\begin{proof}[Proof of Lemma~\ref{lem:ridge-signal-equivalence}]
Fix $\zeta,\xi\in\mathbb C\setminus\mathbb R_{\geq0}$, and define
\[
        \Delta_n(\zeta,\xi)
        \coloneqq
        \{-\zeta v_n(\zeta)\}\{-\xi v_n(\xi)\}
        \tr\!\left[
        \Theta_nR_n(\zeta)\Sigma_nR_n(\xi)
        \right]
        -
        \tr\!\left[
        \Theta_nR_n(\zeta)\hat\Sigma_nR_n(\xi)
        \right].
\]
Choose a compact set $K\subseteq\mathbb C\setminus\mathbb R_{\geq0}$
containing $\zeta$ and $\xi$. By
Lemma~\ref{lem:ridge-signal-equivalence:rates},
$\sup_n\|\Theta_n\|_F<\infty$, and \eqref{gamma_n:cond},
\begin{align}\label{eq:signal-in-probability}
        \Delta_n(\zeta,\xi)
        \overset{\mathbb P}{\longrightarrow}0.
\end{align}
We shall prove
\begin{align}\label{eq:signal-centered-as}
        \Delta_n(\zeta,\xi)-\E\{\Delta_n(\zeta,\xi)\}
        \longrightarrow0
        \qquad\text{almost surely}.
\end{align}

First, we simplify the expression whose martingale increments must
be estimated. Define the deterministic matrix
\begin{align}\label{eq:signal-reduced-matrix}
        C_n
        \coloneqq
        \{-\zeta v_n(\zeta)\}\{-\xi v_n(\xi)\}\Sigma_n
        -\zeta I_p.
\end{align}
By Assumption~\ref{A2} and
\eqref{eq:signal-uniform-scalar-bounds},
$\sup_n\|C_n\|_2\leq C_{\zeta,\xi}$.
Since $R_n(\zeta)\hat\Sigma_n=I_p+\zeta R_n(\zeta)$, we have
\begin{align}\label{eq:signal-reduced-representation}
        \Delta_n(\zeta,\xi)
        =
        \tr\!\left[\Theta_nR_n(\zeta)C_nR_n(\xi)\right]
        -\tr\!\left[\Theta_nR_n(\xi)\right].
\end{align}
The corresponding identity holds for the leave-one-out matrices.
Accordingly, define
\begin{align}\label{eq:signal-reduced-leave-one-out}
        \Delta_{n,-j}(\zeta,\xi)
        \coloneqq
        \tr\!\left[\Theta_nR_{n,-j}(\zeta)C_nR_{n,-j}(\xi)\right]
        -\tr\!\left[\Theta_nR_{n,-j}(\xi)\right].
\end{align}

Let
\[
        \mathcal F_j=\sigma(x_1,\ldots,x_j),
        \qquad
        \E_j(\,\cdot\,)=\E(\,\cdot\,\mid\mathcal F_j),
        \qquad j=0,\ldots,n,
\]
where $\mathcal F_0$ is the trivial sigma-field.
Since $\Delta_{n,-j}(\zeta,\xi)$ is
$\mathcal F_{n,-j}$-measurable, independence of the rows implies
$(\E_j-\E_{j-1})\Delta_{n,-j}(\zeta,\xi)=0$.
Thus, with
\[
        d_{n,j}
        \coloneqq
        (\E_j-\E_{j-1})
        \{\Delta_n(\zeta,\xi)-\Delta_{n,-j}(\zeta,\xi)\},
\]
we obtain
\begin{align}\label{eq:signal-martingale-decomposition}
        \Delta_n(\zeta,\xi)-\E\{\Delta_n(\zeta,\xi)\}
        =\sum_{j=1}^n d_{n,j}.
\end{align}
We next establish
\begin{align}\label{eq:signal-martingale-increment-bound}
        \E|d_{n,j}|^4\leq\frac{C_{\zeta,\xi}}{n^4}.
\end{align}

For this purpose, define the following
$\mathcal F_{n,-j}$-measurable matrices:
\begin{align}\label{eq:signal-deletion-matrices}
\begin{split}
        A_{n,j}
        &\coloneqq
        R_{n,-j}(\zeta)C_nR_{n,-j}(\xi),
        \\
        M_{n,j}^{(1)}
        &\coloneqq
        R_{n,-j}(\zeta)C_nR_{n,-j}(\xi)
        \Theta_nR_{n,-j}(\zeta),
        \\
        M_{n,j}^{(2)}
        &\coloneqq
        R_{n,-j}(\xi)\Theta_n
        \{I_p-R_{n,-j}(\zeta)C_n\}R_{n,-j}(\xi),
        \\
        M_{n,j}^{(3)}
        &\coloneqq
        R_{n,-j}(\xi)\Theta_nR_{n,-j}(\zeta).
\end{split}
\end{align}
By the Sherman--Morrison formula,
\begin{align}\label{eq:signal-SM-explicit}
        R_n(w)
        =R_{n,-j}(w)
        -\frac{\beta_j(w)}n
        R_{n,-j}(w)x_jx_j^\top R_{n,-j}(w),
        \qquad w\in\{\zeta,\xi\}.
\end{align}
Substituting \eqref{eq:signal-SM-explicit} into
\eqref{eq:signal-reduced-representation} and subtracting
\eqref{eq:signal-reduced-leave-one-out} gives
\begin{align}\label{eq:signal-three-term-deletion}
\begin{split}
        \Delta_n(\zeta,\xi)-\Delta_{n,-j}(\zeta,\xi)
        &=-\frac{\beta_j(\zeta)}n
        x_j^\top M_{n,j}^{(1)}x_j
        +\frac{\beta_j(\xi)}n
        x_j^\top M_{n,j}^{(2)}x_j
        \\
        &\quad+
        \frac{\beta_j(\zeta)\beta_j(\xi)}{n^2}
        \{x_j^\top A_{n,j}x_j\}
        \{x_j^\top M_{n,j}^{(3)}x_j\}.
\end{split}
\end{align}
By \eqref{eq:general-resolvent-bound},
\eqref{eq:signal-reduced-matrix}, and
\eqref{eq:signal-deletion-matrices},
\begin{align}\label{eq:signal-deletion-matrix-bounds}
        \|A_{n,j}\|_2\leq C_{\zeta,\xi},
        \qquad
        \max_{1\leq r\leq3}\|M_{n,j}^{(r)}\|_F
        \leq C_{\zeta,\xi}\|\Theta_n\|_F.
\end{align}

For any $\mathcal F_{n,-j}$-measurable complex matrix $M$,
Lemma~\ref{lem_b26_bai_silverstein}, applied conditionally on
$\mathcal F_{n,-j}$ with $q=4$, and Assumption~\ref{A2} give
\begin{align}\label{eq:signal-conditional-quadratic-fourth}
\begin{split}
        \E\!\left(
        |x_j^\top Mx_j-\tr(\Sigma_nM)|^4
        \,\middle|\,\mathcal F_{n,-j}\right)
        \leq
        C\|\Sigma_n^{1/2}M\Sigma_n^{1/2}\|_F^4
        \leq C\|M\|_F^4.
\end{split}
\end{align}
We also have
\begin{align}\label{eq:signal-trace-Frobenius-bound}
        |\tr(\Sigma_nM)|
        \leq\sqrt p\,\|\Sigma_n\|_2\|M\|_F.
\end{align}
Recall $b_j(w)$ from \eqref{eq:a_jandb_j}, and note that 
by \eqref{eq:general-denominator-bound} and
\eqref{eq:general-trace-denominator-bound},
\begin{align}\label{eq:signal-beta-b-uniform}
        |\beta_j(w)|+|b_j(w)|\leq C_{\zeta,\xi},
        \qquad w\in\{\zeta,\xi\}.
\end{align}
Moreover, \eqref{eq:beta_b_gamma} gives
\[
\begin{split}
        \beta_j(w)-b_j(w)
        =-\frac{\beta_j(w)b_j(w)}n
        \left\{
        x_j^\top R_{n,-j}(w)x_j
        -\tr[\Sigma_nR_{n,-j}(w)]
        \right\}.
\end{split}
\]
Applying \eqref{eq:signal-conditional-quadratic-fourth} with
$M=R_{n,-j}(w)$, and using
$\|R_{n,-j}(w)\|_F\leq\sqrt p\,\|R_{n,-j}(w)\|_2$, yields
\begin{align}\label{eq:signal-beta-b-fourth}
        \E\!\left(
        |\beta_j(w)-b_j(w)|^4
        \,\middle|\,\mathcal F_{n,-j}\right)
        \leq\frac{C_{\zeta,\xi}p^2}{n^4},
        \qquad w\in\{\zeta,\xi\}.
\end{align}
For the third term in \eqref{eq:signal-three-term-deletion}, define
\begin{align}\label{eq:signal-c-coefficients}
\begin{split}
        c_{n,j}
        &\coloneqq
        \frac{\beta_j(\zeta)\beta_j(\xi)}n
        x_j^\top A_{n,j}x_j,
        \\
        c_{n,j}^{0}
        &\coloneqq
        \frac{b_j(\zeta)b_j(\xi)}n
        \tr(\Sigma_nA_{n,j}).
\end{split}
\end{align}
Note that the coefficient $c_{n,j}^{0}$ is
$\mathcal F_{n,-j}$-measurable. We show both that $c_{n,j}$ is
uniformly bounded and that $c_{n,j}-c_{n,j}^{0}$ satisfies the
same fourth-moment bound as in \eqref{eq:signal-beta-b-fourth}.

Since $R_{n,-j}(w)^\top=R_{n,-j}(w)$, the rank-one operator-norm
identity and \eqref{eq:signal-SM-explicit} give
\begin{align}\label{eq:signal-rank-one-norm-bound}
\begin{split}
        \frac{|\beta_j(w)|}{n}\|R_{n,-j}(w)x_j\|_2^2
        &=\|R_{n,-j}(w)-R_n(w)\|_2
        \\
        &\leq\|R_{n,-j}(w)\|_2+\|R_n(w)\|_2
        \leq\frac{2}{\operatorname{dist}(w,\mathbb R_{\geq0})}.
\end{split}
\end{align}
Thus, by Cauchy--Schwarz,
\begin{align}\label{eq:signal-c-uniform-bound}
\begin{split}
        |c_{n,j}|
        &\leq
        \|C_n\|_2\sqrt{|\beta_j(\zeta)\beta_j(\xi)|}
        \left\{
        \frac{|\beta_j(\zeta)|}{n}
        \|R_{n,-j}(\zeta)x_j\|_2^2
        \right\}^{1/2}
        \\
        &\qquad\times
        \left\{
        \frac{|\beta_j(\xi)|}{n}
        \|R_{n,-j}(\xi)x_j\|_2^2
        \right\}^{1/2}
        \leq C_{\zeta,\xi},
\end{split}
\end{align}
where the last inequality uses
\eqref{eq:signal-reduced-matrix},
\eqref{eq:signal-beta-b-uniform}, and
\eqref{eq:signal-rank-one-norm-bound}.

Next, expanding only the scalar difference gives
\begin{align}\label{eq:signal-c-difference}
\begin{split}
        c_{n,j}-c_{n,j}^{0}
        &=\frac{\beta_j(\zeta)\beta_j(\xi)}n
        \{x_j^\top A_{n,j}x_j-\tr(\Sigma_nA_{n,j})\}
        \\
        &\quad+
        \frac{\tr(\Sigma_nA_{n,j})}{n}
        \left[
        \beta_j(\zeta)\{\beta_j(\xi)-b_j(\xi)\}
        +b_j(\xi)\{\beta_j(\zeta)-b_j(\zeta)\}
        \right].
\end{split}
\end{align}
By \eqref{eq:signal-deletion-matrix-bounds},
$\|A_{n,j}\|_F\leq C_{\zeta,\xi}\sqrt p$ and
$|\tr(\Sigma_nA_{n,j})|\leq C_{\zeta,\xi}p$.
Consequently, \eqref{eq:signal-c-difference},
\eqref{eq:signal-conditional-quadratic-fourth},
\eqref{eq:signal-beta-b-uniform}, and
\eqref{eq:signal-beta-b-fourth} imply
\begin{align}\label{eq:signal-c-comparison-fourth}
\begin{split}
        \E\!\left(
        |c_{n,j}-c_{n,j}^{0}|^4
        \,\middle|\,\mathcal F_{n,-j}\right)
        &\leq
        C_{\zeta,\xi}\left\{
        \frac{p^2}{n^4}+\frac{p^6}{n^8}
        \right\}
        \leq\frac{C_{\zeta,\xi}p^2}{n^4}.
\end{split}
\end{align}
The last inequality uses $p/n=O(1)$ from \eqref{gamma_n:cond}.

We can now estimate all three terms of
\eqref{eq:signal-three-term-deletion}. Let $(a,a^0,M)$ be any
one of the three triples
\[
        (\beta_j(\zeta),b_j(\zeta),M_{n,j}^{(1)}),
        \qquad
        (\beta_j(\xi),b_j(\xi),M_{n,j}^{(2)}),
        \qquad
        (c_{n,j},c_{n,j}^{0},M_{n,j}^{(3)}).
\]
In each case, $a^0$ and $M$ are
$\mathcal F_{n,-j}$-measurable, $|a|\leq C_{\zeta,\xi}$,
and the conditional fourth moment of $a-a^0$ is at most
$C_{\zeta,\xi}p^2/n^4$, by
\eqref{eq:signal-beta-b-uniform},
\eqref{eq:signal-beta-b-fourth},
\eqref{eq:signal-c-uniform-bound}, and
\eqref{eq:signal-c-comparison-fourth}.
The exact identity
\[
        \frac a n x_j^\top Mx_j
        -\frac{a^0}{n}\tr(\Sigma_nM)
        =\frac a n\{x_j^\top Mx_j-\tr(\Sigma_nM)\}
        +\frac{a-a^0}{n}\tr(\Sigma_nM)
\]
therefore yields
\begin{align}\label{eq:signal-weighted-quadratic-fourth}
\begin{split}
        \E\!\left(
        \left|
        \frac a n x_j^\top Mx_j
        -\frac{a^0}{n}\tr(\Sigma_nM)
        \right|^4
        \,\middle|\,\mathcal F_{n,-j}\right)
        &\leq
        \frac{C_{\zeta,\xi}}{n^4}\|M\|_F^4
        +\frac{C_{\zeta,\xi}p^2}{n^8}
        |\tr(\Sigma_nM)|^4
        \\
        &\qquad\leq
        \frac{C_{\zeta,\xi}}{n^4}
        \left(1+\frac{p^4}{n^4}\right)\|M\|_F^4
        \leq
        \frac{C_{\zeta,\xi}}{n^4}\|\Theta_n\|_F^4.
\end{split}
\end{align}
Here the first inequality uses
$|u+v|^4\leq8(|u|^4+|v|^4)$ and
\eqref{eq:signal-conditional-quadratic-fourth}; the second uses
\eqref{eq:signal-trace-Frobenius-bound}; and the last uses
\eqref{eq:signal-deletion-matrix-bounds} and
\eqref{gamma_n:cond}.

Define the $\mathcal F_{n,-j}$-measurable quantity
\begin{align}\label{eq:signal-measurable-proxy}
\begin{split}
        \Lambda_{n,j}
        &\coloneqq
        -\frac{b_j(\zeta)}n\tr(\Sigma_nM_{n,j}^{(1)})
        +\frac{b_j(\xi)}n\tr(\Sigma_nM_{n,j}^{(2)})
        +\frac{c_{n,j}^{0}}n\tr(\Sigma_nM_{n,j}^{(3)}).
\end{split}
\end{align}
Using \eqref{eq:signal-c-coefficients} in
\eqref{eq:signal-three-term-deletion} and applying
\eqref{eq:signal-weighted-quadratic-fourth} to its three terms,
we obtain
\begin{align}\label{eq:signal-deletion-proxy-fourth}
        \E\left(\left|
        \Delta_n(\zeta,\xi)-\Delta_{n,-j}(\zeta,\xi)
        -\Lambda_{n,j}
        \right|^4\right)
        \leq\frac{C_{\zeta,\xi}}{n^4}\|\Theta_n\|_F^4
        \leq\frac{C_{\zeta,\xi}}{n^4}.
\end{align}
Independence of the rows gives
$(\E_j-\E_{j-1})\Lambda_{n,j}=0$. Hence
\[
        d_{n,j}
        =(\E_j-\E_{j-1})
        \{\Delta_n(\zeta,\xi)-\Delta_{n,-j}(\zeta,\xi)
        -\Lambda_{n,j}\}.
\]
By conditional Jensen's inequality,
\[
        \E\left(|d_{n,j}|^4\right)
        \leq16\E\left(\left|
        \Delta_n(\zeta,\xi)-\Delta_{n,-j}(\zeta,\xi)
        -\Lambda_{n,j}
        \right|^4\right)
        \leq\frac{C_{\zeta,\xi}}{n^4},
\]
which proves \eqref{eq:signal-martingale-increment-bound}.

Applying the Burkholder inequality
\citep[Lemma~2.12]{Silversteinbook2010} to
\eqref{eq:signal-martingale-decomposition} gives
\[
        \E\left(\left|
        \Delta_n(\zeta,\xi)-\E\{\Delta_n(\zeta,\xi)\}
        \right|^4\right)\leq
        C\E\left(\sum_{j=1}^n\E_{j-1}|d_{n,j}|^2\right)^2
        +C\sum_{j=1}^n\E|d_{n,j}|^4.
\]
By Cauchy--Schwarz and conditional Jensen,
\[
        \E\left(\sum_{j=1}^n\E_{j-1}|d_{n,j}|^2\right)^2
        \leq n\sum_{j=1}^n
        \E(\E_{j-1}|d_{n,j}|^2)^2
        \leq n\sum_{j=1}^n\E|d_{n,j}|^4.
\]
Together with \eqref{eq:signal-martingale-increment-bound},
these inequalities imply
\begin{align}\label{eq:signal-centered-fourth}
        \E\left(\left|
        \Delta_n(\zeta,\xi)-\E\{\Delta_n(\zeta,\xi)\}
        \right|^4\right)
        \leq\frac{C_{\zeta,\xi}}{n^2}.
\end{align}
Consequently, for every $\varepsilon>0$, Markov's inequality gives
\[
        \sum_{n=1}^\infty
        \mathbb P\left(
        |\Delta_n(\zeta,\xi)-\E\{\Delta_n(\zeta,\xi)\}|
        >\varepsilon\right)
        \leq\frac{C_{\zeta,\xi}}{\varepsilon^4}
        \sum_{n=1}^\infty\frac1{n^2}<\infty.
\]

Finally, \eqref{eq:signal-in-probability} and
\eqref{eq:signal-centered-as} imply
\[
\begin{split}
        \E\{\Delta_n(\zeta,\xi)\}
        &=\Delta_n(\zeta,\xi)
        -\bigl[\Delta_n(\zeta,\xi)-\E\{\Delta_n(\zeta,\xi)\}\bigr]
        \overset{\mathbb P}{\longrightarrow}0.
\end{split}
\]
Since this quantity is deterministic, convergence in probability
is equivalent to ordinary convergence. Therefore,
\[
        \Delta_n(\zeta,\xi)
        =\bigl[\Delta_n(\zeta,\xi)-\E\{\Delta_n(\zeta,\xi)\}\bigr]
        +\E\{\Delta_n(\zeta,\xi)\}
        \longrightarrow0
        \qquad\text{almost surely}.
\]
This proves \eqref{eq:ridge-signal-equivalence}.
\end{proof}

\subsection{Proof of Theorem~\ref{thm:double-resolvent-limit}}

\begin{proof}[Proof of Theorem~\ref{thm:double-resolvent-limit}] 
Fix a compact set $K\subseteq\mathbb C\setminus\mathbb R_{\geq0}$.
Choose $0<\rho<\operatorname{dist}(K,\mathbb R_{\geq0})/2$ and set
\[
    K_1=\{w\in\mathbb C:\operatorname{dist}(w,K)\leq\rho\}.
\]

We first show that $|-wv_n(w)|$ is uniformly
bounded away from zero for $w\in K_1$ and $n\geq1$.
By \eqref{stielt}, with $(\gamma_0,H_0)=(\gamma_n,H_n)$, and the definition of $R_n^{\mathrm{det}}$ in \eqref{eq:def_R},
\[
    \frac{1}{-wv_n(w)}
    =
    1+\frac{1}{n}
    \tr[\Sigma_nR_n^{\mathrm{det}}(w)].
\]
Consequently, Assumption~\ref{A2},
\eqref{gamma_n:cond}, and
\eqref{eq:deterministic-resolvent-bound} imply
\begin{equation}\label{eq:double-scalar-inverse-bound}
    \sup_{n\geq1}\sup_{w\in K_1}
    \left|\frac{1}{-wv_n(w)}\right|
    \leq
    1+
    \frac{
        \sup_{n\geq1}\gamma_n\|\Sigma_n\|_2
    }{
        \operatorname{dist}(K_1,\mathbb R_{\geq0})
    }
    <\infty.
\end{equation}

Next, define
\[
    f_n(w)
    \coloneqq
    w\tr\!\left[
        \Theta_n
        \{R_n(w)-R_n^{\mathrm{det}}(w)\}
    \right],
    \qquad
    w\in\mathbb C\setminus\mathbb R_{\geq0}.
\]
Theorem~\ref{Lemma:GreenRom},
Corollary~\ref{cor:GreenRom-derivatives:rates},
and the product rule give
\begin{equation}\label{eq:double-single-error-rate}
    \sup_{w\in K_1}
    \left\{
        |f_n(w)|
        +
        \left|\frac{d}{dw}f_n(w)\right|
    \right\}
    =
    O_{\mathbb P}\!\left(
        \frac{\|\Theta_n\|_F}{\sqrt n}
    \right).
\end{equation}
For $\zeta\neq\xi$, the resolvent identity gives
\[
    R_n(\zeta)\widehat\Sigma_nR_n(\xi)
    =
    \frac{
        \zeta R_n(\zeta)-\xi R_n(\xi)
    }{\zeta-\xi},
\]
and for $\zeta=\xi$, we have
\[
    R_n(\zeta)\widehat\Sigma_nR_n(\zeta)
    =
    R_n(\zeta)+\zeta R_n(\zeta)^2
    =
    \frac{d}{d\zeta}[\zeta R_n(\zeta)].
\]
Thus, by the definition of $K_n^\Theta(\zeta,\xi)$ in \eqref{def_K_theta},
\begin{equation}\label{eq:double-empirical-error}
\begin{split}
    &\tr\!\left[
        \Theta_nR_n(\zeta)\widehat\Sigma_nR_n(\xi)
    \right]
    -
    \{-\zeta v_n(\zeta)\}\{-\xi v_n(\xi)\}
    K_n^\Theta(\zeta,\xi)\\
    &\qquad=
    \begin{cases}
        \displaystyle
        \frac{f_n(\zeta)-f_n(\xi)}{\zeta-\xi},
        & \zeta\neq\xi,\\[3mm]
        \displaystyle
        \frac{d}{d\zeta}f_n(\zeta),
        & \zeta=\xi.
    \end{cases}
\end{split}
\end{equation}

Next, we bound this expression uniformly over $K\times K$.
If $0<|\zeta-\xi|\leq\rho$, the line segment joining
$\xi$ and $\zeta$ lies in $K_1$, since
\[
    \operatorname{dist}\bigl(\xi+t(\zeta-\xi),K\bigr)
    \leq t|\zeta-\xi|
    \leq\rho,
    \qquad t\in[0,1].
\]
Integrating the derivative along this segment gives
\[
    \frac{f_n(\zeta)-f_n(\xi)}{\zeta-\xi}
    =
    \int_0^1
    \left.
        \frac{d}{dw}f_n(w)
    \right|_{w=\xi+t(\zeta-\xi)}
    \,dt,
\]
and hence
\[
    \left|
        \frac{f_n(\zeta)-f_n(\xi)}{\zeta-\xi}
    \right|
    \leq
    \sup_{w\in K_1}
    \left|\frac{d}{dw}f_n(w)\right|.
\]
The same bound applies to the diagonal expression
in \eqref{eq:double-empirical-error}.
If $|\zeta-\xi|>\rho$, then
\[
    \left|
        \frac{f_n(\zeta)-f_n(\xi)}{\zeta-\xi}
    \right|
    \leq
    \frac{2}{\rho}\sup_{w\in K_1}|f_n(w)|.
\]
Together with \eqref{eq:double-single-error-rate},
these estimates imply
\begin{equation}\label{eq:double-empirical-error-rate}
\sup_{\zeta,\xi\in K}
    \left|
        \tr\!\left[
            \Theta_nR_n(\zeta)\widehat\Sigma_nR_n(\xi)
        \right]
        -
        \{-\zeta v_n(\zeta)\}\{-\xi v_n(\xi)\}
        K_n^\Theta(\zeta,\xi)
    \right|=
    O_{\mathbb P}\!\left(
        \frac{\|\Theta_n\|_F}{\sqrt n}
    \right).
\end{equation}
Moreover, Lemma~\ref{lem:ridge-signal-equivalence:rates}
and \eqref{gamma_n:cond} give
\begin{equation}\label{eq:double-signal-error-rate}
    \sup_{\zeta,\xi\in K}
    \left|
        \{-\zeta v_n(\zeta)\}\{-\xi v_n(\xi)\}
        \tr\!\left[
            \Theta_nR_n(\zeta)\Sigma_nR_n(\xi)
        \right]
        -
        \tr\!\left[
            \Theta_nR_n(\zeta)\widehat\Sigma_nR_n(\xi)
        \right]
    \right|=
    O_{\mathbb P}\!\left(
        \frac{\|\Theta_n\|_F}{\sqrt n}
    \right).
\end{equation}
Adding the expressions in
\eqref{eq:double-empirical-error-rate}
and~\eqref{eq:double-signal-error-rate} yields
\[
    \sup_{\zeta,\xi\in K}
    \left|
        \{-\zeta v_n(\zeta)\}\{-\xi v_n(\xi)\}
        \left\{
            \tr[
                \Theta_nR_n(\zeta)\Sigma_nR_n(\xi)]
            -
            K_n^\Theta(\zeta,\xi)
        \right\}
    \right|=
    O_{\mathbb P}\!\left(
        \frac{\|\Theta_n\|_F}{\sqrt n}
    \right).
\]
Finally, \eqref{eq:double-scalar-inverse-bound}
allows us to divide by
$\{-\zeta v_n(\zeta)\}\{-\xi v_n(\xi)\}$
uniformly over $(\zeta, \xi)\in K\times K$, proving the desired assertion.
\end{proof}

\subsection{Proof of Theorem~\ref{thm:double-resolvent-limit2}}

\begin{proof}[Proof of Theorem~\ref{thm:double-resolvent-limit2}]
Write $X\coloneqq\mathbb C\setminus\mathbb R_{\geq0}$ and define
\[
        F_n(\zeta,\xi)
        \coloneqq
        \tr\!\left[
        \Theta_nR_n(\zeta)\Sigma_nR_n(\xi)
        \right]
        -
        K_n^\Theta(\zeta,\xi),
        \qquad (\zeta,\xi)\in X^2.
\]
We suppress the dependence on $\omega$ unless needed.
Let $E$ be an event of probability one on which all sample
covariance matrices are positive semidefinite. All resolvent
bounds below hold uniformly for $\omega\in E$.

We first consider a fixed pair $(\zeta,\xi)\in X^2$ with
$\zeta\neq\xi$.
Proposition~\ref{prop:RM-recovered} gives, for every fixed $w\in X$,
\[
        f_n(w)
        \coloneqq
        w\tr\!\left[
        \Theta_n
        \{R_n(w)-R_n^{\operatorname{det}}(w)\}
        \right]
        \longrightarrow0
        \qquad\text{almost surely}.
\]
In particular, this finite-$n$ approximation does not require
convergence of $\gamma_n$ or $H_n$.

Note that
\[
        R_n(\zeta)\hat\Sigma_nR_n(\xi)
        =
        \frac{\zeta R_n(\zeta)-\xi R_n(\xi)}
             {\zeta-\xi}.
\]
Moreover, by \eqref{def_K_theta},
\[
        \{-\zeta v_n(\zeta)\}\{-\xi v_n(\xi)\}
        K_n^\Theta(\zeta,\xi)
        =
        \frac{
        \zeta\tr\{\Theta_nR_n^{\operatorname{det}}(\zeta)\}
        -
        \xi\tr\{\Theta_nR_n^{\operatorname{det}}(\xi)\}
        }{\zeta-\xi}.
\]
Subtracting these trace expressions yields
\begin{align}\label{eq:double-empirical-error-as}
\begin{split}
        &\tr\!\left[
        \Theta_nR_n(\zeta)\hat\Sigma_nR_n(\xi)
        \right]
        -
        \{-\zeta v_n(\zeta)\}\{-\xi v_n(\xi)\}
        K_n^\Theta(\zeta,\xi)
        \\
        &\qquad=
        \frac{f_n(\zeta)-f_n(\xi)}{\zeta-\xi}
        \longrightarrow0
        \qquad\text{almost surely}.
\end{split}
\end{align}
Since $\|\Theta_n\|_F\leq\|\Theta_n\|_{1}$,
Lemma~\ref{lem:ridge-signal-equivalence} applies. Combining
\eqref{eq:ridge-signal-equivalence} with
\eqref{eq:double-empirical-error-as} gives
\[
        \{-\zeta v_n(\zeta)\}\{-\xi v_n(\xi)\}
        F_n(\zeta,\xi)
        \longrightarrow0
        \qquad\text{almost surely}.
\]

By \eqref{stielt}, applied with
$(\gamma_0,H_0)=(\gamma_n,H_n)$, and the definition of
$R_n^{\operatorname{det}}$,
\begin{align}\label{eq:double-inverse-scalar-identity}
\begin{split}
        \{-wv_n(w)\}^{-1}
        &=
        1-\frac{\gamma_n}{w}
        \int_0^\infty
        \frac{t}{1+tv_n(w)}\,dH_n(t)
        \\
        &=
        1+\frac1n
        \tr\!\left[
        \Sigma_nR_n^{\operatorname{det}}(w)
        \right],
        \qquad w\in X.
\end{split}
\end{align}
Thus, by \eqref{eq:deterministic-resolvent-bound},
Assumption~\ref{A2}, and \eqref{gamma_n:cond},
\[
        \sup_{n\geq1}|\{-wv_n(w)\}^{-1}|
        \leq
        1+
        \frac{
        \sup_{n\geq1}\{\gamma_n\|\Sigma_n\|_2\}
        }{
        \operatorname{dist}(w,\mathbb R_{\geq0})
        }
        <\infty.
\]
It follows that
\begin{align}\label{eq:double-off-diagonal-as}
        F_n(\zeta,\xi)\longrightarrow0
        \qquad\text{almost surely}
\end{align}
for every fixed pair $(\zeta,\xi)\in X^2$ with $\zeta\neq\xi$.

We next establish local uniform bounds.
Let $K\subset X$ be compact, choose
\[
        0<\rho<\operatorname{dist}(K,\mathbb R_{\geq0}),
\]
and set
\[
        K_1
        \coloneqq
        \{w\in\mathbb C:\operatorname{dist}(w,K)\leq\rho\}.
\]
Choose $0<\rho_1<\operatorname{dist}(K_1,\mathbb R_{\geq0})$
and define
\[
        K_2
        \coloneqq
        \{w\in\mathbb C:\operatorname{dist}(w,K_1)\leq\rho_1\}.
\]
Then $K\subseteq K_1\subseteq K_2\subset X$, and all three
sets are compact.

Define the holomorphic function
\[
        g_n(w)
        \coloneqq
        w\tr\!\left[
        \Theta_nR_n^{\operatorname{det}}(w)
        \right],
        \qquad w\in X.
\]
By \eqref{eq:deterministic-resolvent-bound} and
$\sup_n\|\Theta_n\|_{1}<\infty$,

\begin{align}\label{eq:double-g-uniform-bound}
        \sup_{n\geq1}\sup_{w\in K_2}|g_n(w)|
        \leq
        \frac{
        \sup_{w\in K_2}|w|\,
        \sup_{n\geq1}\|\Theta_n\|_{1}
        }{
        \operatorname{dist}(K_2,\mathbb R_{\geq0})
        }
        \leq C_{K_2}<\infty.
\end{align}
Moreover, \eqref{eq:double-inverse-scalar-identity} gives
\begin{align}\label{eq:double-inverse-scalar-compact}
        \sup_{n\geq1}\sup_{w\in K_2}
        |\{-wv_n(w)\}^{-1}|
        \leq
        1+
        \frac{
        \sup_{n\geq1}\{\gamma_n\|\Sigma_n\|_2\}
        }{
        \operatorname{dist}(K_2,\mathbb R_{\geq0})
        }
        \leq C_{K_2}.
\end{align}
Here and below, constants are independent of $n$ and $\omega$.

For $w\in K_1$, the closed disc of radius $\rho_1$ centered
at $w$ is contained in $K_2$. Hence Cauchy's estimate yields
\begin{align}\label{eq:double-g-derivative-bound}
        \sup_{n\geq1}\sup_{w\in K_1}|g_n'(w)|
        \leq
        \frac1{\rho_1}
        \sup_{n\geq1}\sup_{u\in K_2}|g_n(u)|
        \leq
        \frac{C_{K_2}}{\rho_1}.
\end{align}
For $z,w\in K$ with $0<|z-w|\leq\rho$, the line segment
joining $w$ to $z$ lies in $K_1$. Therefore,
\[
        \left|\frac{g_n(z)-g_n(w)}{z-w}\right|
        =
        \left|
        \int_0^1 g_n'\bigl(w+t(z-w)\bigr)\,dt
        \right|
        \leq
        \frac{C_{K_2}}{\rho_1}.
\]
If $|z-w|>\rho$, then
\[
        \left|\frac{g_n(z)-g_n(w)}{z-w}\right|
        \leq
        \frac{|g_n(z)|+|g_n(w)|}{\rho}
        \leq
        \frac{2C_{K_2}}{\rho}.
\]
For $z=w$, the corresponding divided difference is
$g_n'(z)$ and is bounded by
\eqref{eq:double-g-derivative-bound}.
Using these bounds in \eqref{def_K_theta}, together with
\eqref{eq:double-inverse-scalar-compact}, gives
\[
        \sup_{n\geq1}\sup_{z,w\in K}|K_n^\Theta(z,w)|
        \leq
        C_{K_2}^3
        \max\left\{\frac1{\rho_1},\frac2\rho\right\}.
\]
Furthermore, Assumption~\ref{A2} and
\eqref{eq:general-resolvent-bound} imply
\[
        \sup_{n\geq1}\sup_{z,w\in K}
        \left|
        \tr\!\left[
        \Theta_nR_n(z)\Sigma_nR_n(w)
        \right]
        \right|
        \leq
        \frac{
        \sup_{n\geq1}
        \{\|\Theta_n\|_{1}\|\Sigma_n\|_2\}
        }{
        \operatorname{dist}(K,\mathbb R_{\geq0})^2
        }.
\]
Consequently,
\begin{align}\label{eq:double-local-uniform-bound}
        \sup_{n\geq1}\sup_{z,w\in K}|F_n(z,w,\omega)|
        \leq C_K<\infty,
        \qquad \omega\in E,
\end{align}
where $C_K$ may also depend on the chosen radii.

The functions $R_n$ and $R_n^{\operatorname{det}}$ are
holomorphic on $X$, and
\eqref{eq:double-inverse-scalar-identity} shows that
$w\mapsto\{-wv_n(w)\}^{-1}$ is holomorphic on $X$.
The divided difference of $g_n$, with value $g_n'(w)$ at
$z=w$, is holomorphic in each variable.
Thus, by \eqref{def_K_theta}, $F_n$ is holomorphic in each
variable on $X^2$ for every $\omega\in E$.

Fix $(\zeta_0,\xi_0)\in X^2$ and choose $r>0$ such that
\[
        L
        \coloneqq
        \{z:|z-\zeta_0|\leq2r\}
        \cup
        \{w:|w-\xi_0|\leq2r\}
        \subset X.
\]
Applying \eqref{eq:double-local-uniform-bound} with $K=L$,
Cauchy's estimates in each variable give
\[
        \sup_{n\geq1}
        \sup_{\substack{|z-\zeta_0|\leq r\\|w-\xi_0|\leq r}}
        \left\{
        |\partial_zF_n(z,w,\omega)|
        +
        |\partial_wF_n(z,w,\omega)|
        \right\}
        \leq
        \frac{2C_L}{r},
        \qquad \omega\in E.
\]
Integrating along line segments in the two discs therefore yields
\begin{align}\label{eq:double-local-lipschitz}
|F_n(z_1,w_1,\omega)-F_n(z_2,w_2,\omega)|
        \leq
        \frac{2C_L}{r}
        \bigl(|z_1-z_2|+|w_1-w_2|\bigr),
\end{align}
whenever $|z_k-\zeta_0|\leq r$ and $|w_k-\xi_0|\leq r$,
$k=1,2$, uniformly in $n$ and $\omega\in E$.
In particular, by a finite covering argument, the family
$\{F_n(\cdot,\cdot,\omega)\}_{n\geq1}$ is equicontinuous
on every compact subset of $X^2$ for every $\omega\in E$.

We now extend the pointwise almost sure convergence to the
diagonal.
Fix $\zeta\in X$ and choose a deterministic sequence
$\xi_k\in X\setminus\{\zeta\}$ converging to $\zeta$.
By \eqref{eq:double-off-diagonal-as}, there exists an event
$A_\zeta\subseteq E$ with $\mathbb P(A_\zeta)=1$ such that
\[
        F_n(\zeta,\xi_k,\omega)\longrightarrow0
        \qquad\text{for every }k\geq1
        \text{ and }\omega\in A_\zeta.
\]
Applying \eqref{eq:double-local-lipschitz} in a neighborhood
of $(\zeta,\zeta)$ gives, for all sufficiently large $k$,
\[
        \limsup_{n\to\infty}|F_n(\zeta,\zeta,\omega)|
        \leq
        L_\zeta|\zeta-\xi_k|,
        \qquad \omega\in A_\zeta,
\]
where $L_\zeta$ is independent of $n$, $k$, and $\omega$.
Letting $k\to\infty$ proves that
$F_n(\zeta,\zeta)\to0$ almost surely.
Together with \eqref{eq:double-off-diagonal-as}, this establishes
\[
        F_n(\zeta,\xi)\longrightarrow0
        \qquad\text{almost surely}
\]
for every fixed pair $(\zeta,\xi)\in X^2$.

It remains to obtain a common event of probability one.
For $m\geq1$, define
\[
        \mathcal K_m
        \coloneqq
        \left\{
        w\in\mathbb C:
        |w|\leq m,\ 
        \operatorname{dist}(w,\mathbb R_{\geq0})\geq m^{-1}
        \right\}.
\]
Each $\mathcal K_m^2$ is compact, and $X^2=\bigcup_{m\geq1}\mathcal K_m^2$.
For every fixed pair $(\zeta,\xi)$, the maps
$\omega\mapsto F_n(\zeta,\xi,\omega)$ are measurable.
The pointwise almost sure convergence and the equicontinuity
on $\mathcal K_m^2$ therefore allow us to apply
Lemma~\ref{lem:Georg_global_nullset} with
$X=\mathcal K_m^2$ and $Y=\mathbb C$.
This gives an event $\widetilde C_m$ of probability one on which
$F_n$ converges at every point of $\mathcal K_m^2$ to a continuous
limiting function.

Choose a countable dense subset $D_m\subset\mathcal K_m^2$ and set
\[
        C_m
        \coloneqq
        \widetilde C_m\cap E
        \cap
        \bigcap_{(\zeta,\xi)\in D_m}
        \left\{
        \omega:
        F_n(\zeta,\xi,\omega)\longrightarrow0
        \right\}.
\]
Then $\mathbb P(C_m)=1$.
For every $\omega\in C_m$, the continuous limiting function
vanishes on $D_m$ and hence on all of $\mathcal K_m^2$.

Finally, $\Omega_0\coloneqq\bigcap_{m\geq1}C_m$ has probability one. Since $X^2=\bigcup_{m\geq1}\mathcal K_m^2$,
for every $\omega\in\Omega_0$ we have
\[
        F_n(\zeta,\xi,\omega)\longrightarrow0
        \qquad\text{for every }(\zeta,\xi)\in X^2.
\]
This proves the assertion.
\end{proof}

\section{Proofs of the results in  Section~\ref{sec:applications}}

\subsection{Proof of Proposition~\ref{prop:app-precision-loss}}
\begin{proof}[Proof of Proposition~\ref{prop:app-precision-loss}]
The identities \eqref{eq:fixed_point} and
\eqref{eq:emp_fixed_point}, together with the fixed-point
equation \eqref{stielt}, Theorem~\ref{Lemma:GreenRom},
and Corollary~\ref{cor:GreenRom-derivatives:rates},
applied with $\Theta_n=I_p/p$, yield
\begin{equation}\label{eq:app-precision-companion-rate}
    \sup_{\lambda\in\Lambda}
    \left\{
        |\widehat v_n(-\lambda)-v_n(-\lambda)|
        +
        \left|
            \frac{d}{d\lambda}\widehat v_n(-\lambda)
            -
            \frac{d}{d\lambda}v_n(-\lambda)
        \right|
    \right\}
    =
    O_{\mathbb P}(n^{-1}).
\end{equation}
Moreover, \eqref{stielt}, with
$(\gamma_0,H_0)=(\gamma_n,H_n)$ and $\zeta=-\lambda$,
implies that
\begin{align}\label{eq:vn_lowerbound}
    \frac{1}{v_n(-\lambda)}
    =
    \lambda+
    \frac{\gamma_n}{p}
    \tr\!\left[
        \Sigma_n
        (I_p+v_n(-\lambda)\Sigma_n)^{-1}
    \right]
    \leq
    \lambda+\gamma_n\|\Sigma_n\|_2.
\end{align}
Hence, there exists a constant $c>0$ such that
$v_n(-\lambda)\geq c$ for all $n$ and $\lambda\in\Lambda$.
By \eqref{eq:app-precision-companion-rate},
\[
    \mathbb P\left(
        \inf_{\lambda\in\Lambda}
        \widehat v_n(-\lambda)\geq c/2
    \right)\longrightarrow1.
\]
Furthermore, since $v_n$ and $\widehat v_n$ are
Stieltjes transforms of probability measures on
$[0,\infty)$, we have
\[
    0<v_n(-\lambda),\widehat v_n(-\lambda)
    \leq\lambda_-^{-1}
\quad\text{and} \quad 0\leq
    -\frac{d}{d\lambda}v_n(-\lambda),
    -\frac{d}{d\lambda}\widehat v_n(-\lambda)
    \leq\lambda_-^{-2}.
\]

On the event considered above, the mean value theorem
gives, for $k=1,2,4$,
\[
    \left|
        \frac{1}{\widehat v_n(-\lambda)^k}
        -
        \frac{1}{v_n(-\lambda)^k}
    \right|
    \leq
    k\left(\frac{2}{c}\right)^{k+1}
    |\widehat v_n(-\lambda)-v_n(-\lambda)|.
\]
In particular, using
\eqref{eq:app-precision-det-traces}
and~\eqref{eq:app-precision-empirical-traces-def},
we obtain
\[
    |\widehat T_{1,n}(\lambda)
      -T_{1,n}^{\mathrm{det}}(\lambda)|
    \leq
    C|\widehat v_n(-\lambda)-v_n(-\lambda)|.
\]
For the second trace functional, we similarly have
\[
\begin{split}
    |\widehat T_{2,n}(\lambda)
      -T_{2,n}^{\mathrm{det}}(\lambda)|&\leq
    \frac{1}{\gamma_n\lambda^2}
    \left\{
        \frac{1}{\widehat v_n(-\lambda)^4}
        \left|
            \frac{d}{d\lambda}\widehat v_n(-\lambda)
            -
            \frac{d}{d\lambda}v_n(-\lambda)
        \right|\right.\\
    &\qquad+
        \left|\frac{d}{d\lambda}v_n(-\lambda)\right|
        \left|
            \frac{1}{\widehat v_n(-\lambda)^4}
            -
            \frac{1}{v_n(-\lambda)^4}
        \right|\left.+
        \left|
            \frac{1}{\widehat v_n(-\lambda)^2}
            -
            \frac{1}{v_n(-\lambda)^2}
        \right|
    \right\}\\
    &\leq
    C\left\{
        |\widehat v_n(-\lambda)-v_n(-\lambda)|
        +
        \left|
            \frac{d}{d\lambda}\widehat v_n(-\lambda)
            -
            \frac{d}{d\lambda}v_n(-\lambda)
        \right|
    \right\},
\end{split}
\]
where $C$ is independent of $n$ and $\lambda$.
Here, we used that $\gamma_n$ and $\lambda$ are
uniformly bounded away from zero.
Taking the supremum over $\Lambda$ and using
\eqref{eq:app-precision-companion-rate} yields
\[
    \sup_{\lambda\in\Lambda}
    \max_{j=1,2}
    \left|
        \widehat T_{j,n}(\lambda)
        -
        T_{j,n}^{\mathrm{det}}(\lambda)
    \right|
    =
    O_{\mathbb P}(n^{-1}).
\]
Since $\mathcal A$ is bounded, we obtain
\[
\sup_{(\alpha,\lambda)\in\mathcal A\times\Lambda}
    \left|
        \widehat{\mathcal Q}_n(\alpha,\lambda)
        -
        \mathcal Q_n^{\mathrm{det}}(\alpha,\lambda)
    \right|
    =
    O_{\mathbb P}(n^{-1}).
\]
Moreover, since $\|\Sigma_n/p\|_F=O(n^{-1/2})$, Theorems~\ref{Lemma:GreenRom}
and~\ref{thm:double-resolvent-limit}, applied with
$\Theta_n=\Sigma_n/p$, together with the boundedness of
$\mathcal A$, yield
\[
    \sup_{(\alpha,\lambda)\in\mathcal A\times\Lambda}
    \left|
        \mathcal Q_n(\alpha,\lambda)
        -
        \mathcal Q_n^{\mathrm{det}}(\alpha,\lambda)
    \right|
    =
    O_{\mathbb P}(n^{-1}).
\]
Combining this with the preceding estimate proves
\eqref{eq:app-precision-empirical-de}.

To prove the second assertion, the triangle inequality
and \eqref{eq:app-precision-empirical-de} give
\[
    \sup_{(\alpha,\lambda)\in\mathcal A\times\Lambda}
    \left|
        \widehat{\mathcal Q}_n(\alpha,\lambda)
        -
        \mathcal Q_n(\alpha,\lambda)
    \right|
    =
    O_{\mathbb P}(n^{-1}).
\]
Finally, we have
\[
    0
    \leq
    \mathcal Q_n(\widehat\alpha_n,\widehat\lambda_n)
    -
    \inf_{(\alpha,\lambda)\in\mathcal A\times\Lambda}
    \mathcal Q_n(\alpha,\lambda)\leq
    2\sup_{(\alpha,\lambda)\in\mathcal A\times\Lambda}
    \left|
        \widehat{\mathcal Q}_n(\alpha,\lambda)
        -
        \mathcal Q_n(\alpha,\lambda)
    \right|,
\]
and together with the preceding estimate the second assertion follows.
This completes the proof.
\end{proof}

\subsection{Proof of Proposition~\ref{prop:app-ridge-tuning}}
\begin{proof}[Proof of Proposition~\ref{prop:app-ridge-tuning}]
Using \eqref{eq:app-precision-companion-rate} and the
uniform lower bounds established in the proof of
Proposition~\ref{prop:app-precision-loss}, we obtain
\begin{equation}\label{eq:app-gcv-denominator-rate}
    \sup_{\lambda\in\Lambda}
    \left|
        \frac{1}{\widehat v_n(-\lambda)^2}
        -
        \frac{1}{v_n(-\lambda)^2}
    \right|
    =
    O_{\mathbb P}(n^{-1}).
\end{equation}
We first consider the numerator in
\eqref{eq:app-gcv-companion-form}.
Using \eqref{eq:lin_mod} and the identities
$\underline R_n(-\lambda)X_n=X_nR_n(-\lambda)$ and
$\frac{d}{d\lambda}\widehat v_n(-\lambda)
=-n^{-1}\tr[\underline R_n(-\lambda)^2]$,
we obtain
\begin{align}\label{eq:GCV_dec_1}
\begin{split}
    \frac{1}{n}
    y_n^\top\underline R_n(-\lambda)^2y_n
    &=
    \beta_n^\top
    \widehat\Sigma_nR_n(-\lambda)^2\beta_n
    -
    \sigma_n^2\frac{d}{d\lambda}
    \widehat v_n(-\lambda)\\
    &\quad+
    \frac{2}{n}
    \beta_n^\top R_n(-\lambda)^2X_n^\top u_n\\
    &\quad+
    \frac{1}{n}
    \left\{
        u_n^\top\underline R_n(-\lambda)^2u_n
        -
        \sigma_n^2
        \tr[\underline R_n(-\lambda)^2]
    \right\}.
\end{split}
\end{align}

Conditional on $X_n$, the last two terms have mean zero. To control these terms uniformly over $\Lambda$, let
$\Gamma$ be a fixed rectangular contour enclosing
$\Lambda$ and contained in
$\{z\in\mathbb C:\operatorname{Re}z>\lambda_-/2\}$.
For $z\in\Gamma$, we have the bounds
\[
    \|R_n(-z)\|_2,\,
    \|\underline R_n(-z)\|_2
    \leq\frac{2}{\lambda_-},\qquad\text{and}\qquad
    \|\widehat\Sigma_nR_n(-z)^2
        R_n(-\overline z)^2\|_2
    \leq
    \left(1+|z|\frac{2}{\lambda_-}\right)
    \left(\frac{2}{\lambda_-}\right)^3
    \leq C.
\]
By independence of $u_n$ and $X_n$, we have
\begin{align}\label{eq:sec_first_equal}
    \mathbb E\!\left[
        \left|
            \frac{1}{n}
            \beta_n^\top R_n(-z)^2X_n^\top u_n
        \right|^2
        \,\middle|\,X_n
    \right]=
    \frac{\sigma_n^2}{n}
    \beta_n^\top R_n(-z)^2
    \widehat\Sigma_nR_n(-\overline z)^2\beta_n\leq\frac{C}{n}.
\end{align}
Here and below, $C$ denotes a constant independent of
$n$, $z\in\Gamma$, and the realization of $X_n$.
Similarly, Lemma~B.26 of
\citet{Silversteinbook2010}, applied to the real and
imaginary parts of the quadratic form, yields
\begin{align}\label{eq:sec_sec_equal}
\mathbb E\!\left[
        \left|
            \frac{1}{n}
            \left\{
                u_n^\top\underline R_n(-z)^2u_n
                -
                \sigma_n^2
                \tr[\underline R_n(-z)^2]
            \right\}
        \right|^2
        \,\middle|\,X_n
    \right]\leq
    \frac{C}{n^2}
    \tr\!\left[
        \underline R_n(-z)^2
        \underline R_n(-\overline z)^2
    \right]
    \leq\frac{C}{n}.
\end{align}
Consider the two centered terms as functions of the
complex parameter $z$,
\[
    f_{1,n}(z)
    \coloneqq
    \frac{2}{n}
    \beta_n^\top R_n(-z)^2X_n^\top u_n
\]
and
\[
    f_{2,n}(z)
    \coloneqq
    \frac{1}{n}
    \left\{
        u_n^\top\underline R_n(-z)^2u_n
        -
        \sigma_n^2\tr[\underline R_n(-z)^2]
    \right\}.
\]
Since both resolvents are holomorphic on
$\{z\in\mathbb C:\operatorname{Re}z>0\}$,
the functions $f_{1,n}$ and $f_{2,n}$ are holomorphic
in a neighborhood of the closed region enclosed by
$\Gamma$. Hence, for every $\lambda\in\Lambda$,
Cauchy's integral formula gives
\[
    f_{j,n}(\lambda)
    =
    \frac{1}{2\pi i}
    \int_\Gamma
    \frac{f_{j,n}(z)}{z-\lambda}\,dz,
    \qquad j=1,2.
\]
Since $\operatorname{dist}(\Lambda,\Gamma)>0$, it follows
that
\[
    \sup_{\lambda\in\Lambda}|f_{j,n}(\lambda)|
    \leq
    \frac{1}{2\pi\operatorname{dist}(\Lambda,\Gamma)}
    \int_\Gamma|f_{j,n}(z)|\,|dz|.
\]
Taking conditional expectations and applying
Cauchy--Schwarz, we obtain
\[
    \mathbb E\!\left[
        \sup_{\lambda\in\Lambda}|f_{j,n}(\lambda)|
        \,\middle|\,X_n
    \right]\leq
    \frac{1}{2\pi\operatorname{dist}(\Lambda,\Gamma)}
    \int_\Gamma
    \left(
        \mathbb E[|f_{j,n}(z)|^2\mid X_n]
    \right)^{1/2}|dz|\leq\frac{C}{\sqrt n},
    \qquad j=1,2,
\]
where the last inequality follows from \eqref{eq:sec_first_equal} and \eqref{eq:sec_sec_equal} and the fact that
$\Gamma$ is fixed. Taking expectations in the preceding display and applying
Markov's inequality yields
\[
    \sup_{\lambda\in\Lambda}|f_{j,n}(\lambda)|
    =
    O_{\mathbb P}(n^{-1/2}),
    \qquad j=1,2.
\]
Combining the preceding display with \eqref{eq:GCV_dec_1} and using the triangle
inequality, we conclude that
\begin{equation}\label{eq:app-gcv-noise-rate}
    \sup_{\lambda\in\Lambda}
    \biggl|
        \frac{1}{n}
        y_n^\top\underline R_n(-\lambda)^2y_n
        -
        \beta_n^\top
        \widehat\Sigma_nR_n(-\lambda)^2\beta_n
        +
        \sigma_n^2\frac{d}{d\lambda}
        \widehat v_n(-\lambda)
    \biggr|
    =
    O_{\mathbb P}(n^{-1/2}).
\end{equation}
Moreover, the identity
$\widehat\Sigma_nR_n(-\lambda)^2
=\frac{d}{d\lambda}[\lambda R_n(-\lambda)]$,
Theorem~\ref{Lemma:GreenRom}, and
Corollary~\ref{cor:GreenRom-derivatives:rates},
applied with $\Theta_n=\beta_n\beta_n^\top$, imply
\[
    \sup_{\lambda\in\Lambda}
    \left|
        \beta_n^\top
        \widehat\Sigma_nR_n(-\lambda)^2\beta_n
        -
        \frac{d}{d\lambda}
        \left[
            \lambda\beta_n^\top
            R_n^{\mathrm{det}}(-\lambda)\beta_n
        \right]
    \right|
    =
    O_{\mathbb P}(n^{-1/2}).
\]
Combining this estimate with
\eqref{eq:app-precision-companion-rate}
and~\eqref{eq:app-gcv-noise-rate}, and using
$\sup_n\sigma_n^2<\infty$, gives
\begin{equation}\label{eq:app-gcv-numerator-rate}
    \sup_{\lambda\in\Lambda}
    \biggl|
        \frac{1}{n}
        y_n^\top\underline R_n(-\lambda)^2y_n
        -
        \frac{d}{d\lambda}
        \left[
            \lambda\beta_n^\top
            R_n^{\mathrm{det}}(-\lambda)\beta_n
        \right]+
        \sigma_n^2\frac{d}{d\lambda}v_n(-\lambda)
    \biggr|
    =
    O_{\mathbb P}(n^{-1/2}).
\end{equation}

By \eqref{eq:def_R} and \eqref{def_K_theta},
\[
\begin{split}
    \frac{d}{d\lambda}
    \left[
        \lambda\beta_n^\top
        R_n^{\mathrm{det}}(-\lambda)\beta_n
    \right]&=
    -\frac{d}{d\lambda}v_n(-\lambda)\,
    \beta_n^\top
    (I_p+v_n(-\lambda)\Sigma_n)^{-1}
    \Sigma_n
    (I_p+v_n(-\lambda)\Sigma_n)^{-1}\beta_n\\
    &\quad=
    \lambda^2v_n(-\lambda)^2
    K_n^{\beta_n\beta_n^\top}(-\lambda,-\lambda).
\end{split}
\]
Moreover, differentiating \eqref{stielt}, with
$(\gamma_0,H_0)=(\gamma_n,H_n)$ and $\zeta=-\lambda$,
with respect to $\lambda$ and using the definition of $h_n$
yields
\begin{align}\label{eq:derivative_ident_vn}
    -\frac{
        \frac{d}{d\lambda}v_n(-\lambda)
    }{
        v_n(-\lambda)^2
    }
    =
    1+\gamma_n
    \left\{
        h_n(-\lambda)
        +
        \lambda\frac{d}{d\lambda}h_n(-\lambda)
    \right\}.
\end{align}

Together with \eqref{eq:app-ridge-risk-de}, these
identities give
\begin{equation}\label{eq:app-gcv-det-identity}
    \frac{d}{d\lambda}
    \left[
        \lambda\beta_n^\top
        R_n^{\mathrm{det}}(-\lambda)\beta_n
    \right]
    -
    \sigma_n^2\frac{d}{d\lambda}v_n(-\lambda)
    =
    v_n(-\lambda)^2
    \left\{
        \sigma_n^2+\mathcal R_n^{\mathrm{det}}(\lambda)
    \right\}.
\end{equation}
Moreover, the right-hand side is uniformly bounded
in $n$ and $\lambda\in\Lambda$. Indeed, \eqref{eq:app-gcv-det-identity} and the preceding
expression for the derivative of
$\lambda\beta_n^\top R_n^{\mathrm{det}}(-\lambda)\beta_n$,
together with
$0\leq-\frac{d}{d\lambda}v_n(-\lambda)\leq\lambda_-^{-2}$
and $\|(I_p+v_n(-\lambda)\Sigma_n)^{-1}\|_2\leq1$,
give
\begin{align}\label{eq:upper_det_risk_1}
    \left|
        v_n(-\lambda)^2
        \left\{
            \sigma_n^2+\mathcal R_n^{\mathrm{det}}(\lambda)
        \right\}
    \right|
    \leq
    \lambda_-^{-2}
    \left(
        \|\Sigma_n\|_2\|\beta_n\|_2^2+\sigma_n^2
    \right)
    \leq C.
\end{align}

Combining \eqref{eq:vn_lowerbound} and \eqref{eq:app-gcv-denominator-rate} yields
\begin{align}\label{eq:hatvn_OP1}
    \sup_{\lambda\in\Lambda}
    \frac{1}{\widehat v_n(-\lambda)^2}
    =
    O_{\mathbb P}(1).
\end{align}
Using \eqref{eq:app-gcv-companion-form}, \eqref{eq:app-gcv-numerator-rate}, \eqref{eq:app-gcv-det-identity}, and \eqref{eq:upper_det_risk_1}, we obtain
\[
\begin{split}
    &\sup_{\lambda\in\Lambda}
    \left|
        \operatorname{GCV}_n(\lambda)
        -\sigma_n^2-\mathcal R_n^{\mathrm{det}}(\lambda)
    \right|\\
    &\quad\leq
    \left(
        \sup_{\lambda\in\Lambda}
        \frac{1}{\widehat v_n(-\lambda)^2}
    \right)
    \sup_{\lambda\in\Lambda}
    \left|
        \frac{1}{n}
        y_n^\top\underline R_n(-\lambda)^2y_n
        -
        v_n(-\lambda)^2
        \left\{
            \sigma_n^2+\mathcal R_n^{\mathrm{det}}(\lambda)
        \right\}
    \right|\\
    &\qquad+
    C\sup_{\lambda\in\Lambda}
    \left|
        \frac{1}{\widehat v_n(-\lambda)^2}
        -
        \frac{1}{v_n(-\lambda)^2}
    \right|.
\end{split}
\]
By \eqref{eq:app-gcv-numerator-rate}, \eqref{eq:app-gcv-det-identity} and \eqref{eq:hatvn_OP1}, the first term in the upper bound is
$O_{\mathbb P}(n^{-1/2})$ and 
the second term is $O_{\mathbb P}(n^{-1})$ by
\eqref{eq:app-gcv-denominator-rate}.
Consequently,
\begin{equation}\label{eq:app-gcv-deterministic-rate}
    \sup_{\lambda\in\Lambda}
    \left|
        \operatorname{GCV}_n(\lambda)
        -
        \sigma_n^2
        -
        \mathcal R_n^{\mathrm{det}}(\lambda)
    \right|
    =
    O_{\mathbb P}(n^{-1/2}).
\end{equation}

We next relate this approximation to $\mathcal L_n$.
By the independence of $(x_0,y_0)$ and $(X_n,y_n)$,
\[
    \mathcal L_n(\lambda)
    =
    (\widehat\beta_n(\lambda)-\beta_n)^\top\Sigma_n(\widehat\beta_n(\lambda)-\beta_n)
    +
    \sigma_n^2.
\]
Using
$\widehat\beta_n(\lambda)-\beta_n
=-\lambda R_n(-\lambda)\beta_n
+n^{-1}R_n(-\lambda)X_n^\top u_n$
and \eqref{eq:app-ridge-risk}, we obtain
\begin{equation}\label{eq:app-gcv-error-expansion}
\begin{split}
    &\mathcal L_n(\lambda)-\sigma_n^2
    -\mathcal R_n(\lambda)\\
    &\quad=
    -\frac{2\lambda}{n}
    \beta_n^\top R_n(-\lambda)\Sigma_n
    R_n(-\lambda)X_n^\top u_n\\
    &\qquad+
    \frac{1}{n^2}
    \left\{
        u_n^\top X_nR_n(-\lambda)\Sigma_n
        R_n(-\lambda)X_n^\top u_n
        -
        \sigma_n^2\tr[
            X_nR_n(-\lambda)\Sigma_n
            R_n(-\lambda)X_n^\top]
    \right\}.
\end{split}
\end{equation}
Both terms after the equality have conditional
mean zero given $X_n$.

We apply the same contour argument as above.
For $z\in\Gamma$, we have the following bounds
\[
    \frac{1}{n}
    \|X_nR_n(-z)\Sigma_nR_n(-z)\beta_n\|_2^2
    \leq C
\qquad \text{and} \qquad
    \frac{1}{n^2}
    \|X_nR_n(-z)\Sigma_nR_n(-z)X_n^\top\|_F^2
    \leq Cn.
\]
Indeed, using
$\widehat\Sigma_nR_n(-z)=I_p-zR_n(-z)$,
we obtain
\[
    \frac{\|X_nR_n(-z)\|_2^2}{n}
    =
    \|R_n(-\overline z)
      \widehat\Sigma_nR_n(-z)\|_2\leq
    \|R_n(-z)\|_2
    \left(1+|z|\|R_n(-z)\|_2\right)
    \leq C,
\]
uniformly for $z\in\Gamma$.
Together with the uniform bounds on
$\|\Sigma_n\|_2$, $\|\beta_n\|_2$, and
$\|R_n(-z)\|_2$, this yields the first bound.
For the second bound, we use
\[
    \frac{1}{n}
    \|X_nR_n(-z)\Sigma_nR_n(-z)X_n^\top\|_2
    \leq
    \frac{\|X_nR_n(-z)\|_2^2}{n}\|\Sigma_n\|_2
    \leq C
\]
and $\|A\|_F\leq\sqrt n\|A\|_2$, for every $n\times n$ matrix $A$.

By independence of $u_n$ and $X_n$, we obtain
\[
    \mathbb E\!\left[
        \left|
            \frac{2z}{n}
            \beta_n^\top R_n(-z)\Sigma_n
            R_n(-z)X_n^\top u_n
        \right|^2
        \,\middle|\,X_n
    \right]=
    \frac{4\sigma_n^2|z|^2}{n^2}
    \|X_nR_n(-\overline z)\Sigma_n
    R_n(-\overline z)\beta_n\|_2^2
    \leq\frac{C}{n}.
\]
Similarly, Lemma~B.26 of
\citet{Silversteinbook2010} yields
\[
\begin{split}
    &\mathbb E\!\left[
        \left|
            \frac{1}{n^2}
            \left\{
                u_n^\top X_nR_n(-z)\Sigma_n
                R_n(-z)X_n^\top u_n
                -
                \sigma_n^2\tr[
                    X_nR_n(-z)\Sigma_n
                    R_n(-z)X_n^\top]
            \right\}
        \right|^2
        \,\middle|\,X_n
    \right]\\
    &\quad\leq
    \frac{C}{n^4}
    \|X_nR_n(-z)\Sigma_nR_n(-z)X_n^\top\|_F^2
    \leq\frac{C}{n}.
\end{split}
\]
Both centered terms are holomorphic in a neighborhood
of the region enclosed by $\Gamma$.
Applying Cauchy's integral formula and Cauchy--Schwarz
as above, followed by Markov's inequality, gives
\begin{equation}\label{eq:app-gcv-error-fluctuations}
    \sup_{\lambda\in\Lambda}
    \left|
        \mathcal L_n(\lambda)-\sigma_n^2
        -\mathcal R_n(\lambda)
    \right|
    =
    O_{\mathbb P}(n^{-1/2}).
\end{equation}

Finally, by the triangle inequality,
\[
\begin{split}
    \sup_{\lambda\in\Lambda}
    |\operatorname{GCV}_n(\lambda)-\mathcal L_n(\lambda)|&\leq
    \sup_{\lambda\in\Lambda}
    |\operatorname{GCV}_n(\lambda)
      -\sigma_n^2-\mathcal R_n^{\mathrm{det}}(\lambda)|\\
    &\qquad+
    \sup_{\lambda\in\Lambda}
    |\mathcal R_n^{\mathrm{det}}(\lambda)
      -\mathcal R_n(\lambda)|\\
    &\qquad+
    \sup_{\lambda\in\Lambda}
    |\sigma_n^2+\mathcal R_n(\lambda)
      -\mathcal L_n(\lambda)|.
\end{split}
\]
The first and third terms are $O_{\mathbb P}(n^{-1/2})$
by \eqref{eq:app-gcv-deterministic-rate} and \eqref{eq:app-gcv-error-fluctuations}, respectively.
The second term is of the same order by
Proposition~\ref{prop:app-ridge-risk}, since
$\|\beta_n\|_2$ and $\sigma_n^2$ are uniformly bounded.
This establishes \eqref{eq:app-gcv-empirical-de}.

For the second assertion, fix a realization of the training data.
Since $\widehat\lambda_n$ minimizes $\operatorname{GCV}_n$ over
$\Lambda$, for every $\lambda\in\Lambda$,
\[
\begin{split}
    \mathcal L_n(\widehat\lambda_n)-\mathcal L_n(\lambda)
    &\leq
    |\mathcal L_n(\widehat\lambda_n)
      -\operatorname{GCV}_n(\widehat\lambda_n)|
    +|\operatorname{GCV}_n(\lambda)-\mathcal L_n(\lambda)|\\
    &\leq 2\sup_{t\in\Lambda}
         |\operatorname{GCV}_n(t)-\mathcal L_n(t)|.
\end{split}
\]
Taking the infimum of $\mathcal L_n(\lambda)$ over $\Lambda$
yields the pathwise inequality
\[
    0
    \leq
    \mathcal L_n(\widehat\lambda_n)
    -
    \inf_{\lambda\in\Lambda}\mathcal L_n(\lambda)
    \leq
    2\sup_{\lambda\in\Lambda}
    \left|
        \operatorname{GCV}_n(\lambda)
        -
        \mathcal L_n(\lambda)
    \right|.
\]
The assertion now follows from
\eqref{eq:app-gcv-empirical-de}.
\end{proof}

\section{Proofs of the results in Section~\ref{sec:Rubmes_erroranalysis}}\label{sec:averaging-proofs}

\subsection{Proof of Lemma~\ref{RubMes_Lem4_fix}}

\begin{proof}[Proof of Lemma~\ref{RubMes_Lem4_fix}]

Without loss of generality, we prove the result under the first condition in
\eqref{cond:lem4}, namely 
\[
\sup_{n\ge1}\sup_{\zeta\in K_1}
        \|E_{n}(\zeta)\|_2<\infty,
        \qquad
\sup_{\zeta\in K_1}\|F_{n}(\zeta)\|_F\le C_{K_1}r_n.\]
For each $j \in \{1,\dots,n\}$ and $\zeta\in \mathbb{C}\setminus\mathbb{R}_{\geq 0}$, define
\[
Y_j(\zeta)
\coloneqq x_j^\top E_n R_{n,-j}(\zeta) F_n x_j
-
\tr\bigl(\Sigma_n E_n R_{n,-j}(\zeta) F_n\bigr)
\]
and 
\begin{align} \label{eq_def_S_n_zeta}
S_n(\zeta)\coloneqq \frac{1}{n}\sum_{j=1}^n Y_j(\zeta).
\end{align}
Since $K\subset\operatorname{int}(K_1)$, we can choose
$\rho>0$ such that $\{w\in\mathbb C:\operatorname{dist}(w,K)\leq\rho\}
    \subseteq K_1$.
The proof is divided into two steps. First, we establish
\begin{align}\label{first:step:lem4}
\sup_{\zeta\in K_1}
\mathbb E\bigl[|S_n(\zeta)|^2\bigr]
\le
C_{K_1} r_n^2\left(\frac{1}{n}\vee \frac{p}{n^2}\right).%
\end{align}
In the second step, we transfer the bound in \eqref{first:step:lem4} to a
uniform bound on \(K\). More precisely, we will show that
\begin{align}\label{sec:step:lem4}
\mathbb E\left[
\sup_{\zeta\in K}|S_n(\zeta)|^2
\right]
\le
C_{K_1} r_n^2 \left(\frac{1}{n}\vee \frac{p}{n^2}\right).%
\end{align}
Then, the proof of Lemma \ref{RubMes_Lem4_fix} can be completed using \eqref{first:step:lem4}, \eqref{sec:step:lem4} and Markov's inequality.
Indeed, these results give
\[
\sup_{\zeta\in K}|S_n(\zeta)|
=
O_{\mathbb P}\!\left(r_n\,\biggl\{
n^{-1/2}\vee \frac{\sqrt p}{n}\biggr\}
\right).
\]
Thus, it remains to prove \eqref{first:step:lem4} and \eqref{sec:step:lem4}.

We start by proving \eqref{first:step:lem4}.
When the dependence on $\zeta$ is suppressed below,
$E_n$ and $F_n$ mean $E_n(\zeta)$ and $F_n(\zeta)$. We decompose 
\begin{align}\label{dec1}
\mathbb{E}\!\left[\left|S_n(\zeta)\right|^2\right]=\left|\mathbb{E}\!\left[\left|S_n(\zeta)\right|^2\right]\right|\leq \frac{1}{n^2}\left|\sum_{j,k=1}^n\mathbb{E}\!\left[Y_j\overline{Y}_k\right]\right|
\leq
\frac{1}{n^2}\sum_{j=1}^n\mathbb{E}\bigl(|Y_j|^2\bigr)
+
\frac{1}{n^2}\sum_{\substack{j,k=1,\\
j\neq k}}^n|\mathbb{E}(Y_j\overline{Y}_k)|.
\end{align}
We first bound $\mathbb E|Y_j|^2$. Applying Lemma \ref{lem_b26_bai_silverstein} conditionally on
\(\hat\Sigma_{n,-j}\), and using Assumption~\ref{A3}, we obtain, for every 
\(j\in \{1,\dots,n\}\) and any \(\zeta\in K_1\), almost surely,
\[
\mathbb E\!\left(|Y_j(\zeta)|^2\,\middle|\,\hat\Sigma_{n,-j}\right)
\le C_{K_1}
\operatorname{tr}\!\left(
\Sigma_n^{1/2}E_{n}R_{n,-j}(\zeta)F_{n}\Sigma_nF_{n}^\ast R_{n,-j}(\zeta)^\ast E_{n}^\ast
\Sigma_n^{1/2}\right).
\]
Moreover,
\[
\operatorname{tr}\!\left(\Sigma_n^{1/2}
E_{n}R_{n,-j}(\zeta)F_{n}\Sigma_nF_{n}^\ast R_{n,-j}(\zeta)^\ast E_{n}^\ast
\Sigma_n^{1/2}\right)
=
\|\Sigma_n^{1/2}E_{n}R_{n,-j}(\zeta)F_{n}\Sigma_n^{1/2}\|_F^2.
\]
Using
$
\sup_{\zeta\in K_1}\sup_{1\le j\le n}\|R_{n,-j}(\zeta)\|_2\le C_{K_1}$, Assumption~\ref{A2}, and \eqref{cond:lem4}, we obtain, for every $\zeta\in K_1$
\[
\|\Sigma_n^{1/2}E_{n}R_{n,-j}(\zeta)F_{n}\Sigma_n^{1/2}\|_F^2
\le C_{K_1} \|F_{n}\|_F^2\leq C_{K_1} r_n^2.
\]
Consequently, the first term on the right-hand side of \eqref{dec1}
satisfies 
\begin{align}\label{stat1}
\sup_{\zeta\in K_1}
\frac{1}{n^2}\sum_{j=1}^n
\mathbb E\bigl(|Y_j(\zeta)|^2\bigr)
\le
\frac{C_{K_1} r_n^2}{n}.
\end{align}
We next turn to the second term on the right-hand side of \eqref{dec1}. 
Using the Sherman--Morrison formula, we obtain
\[
R_{n,-j}(\zeta)
=
R_{n,-jk}(\zeta)
-
\frac{\frac1n R_{n,-jk}(\zeta)x_k x_k^\top R_{n,-jk}(\zeta)}
{1+\frac1n x_k^\top R_{n,-jk}(\zeta)x_k}.
\]
Therefore \(Y_j(\zeta)=A_{jk}(\zeta)+\Delta_{jk}(\zeta)\), where
\begin{align*}
A_{jk}(\zeta)
&=
x_j^\top E_{n}(\zeta)R_{n,-jk}(\zeta)F_{n}(\zeta)x_j
-
\operatorname{tr}(\Sigma_n E_{n}(\zeta)R_{n,-jk}(\zeta)F_{n}(\zeta)),\\[2 mm]
\Delta_{jk}(\zeta)
&=
-
\frac{1}
{1+\frac1n x_k^\top R_{n,-jk}(\zeta)x_k}
\Bigl\{
\frac1n x_j^\top E_{n}(\zeta)R_{n,-jk}(\zeta)x_k x_k^\top R_{n,-jk}(\zeta)F_{n}(\zeta)x_j
\\[2 mm]
&\qquad-
\frac1n\operatorname{tr}
\left(\Sigma_n E_{n}(\zeta)R_{n,-jk}(\zeta)x_k x_k^\top R_{n,-jk}(\zeta)F_{n}(\zeta)\right)
\Bigr\}.
\end{align*}
Analogously, interchanging the roles of \(j\) and \(k\) yields the decomposition \(Y_k=A_{kj}+\Delta_{kj}\).
Conditional on $\mathcal F_{n,-jk}$, the random variables \(A_{jk}\) and
\(A_{kj}\) are independent and have mean zero. Hence
\[
\mathbb{E}\bigl[A_{jk}(\zeta)\overline{A_{kj}(\zeta)}\bigr]=\mathbb{E}\bigl[\mathbb E\bigl(A_{jk}(\zeta)\overline{A_{kj}(\zeta)}\mid \mathcal F_{n,-jk}\bigr)\bigr]=0.
\]
Moreover, \(\overline{A_{kj}(\zeta)}\) is \(\mathcal F_{n,-j}\)-measurable, while
$\E(\Delta_{jk}(\zeta)\mid \mathcal F_{n,-j})=0$ and therefore
\[
\E\left(\Delta_{jk}(\zeta)\overline{A_{kj}(\zeta)}\right)=0.
\]
Similarly, by symmetry
$\E\left(A_{jk}(\zeta)\overline{\Delta_{kj}(\zeta)}\right)=0$.
Consequently,
\[
\left|\E\left(Y_j(\zeta)\overline{Y_k(\zeta)}\right)\right|
=
\left|\E\left(\Delta_{jk}(\zeta)\overline{\Delta_{kj}(\zeta)}\right)\right|\leq \E\left(\left|\Delta_{jk}(\zeta)\overline{\Delta_{kj}(\zeta)}\right|\right).
\]
Using the Cauchy-Schwarz inequality, we get
\begin{align*}
\E(|\Delta_{jk}(\zeta)\overline{\Delta_{kj}(\zeta)}|)&\leq\E\bigl(|\Delta_{jk}(\zeta)|^2\bigr)^{1/2}\,\E\bigl(|\Delta_{kj}(\zeta)|^2\bigr)^{1/2}.
\end{align*}
Using \eqref{eq:general-denominator-bound}, we have,  for any $\zeta \in \mathbb{C}\setminus \mathbb{R}_{\geq 0}$,
\begin{align}\label{step0}
\begin{split}
\mathbb{E}(|\Delta_{jk}(\zeta)|^2)&\leq C \mathbb{E}\biggl[\Bigl|\frac{1}{n}x_j^\top E_{n}(\zeta) R_{n,-jk}(\zeta) x_k x_k^\top R_{n,-jk}(\zeta) F_{n}(\zeta)x_j\\
&\qquad-\frac{1}{n}\tr\bigl(\Sigma_n E_{n}(\zeta) R_{n,-jk}(\zeta) x_k x_k^\top R_{n,-jk}(\zeta)F_{n}(\zeta)\bigr)\Bigr|^2\biggr].
\end{split}
\end{align}
Conditionally on $\mathcal F_{n,-j}$, Lemma~\ref{lem_b26_bai_silverstein}
gives, almost surely,
\begin{equation*}%
\begin{split}
 &\mathbb E\left[
  \left|\frac1n x_j^\top E_nR_{n,-jk}x_kx_k^\top R_{n,-jk}F_nx_j
     -\frac1n\tr(\Sigma_nE_nR_{n,-jk}x_kx_k^\top R_{n,-jk}F_n)
  \right|^2\,\middle|\,\mathcal F_{n,-j}\right]\\
 &\qquad\leq\frac C{n^2}
   \|\Sigma_n^{1/2}E_nR_{n,-jk}x_kx_k^\top
      R_{n,-jk}F_n\Sigma_n^{1/2}\|_F^2
   =\frac C{n^2}\Delta_{jk,1}(\zeta)\Delta_{jk,2}(\zeta),
\end{split}
\end{equation*}
where every resolvent and deterministic matrix in this display
is evaluated at $\zeta$, and
\[
\begin{split}
    \Delta_{jk,1}(\zeta)
      &=x_k^\top R_{n,-jk}(\zeta)F_n(\zeta)\Sigma_n
          F_n(\zeta)^*R_{n,-jk}(\zeta)^*x_k,\\
    \Delta_{jk,2}(\zeta)
      &=x_k^\top R_{n,-jk}(\zeta)^*E_n(\zeta)^*\Sigma_n
          E_n(\zeta)R_{n,-jk}(\zeta)x_k.
\end{split}
\]
Using the Cauchy-Schwarz inequality for the conditional expectation, we obtain, for any $\zeta\in \mathbb{C}\setminus\mathbb{R}_{\geq 0}$, almost surely,
\begin{align}
    \label{step3}
\mathbb{E}\bigl[|\Delta_{jk,1} \Delta_{jk,2}|\big|\mathcal F_{n,-jk}\bigr]\leq \left\{ \mathbb{E}\bigl[|\Delta_{jk,1}|^2\big|\mathcal F_{n,-jk}\bigr] \mathbb{E}\bigl[|\Delta_{jk,2}|^2\big|\mathcal F_{n,-jk}\bigr] \right\}^{1/2}.
\end{align}
Using $\sup_{\zeta\in K_1}\sup_{1\le j<k\le n}
\|R_{n,-jk}(\zeta)\|_2\le C$, Assumption~\ref{A2} and \eqref{cond:lem4},
we obtain, conditionally on \(\hat\Sigma_{n,-jk}\), for every
\(\zeta\in K_1\), almost surely,
\[
\begin{aligned}
\mathbb E\!\left[
|\Delta_{jk,2}(\zeta)|^2
\,\middle|\,\mathcal F_{n,-jk}
\right]
&=
\mathbb E\!\left[
\left|
x_k^\top R_{n,-jk}(\zeta)^*
E_n(\zeta)^*\Sigma_n E_n(\zeta)
R_{n,-jk}(\zeta)x_k
\right|^2
\,\middle|\,\mathcal F_{n,-jk}
\right] \\
&\le
\mathbb E\!\left[
\|R_{n,-jk}(\zeta)^*E_n(\zeta)^*\Sigma_n E_n(\zeta)R_{n,-jk}(\zeta)\|_2^2
\|x_k\|_2^4
\,\middle|\,\mathcal F_{n,-jk}
\right] \\
&\le
C_{K_1}\,\mathbb E\!\left[
\|z_k\|_2^4
\,\middle|\,\mathcal F_{n,-jk}
\right] \\
&=
C_{K_1}\,\mathbb E\bigl(\|z_k\|_2^4\bigr) .
\end{aligned}
\]
By Assumption~\ref{A3},
\[
\mathbb E\bigl(\|z_k\|_2^4\bigr)
=
\mathbb E\left[\biggl(\sum_{\ell=1}^p z_{k,\ell}^2\biggr)^2\right] \leq
Cp^2 .
\]
Therefore, we have, for every $\zeta\in K_1$,
\[
\mathbb E\!\left[
|\Delta_{jk,2}(\zeta)|^2\right]
\le C_{K_1}p^2.
\]
Define 
\[
B_n=B_n(\zeta)\coloneqq \Sigma_n^{1/2}R_{n,-jk}(\zeta)F_n(\zeta)\Sigma_nF_n(\zeta)^\ast R_{n,-jk}(\zeta)^*\Sigma_n^{1/2}.
\]
The matrix $B_n=B_n(\zeta)$ is nonnegative definite and
$\mathcal F_{n,-jk}$-measurable. Using $|a+b|^2\leq2|a|^2+2|b|^2$
and Lemma~\ref{lem_b26_bai_silverstein} conditionally gives
\begin{align}\label{step5}
\begin{split}
    &\mathbb E[|\Delta_{jk,1}(\zeta)|^2\mid\mathcal F_{n,-jk}]\\
    &\quad\leq
    2\mathbb E[|z_k^\top B_n z_k-\tr B_n|^2\mid\mathcal F_{n,-jk}]
    +2|\tr B_n|^2\\
    &\quad\leq C_{K_1}\{\tr(B_nB_n^*)+|\tr B_n|^2\}\\
    &\quad\leq C_{K_1}\left\{
       \tr[(F_n(\zeta)F_n(\zeta)^*)^2]
       +\bigl(\tr[F_n(\zeta)F_n(\zeta)^*]\bigr)^2\right\}\\
    &\quad\leq C_{K_1}\|F_n(\zeta)\|_F^4.
\end{split}
\end{align}
Combining \eqref{step0}--\eqref{step5}, we obtain,
\begin{align*}
    \sup_{\zeta\in K_1}\mathbb E(|\Delta_{jk}(\zeta)|^2)
    &\leq\frac{C_{K_1}}{n^2}\sup_{\zeta\in K_1}
       \left\{\mathbb E(|\Delta_{jk,1}(\zeta)|^2)
              \mathbb E(|\Delta_{jk,2}(\zeta)|^2)\right\}^{1/2}\\
    &\leq C_{K_1}\frac p{n^2}
       \sup_{\zeta\in K_1}\|F_n(\zeta)\|_F^2
    \leq C_{K_1}\frac p{n^2}r_n^2.
\end{align*}

Similarly, by repeating the argument in \eqref{step0}--\eqref{step5} with the roles of
$j$ and $k$ interchanged, we obtain
\[
\sup_{\zeta\in K_1}\mathbb{E}\!\left(|\Delta_{kj}(\zeta)|^2\right)
\leq C_{K_1} \frac{p}{n^2}\sup_{\zeta\in K_1}\|F_n(\zeta)\|_F^2\leq C_{K_1}\frac{p}{n^2}r_n^2.
\]
Therefore,
\begin{align}\label{sup:K1:bound}
    \sup_{\zeta\in K_1}\mathbb{E}\left[\bigl|S_n(\zeta)\bigr|^2\right]\leq C_{K_1} r_n^2\biggl(n^{-1}\vee\frac{p}{n^2}\biggr).
\end{align}

This completes the proof of \eqref{first:step:lem4}, and it remains to establish \eqref{sec:step:lem4}.
Put $K_\rho\coloneqq\{w\in\mathbb C:\operatorname{dist}(w,K)\leq\rho\}\subseteq K_1$.
We apply Lemma~\ref{lem:cauchy_reduction_random} with $K_\rho$ as its enlarged compact set.
We verify that $S_n$ satisfies its assumptions. For each fixed
\(\zeta\in\mathbb C\setminus\mathbb R_{\ge0}\), the matrix
\(R_{n,-j}(\zeta)\) is a measurable function of \(\{z_i:i\ne j\}\). Hence
\(Y_j(\zeta)\) is \(\mathcal F\)-measurable, and so is
\(S_n(\zeta)=n^{-1}\sum_{j=1}^nY_j(\zeta)\). Moreover, for each
\(\omega\in\Omega\), the map
$\zeta\mapsto R_{n,-j}(\zeta,\omega)$
is holomorphic on \(\mathbb C\setminus\mathbb R_{\ge0}\), since the spectrum
of the corresponding leave-one-out sample covariance matrix is contained in
\(\mathbb R_{\ge0}\). Therefore, since \(E_n(\zeta)\) and \(F_n(\zeta)\) are deterministic and
holomorphic on \(\mathbb C\setminus\mathbb R_{\ge0}\) by assumption, the maps
\(\zeta\mapsto Y_j(\zeta,\omega)\) and hence
\(\zeta\mapsto S_n(\zeta,\omega)\) are holomorphic on
\(\mathbb C\setminus\mathbb R_{\ge0}\). Thus \(S_n\) satisfies the assumptions
of Lemma~\ref{lem:cauchy_reduction_random}. Together with
\eqref{sup:K1:bound}, this gives

\begin{align*}
\mathbb E\left(\sup_{\zeta\in K}|S_n(\zeta)|^2\right)
&\le
C_{K_\rho}\sup_{\zeta\in K_\rho}\mathbb E\bigl(|S_n(\zeta)|^2\bigr)\le
C_{K_1}\sup_{\zeta\in K_1}\mathbb E\bigl(|S_n(\zeta)|^2\bigr) \le
C_{K_1} r_n^2\left(n^{-1}\vee\frac{p}{n^2}\right).
\end{align*}
Thus, \eqref{sec:step:lem4} holds and the proof concludes.

\end{proof}

\subsection{Proof of Lemma~\ref{RubMes_Lem4_fix_2}}

\textit{Proof outline of  Lemma \ref{RubMes_Lem4_fix_2}}
We first explain the main idea of the proof before providing the details. The proof itself is rather technical and requires some additional notation. We therefore present the main steps using only the notation necessary to convey the underlying argument. Consider the same notation as in Lemma~\ref{RubMes_Lem4_fix}.
The goal is to show that, for $\zeta\in \mathbb C\setminus\mathbb R_{\geq 0}$,
\begin{align}
\label{eq:goal_lem}
\E \left| \frac1n\sum_{i=1}^n Y_i \right|^4 \leq C_\zeta n^{-2},
\end{align}

which in turn requires the analysis of terms of the form
\begin{align}
\label{eq:Y_terms}
\E \left[ Y_{a_1} Y_{a_2} \overline{Y_{a_3}} \overline{Y_{a_4}} \right], \quad 1 \leq a_1, a_2, a_3, a_4 \leq n.
\end{align}
Although $\E [ Y_a] =0$, this property alone is not immediately helpful for analyzing the mixed moment in \eqref{eq:Y_terms}. For instance, for $a_1 \neq a_2$, the random variables $Y_{a_1}$ and $Y_{a_2}$ are generally dependent since the resolvent term in $Y_{a_1}$ depends on $x_{a_2}$.
This motivates the decomposition
\begin{align*}
Y_{a_k} = \sum_{H \subset \{a_1, a_2, a_3, a_4\} \setminus { \{a_k \} }} X_{k,H}.
\end{align*}
We refrain from defining $X_{k,H}$ for this proof sketch and instead present two important properties.
First, for $q\in\{a_1,a_2,a_3,a_4\}\setminus(H\cup\{a_k\})$, $X_{k,H}$ is independent of $x_q$.
Second, the variables $X_{k,H}$ satisfy the moment inequality
\begin{align} \label{eq:mom_bound_X_kH}
\E \left[ |X_{k,H}|^4 \right] \leq C_\zeta n^{-2|H|}.
\end{align}
The second property builds on the Sherman-Morrison-Woodbury formula, which is routinely used for the analysis of resolvent terms (see, e.g., \cite{Silversteinbook2010,zheng_et_al_2015,zhang2022asymptotic}). Here, every additional observation $x_a$, $a \in H$, generates another rank-one resolvent difference and thus an additional factor $n^{-1/2}$.
To transfer the moment bound \eqref{eq:mom_bound_X_kH} to the term \eqref{eq:Y_terms}, we need to study
\begin{align} \label{eq:mixed_mom_X_kH}
\E \left[ X_{1, H_1} X_{2, H_2}
\overline{X_{3, H_3}}
\overline{X_{4, H_4}}
\right].
\end{align}
Let $Q$ denote the set of indices that appear exactly once in $(a_1, a_2, a_3, a_4)$. Suppose that $q\in Q$ does not belong to any $H_k$. Then, the moment \eqref{eq:mixed_mom_X_kH} is zero, since three of the random variables are measurable with respect to $\mathcal{F}_{n,-q}$ and one is conditionally centered with respect to $\mathcal{F}_{n,-q}$.
Therefore, only terms in \eqref{eq:mixed_mom_X_kH} with $Q\subset H_1 \cup H_2 \cup H_3 \cup H_4$ have a nontrivial contribution.
Thus, by Hölder's inequality and \eqref{eq:mom_bound_X_kH}, we obtain
\begin{align*}
\left| \E \left[ X_{1, H_1} X_{2, H_2}
\overline{X_{3, H_3}}
\overline{X_{4, H_4}}
\right] \right|
\leq C_\zeta n^{-|Q|/2}.
\end{align*}
Now the proof concludes by counting the number of quadruples $(a_1, a_2, a_3, a_4)$ that yield a non-zero term \eqref{eq:mixed_mom_X_kH}, which then shows \eqref{eq:goal_lem}.

\begin{proof}[Proof of Lemma~\ref{RubMes_Lem4_fix_2}]
The proof is divided into several parts. We start by fixing a \(\zeta\in\mathbb C\setminus\mathbb R_{\geq0}\), and suppress the
dependence on \(\zeta\) in the notation.

We first consider the case
\[
        \sup_{n\geq1}\|E_n\|_2<\infty,
        \qquad
        \sup_{n\geq1}\|F_n\|_F<\infty.
\]
The second case will be reduced to this one at the end of the proof.
 We shall prove that %
\begin{equation}
\label{eq:proof-fourth-moment-target}
        \E(|S_n|^4)\leq\frac{C_\zeta}{n^2}
\end{equation}
and conclude by Markov's inequality and the Borel--Cantelli lemma. Indeed,
for every \(\varepsilon>0\), Markov's inequality gives
\[
        \sum_{n=1}^\infty
        \mathbb P(|S_n|>\varepsilon)
        \leq
        \frac1{\varepsilon^4}
        \sum_{n=1}^\infty\E(|S_n|^4)
        \leq
        \frac{C_\zeta}{\varepsilon^4}
        \sum_{n=1}^\infty\frac1{n^2}
        <\infty.
\]
Hence, by the Borel--Cantelli lemma,
\[
        S_n\longrightarrow0
        \qquad\text{almost surely}.
\]
It therefore remains to prove \eqref{eq:proof-fourth-moment-target}.

\bigskip

Expanding the fourth power gives
\begin{equation}
\label{eq:proof-expand-fourth-power}
\begin{split}
        \E\left(\left|\sum_{a=1}^nY_a\right|^4\right)
        &=
        \sum_{a_1,a_2,a_3,a_4=1}^n
        \E\left[
        Y_{a_1}Y_{a_2}
        \overline{Y_{a_3}}\,
        \overline{Y_{a_4}}
        \right].
\end{split}
\end{equation}
Fix one ordered quadruple $(a_1,a_2,a_3,a_4)\in\{1,\ldots,n\}^4$ and set
\[
        I\coloneqq\{a_1,a_2,a_3,a_4\}.
\]

For \(k\in\{1,\ldots,4\}\) and
\(H\subseteq I\setminus\{a_k\}\), we define
\begin{equation}
\label{eq:proof-X-kH-definition}
        X_{k,H}
        \coloneqq
        \sum_{L\subseteq H}
        (-1)^{|H|-|L|}
        \left[
        x_{a_k}^\top
        E_nR_{n,-(I\setminus L)}F_n x_{a_k}
        -
        \tr\left\{\Sigma_n
        E_nR_{n,-(I\setminus L)}F_n
        \right\}
        \right].
\end{equation}
Note that the dependence of $X_{k,H}$ on $I$ is suppressed in the notation.
Next, we establish that the variables in \eqref{eq:proof-expand-fourth-power} admit the decomposition
\begin{equation}
\label{eq:proof-Y-component-decomposition}
        Y_{a_k}
        =
        \sum_{H\subseteq I\setminus\{a_k\}}X_{k,H}.
\end{equation}
Indeed, inserting \eqref{eq:proof-X-kH-definition} and
interchanging the finite sums gives
\begin{align*}
    \sum_{H\subseteq I\setminus\{a_k\}}X_{k,H}&=\sum_{L\subseteq I\setminus\{a_k\}}
    \left\{\sum_{H:\,L\subseteq H\subseteq I\setminus\{a_k\}}
             (-1)^{|H|-|L|}\right\}\\
    &\qquad\quad\times
    \left[x_{a_k}^\top E_nR_{n,-(I\setminus L)}F_nx_{a_k}
       -\tr\{\Sigma_nE_nR_{n,-(I\setminus L)}F_n\}\right].
\end{align*}

For fixed \(L\subseteq I\setminus\{a_k\}\), every set \(H\) satisfying $L\subseteq H\subseteq I\setminus\{a_k\}$ can be written uniquely as the union of disjoint sets $H=L\cup K$ with $K\subseteq I\setminus(\{a_k\}\cup L)$.
Since \(K\cap L=\varnothing\), we have $|H|-|L|=|K|$.
Consequently,
\[
\begin{split}
\sum_{\substack{
H:\,
L\subseteq H\subseteq I\setminus\{a_k\}
}}
(-1)^{|H|-|L|}
&=
\sum_{K\subseteq I\setminus(\{a_k\}\cup L)}
(-1)^{|K|}
\\
&=
\sum_{j=0}^{|I\setminus(\{a_k\}\cup L)|}
\binom{|I\setminus(\{a_k\}\cup L)|}{j}(-1)^j
\\
&=
(1-1)^{|I\setminus(\{a_k\}\cup L)|}.
\end{split}
\]
The last expression in the previous display is zero unless $I\setminus(\{a_k\}\cup L)=\varnothing$,
or equivalently, $L=I\setminus\{a_k\}$.
In that case, the coefficient equals one. The only surviving term is
therefore
\[
        x_{a_k}^\top E_nR_{n,-a_k}F_n x_{a_k}
        -
        \tr(\Sigma_n E_nR_{n,-a_k}F_n)
        =
        Y_{a_k},
\]
which proves \eqref{eq:proof-Y-component-decomposition}.
\bigskip

Recall that $\mathcal F_{n,-a}=\sigma(x_b:b\in \{1,\dots,n\}\setminus\{a\})$. Since \(H\subseteq I\setminus\{a_k\}\), we have $a_k\in I\setminus L$ for every $L\subseteq H$.
Thus every resolvent in \eqref{eq:proof-X-kH-definition} leaves out
\(a_k\) and is independent of \(x_{a_k}\). Therefore,
\begin{equation}
\label{eq:proof-component-centered}
        \E\left(
        X_{k,H}
        \,\middle|\,
        \mathcal F_{n,-a_k}
        \right)=0
\end{equation}
for every \(H\subseteq I\setminus\{a_k\}\). The same assertion holds after
complex conjugation.

We next prove that
\begin{equation}
\label{eq:proof-component-L4-bound}
        \E (|X_{k,H}|^4)^{1/4}
        \leq
        C n^{-|H|/2}
\end{equation}
for every \(H\subseteq I\setminus\{a_k\}\). Since
\(|I|\leq4\), we only have to consider the cases $|H|\in\{0,1,2,3\}$.

We first consider $H=\varnothing$. Conditionally on $\mathcal F_{n,-a_k}$ , the matrix $E_nR_{n,-I}F_n$ is deterministic and
independent of $x_{a_k}$. Hence,
Lemma~\ref{lem_b26_bai_silverstein} applied conditionally,  Assumption~\ref{A2} and \eqref{eq:general-resolvent-bound}, give
\begin{align}\label{eq:case:empty}
\begin{split}
        \E\left(
        |X_{k,\varnothing}|^4
        \,\middle|\,
        \mathcal F_{n,-a_k}
        \right)&=\E\left( | x_{a_k}^{\top}E_nR_{n,-I}F_n x_{a_k} - \tr(\Sigma_n E_nR_{n,-I}F_n) |^4 \,\middle|\, \mathcal F_{n,-a_k} \right)\\
        &\leq
        C\|\Sigma_n^{1/2}E_nR_{n,-I}F_n\Sigma_n^{1/2}\|_F^4\\
        &\leq
        C\|E_n\|_2^4
        \|R_{n,-I}\|_2^4
        \|F_n\|_F^4
        \leq C.
\end{split}        
\end{align}
Taking expectations yields $\{\mathbb E|X_{k,\varnothing}|^4\}^{1/4}\leq C$.

Now let \(H\neq\varnothing\), put \(h\coloneqq|H|\), and write
\begin{equation*}
        \mathcal R_H
        \coloneqq
        \sum_{L\subseteq H}
        (-1)^{h-|L|}
        R_{n,-(I\setminus L)}.
\end{equation*}
Then
\begin{equation*}
        X_{k,H}
        =
        x_{a_k}^\top E_n\mathcal R_HF_n x_{a_k}
        -
        \tr(\Sigma_n E_n\mathcal R_HF_n).
\end{equation*}
For any deletion set \(D\subseteq\{1,\dots,n\}\) and any \(q\notin D\), define
\[
        \beta_{D,q}
        \coloneqq
        \left\{
        1+\frac1n
        x_q^\top R_{n,-(D\cup\{q\})}x_q
        \right\}^{-1}.
\]

For each nonempty $H$, put $D=I\setminus H$.
We now bound $E_n\mathcal R_HF_n$ in the cases $|H|\in\{1,2,3\}$.
If $H=\{q\}$, using \eqref{eq:conditional-mixed-matrix} from Lemma~\ref{lem:conditional-moment-estimates}, \eqref{eq:general-denominator-bound} from Lemma~\ref{lem:resolvent-denominator-bounds}, \eqref{eq:proof-RH-one} from Lemma~\ref{lem:finite-resolvent-difference-expansions}, and Assumption~\ref{A2}, we have

\[
\E\left(\|E_n\mathcal R_HF_n\|_F^4\right)
\leq
\frac{C}{n^4}
\E\left[
\|E_nR_{n,-I}x_q\|_2^4
\|x_q^\top R_{n,-I}F_n\|_2^4
\right]
\leq
C\frac{p^2}{n^4}
\leq
C n^{-2}.
\]

Suppose that \(|H|=2\), with \(H=\{q,r\}\). 
Multiplying the first summand in \eqref{eq:proof-RH-two} from left by $E_n$ and from right with $F_n$, the Frobenius norm of this expression can be bounded by
\begin{align*}
&
\left\|
\frac1{n^2}
\beta_{D,q}\beta_{D\cup\{q\},r}
E_nR_{n,-(D\cup\{q\})}
x_qx_q^\top
R_{n,-I}x_rx_r^\top
R_{n,-I}F_n
\right\|_F
\\
&=
\frac1{n^2}
\left|
\beta_{D,q}\beta_{D\cup\{q\},r}
\right|
\left\|
\bigl(E_nR_{n,-(D\cup\{q\})}x_q\bigr)
\bigl(x_q^\top R_{n,-I}x_r\bigr)
\bigl(x_r^\top R_{n,-I}F_n\bigr)
\right\|_F
\\
&=
\frac1{n^2}
\left|
\beta_{D,q}\beta_{D\cup\{q\},r}
\right|
\|E_nR_{n,-(D\cup\{q\})}x_q\|_2
\left|
x_q^\top R_{n,-I}x_r
\right|
\|x_r^\top R_{n,-I}F_n\|_2
\\
&\leq
\frac{C}{n^2}
\|E_nR_{n,-(D\cup\{q\})}x_q\|_2
\left|
x_q^\top R_{n,-I}x_r
\right|
\|x_r^\top R_{n,-I}F_n\|_2.
\end{align*}

Using
\eqref{eq:conditional-mixed-vector}, \eqref{eq:general-denominator-bound}, and $\sup_{n\in \mathbb{N}}\|E_n\|_2<\infty,$ gives
\[
\begin{aligned}
    &\E\left[
       \|E_nR_{n,-(D\cup\{q\})}x_q\|_2^4
       |x_q^\top R_{n,-I}x_r|^4
       \,\middle|\,\mathcal F_{n,-q}\right]\\
    &\quad\leq Cp^2\|E_nR_{n,-(D\cup\{q\})}\|_2^4
                     \|R_{n,-I}x_r\|_2^4
    \leq Cp^2\|R_{n,-I}x_r\|_2^4.
\end{aligned}
\]
By \eqref{eq:conditional-mixed-matrix}, \eqref{eq:general-denominator-bound} and the Frobenius bound for $F_n$, we have
\[
\E\left[
\|R_{n,-I}x_r\|_2^4
\|x_r^\top R_{n,-I}F_n\|_2^4
\middle|
\mathcal F_{n,-r}\right]
\leq
Cp^2
\|R_{n,-I}\|_2^4
\|R_{n,-I}F_n\|_F^4\leq
C p^2.
\]
Combining the preceding three displays, we conclude that
\begin{align}
&\E\left(
\left\|
\frac1{n^2}
\beta_{D,q}\beta_{D\cup\{q\},r}
E_nR_{n,-(D\cup\{q\})}
x_qx_q^\top
R_{n,-I}x_rx_r^\top
R_{n,-I}F_n
\right\|_F^4
\right)
\nonumber\\[3mm]
&\qquad\leq
\frac{C}{n^8}
\E\left[
\|E_nR_{n,-(D\cup\{q\})}x_q\|_2^4
\left|
x_q^\top R_{n,-I}x_r
\right|^4
\|x_r^\top R_{n,-I}F_n\|_2^4
\right]
\nonumber\\[3mm]
&\qquad=
\frac{C}{n^8}
\E\left[
\E\left(
\|E_nR_{n,-(D\cup\{q\})}x_q\|_2^4
\left|
x_q^\top R_{n,-I}x_r
\right|^4
\|x_r^\top R_{n,-I}F_n\|_2^4
\,\middle|\,
\mathcal F_{n,-q}
\right)
\right]
\nonumber\\[3mm]
&\qquad\leq
\frac{Cp^2}{n^8}
\E\left[
\|R_{n,-I}x_r\|_2^4
\|x_r^\top R_{n,-I}F_n\|_2^4
\right]
\nonumber\\[3mm]
&\qquad=
\frac{Cp^2}{n^8}
\E\left[
\E\left(
\|R_{n,-I}x_r\|_2^4
\|x_r^\top R_{n,-I}F_n\|_2^4
\,\middle|\,
\mathcal F_{n,-r}
\right)
\right]
\nonumber\\[3mm]
&\qquad\leq
\frac{Cp^4}{n^8}.
\end{align}
The second summand in \eqref{eq:proof-RH-two} is treated identically after interchanging the roles of $q$ and $r$, and therefore satisfies the same bound. Since $\|A+B\|_F^4\leq 8\{\|A\|_F^4+\|B\|_F^4\}$ and $\gamma_n=O(1)$,
we obtain
\[
        \E\left(\|E_n\mathcal R_HF_n\|_F^4\right)
        \leq
        C\frac{p^4}{n^8}
        \leq
        C n^{-4}.
\]

Finally, suppose that \(|H|=3\), with \(H=\{q,r,s\}\). Multiplying the
first summand in \eqref{eq:proof-RH-three} from the left by \(E_n\) and
from the right by \(F_n\), we obtain
\[
\begin{aligned}
    &\left\|\frac1{n^3}
       \beta_{D,q}\beta_{D\cup\{q\},r}\beta_{D\cup\{q,r\},s}
       E_nR_{n,-(D\cup\{q\})}x_qx_q^\top
       R_{n,-(D\cup\{q,r\})}x_rx_r^\top
       R_{n,-I}x_sx_s^\top R_{n,-I}F_n\right\|_F\\
    &\quad=\frac1{n^3}
       |\beta_{D,q}\beta_{D\cup\{q\},r}\beta_{D\cup\{q,r\},s}|
       \|E_nR_{n,-(D\cup\{q\})}x_q\|_2\\
    &\qquad\quad\times
       |x_q^\top R_{n,-(D\cup\{q,r\})}x_r|
       |x_r^\top R_{n,-I}x_s|\,
       \|x_s^\top R_{n,-I}F_n\|_2\\
    &\quad\leq\frac C{n^3}\|E_nR_{n,-(D\cup\{q\})}x_q\|_2
       |x_q^\top R_{n,-(D\cup\{q,r\})}x_r|
       |x_r^\top R_{n,-I}x_s|\,
       \|x_s^\top R_{n,-I}F_n\|_2.
\end{aligned}
\]

Recall that in the case $|H|=3$, we have $D=I\setminus\{q,r,s\}$. Since $q\in D\cup\{q\}$, $q\in D\cup\{q,r\}$, $q\in I$,
all resolvents occurring in the preceding expression leave out \(q\).
Therefore, conditionally on \(\mathcal F_{n,-q}\), the $E_nR_{n,-(D\cup\{q\})}$ and $R_{n,-(D\cup\{q,r\})}x_r$ are independent of \(x_q\). Hence, by
\eqref{eq:conditional-mixed-vector},
\begin{align*}
\E\left[
\|E_nR_{n,-(D\cup\{q\})}x_q\|_2^4
\left|
x_q^\top R_{n,-(D\cup\{q,r\})}x_r
\right|^4
\,\middle|\,
\mathcal F_{n,-q}
\right]
&\leq
Cp^2
\|E_nR_{n,-(D\cup\{q\})}\|_2^4
\|R_{n,-(D\cup\{q,r\})}x_r\|_2^4
\\[2 mm]
&\leq
C p^2
\|R_{n,-(D\cup\{q,r\})}x_r\|_2^4.
\end{align*}

Similarly, since $r\in D\cup\{q,r\}$ and $r\in I$,
the resolvents \(R_{n,-(D\cup\{q,r\})}\) and \(R_{n,-I}\) leave out \(r\).
Thus, conditionally on \(\mathcal F_{n,-r}\),
\eqref{eq:conditional-mixed-vector} gives
\[
\E\left[
\|R_{n,-(D\cup\{q,r\})}x_r\|_2^4
\left|
x_r^\top R_{n,-I}x_s
\right|^4
\,\middle|\,
\mathcal F_{n,-r}
\right]\leq
Cp^2
\|R_{n,-(D\cup\{q,r\})}\|_2^4
\|R_{n,-I}x_s\|_2^4
\\[2 mm]
\leq
C p^2
\|R_{n,-I}x_s\|_2^4.
\]

Finally, since \(s\in I\), the resolvent \(R_{n,-I}\) leaves out \(s\). We may therefore apply
\eqref{eq:conditional-mixed-matrix},
to obtain
\[
\E\left[
\|R_{n,-I}x_s\|_2^4
\|x_s^\top R_{n,-I}F_n\|_2^4
\,\middle|\,
\mathcal F_{n,-s}
\right]\leq
Cp^2
\|R_{n,-I}\|_2^4
\|R_{n,-I}F_n\|_F^4
\leq
C p^2.
\]

Combining the preceding estimates and applying the tower property
successively, we obtain
\begin{align*}
&\E\left(
\left\|
-\frac1{n^3}
\beta_{D,q}
\beta_{D\cup\{q\},r}
\beta_{D\cup\{q,r\},s}
E_nR_{n,-(D\cup\{q\})}
x_qx_q^\top
R_{n,-(D\cup\{q,r\})}
x_rx_r^\top
R_{n,-I}
x_sx_s^\top
R_{n,-I}F_n
\right\|_F^4
\right)
\\
&\quad\leq
\frac{C}{n^{12}}
\E\left[
\|E_nR_{n,-(D\cup\{q\})}x_q\|_2^4
\left|
x_q^\top R_{n,-(D\cup\{q,r\})}x_r
\right|^4
\left|
x_r^\top R_{n,-I}x_s
\right|^4
\|x_s^\top R_{n,-I}F_n\|_2^4
\right]
\\
&\quad=
\frac{C}{n^{12}}
\E\left[
\left|
x_r^\top R_{n,-I}x_s
\right|^4
\|x_s^\top R_{n,-I}F_n\|_2^4\,
\E\left(
\|E_nR_{n,-(D\cup\{q\})}x_q\|_2^4
\left|
x_q^\top R_{n,-(D\cup\{q,r\})}x_r
\right|^4
\,\middle|\,
\mathcal F_{n,-q}
\right)
\right]
\\
&\quad\leq
\frac{Cp^2}{n^{12}}
\E\left[
\|R_{n,-(D\cup\{q,r\})}x_r\|_2^4
\left|
x_r^\top R_{n,-I}x_s
\right|^4
\|x_s^\top R_{n,-I}F_n\|_2^4
\right]
\\
&\quad=
\frac{Cp^2}{n^{12}}
\E\left[
\|x_s^\top R_{n,-I}F_n\|_2^4\,
\E\left(
\|R_{n,-(D\cup\{q,r\})}x_r\|_2^4
\left|
x_r^\top R_{n,-I}x_s
\right|^4
\,\middle|\,
\mathcal F_{n,-r}
\right)
\right]
\\
&\quad\leq
\frac{Cp^4}{n^{12}}
\E\left[
\|R_{n,-I}x_s\|_2^4\,
\|x_s^\top R_{n,-I}F_n\|_2^4
\right]
\\
&\quad=
\frac{Cp^4}{n^{12}}
\E\left[
\E\left(
\|R_{n,-I}x_s\|_2^4
\|x_s^\top R_{n,-I}F_n\|_2^4
\,\middle|\,
\mathcal F_{n,-s}
\right)
\right]
\\
&\quad\leq
\frac{Cp^6}{n^{12}}.
\end{align*}

The other five summands in \eqref{eq:proof-RH-three} are treated
identically after permuting the roles of \(q\), \(r\), and \(s\), and
therefore satisfy the same bound. 
By the triangle inequality and Hölder's inequality,
\[\left\|\sum_{j=1}^6A_j\right\|_F^4\leq
        \left(\sum_{j=1}^6\|A_j\|_F\right)^4\leq
        \left(\sum_{j=1}^6 1^{4/3}\right)^3
        \sum_{j=1}^6\|A_j\|_F^4
        =
        6^3\sum_{j=1}^6\|A_j\|_F^4.
\]
Since $p/n=O(1)$, we conclude that
\begin{align*}
        \E\left(\| E_n\mathcal R_HF_n\|_F^4\right)
        &\leq
        C \frac{p^6}{n^{12}}
        =
        C
        \left(\frac pn\right)^6
        \frac1{n^6}
        \leq
        C n^{-6}.
\end{align*}
Combining the three cases of $|H|$, we have proved that
\begin{equation}
\label{eq:proof-resolvent-component-bound}
\E\left(\|E_n\mathcal R_HF_n\|_F^4\right)
        \leq
        C n^{-2|H|}
\end{equation}
for every nonempty \(H\subseteq I\setminus\{a_k\}\).
Conditionally on \(\mathcal F_{n,-a_k}\), the matrix
\(E_n\mathcal R_HF_n\) is deterministic and independent of \(z_{a_k}\).
Therefore, Lemma~\ref{lem_b26_bai_silverstein} and
\eqref{eq:proof-resolvent-component-bound} yield
\[
        \E(|X_{k,H}|^4)
        \leq
        C\E\left(\|\Sigma_n^{1/2}E_n\mathcal R_HF_n\Sigma_n^{1/2}\|_F^4\right)\leq C\E\left(\|E_n\mathcal R_HF_n\|_F^4\right)
        \leq
        C n^{-2|H|}.
\]
Taking fourth roots, we have
\[
       \E(|X_{k,H}|^4)^{1/4}
        \leq
        C n^{-|H|/2}, \qquad\text{for }|H|\in \{1,2,3\}.
\]
Together with the case \(H=\varnothing\) in \eqref{eq:case:empty}, this proves
\eqref{eq:proof-component-L4-bound}.

We now return to the mixed fourth moment. By
\eqref{eq:proof-Y-component-decomposition},
\begin{equation}
\label{eq:proof-component-product-expansion}
\begin{split}
Y_{a_1}Y_{a_2}
\overline{Y_{a_3}}\,
\overline{Y_{a_4}}
=
\sum_{H_1\subseteq I\setminus\{a_1\}}
\sum_{H_2\subseteq I\setminus\{a_2\}}
\sum_{H_3\subseteq I\setminus\{a_3\}}
\sum_{H_4\subseteq I\setminus\{a_4\}}
X_{1,H_1}X_{2,H_2}
\overline{X_{3,H_3}}\,
\overline{X_{4,H_4}}.
\end{split}
\end{equation}
Since each set \(I\setminus\{a_k\}\) contains at most three elements, the
right-hand side contains only a constant number of terms.

Let \(Q\subseteq I\) denote the set of indices which occur exactly once among
\(a_1,a_2,a_3,a_4\). Fix one term in
\eqref{eq:proof-component-product-expansion}, and suppose that there exists a \(q\in Q\) such that
\begin{equation}
\label{eq:proof-uncovered-singleton}
        q\notin H_1\cup H_2\cup H_3\cup H_4.
\end{equation}
Let \(k_0\) be the unique index such that \(a_{k_0}=q\). By
\eqref{eq:proof-component-centered}, \(X_{k_0,H_{k_0}}\) is conditionally
centered with respect to \(\mathcal F_{n,-q}\), and the same holds for its
complex conjugate if \(k_0\in\{3,4\}\).

For every \(k\in \{1,\dots,4\}\setminus\{k_0\}\), we have $q\neq a_k$ by the definition of $Q$ and $q\notin H_k$ by \eqref{eq:proof-uncovered-singleton}.
Moreover, for every \(L\subseteq H_k\), the relation \(q\notin H_k\)
implies \(q\notin L\). Since \(q\in I\), it follows that $q\in I\setminus L$.
Now, recall the definition of $X_{k,H_k}$,
\[
\begin{split}
        X_{k,H_k}
        &=
        \sum_{L\subseteq H_k}
        (-1)^{|H_k|-|L|}
        \left[
        x_{a_k}^\top
        E_nR_{n,-(I\setminus L)}F_n x_{a_k}
        -
        \tr\left\{\Sigma_n
        E_nR_{n,-(I\setminus L)}F_n
        \right\}
        \right].
\end{split}
\]
For every \(L\subseteq H_k\), the inclusion \(q\in I\setminus L\)
shows that the resolvent \(R_{n,-(I\setminus L)}\) leaves out \(x_q\). Therefore, \(R_{n,-(I\setminus L)}\) is
\(\mathcal F_{n,-q}\)-measurable. Since \(a_k\neq q\), the vector
\(x_{a_k}\) is also \(\mathcal F_{n,-q}\)-measurable. Consequently, every
summand in the preceding display, and therefore \(X_{k,H_k}\), is
\(\mathcal F_{n,-q}\)-measurable.
The same is true after
complex conjugation. Consequently,
\[
\E\left[
X_{1,H_1}X_{2,H_2}
\overline{X_{3,H_3}}\,
\overline{X_{4,H_4}}
\right]=
\E\left[
\E\left(
X_{1,H_1}X_{2,H_2}
\overline{X_{3,H_3}}\,
\overline{X_{4,H_4}}
\,\middle|\,
\mathcal F_{n,-q}
\right)
\right]
=0.
\]
Therefore, a term can have nonzero expectation only if, for every
\(q\in Q\),
\begin{equation}
\label{eq:proof-covering-condition}
        q\in H_1\cup H_2\cup H_3\cup H_4.
\end{equation}
In particular, it follows that,
\begin{align}
\label{eq:proof-total-component-order}
        |H_1|+|H_2|+|H_3|+|H_4|
        \geq
        |Q|.
\end{align}
For every term satisfying \eqref{eq:proof-covering-condition}, Hölder's
inequality and \eqref{eq:proof-component-L4-bound} give
\begin{align}
\begin{split}
\left|
\E\left[
X_{1,H_1}X_{2,H_2}
\overline{X_{3,H_3}}\,
\overline{X_{4,H_4}}
\right]
\right|
&\leq
\prod_{k=1}^4\{\mathbb E|X_{k,H_k}|^4\}^{1/4}\\
&\leq
C
n^{-\{
|H_1|+|H_2|+|H_3|+|H_4|
\}/2}\\
&\leq
C n^{-|Q|/2}.
\label{eq:proof-one-component-product-bound}
\end{split}
\end{align}

Since the number of terms in
\eqref{eq:proof-component-product-expansion} is at most $2^{12}$, we conclude that
\begin{equation}
\label{eq:proof-mixed-fourth-bound}
        \left|
        \E\left[
        Y_{a_1}Y_{a_2}
        \overline{Y_{a_3}}\,
        \overline{Y_{a_4}}
        \right]
        \right|
        \leq
        C n^{-|Q|/2}.
\end{equation}

It remains to sum \eqref{eq:proof-mixed-fourth-bound} over the ordered
quadruples. If \(|Q|=0\), then either all four indices coincide or two
distinct indices occur twice each. There are \(O(n^2)\) such ordered
quadruples, and every mixed moment is bounded by \(C\).

If \(|Q|=1\), the multiplicity pattern is \(3+1\). There are \(O(n^2)\)
such ordered quadruples, and every mixed moment is bounded by
\(C n^{-1/2}\).

If \(|Q|=2\), the multiplicity pattern is \(2+1+1\). There are \(O(n^3)\)
such ordered quadruples, and every mixed moment is bounded by
\(C n^{-1}\).

The case \(|Q|=3\) is impossible.

Finally, if \(|Q|=4\), all four indices are distinct. There are \(O(n^4)\)
such ordered quadruples, and every mixed moment is bounded by
\(C n^{-2}\).

Using the triangle inequality in
\eqref{eq:proof-expand-fourth-power}, we obtain
\[
        \E\left(
        \left|\sum_{a=1}^nY_a\right|^4
        \right)
        \leq
        C
        \left(
        n^2
        +
        n^2n^{-1/2}
        +
        n^3n^{-1}
        +
        n^4n^{-2}
        \right)
        \leq
        C n^2.
\]
This proves \eqref{eq:proof-fourth-moment-target} under the first
condition. The same conclusion holds under the second condition by applying the preceding argument 
 with \(F_n^\top\) in place
of \(E_n\) and \(E_n^\top\) in place of \(F_n\).
\end{proof}

\subsection{Proof of Lemma~\ref{lem:ridge_sandwich}}

\begin{proof}[Proof of Lemma~\ref{lem:ridge_sandwich}]
Fix a compact set
$K\subseteq\mathbb C\setminus\mathbb R_{\geq0}$.
Choose $K_1$ as in Lemma~\ref{lem:cauchy_reduction_random}. For
$\zeta,\xi\in\mathbb C\setminus\mathbb R_{\geq0}$, define
\[
S_n(\zeta,\xi)
\coloneqq
\frac1n\sum_{j=1}^n Y_j(\zeta,\xi),
\]
where
\begin{align*}
Y_j(\zeta,\xi)
\coloneqq
z_j^\top\Sigma_n^{1/2} R_{n,-j}(\zeta)\Theta_nR_{n,-j}(\xi)\Sigma_n^{1/2}z_j-
\tr\!\left\{
\Sigma_nR_{n,-j}(\zeta)\Theta_nR_{n,-j}(\xi)
\right\}.
\end{align*}

First, we show that
\begin{align}\label{eq:dres:part1}
\begin{split}
        \sup_{\zeta,\xi\in K_1}
        \mathbb E(|S_n(\zeta,\xi)|^2)&=\sup_{\zeta,\xi\in K_1}
        \left| \E(|S_n(\zeta,\xi)|^2) \right|\\
        &\le
        \frac1{n^2}\sum_{j=1}^n
        \sup_{\zeta,\xi\in K_1}
        \mathbb E (|Y_j(\zeta,\xi)|^2)
        +
        \frac1{n^2}
        \sum_{\substack{j,k=1\\ j\neq k}}^n
        \sup_{\zeta,\xi\in K_1}
        \left|\mathbb E\left[
        Y_j(\zeta,\xi)\overline{Y_k(\zeta,\xi)}
        \right]\right|
        \\
        &\le
        C_K\|\Theta_n\|_F^2
        \left(
        \frac1n+\frac{p}{n^2}
        \right).
        \end{split}
\end{align}

For the diagonal part of $\E|S_n(\zeta,\xi)|^2$, we condition on
$\mathcal F_{n,-j}$. The matrix
\[
    \Sigma_n^{1/2}R_{n,-j}(\zeta)\Theta_n
    R_{n,-j}(\xi)\Sigma_n^{1/2}
\]
is measurable with respect to this sigma-algebra and independent of $z_j$.
Assumption~\ref{A2} and \eqref{eq:general-resolvent-bound} give
\[
\left\|
        \Sigma_n^{1/2}R_{n,-j}(\zeta)\Theta_nR_{n,-j}(\xi)\Sigma_n^{1/2}
\right\|_F
\leq
C_K\|\Theta_n\|_F .
\]
Hence, Lemma~\ref{lem_b26_bai_silverstein}, applied conditionally on
\(\mathcal F_{n,-j}\), gives
\[
        \mathbb E\!\left(
        |Y_j(\zeta,\xi)|^2
        \,\middle|\,
        \mathcal F_{n,-j}
        \right)
        \leq
        C_K\|\Theta_n\|_F^2 .
\]
Therefore,
\begin{align}\label{eq:diag:dres:bound}
        \frac1{n^2}
        \sum_{j=1}^n
        \mathbb E(|Y_j(\zeta,\xi)|^2)
        \leq
        \frac{C_K\|\Theta_n\|_F^2}{n}.
\end{align}

We next bound the off-diagonal terms. Fix $j\neq k$ and $w\in K_1$.
The Sherman--Morrison formula gives
\[
        R_{n,-j}(w)
        =
        R_{n,-jk}(w)
        -
        U_k(w),
\]
where
\[
        U_k(w)
        \coloneqq
        \frac{
        n^{-1}
        R_{n,-jk}(w)\Sigma_n^{1/2}z_kz_k^\top\Sigma_n^{1/2}
        R_{n,-jk}(w)}
        {1+n^{-1}z_k^\top\Sigma_n^{1/2}R_{n,-jk}(w)\Sigma_n^{1/2}z_k}.
\]
Define
\[
\begin{split}
A_{jk}(\zeta,\xi)
\coloneqq
z_j^\top\Sigma_n^{1/2}
R_{n,-jk}(\zeta)\Theta_nR_{n,-jk}(\xi)
\Sigma_n^{1/2}z_j
-
\tr\!\left(
\Sigma_nR_{n,-jk}(\zeta)\Theta_nR_{n,-jk}(\xi)
\right),
\end{split}
\]
and write
\[
        Y_j(\zeta,\xi)
        =
        A_{jk}(\zeta,\xi)
        +
        \Delta_{jk}(\zeta,\xi),
\]
where
\[
\begin{split}
\Delta_{jk}(\zeta,\xi)
&\coloneqq
z_j^\top\Sigma_n^{1/2}
\Big[
R_{n,-j}(\zeta)\Theta_nR_{n,-j}(\xi)
-
R_{n,-jk}(\zeta)\Theta_nR_{n,-jk}(\xi)
\Big]
\Sigma_n^{1/2}z_j
\\
&\quad-
\tr\!\left(
\Sigma_n
\Big[
R_{n,-j}(\zeta)\Theta_nR_{n,-j}(\xi)
-
R_{n,-jk}(\zeta)\Theta_nR_{n,-jk}(\xi)
\Big]
\right).
\end{split}
\]
Using \(R_{n,-j}(w)=R_{n,-jk}(w)-U_k(w)\), this becomes
\[
\begin{split}
\Delta_{jk}(\zeta,\xi)
={}&
z_j^\top\Sigma_n^{1/2}
\Big[
-U_k(\zeta)\Theta_nR_{n,-jk}(\xi)
-
R_{n,-jk}(\zeta)\Theta_nU_k(\xi)
+
U_k(\zeta)\Theta_nU_k(\xi)
\Big]
\Sigma_n^{1/2}z_j
\\
&-
\tr\!\left(
\Sigma_n
\Big[
-U_k(\zeta)\Theta_nR_{n,-jk}(\xi)
-
R_{n,-jk}(\zeta)\Theta_nU_k(\xi)
+
U_k(\zeta)\Theta_nU_k(\xi)
\Big]
\right).
\end{split}
\]
Similarly,
\[
        Y_k(\zeta,\xi)
        =
        A_{kj}(\zeta,\xi)
        +
        \Delta_{kj}(\zeta,\xi),
\]
where
\[
A_{kj}(\zeta,\xi)=
z_k^\top\Sigma_n^{1/2}
R_{n,-jk}(\zeta)\Theta_nR_{n,-jk}(\xi)
\Sigma_n^{1/2}z_k-
\tr\!\left(
\Sigma_nR_{n,-jk}(\zeta)\Theta_nR_{n,-jk}(\xi)
\right).
\]
Conditional on $\mathcal F_{n,-jk}$,
the matrix
$\Sigma_n^{1/2}
        R_{n,-jk}(\zeta)\Theta_nR_{n,-jk}(\xi)
        \Sigma_n^{1/2}$
is deterministic and independent of both \(z_j\) and \(z_k\). Hence,
\[
        \mathbb E\!\left(
        A_{jk}(\zeta,\xi)
        \,\middle|\,
        \mathcal F_{n,-jk}
        \right)=0,
        \qquad
        \mathbb E\!\left(
        \overline{A_{kj}(\zeta,\xi)}
        \,\middle|\,
        \mathcal F_{n,-jk}
        \right)=0.
\]
Moreover, conditionally on \(\mathcal F_{n,-jk}\), the random variables
\(A_{jk}(\zeta,\xi)\) and \(A_{kj}(\zeta,\xi)\) are independent. Therefore,
\[
\begin{split}
\mathbb E\!\left[
A_{jk}(\zeta,\xi)\overline{A_{kj}(\zeta,\xi)}
\right]
&=
\mathbb E\!\left[
\mathbb E\!\left(
A_{jk}(\zeta,\xi)\overline{A_{kj}(\zeta,\xi)}
\,\middle|\,
\mathcal F_{n,-jk}
\right)
\right]
\\[2 mm]
&=
\mathbb E\!\left[
\mathbb E\!\left(
A_{jk}(\zeta,\xi)
\,\middle|\,
\mathcal F_{n,-jk}
\right)
\mathbb E\!\left(
\overline{A_{kj}(\zeta,\xi)}
\,\middle|\,
\mathcal F_{n,-jk}
\right)
\right]
=0 .
\end{split}
\]
Furthermore, \(\overline{A_{kj}(\zeta,\xi)}\) is
\(\mathcal F_{n,-j}\)-measurable, while
\[
        \mathbb E\!\left(
        \Delta_{jk}(\zeta,\xi)
        \,\middle|\,
        \mathcal F_{n,-j}
        \right)=0,
\]
because \(\Delta_{jk}(\zeta,\xi)\) is a centered quadratic form in \(z_j\)
with a matrix independent of \(z_j\). Hence
\[
        \mathbb E\!\left[
        \Delta_{jk}(\zeta,\xi)\overline{A_{kj}(\zeta,\xi)}
        \right]=0.
\]
Similarly,
\[
        \mathbb E\!\left[
        A_{jk}(\zeta,\xi)\overline{\Delta_{kj}(\zeta,\xi)}
        \right]=0.
\]
Consequently, by the Cauchy--Schwarz inequality
\begin{align}\label{eq:dres:delta:bound}
\begin{split}
        \left|
        \mathbb E\left[
        Y_j(\zeta,\xi)\overline{Y_k(\zeta,\xi)}
        \right]
        \right|
        &=
        \bigl|
        \mathbb E\left[
        \Delta_{jk}(\zeta,\xi)
        \overline{\Delta_{kj}(\zeta,\xi)}
        \right]
        \bigr|  \\
        &\le
        \mathbb E\bigl[
        |\Delta_{jk}(\zeta,\xi)|
        |\Delta_{kj}(\zeta,\xi)|
        \bigr]\\
        &\le
        \left\{
        \mathbb E(|\Delta_{jk}(\zeta,\xi)|^2)
        \right\}^{1/2}
        \left\{
        \mathbb E(|\Delta_{kj}(\zeta,\xi)|^2)
        \right\}^{1/2}.
\end{split}
\end{align}

We now bound \(\E(|\Delta_{jk}(\zeta,\xi)|^2)\). Conditionally on
\(\mathcal F_{n,-j}\), it is a centered quadratic form in \(z_j\) with
coefficient matrix
\begin{align*}
B_{jk}(\zeta,\xi)
\coloneqq
\Sigma_n^{1/2}
\Big[
-U_k(\zeta)\Theta_nR_{n,-jk}(\xi)
-
R_{n,-jk}(\zeta)\Theta_nU_k(\xi)+
U_k(\zeta)\Theta_nU_k(\xi)
\Big]
\Sigma_n^{1/2}.
\end{align*}
Thus, by Lemma~\ref{lem_b26_bai_silverstein}, applied conditionally on
\(\mathcal F_{n,-j}\),
\[
        \mathbb E\!\left(
        |\Delta_{jk}(\zeta,\xi)|^2
        \,\middle|\,
        \mathcal F_{n,-j}
        \right)
        \leq
        C\|B_{jk}(\zeta,\xi)\|_F^2 .
\]
We now bound the terms in \(B_{jk}\).
Decompose
\[
        B_{jk}(\zeta,\xi)
        =
        B_{jk}^{(1)}(\zeta,\xi)
        +
        B_{jk}^{(2)}(\zeta,\xi)
        +
        B_{jk}^{(3)}(\zeta,\xi),
\]
where
\[
        B_{jk}^{(1)}(\zeta,\xi)
        \coloneqq
        -\Sigma_n^{1/2}U_k(\zeta)\Theta_nR_{n,-jk}(\xi)\Sigma_n^{1/2},
\]
\[
        B_{jk}^{(2)}(\zeta,\xi)
        :=
        -\Sigma_n^{1/2}R_{n,-jk}(\zeta)\Theta_nU_k(\xi)\Sigma_n^{1/2},
\]
and
\[
        B_{jk}^{(3)}(\zeta,\xi)
        :=
        \Sigma_n^{1/2}U_k(\zeta)\Theta_nU_k(\xi)\Sigma_n^{1/2}.
\]
We bound the three terms separately.

For the first term, using the rank-one form of \(U_k(\zeta)\) and \eqref{eq:general-denominator-bound}, we get
\[
\left\|B_{jk}^{(1)}(\zeta,\xi)\right\|_F^2
\le
\frac{C_K}{n^2}
\left\|
\Sigma_n^{1/2}R_{n,-jk}(\zeta) x_k
\right\|_2^2
\left\|
x_k^\top R_{n,-jk}(\zeta)\Theta_nR_{n,-jk}(\xi)\Sigma_n^{1/2}
\right\|_2^2.
\]
Conditioning on \(\mathcal F_{n,-jk}\), the matrices \(R_{n,-jk}(\zeta)\) and
\(R_{n,-jk}(\xi)\) are deterministic and \(x_k=\Sigma_n^{1/2}z_k\) is independent of
\(\mathcal F_{n,-jk}\). Hence, by the Cauchy--Schwarz inequality,
\begin{align}\label{eq:dres:B1}
\begin{split}
\mathbb E\left[
\left\|B_{jk}^{(1)}(\zeta,\xi)\right\|_F^2
\,\middle|\,
\mathcal F_{n,-jk}
\right]
&\le
\frac{C_K}{n^2}
\left\{
\mathbb E\left[
\left\|
\Sigma_n^{1/2}R_{n,-jk}(\zeta) x_k
\right\|_2^4
\,\middle|\,
\mathcal F_{n,-jk}
\right]
\right\}^{1/2}
\\
&\qquad\times
\left\{
\mathbb E\left[
\left\|
x_k^\top R_{n,-jk}(\zeta)\Theta_nR_{n,-jk}(\xi)\Sigma_n^{1/2}
\right\|_2^4
\,\middle|\,
\mathcal F_{n,-jk}
\right]
\right\}^{1/2}.
\end{split}
\end{align}
For the first factor on the right-hand side of \eqref{eq:dres:B1},
Assumption~\ref{A2} and \eqref{eq:general-resolvent-bound} give
\[
        \left\|
        \Sigma_n^{1/2}R_{n,-jk}(\zeta) x_k
        \right\|_2
        \le C_K\|z_k\|_2 .
\]
Hence, by Assumption~\ref{A3} and the preceding display,
\[
\mathbb E\left[
\left\|
\Sigma_n^{1/2}R_{n,-jk}(\zeta) x_k
\right\|_2^4
\,\middle|\,
\mathcal F_{n,-jk}
\right]
\le
C_K\mathbb E(\|z_k\|_2^4)
\le
C_Kp^2 .
\]
For the second factor on the right-hand side of \eqref{eq:dres:B1}, set
\[
        M_{1}
        \coloneqq
        R_{n,-jk}(\zeta)
        \Theta_n
        R_{n,-jk}(\xi)
        \Sigma_n
        R_{n,-jk}(\xi)^*
        \Theta_n^*
        R_{n,-jk}(\zeta)^*.
\]
Then,
\[
        \left\|
        x_k^\top R_{n,-jk}(\zeta)\Theta_nR_{n,-jk}(\xi)\Sigma_n^{1/2}
        \right\|_2^2
        =
        x_k^\top M_1x_k
        =
        z_k^\top \Sigma_n^{1/2}M_1\Sigma_n^{1/2}z_k .
\]
The matrix $M_1$ is Hermitian and nonnegative definite, and
\[
        \operatorname{tr}(\Sigma_n^{1/2}M_1\Sigma_n^{1/2})
        =
        \left\|
        \Sigma_n^{1/2}R_{n,-jk}(\zeta)\Theta_n
        R_{n,-jk}(\xi)\Sigma_n^{1/2}
        \right\|_F^2
        \le C_K\|\Theta_n\|_F^2.
\]
Consequently,
\[
        \operatorname{tr}\!\left[
        \left(\Sigma_n^{1/2}M_1\Sigma_n^{1/2}\right)
        \left(\Sigma_n^{1/2}M_1\Sigma_n^{1/2}\right)^*
        \right]
        \le
        \left\{
        \operatorname{tr}(\Sigma_n^{1/2}M_1\Sigma_n^{1/2})
        \right\}^2
        \le C_K\|\Theta_n\|_F^4 .
\]
Applying Lemma~\ref{lem_b26_bai_silverstein} conditionally on
$\mathcal F_{n,-jk}$ with $q=2$, and using
$|a+b|^2\leq2|a|^2+2|b|^2$, gives
\begin{align*}
&\mathbb E\left[
\left\|
x_k^\top R_{n,-jk}(\zeta)\Theta_nR_{n,-jk}(\xi)\Sigma_n^{1/2}
\right\|_2^4
\,\middle|\,
\mathcal F_{n,-jk}
\right]
\\
&\quad\leq
2\mathbb E\left[
\left|
z_k^\top\Sigma_n^{1/2}M_1\Sigma_n^{1/2}z_k
-
\tr(\Sigma_n^{1/2}M_1\Sigma_n^{1/2})
\right|^2
\,\middle|\,
\mathcal F_{n,-jk}
\right]
\\
&\qquad+
2\left|
\tr(\Sigma_n^{1/2}M_1\Sigma_n^{1/2})
\right|^2\\
&\quad\leq
C_K\|\Theta_n\|_F^4 .
\end{align*}
Combining the preceding arguments yields
\begin{align}\label{eq:expB1:bound}
        \mathbb E\left[
        \left\|B_{jk}^{(1)}(\zeta,\xi)\right\|_F^2
        \,\middle|\,
        \mathcal F_{n,-jk}
        \right]
        \le
        \frac{C_Kp}{n^2}\|\Theta_n\|_F^2.
\end{align}

For $B_{jk}^{(2)}(\zeta,\xi)$, we apply the preceding argument with
$\zeta$ and $\xi$ interchanged and $\Theta_n^\top$ in place of
$\Theta_n$, and transpose the resulting matrix. This is justified by
$R_{n,-jk}(w)^\top=R_{n,-jk}(w)$, $U_k(w)^\top=U_k(w)$, and
$\|\Theta_n^\top\|_F=\|\Theta_n\|_F$; no symmetry of $\Theta_n$
is required. More explicitly,
\[
\left\|B_{jk}^{(2)}(\zeta,\xi)\right\|_F^2
\le
\frac{C_K}{n^2}
\left\|
\Sigma_n^{1/2}R_{n,-jk}(\zeta)\Theta_nR_{n,-jk}(\xi) x_k
\right\|_2^2
\left\|
x_k^\top R_{n,-jk}(\xi)\Sigma_n^{1/2}
\right\|_2^2 .
\]
Here,
\[
\mathbb E\left[
\left\|
x_k^\top R_{n,-jk}(\xi)\Sigma_n^{1/2}
\right\|_2^4
\,\middle|\,
\mathcal F_{n,-jk}
\right]
\le C_Kp^2,
\]
while
\[
\mathbb E\left[
\left\|
\Sigma_n^{1/2}R_{n,-jk}(\zeta)\Theta_nR_{n,-jk}(\xi) x_k
\right\|_2^4
\,\middle|\,
\mathcal F_{n,-jk}
\right]
\le C_K\|\Theta_n\|_F^4 .
\]
Thus, again by the Cauchy--Schwarz inequality,
\begin{align}\label{eq:expB2:bound}
        \mathbb E\left[
        \left\|B_{jk}^{(2)}(\zeta,\xi)\right\|_F^2
        \,\middle|\,
        \mathcal F_{n,-jk}
        \right]
        \le
        \frac{C_Kp}{n^2}\|\Theta_n\|_F^2.
\end{align}        

For $B_{jk}^{(3)}$, first recall that
$U_k(\xi)=R_{n,-jk}(\xi)-R_{n,-j}(\xi)$,
and hence, by \eqref{eq:general-resolvent-bound}, we have
$\|U_k(\xi)\|_2\leq2/\operatorname{dist}(K_1,\mathbb R_{\geq0})\leq C_K$.
Therefore,
\[
\left\|B_{jk}^{(3)}(\zeta,\xi)\right\|_F^2
=
\left\|
\Sigma_n^{1/2}U_k(\zeta)\Theta_nU_k(\xi)\Sigma_n^{1/2}
\right\|_F^2\le
C_K
\left\|
\Sigma_n^{1/2}U_k(\zeta)\Theta_n
\right\|_F^2.
\]
Using \eqref{eq:general-denominator-bound}, we obtain
\[
\begin{split}
\left\|
\Sigma_n^{1/2}U_k(\zeta)\Theta_n
\right\|_F^2
&\le
\frac{C_K}{n^2}
\left\|
\Sigma_n^{1/2}R_{n,-jk}(\zeta) x_k
\right\|_2^2
\left\|
x_k^\top R_{n,-jk}(\zeta)\Theta_n
\right\|_2^2 .
\end{split}
\]
As above,
\[
\mathbb E\left[
\left\|
\Sigma_n^{1/2}R_{n,-jk}(\zeta) x_k
\right\|_2^4
\,\middle|\,
\mathcal F_{n,-jk}
\right]
\le C_Kp^2.
\]
Moreover,
\[
\mathbb E\left[
\left\|
x_k^\top R_{n,-jk}(\zeta)\Theta_n
\right\|_2^4
\,\middle|\,
\mathcal F_{n,-jk}
\right]
\le C_K\|\Theta_n\|_F^4 .
\]
Hence,
\begin{align}\label{eq:expB3:bound}
        \mathbb E\left[
        \left\|B_{jk}^{(3)}(\zeta,\xi)\right\|_F^2
        \,\middle|\,
        \mathcal F_{n,-jk}
        \right]
        \le
        \frac{C_Kp}{n^2}\|\Theta_n\|_F^2 .
\end{align}

Using the triangle inequality for the Frobenius norm and the elementary bound
$(a_1+a_2+a_3)^2\le 3(a_1^2+a_2^2+a_3^2)$, for $a_1,a_2,a_3\ge0$,
with \(a_\ell=\|B_{jk}^{(\ell)}\|_F\), \(\ell=1,2,3\), we obtain
\[
        \|B_{jk}\|_F^2
        =
        \left\|B_{jk}^{(1)}+B_{jk}^{(2)}+B_{jk}^{(3)}\right\|_F^2
        \le
        3\left(
        \|B_{jk}^{(1)}\|_F^2
        +
        \|B_{jk}^{(2)}\|_F^2
        +
        \|B_{jk}^{(3)}\|_F^2
        \right).
\]
Taking expectations and combining the bounds in \eqref{eq:expB1:bound}, \eqref{eq:expB2:bound} and \eqref{eq:expB3:bound} gives
\[
        \mathbb E(\|B_{jk}(\zeta,\xi)\|_F^2)
        \le
        \frac{C_Kp}{n^2}\|\Theta_n\|_F^2 .
\]
Since, conditionally on \(\mathcal F_{n,-j}\), \(\Delta_{jk}(\zeta,\xi)\) is
a centered quadratic form in \(z_j\) with coefficient matrix
\(B_{jk}(\zeta,\xi)\), Lemma~\ref{lem_b26_bai_silverstein} gives
\[
        \mathbb E(|\Delta_{jk}(\zeta,\xi)|^2)
        \le
        C\,\mathbb E(\|B_{jk}(\zeta,\xi)\|_F^2)
        \le
        \frac{C_Kp}{n^2}\|\Theta_n\|_F^2 .
\]
The same estimate holds for \(\Delta_{kj}(\zeta,\xi)\). Therefore, using \eqref{eq:dres:delta:bound} and the preceding arguments, we obtain
\[
\left|
        \mathbb E\!\left[
        Y_j(\zeta,\xi)\overline{Y_k(\zeta,\xi)}
        \right]
\right|
\leq
        \frac{C_Kp}{n^2}\|\Theta_n\|_F^2.
\]
Combining this with the diagonal bound \eqref{eq:diag:dres:bound}, and using
that all bounds hold uniformly over \(K_1\), we obtain the upper bound in
\eqref{eq:dres:part1}.

For every fixed $(\zeta,\xi)$, the variable $S_n(\zeta,\xi)$ is
$\mathcal F$-measurable. Moreover, for every $\omega\in\Omega$, the map
$(\zeta,\xi)\mapsto S_n(\zeta,\xi,\omega)$ is continuous on
$(\mathbb C\setminus\mathbb R_{\geq0})^2$ and holomorphic in each
variable separately, since each leave-one-out resolvent has these
properties. Thus, the two-dimensional part of
Lemma~\ref{lem:cauchy_reduction_random} applies and gives
\[\mathbb E\left(
        \sup_{\zeta,\xi\in K}|S_n(\zeta,\xi)|^2\right)
        \leq
        C_K\sup_{\zeta,\xi\in K_1} \mathbb E\left(
        |S_n(\zeta,\xi)|^2\right)
        \leq
        C_K\|\Theta_n\|_F^2
        \left(\frac1n+\frac{p}{n^2}\right).
\]
Markov's inequality therefore yields
\[
        \sup_{\zeta,\xi\in K}|S_n(\zeta,\xi)|
        =
        O_{\mathbb P}\left(
        \|\Theta_n\|_F
        \left\{\frac1{\sqrt n}\vee\frac{\sqrt p}{n}\right\}
        \right),
\]
which proves the assertion.
\end{proof}

\subsection{Proof of Lemma~\ref{lem:ridge_sandwich2}}

\begin{proof}[Proof of Lemma~\ref{lem:ridge_sandwich2}]
Fix $\zeta,\xi\in\mathbb C\setminus\mathbb R_{\geq0}$. Using $x_a=\Sigma_n^{1/2}z_a$, write
\[
        Y_a
        \coloneqq
        x_a^\top R_{n,-a}(\zeta)\Theta_nR_{n,-a}(\xi)x_a
        -
        \tr\!\left[
        \Sigma_nR_{n,-a}(\zeta)\Theta_nR_{n,-a}(\xi)
        \right],
        \qquad
        S_n
        \coloneqq
        \frac1n\sum_{a=1}^nY_a.
\]
As in the proof of Lemma~\ref{RubMes_Lem4_fix_2}, we shall prove that
\begin{align}\label{eq:sandwich-fourth-moment-target}
        \E(|S_n|^4)
        \leq
        \frac{C_{\zeta,\xi}}{n^2}
\end{align}
and conclude by Markov's inequality and the Borel--Cantelli lemma.
All constants below are independent of $n$ and of the deletion
indices.

Expanding the fourth power gives
\begin{align}\label{eq:sandwich-expand-fourth-power}
        \E\left|\sum_{a=1}^nY_a\right|^4
        =
        \sum_{a_1,a_2,a_3,a_4=1}^n
        \E\left[
        Y_{a_1}Y_{a_2}
        \overline{Y_{a_3}}\overline{Y_{a_4}}
        \right].
\end{align}
Fix one ordered quadruple $(a_1,a_2,a_3,a_4)$ and set
\[
        I\coloneqq\{a_1,a_2,a_3,a_4\}.
\]
For $k\in\{1,\ldots,4\}$ and $H\subseteq I\setminus\{a_k\}$,
define, as in \eqref{eq:proof-X-kH-definition},
\begin{align}\label{eq:sandwich-X-kH-definition}
\begin{split}
        X_{k,H}
        \coloneqq
        \sum_{L\subseteq H}
        (-1)^{|H|-|L|}
        \Bigl[
        &x_{a_k}^\top
        R_{n,-(I\setminus L)}(\zeta)\Theta_n
        R_{n,-(I\setminus L)}(\xi)x_{a_k}
        \\
        &-
        \tr\!\left\{
        \Sigma_n
        R_{n,-(I\setminus L)}(\zeta)\Theta_n
        R_{n,-(I\setminus L)}(\xi)
        \right\}
        \Bigr].
\end{split}
\end{align}
The same inclusion--exclusion calculation as in the proof of
\eqref{eq:proof-Y-component-decomposition} applies here. Indeed, the
argument only uses the dependence of the quadratic form on the deletion
set $I\setminus L$ and does not depend on the number of resolvent
factors. In the present setting, both resolvents are indexed by the same
deletion set $I\setminus L$. Thus, inserting the definition of
$X_{k,H}$ and interchanging the order of summation gives
\[
\begin{split}
\sum_{H\subseteq I\setminus\{a_k\}}X_{k,H}
&=
\sum_{L\subseteq I\setminus\{a_k\}}
\left\{
\sum_{\substack{H:\,
L\subseteq H\subseteq I\setminus\{a_k\}}}
(-1)^{|H|-|L|}
\right\}                                      \\
&\qquad\times
\Bigl[
x_{a_k}^\top
R_{n,-(I\setminus L)}(\zeta)\Theta_n
R_{n,-(I\setminus L)}(\xi)x_{a_k}
-
\tr\!\left\{
\Sigma_nR_{n,-(I\setminus L)}(\zeta)\Theta_n
R_{n,-(I\setminus L)}(\xi)
\right\}
\Bigr].
\end{split}
\]
For each fixed $L\subseteq I\setminus\{a_k\}$, the coefficient in
braces is
\[
        \sum_{\substack{H:\,
        L\subseteq H\subseteq I\setminus\{a_k\}}}
        (-1)^{|H|-|L|}
        =
        \sum_{V\subseteq I\setminus(\{a_k\}\cup L)}(-1)^{|V|}
        =
        \begin{cases}
        1, & L=I\setminus\{a_k\},\\
        0, & L\neq I\setminus\{a_k\}.
        \end{cases}
\]
Thus only the term with $I\setminus L=\{a_k\}$ survives. Its two
resolvents are $R_{n,-a_k}(\zeta)$ and $R_{n,-a_k}(\xi)$, and the
resulting centered quadratic form is precisely $Y_{a_k}$. Hence,
\begin{align}\label{eq:sandwich-Y-component-decomposition}
        Y_{a_k}
        =
        \sum_{H\subseteq I\setminus\{a_k\}}X_{k,H}.
\end{align}
Moreover, $a_k\in I\setminus L$ for every $L\subseteq H$.
Hence both resolvents in every summand of
\eqref{eq:sandwich-X-kH-definition} are
$\mathcal F_{n,-a_k}$-measurable. Their product with $\Theta_n$
is therefore independent of $x_{a_k}$. Since
\[
        \E(x_{a_k}x_{a_k}^\top\mid\mathcal F_{n,-a_k})=\Sigma_n,
\]
each quadratic form has conditional expectation equal to its
corresponding trace term. Consequently,
\begin{align}\label{eq:sandwich-component-centered}
        \E(X_{k,H}\mid\mathcal F_{n,-a_k})=0.
\end{align}
If $q\in I\setminus(\{a_k\}\cup H)$, then $q\in I\setminus L$
for every $L\subseteq H$. Both resolvents therefore also leave out
$q$, and $x_{a_k}$ is $\mathcal F_{n,-q}$-measurable because
$a_k\neq q$. Thus $X_{k,H}$ is $\mathcal F_{n,-q}$-measurable.
These assertions remain valid after complex conjugation.

It remains to establish the counterpart of
\eqref{eq:proof-component-L4-bound}, namely
\begin{align}\label{eq:sandwich-component-L4-bound}
        \E(|X_{k,H}|^4)^{1/4}
        \leq
        C_{\zeta,\xi}n^{-|H|/2}.
\end{align}
For $H=\varnothing$, Lemma~\ref{lem_b26_bai_silverstein},
applied conditionally on $\mathcal F_{n,-a_k}$, together with
Assumption~\ref{A2} and \eqref{eq:general-resolvent-bound}, gives
\begin{align*}
\E(|X_{k,\varnothing}|^4)
        \leq
        C\E\left\|
        \Sigma_n^{1/2}
        R_{n,-I}(\zeta)\Theta_nR_{n,-I}(\xi)
        \Sigma_n^{1/2}
        \right\|_F^4
        \leq
        C_{\zeta,\xi}\|\Theta_n\|_F^4
        \leq
        C_{\zeta,\xi}.
\end{align*}
For nonempty $H\subseteq I\setminus\{a_k\}$, write
\begin{align}\label{eq:sandwich-B-component}
        B_H
        \coloneqq
        \sum_{L\subseteq H}
        (-1)^{|H|-|L|}
        R_{n,-(I\setminus L)}(\zeta)\Theta_n
        R_{n,-(I\setminus L)}(\xi).
\end{align}
Then
\begin{align}\label{eq:sandwich-X-quadratic-form}
        X_{k,H}
        =
        x_{a_k}^\top B_Hx_{a_k}
        -
        \tr(\Sigma_nB_H).
\end{align}

To estimate $B_H$, we condition on the rows outside the set of indices
being added back. Under this conditioning, the test matrix appearing in
the one-resolvent term is fixed, and the moment calculation from
Lemma~\ref{RubMes_Lem4_fix_2} can be repeated conditionally. For $U\subseteq J\subseteq\{1,\ldots,n\}$ with $|U|\leq3$, write
\[
        \mathcal R_U^{\,J}(w)
        \coloneqq
        \sum_{L\subseteq U}
        (-1)^{|U|-|L|}
        R_{n,-(J\setminus L)}(w),
        \qquad
        w\in\{\zeta,\xi\}.
\]
Thus $\mathcal R_U^{\,J}$ is the finite difference from
Lemma~\ref{lem:finite-resolvent-difference-expansions}, with the
initial deletion set $J$ displayed explicitly. In particular,
$\mathcal R_\varnothing^{\,J}(w)=R_{n,-J}(w)$.

We claim that, for every $\mathcal F_{n,-U}$-measurable complex
$p\times p$ matrix $A$,
\begin{align}\label{eq:sandwich-conditional-difference-bound}
        \E\left(
        \|\mathcal R_U^{\,J}(w)A\|_F^4
        +
        \|A\mathcal R_U^{\,J}(w)\|_F^4
        \,\middle|\,
        \mathcal F_{n,-U}
        \right)
        \leq
        C_{\zeta,\xi}
        \frac{p^{2|U|}}{n^{4|U|}}
        \|A\|_F^4.
\end{align}
Conditionally on $\mathcal F_{n,-U}$, the matrix $A$ is fixed,
while the rows indexed by $U$ remain independent and satisfy the
same moment bounds. In particular, for every $q\in U$, $\mathcal F_{n,-U}\subseteq\mathcal F_{n,-q}$,
so $A$ remains measurable at every single-row conditioning step.
All conditional inequalities below hold almost surely.

For $U=\varnothing$, the claim follows from
\eqref{eq:general-resolvent-bound}. For nonempty $U$, set
$D=J\setminus U$, and recall the definition
from Lemma~\ref{lem:finite-resolvent-difference-expansions},
\[
        \beta_{G,q}(w)
        =
        \left\{
        1+\frac1n x_q^\top R_{n,-(G\cup\{q\})}(w)x_q
        \right\}^{-1},
        \qquad q\notin G.
\]
These factors are bounded uniformly in the deletion indices by
\eqref{eq:general-denominator-bound}.

First, let $U=\{q\}$. By \eqref{eq:proof-RH-one},
\[
        \mathcal R_{\{q\}}^{\,J}(w)A
        =
        -\frac{\beta_{D,q}(w)}n
        R_{n,-J}(w)x_qx_q^\top R_{n,-J}(w)A.
\]
Both $R_{n,-J}(w)$ and $A$ are
$\mathcal F_{n,-q}$-measurable. Therefore, the rank-one form of
this expression and \eqref{eq:conditional-mixed-matrix} give
\begin{align*}
        \E\left(
        \|\mathcal R_{\{q\}}^{\,J}(w)A\|_F^4
        \,\middle|\,
        \mathcal F_{n,-q}
        \right)
        &\leq
        \frac{C_{\zeta,\xi}}{n^4}
        \E\left(
        \|R_{n,-J}(w)x_q\|_2^4
        \|x_q^\top R_{n,-J}(w)A\|_2^4
        \,\middle|\,
        \mathcal F_{n,-q}
        \right)
        \\
        &\leq
        C_{\zeta,\xi}\frac{p^2}{n^4}\|A\|_F^4.
\end{align*}
Here we used $x_q=\Sigma_n^{1/2}z_q$, Assumption~\ref{A2},
and \eqref{eq:general-resolvent-bound}.

Next, let $U=\{q,r\}$. By \eqref{eq:proof-RH-two}, the first
of the two summands in $\mathcal R_U^{\,J}(w)A$ is
\[
        \frac1{n^2}
        \beta_{D,q}(w)\beta_{D\cup\{q\},r}(w)
        R_{n,-(D\cup\{q\})}(w)x_qx_q^\top
        R_{n,-J}(w)x_rx_r^\top R_{n,-J}(w)A.
\]
Using \eqref{eq:general-denominator-bound} to bound
$\beta_{D,q}(w)$ and $\beta_{D\cup\{q\},r}(w)$,
its conditional fourth moment is bounded by
\[
\frac{C_{\zeta,\xi}}{n^8}
        \E\Bigl[
        \|R_{n,-(D\cup\{q\})}(w)x_q\|_2^4
        |x_q^\top R_{n,-J}(w)x_r|^4
        \|x_r^\top R_{n,-J}(w)A\|_2^4
        \,\Bigm|\,
        \mathcal F_{n,-U}
        \Bigr].
\]
To integrate out $x_q$, note that both resolvents in this
expression leave out $q$. Hence
\eqref{eq:conditional-mixed-vector} gives
\begin{align*}
        \E\left[
        \|R_{n,-(D\cup\{q\})}(w)x_q\|_2^4
        |x_q^\top R_{n,-J}(w)x_r|^4
        \,\middle|\,
        \mathcal F_{n,-q}
        \right]
        \leq
        C_{\zeta,\xi}p^2\|R_{n,-J}(w)x_r\|_2^4.
\end{align*}
The remaining row-vector norm is
$\mathcal F_{n,-q}$-measurable. Moreover, since $R_{n,-J}(w)$
and $A$ also leave out $r$, \eqref{eq:conditional-mixed-matrix}
gives
\[
        \E\left[
        \|R_{n,-J}(w)x_r\|_2^4
        \|x_r^\top R_{n,-J}(w)A\|_2^4
        \,\middle|\,
        \mathcal F_{n,-r}
        \right]
        \leq C_{\zeta,\xi}p^2\|A\|_F^4.
\]
Using the tower property with
$\mathcal F_{n,-U}\subseteq\mathcal F_{n,-q}$ first, and then
with $\mathcal F_{n,-U}\subseteq\mathcal F_{n,-r}$, bounds the
conditional fourth moment of the displayed summand by
\[
        \frac{C_{\zeta,\xi}p^2}{n^8}
        \E\left[
        \|R_{n,-J}(w)x_r\|_2^4
        \|x_r^\top R_{n,-J}(w)A\|_2^4
        \,\middle|\,
        \mathcal F_{n,-U}
        \right]\leq
        C_{\zeta,\xi}\frac{p^4}{n^8}\|A\|_F^4.
\]
The second summand satisfies the same bound after interchanging
$q$ and $r$. Thus
\[
    \E[\|\mathcal R_U^{\,J}(w)A\|_F^4\mid\mathcal F_{n,-U}]
    \leq C_{\zeta,\xi}\frac{p^4}{n^8}\|A\|_F^4.
\]

Finally, let $U=\{q,r,s\}$. The first of the six summands in
\eqref{eq:proof-RH-three}, multiplied on the right by $A$, is
\[
\begin{split}
        -\frac1{n^3}
        &\beta_{D,q}(w)
        \beta_{D\cup\{q\},r}(w)
        \beta_{D\cup\{q,r\},s}(w)
        \\
        &\qquad\times
        R_{n,-(D\cup\{q\})}(w)x_qx_q^\top
        R_{n,-(D\cup\{q,r\})}(w)x_rx_r^\top
        R_{n,-J}(w)x_sx_s^\top R_{n,-J}(w)A.
\end{split}
\]
Its conditional fourth moment is bounded by
\[
\begin{split}
        \frac{C_{\zeta,\xi}}{n^{12}}
        \E\Bigl[
        &\|R_{n,-(D\cup\{q\})}(w)x_q\|_2^4
        |x_q^\top R_{n,-(D\cup\{q,r\})}(w)x_r|^4
        \\
        &\quad\times
        |x_r^\top R_{n,-J}(w)x_s|^4
        \|x_s^\top R_{n,-J}(w)A\|_2^4
        \,\Bigm|\,
        \mathcal F_{n,-U}
        \Bigr].
\end{split}
\]
All three resolvents leave out $q$. After integrating out $x_q$,
the remaining resolvents leave out $r$, and after integrating out
$x_r$, the remaining resolvent leaves out $s$. More explicitly,
\eqref{eq:conditional-mixed-vector} gives
\begin{align*}
        &\E\left[
        \|R_{n,-(D\cup\{q\})}(w)x_q\|_2^4
        |x_q^\top R_{n,-(D\cup\{q,r\})}(w)x_r|^4
        \,\middle|\,
        \mathcal F_{n,-q}
        \right]\leq
        C_{\zeta,\xi}p^2
        \|R_{n,-(D\cup\{q,r\})}(w)x_r\|_2^4,
        \\
        &\E\left[
        \|R_{n,-(D\cup\{q,r\})}(w)x_r\|_2^4
        |x_r^\top R_{n,-J}(w)x_s|^4
        \,\middle|\,
        \mathcal F_{n,-r}
        \right]\leq
        C_{\zeta,\xi}p^2\|R_{n,-J}(w)x_s\|_2^4.
\end{align*}
At the final step, \eqref{eq:conditional-mixed-matrix} gives
\[
        \E\left[
        \|R_{n,-J}(w)x_s\|_2^4
        \|x_s^\top R_{n,-J}(w)A\|_2^4
        \,\middle|\,
        \mathcal F_{n,-s}
        \right]
        \leq C_{\zeta,\xi}p^2\|A\|_F^4.
\]
At each step, the factors not being integrated are measurable
with respect to the corresponding sigma-field. Since
$\mathcal F_{n,-U}$ is contained in each of these sigma-fields,
the tower property yields the bound
\[
        C_{\zeta,\xi}\frac{p^6}{n^{12}}\|A\|_F^4
\]
for the conditional fourth moment of the displayed summand.
The other five summands are handled by integrating the three rows
in the order in which they occur in the corresponding expansion.
The triangle inequality and H\"older's inequality give the same
bound for their sum, after adjusting the constant.

This proves the bound for $\mathcal R_U^{\,J}(w)A$ in all cases.
Each resolvent is symmetric under transposition, and therefore
\[
        \mathcal R_U^{\,J}(w)^\top=\mathcal R_U^{\,J}(w),
        \qquad
        \|A\mathcal R_U^{\,J}(w)\|_F
        =
        \|\mathcal R_U^{\,J}(w)A^\top\|_F.
\]
Applying the preceding bound to $A^\top$ proves the bound for
$A\mathcal R_U^{\,J}(w)$, since
$\|A^\top\|_F=\|A\|_F$. No symmetry of $A$ is required.
Thus \eqref{eq:sandwich-conditional-difference-bound} holds.

We now use the following product identity:
\begin{align}\label{eq:sandwich-discrete-product}
        B_H
        =
        \sum_{U\subseteq H}
        \mathcal R_U^{\,I}(\zeta)\Theta_n
        \mathcal R_{H\setminus U}^{\,I\setminus U}(\xi).
\end{align}
To verify this identity, first note that
\[
        \mathcal R_{H\setminus U}^{\,I\setminus U}(\xi)
        =
        \sum_{U\subseteq V\subseteq H}
        (-1)^{|H|-|V|}
        R_{n,-(I\setminus V)}(\xi).
\]
Indeed, this follows from its definition by writing $V=U\cup W$
with $W\subseteq H\setminus U$. Expanding also
$\mathcal R_U^{\,I}(\zeta)$, the coefficient of
\[
        R_{n,-(I\setminus L)}(\zeta)\Theta_n
        R_{n,-(I\setminus V)}(\xi),
        \qquad
        L\subseteq V\subseteq H,
\]
is
\[
        (-1)^{|H|-|V|}
        \sum_{U:\,L\subseteq U\subseteq V}
        (-1)^{|U|-|L|}
        =
        \begin{cases}
        (-1)^{|H|-|L|}, & L=V,\\
        0, & L\neq V.
        \end{cases}
\]
Thus the expansion reduces to \eqref{eq:sandwich-B-component}.

We now apply the conditional estimate to each summand of
\eqref{eq:sandwich-discrete-product}. Fix $U\subseteq H$.
Every resolvent in $\mathcal R_U^{\,I}(\zeta)$ leaves out all
indices in $H\setminus U$, because its deletion set is
$I\setminus L$ with $L\subseteq U$. Consequently,
\[
        \mathcal R_U^{\,I}(\zeta)\Theta_n
        \quad\text{is }\mathcal F_{n,-(H\setminus U)}
        \text{-measurable}.
\]
The two finite differences in the product need not be independent.
Instead, condition on $\mathcal F_{n,-(H\setminus U)}$ and apply
\eqref{eq:sandwich-conditional-difference-bound} to the right-hand
finite difference, with deletion sets $H\setminus U\subseteq
I\setminus U$ and with the left-hand matrix
$A=\mathcal R_U^{\,I}(\zeta)\Theta_n$ held fixed. This gives
\begin{align*}
        \E\left[
        \left\|
        \mathcal R_U^{\,I}(\zeta)\Theta_n
        \mathcal R_{H\setminus U}^{\,I\setminus U}(\xi)
        \right\|_F^4
        \,\middle|\,
        \mathcal F_{n,-(H\setminus U)}
        \right]\leq
        C_{\zeta,\xi}
        \frac{p^{2|H\setminus U|}}{n^{4|H\setminus U|}}
        \|\mathcal R_U^{\,I}(\zeta)\Theta_n\|_F^4.
\end{align*}
Taking expectations and applying
\eqref{eq:sandwich-conditional-difference-bound} once more,
now to $\mathcal R_U^{\,I}(\zeta)\Theta_n$ with the deterministic
matrix $A=\Theta_n$, yields
\begin{align*}
        \E\left\|
        \mathcal R_U^{\,I}(\zeta)\Theta_n
        \mathcal R_{H\setminus U}^{\,I\setminus U}(\xi)
        \right\|_F^4
        &\leq
        C_{\zeta,\xi}
        \frac{p^{2|H\setminus U|}}{n^{4|H\setminus U|}}
        \E\|\mathcal R_U^{\,I}(\zeta)\Theta_n\|_F^4
        \\
        &\leq
        C_{\zeta,\xi}
        \frac{p^{2|H\setminus U|}}{n^{4|H\setminus U|}}
        \frac{p^{2|U|}}{n^{4|U|}}\|\Theta_n\|_F^4
        \\
        &=
        C_{\zeta,\xi}
        \frac{p^{2|H|}}{n^{4|H|}}\|\Theta_n\|_F^4
        \leq C_{\zeta,\xi}n^{-2|H|}.
\end{align*}
The last inequality uses $\sup_n\|\Theta_n\|_F<\infty$ and
$p/n=O(1)$, as assumed in the lemma.
This is the step that transfers the one-resolvent component
estimate to the product of two resolvents: the sets $U$ and
$H\setminus U$ are treated successively, using conditional
fourth moments rather than an independence assumption or
higher moments of the two factors. Since \eqref{eq:sandwich-discrete-product} contains at most
$2^{|H|}\leq8$ summands, the triangle inequality and
H\"older's inequality yield
\begin{align}\label{eq:sandwich-BH-fourth-bound}
        \E\|B_H\|_F^4
        \leq
        C_{\zeta,\xi}n^{-2|H|}.
\end{align}

Finally, $B_H$ is $\mathcal F_{n,-a_k}$-measurable.
Applying Lemma~\ref{lem_b26_bai_silverstein} conditionally on
$\mathcal F_{n,-a_k}$ to \eqref{eq:sandwich-X-quadratic-form},
and using Assumption~\ref{A2}, gives
\[
        \E(|X_{k,H}|^4)
        \leq
        C\E\|\Sigma_n^{1/2}B_H\Sigma_n^{1/2}\|_F^4
        \leq
        C_{\zeta,\xi}n^{-2|H|}.
\]
Together with the case $H=\varnothing$, this proves
\eqref{eq:sandwich-component-L4-bound}.

Combining \eqref{eq:sandwich-Y-component-decomposition},
\eqref{eq:sandwich-component-centered}, and
\eqref{eq:sandwich-component-L4-bound}, the argument from
\eqref{eq:proof-component-product-expansion} to
\eqref{eq:proof-mixed-fourth-bound} applies without modification.
Hence, if $Q$ denotes the set of indices which occur exactly once among
$a_1,a_2,a_3,a_4$, then
\[
        \left|
        \E\left[
        Y_{a_1}Y_{a_2}
        \overline{Y_{a_3}}\overline{Y_{a_4}}
        \right]
        \right|
        \leq
        C_{\zeta,\xi}n^{-|Q|/2}.
\]
The counting argument following
\eqref{eq:proof-mixed-fourth-bound} therefore gives
\[
        \E(|S_n|^4)
        \leq
        \frac{C_{\zeta,\xi}}{n^2}.
\]
This proves \eqref{eq:sandwich-fourth-moment-target}.
Consequently, Markov's inequality and the Borel--Cantelli lemma yield
\[
        S_n\longrightarrow0
        \qquad\text{almost surely}.
\]
This proves the assertion.
\end{proof}

\section{Auxiliary results}

\subsection{Measurability and Common Probability Null Sets}

The following elementary lemma ensures measurability of the random variables appearing in the proof of Lemma \ref{RubMes_Lem4_fix}, which are defined as suprema over possibly uncountable index sets.

\begin{lemma}\label{lem_sup_measurable}
Let $K$ be a nonempty compact metric space, and let
$\{S(t):t\in K\}$ be a family of $\mathcal F$-measurable
complex-valued random variables on $(\Omega,\mathcal F,\mathbb P)$.
Suppose that, for every $\omega\in\Omega$, the map
$t\mapsto S(t,\omega)$ is continuous on $K$.
Then $\sup_{t\in K}|S(t)|$ is finite and
$\mathcal F$-measurable.

In particular, the statement applies when the index space
is a nonempty compact subset of $\mathbb C$ or a Cartesian
product of such sets.
\end{lemma}

\begin{proof}
Since $K$ is a compact metric space, it has a countable
dense subset $D$.
For every $\omega\in\Omega$, we have
\[
    \sup_{t\in K}|S(t,\omega)|
    =
    \sup_{t\in D}|S(t,\omega)|.
\]
Indeed, for each $t\in K$, choose a sequence
$\{t_m\}_{m\geq1}\subseteq D$ with $t_m\to t$.
Continuity gives
\[
    |S(t,\omega)|
    =
    \lim_{m\to\infty}|S(t_m,\omega)|
    \leq
    \sup_{u\in D}|S(u,\omega)|.
\]
Taking the supremum over $t\in K$ proves one inequality;
the reverse follows from $D\subseteq K$.

Since $D$ is countable and $|S(t)|$ is
$\mathcal F$-measurable for every $t\in D$,
the supremum over $D$ is $\mathcal F$-measurable.
Finiteness follows from continuity and compactness.
\end{proof}

\begin{lemma}\label{lem:Georg_global_nullset}
Let $X$ be a separable topological space, let $(Y,d)$ be a
complete metric space, and let $(\Omega,\mathcal F,\mathbb P)$
be a probability space. Suppose that $f_n:\Omega\times X\to Y$
satisfies the following conditions:
\begin{enumerate}
    \item For every $x\in X$ and $n\geq1$,
    $\omega\mapsto f_n(\omega,x)$ is measurable from
    $(\Omega,\mathcal F)$ to $(Y,\mathcal B(Y))$.
    \item There is an event $E\in\mathcal F$ with
    $\mathbb P(E)=1$ such that, for every $\omega\in E$,
    the family $\{f_n(\omega,\cdot)\}_{n\geq1}$ is
    equicontinuous with respect to $d$.
\end{enumerate}
Suppose also that, for every $x\in X$, there is an event
$C_x\in\mathcal F$ with $\mathbb P(C_x)=1$ on which
$f_n(\cdot,x)$ converges.
Then there is an event $C\in\mathcal F$ with $\mathbb P(C)=1$
such that, for every $\omega\in C$, $f_n(\omega,x)$ converges
for all $x\in X$, and the resulting function of $x$ is continuous.
\end{lemma}

\begin{proof}
Choose a countable dense subset $D=\{x_1,x_2,\ldots\}$ of $X$
and set
\[
    C=E\cap\bigcap_{k\geq1}C_{x_k}.
\]
Then $C\in\mathcal F$ and $\mathbb P(C)=1$.
Fix $\omega\in C$ and write $g_n(x)=f_n(\omega,x)$.
To prove convergence at an arbitrary $x\in X$, let $\varepsilon>0$.
Equicontinuity gives a neighborhood $U$ of $x$ such that
\[
    d(g_n(x),g_n(y))<\varepsilon/3
    \qquad \text{for every } y\in U \text{and }\ n\geq1.
\]
Choose $x_k\in D\cap U$.
Since $g_n(x_k)$ converges, for all sufficiently large $n,m$,
\[
\begin{split}
    d(g_n(x),g_m(x))
    &\leq d(g_n(x),g_n(x_k))
       +d(g_n(x_k),g_m(x_k))
       +d(g_m(x_k),g_m(x))\\
    &<\varepsilon.
\end{split}
\]
Thus $g_n(x)$ is Cauchy in the specified complete metric $d$
and converges to a point $g(x)\in Y$.
To prove continuity of $g$, fix $x\in X$ and $\varepsilon>0$.
Choose a neighborhood $V$ of $x$ such that
$d(g_n(x),g_n(y))<\varepsilon/2$ for every $y\in V$ and $n\geq1$.
For every $y\in V$, passing to the limit in $n$ gives
\[
    d(g(x),g(y))\leq\varepsilon/2<\varepsilon.
\]
Hence $g$ is continuous.
\end{proof}

\subsection{Auxiliary resolvent and moment bounds}\label{sec:auxilliary_res}

\begin{lemma}[Lemma B.26 of \cite{Silversteinbook2010}] \label{lem_b26_bai_silverstein}
	Let $q\geq 2$, let $A$ be an $n \times n $ non-random matrix and $x=(x_1,...,x_n)^\top$ be a random vector with independent entries. Assume that $\mathbb{E}(x_i)=0$, $\mathbb{E}(x_i^2)=1$, $\mathbb E|x_i|^4\leq\nu_4<\infty$, and
$\mathbb E|x_i|^{2q}\leq\nu_{2q}<\infty$ for every $1 \leq i \leq n$. Then ,
	\begin{align*}
		\mathbb{E}(|x^\top A x-\tr(A)|^q)\leq C_q \bigg((\nu_4 \tr(AA^*))^{q/2}+ \nu_{2q} \tr((AA^*)^{q/2})\bigg),
	\end{align*}
	where $C_q$ is a constant depending on $q$ only.
\end{lemma}

\begin{lemma}
\label{lem:cauchy_reduction_random}
Let \(K\subset \mathbb C\setminus\mathbb R_{\ge0}\) be compact. Choose
\(0<\rho<\operatorname{dist}(K,\mathbb R_{\ge0})\) and define
\[
        K_1:=\{\zeta\in\mathbb C:\operatorname{dist}(\zeta,K)\le \rho\}.
\]
Suppose first that, for each \(n\in\mathbb N\),
$\{S_n(\zeta):\zeta\in\mathbb C\setminus\mathbb R_{\ge0}\}$
is a family of complex-valued random variables such that \(S_n(\zeta)\) is
\(\mathcal F\)-measurable for each
\(\zeta\in\mathbb C\setminus\mathbb R_{\ge0}\), and such that, for every
\(\omega\in\Omega\), the map $\zeta\mapsto S_n(\zeta,\omega)$
is holomorphic on \(\mathbb C\setminus\mathbb R_{\ge0}\). Then $\sup_{\zeta\in K}|S_n(\zeta)|$
is \(\mathcal F\)-measurable, and there exists a constant
\(C_{K_1}<\infty\) such that, for every \(n\in\mathbb N\),
\[
        \mathbb E\left[
        \sup_{\zeta\in K}|S_n(\zeta)|^2
        \right]
        \le
        C_{K_1}
        \sup_{\zeta\in K_1}
        \mathbb E\left(|S_n(\zeta)|^2\right).
\]
Moreover, suppose that, for each \(n\in\mathbb N\), $\{S_n(\zeta,\xi):(\zeta,\xi)\in (\mathbb C\setminus\mathbb R_{\ge0})^2\}$
is a family of complex-valued random variables such that \(S_n(\zeta,\xi)\)
is \(\mathcal F\)-measurable for each
\((\zeta,\xi)\in(\mathbb C\setminus\mathbb R_{\ge0})^2\), and such that, for
every \(\omega\in\Omega\), the map $(\zeta,\xi)\mapsto S_n(\zeta,\xi,\omega)$
is continuous on $(\mathbb C\setminus\mathbb R_{\ge0})^2$ and holomorphic in each variable separately. Then $\sup_{\zeta\in K}\sup_{\xi\in K}|S_n(\zeta,\xi)|$
is \(\mathcal F\)-measurable, and there exists a constant
\(C_{K_1}<\infty\) such that, for every \(n\in\mathbb N\),
\[
        \mathbb E\left[
        \sup_{\zeta\in K}\sup_{\xi\in K}|S_n(\zeta,\xi)|^2
        \right]
        \le
        C_{K_1}
        \sup_{\zeta\in K_1}\sup_{\xi\in K_1}
        \mathbb E\left(|S_n(\zeta,\xi)|^2\right).
\]
\end{lemma}

\begin{proof}
By construction, the set \(K_1\) is compact and satisfies \begin{align}\label{compact:stat:cauchy} K\subseteq K_1\subseteq \mathbb C\setminus\mathbb R_{\ge0}. 
\end{align}
Fix an $\omega\in \Omega$. 
Since \(K\) is
compact, there exist finitely many points \(a_1,\ldots,a_m\in K\) such that
\[
K\subseteq \bigcup_{\ell=1}^m B(a_\ell,\rho/4).
\]
For each \(\ell \in \{1, \ldots, m\}\), let
\[
\Gamma_\ell:=\{w\in\mathbb C: |w-a_\ell|=\rho/2\}.
\]
Recall the definition of $K_1$. Then \(\Gamma_\ell\subset K_1\), because for \(w\in\Gamma_\ell\),
\[
\operatorname{dist}(w,K)\le |w-a_\ell|=\rho/2<\rho.
\]
Moreover, the closed disk \(\overline{B(a_\ell,\rho/2)}\coloneqq\{w\in\mathbb C: |w-a_\ell|\leq\rho/2\}\) is contained in
\(K_1\subseteq \mathbb{C}\setminus\mathbb{R}_{\geq 0}\). In particular, \(S_n(\cdot,\omega)\) is holomorphic on an open neighborhood of
\(\overline{B(a_\ell,\rho/2)}\).
Now fix \(\zeta\in K\), and choose \(\ell\) such that
\(\zeta\in B(a_\ell,\rho/4)\). By Cauchy's integral formula applied on
\(\Gamma_\ell\),
\[
S_n(\zeta)
=
\frac{1}{2\pi i}
\int_{\Gamma_\ell}
\frac{S_n(w)}{w-\zeta}\,dw .
\]
For \(w\in\Gamma_\ell\), we have
\[
|w-\zeta|
\ge
|w-a_\ell|-|\zeta-a_\ell|
\ge
\rho/2-\rho/4
=
\rho/4.
\]
Hence, combining the previous two displays yields the following
\[
|S_n(\zeta)|
\le
\frac{2}{\pi\rho}
\int_{\Gamma_\ell}|S_n(w)|\,|dw|.
\]
Hence,
\[
\sup_{\zeta\in K}|S_n(\zeta)|
\le
C
\max_{1\le \ell\le m}
\int_{\Gamma_\ell}|S_n(w)|\,|dw|.
\]
Consequently,
\begin{align} \label{c_r1}
\sup_{\zeta\in K}|S_n(\zeta)|^2
\le
C
\max_{1\le \ell\le m}
\left(
\int_{\Gamma_\ell}|S_n(w)|\,|dw|
\right)^2
\le
C
\sum_{\ell=1}^m
\left(
\int_{\Gamma_\ell}|S_n(w)|\,|dw|
\right)^2 .
\end{align}
By Cauchy--Schwarz,
\begin{align} \label{c_r2}
\left(
\int_{\Gamma_\ell}|S_n(w)|\,|dw|
\right)^2
\le
\left(\int_{\Gamma_\ell}|S_n(w)|^2\,|dw|\right)
\left(\int_{\Gamma_\ell}|dw|\right)
=
|\Gamma_\ell|
\int_{\Gamma_\ell}|S_n(w)|^2\,|dw|.
\end{align}
Since \(|\Gamma_\ell|=2\pi(\rho/2)=\pi\rho\) for every
\(\ell\in\{1,\dots,m\}\), \eqref{c_r1} and \eqref{c_r2} imply
\[
\sup_{\zeta\in K}|S_n(\zeta)|^2
\le
C
\sum_{\ell=1}^m
\int_{\Gamma_\ell}|S_n(w)|^2\,|dw|.
\]
We point out that the preceding argument is deterministic and holds for every
\(\omega\in\Omega\). Consequently, we conclude that for every \(n\) and \(\omega\in\Omega\),
\begin{align}\label{contour:bound}
\sup_{\zeta\in K}|S_n(\zeta,\omega)|^2
\le
C
\sum_{\ell=1}^m
\int_{\Gamma_\ell}|S_n(w,\omega)|^2\,|dw|.
\end{align}
We are now going to take expectations in \eqref{contour:bound}. Using Lemma~\ref{lem_sup_measurable} it is straightforward to verify that $\sup_{\zeta\in K}|S_n(\zeta)|$ is $\mathcal F$-measurable, since for each $\omega\in \Omega$ the map $\zeta\mapsto S_{n}(\zeta,\omega)$ is holomorphic on $\mathbb{C}\setminus\mathbb{R}_{\geq 0}$. To argue that the expression on the right hand side of \eqref{contour:bound} is $\mathcal F$-measurable, we may use the following parameterization, for arbitrary $\ell\in \{1,\dots,m\}$,
\[
\gamma_\ell(t)=a_\ell+\frac{\rho}{2}e^{it},
\qquad t\in[0,2\pi], \qquad \text{with}\qquad |\gamma_\ell'(t)|=\frac{\rho}{2}. 
\]
Then, for every $\omega\in \Omega$,
\[
\int_{\Gamma_\ell}|S_n(w,\omega)|^2\,|dw|
= \int_0^{2\pi} |S_n(\gamma_{\ell}(t),\omega)|^2\,|\gamma_{\ell}'(t)|dt=
\frac{\rho}{2}
\int_0^{2\pi}|S_n(\gamma_\ell(t),\omega)|^2\,dt.
\] 
It is straightforward to verify that the map \((\omega,t)\mapsto |S_n(\gamma_\ell(t),\omega)|^2\) is measurable in
\(\omega\) for every \(t\), and continuous in \(t\) for every \(\omega\).
Thus it is jointly measurable with respect to $\mathcal F\otimes\mathcal{B}([0,2\pi])$, by \citet[Lemma~4.51]{Aliprantis}. Hence, taking expectations in \eqref{contour:bound}, applying Tonelli's theorem and using $\Gamma_\ell\subset K_1$, gives
\begin{align*}
\mathbb E\left(\sup_{\zeta\in K}|S_n(\zeta)|^2\right)
&\le C
\sum_{\ell=1}^m
\mathbb{E}\left(\int_{0}^{2\pi}|S_n(\gamma_\ell(t))|^2 dt\right)=
C
\sum_{\ell=1}^m
\int_{0}^{2\pi}\mathbb{E}\bigl(|S_n(\gamma_\ell(t))|^2\bigr) dt\\[2 mm]
&\le C \sup_{\zeta\in K_1}\mathbb{E}\bigl(|S_n(\zeta)|^2\bigr)
\end{align*}
We now prove the two-dimensional assertion. The construction of \(K_1\), the
covering points \(a_1,\ldots,a_m\), and the contours $\Gamma_\ell=\{w\in\mathbb C:|w-a_\ell|=\rho/2\}$, for $\ell=1,\ldots,m$, are the same as above. In particular, \(\Gamma_\ell\subset K_1\), and the
closed disk \(\overline{B(a_\ell,\rho/2)}\) is contained in
\(K_1\subset\mathbb C\setminus\mathbb R_{\ge0}\).

Fix \(\omega\in\Omega\). Let \(\zeta,\xi\in K\). Choose
\(\ell,r\in\{1,\ldots,m\}\) such that
\[
        \zeta\in B(a_\ell,\rho/4),
        \qquad
        \xi\in B(a_r,\rho/4).
\]
Since \(S_n(\cdot,\cdot,\omega)\) is holomorphic in each variable separately,
Cauchy's integral formula applied first in the second variable and then in the
first variable gives
\[
\begin{split}
        S_n(\zeta,\xi)
        &=
        \frac{1}{2\pi i}
        \int_{\Gamma_r}
        \frac{S_n(\zeta,u)}{u-\xi}\,du
        =
        \frac{1}{2\pi i}
        \int_{\Gamma_r}
        \frac{1}{u-\xi}
        \left(
        \frac{1}{2\pi i}
        \int_{\Gamma_\ell}
        \frac{S_n(w,u)}{w-\zeta}\,dw
        \right)
        du
        \\
        &=
        \frac{1}{(2\pi i)^2}
        \int_{\Gamma_\ell}
        \int_{\Gamma_r}
        \frac{S_n(w,u)}
        {(w-\zeta)(u-\xi)}
        \,du\,dw .
\end{split}
\]
For \(w\in\Gamma_\ell\) and \(u\in\Gamma_r\), we have
\[
        |w-\zeta|
        \ge
        |w-a_\ell|-|\zeta-a_\ell|
        \ge
        \rho/2-\rho/4
        =
        \rho/4
\]
and similarly $|u-\xi|\ge \rho/4$.
Consequently,
\[
        |S_n(\zeta,\xi)|
        \le
        C_{K_1}
        \int_{\Gamma_\ell}
        \int_{\Gamma_r}
        |S_n(w,u)|\,|du|\,|dw|.
\]
Taking the supremum over \(\zeta,\xi\in K\) yields
\[
        \sup_{\zeta,\xi\in K}|S_n(\zeta,\xi)|
        \le
        C_{K_1}
        \max_{1\le \ell,r\le m}
        \int_{\Gamma_\ell}
        \int_{\Gamma_r}
        |S_n(w,u)|\,|du|\,|dw|.
\]
Therefore,
\[
\begin{split}
        \sup_{\zeta,\xi\in K}|S_n(\zeta,\xi)|^2
        &\le
        C_{K_1}
        \sum_{\ell,r=1}^m
        \left(
        \int_{\Gamma_\ell}
        \int_{\Gamma_r}
        |S_n(w,u)|\,|du|\,|dw|
        \right)^2 .
\end{split}
\]
By Cauchy--Schwarz on the product contour \(\Gamma_\ell\times\Gamma_r\),
\[
\begin{split}
\left(
        \int_{\Gamma_\ell}
        \int_{\Gamma_r}
        |S_n(w,u)|\,|du|\,|dw|
\right)^2 
&\le
\left(
        \int_{\Gamma_\ell}
        \int_{\Gamma_r}
        |S_n(w,u)|^2\,|du|\,|dw|
\right)
\left(
        \int_{\Gamma_\ell}
        \int_{\Gamma_r}
        |du|\,|dw|
\right)
\\
&\quad=
|\Gamma_\ell|\,|\Gamma_r|
        \int_{\Gamma_\ell}
        \int_{\Gamma_r}
        |S_n(w,u)|^2\,|du|\,|dw| .
\end{split}
\]
Since \(|\Gamma_\ell|=\pi\rho\) for every \(\ell\), it follows that
\[
        \sup_{\zeta,\xi\in K}|S_n(\zeta,\xi)|^2
        \le
        C_{K_1}
        \sum_{\ell,r=1}^m
        \int_{\Gamma_\ell}
        \int_{\Gamma_r}
        |S_n(w,u)|^2\,|du|\,|dw|.
\]
This deterministic bound holds for every \(n\) and every \(\omega\in\Omega\).

By Lemma~\ref{lem_sup_measurable}, applied with the
compact metric index space $K\times K$,
the random variable
$\sup_{\zeta,\xi\in K}|S_n(\zeta,\xi)|$
is $\mathcal F$-measurable.

It remains to take expectations in the deterministic contour bound. We
parameterize the contours by
\[
        \gamma_\ell(t)=a_\ell+\frac{\rho}{2}e^{it},
        \qquad t\in[0,2\pi],
        \qquad |\gamma_\ell'(t)|=\frac{\rho}{2}.
\]
Then
\[
\int_{\Gamma_\ell}
        \int_{\Gamma_r}
        |S_n(w,u,\omega)|^2\,|du|\,|dw|
=
        \frac{\rho^2}{4}
        \int_0^{2\pi}\int_0^{2\pi}
        |S_n(\gamma_\ell(t),\gamma_r(s),\omega)|^2
        \,ds\,dt .
\]
For fixed \(s,t\), the map $\omega\mapsto |S_n(\gamma_\ell(t),\gamma_r(s),\omega)|^2$
is \(\mathcal F\)-measurable. For fixed \(\omega\), the map
$(t,s)\mapsto
|S_n(\gamma_\ell(t),\gamma_r(s),\omega)|^2$ is continuous. Hence the map
\[
        (\omega,t,s)\mapsto
        |S_n(\gamma_\ell(t),\gamma_r(s),\omega)|^2
\]
is jointly measurable with respect to
\(\mathcal F\otimes\mathcal B([0,2\pi]^2)\), by
\citet[Lemma~4.51]{Aliprantis}. Therefore the contour integrals are
\(\mathcal F\)-measurable, and Tonelli's theorem applies.

Taking expectations in the deterministic contour bound gives
\[
\mathbb E\left[
        \sup_{\zeta,\xi\in K}|S_n(\zeta,\xi)|^2
\right]
\leq
C_{K_1}
\sum_{\ell,r=1}^m
\int_{\Gamma_\ell}
\int_{\Gamma_r}
\mathbb E\left(|S_n(w,u)|^2\right)
\,|du|\,|dw|.
\]
Since \(\Gamma_\ell,\Gamma_r\subseteq K_1\), we obtain
\[
\begin{split}
\mathbb E\left[
        \sup_{\zeta,\xi\in K}|S_n(\zeta,\xi)|^2
\right]
&\leq
C_{K_1}
\sum_{\ell,r=1}^m
\int_{\Gamma_\ell}
\int_{\Gamma_r}
\sup_{\zeta,\xi\in K_1}
\mathbb E\left(|S_n(\zeta,\xi)|^2\right)
\,|du|\,|dw|
\\
&\leq
C_{K_1}
\sup_{\zeta,\xi\in K_1}
\mathbb E\left(|S_n(\zeta,\xi)|^2\right).
\end{split}
\]
This proves the two-dimensional assertion.
\end{proof}

\begin{lemma}
\label{lem:conditional-moment-estimates}
Let \(z=(z_1,\ldots,z_p)^\top\) be a real random vector with independent
entries satisfying
\[
        \E(z_j)=0,
        \qquad
        \E(z_j^2)=1,
        \qquad
        \max_{1\leq j\leq p}\E(|z_j|^8)\leq C,
\]
and let \(\mathcal G\) be a sigma-field independent of \(z\). Let
\(M,N\) be \(\mathcal G\)-measurable complex \(p\times p\) matrices and
let \(v\) be a \(\mathcal G\)-measurable vector in \(\mathbb C^p\).
Then, almost surely,
\begin{align}
        \E\left[
        \|Mz\|_2^8
        \,\middle|\,
        \mathcal G
        \right]
        &\leq
        C\|M\|_F^8
        \leq
        Cp^4\|M\|_2^8,
        \label{eq:conditional-Mz-eighth}
        \\
        \E\left[
        |z^\top v|^8
        \,\middle|\,
        \mathcal G
        \right]
        &\leq
        C\|v\|_2^8,
        \label{eq:conditional-linear-eighth}
        \\
        \E\left[
        \|z^\top N\|_2^8
        \,\middle|\,
        \mathcal G
        \right]
        &\leq
        C\|N\|_F^8.
        \label{eq:conditional-row-eighth}
\end{align}
Consequently,
\begin{align}
\E\left[
        \|Mz\|_2^4|z^\top v|^4
        \,\middle|\,
        \mathcal G
\right]
&\leq
C\|M\|_F^4\|v\|_2^4
\leq
Cp^2\|M\|_2^4\|v\|_2^4,
\label{eq:conditional-mixed-vector}
\\
\E\left[
        \|Mz\|_2^4\|z^\top N\|_2^4
        \,\middle|\,
        \mathcal G
\right]
&\leq
C\|M\|_F^4\|N\|_F^4
\leq
Cp^2\|M\|_2^4\|N\|_F^4.
\label{eq:conditional-mixed-matrix}
\end{align}
The constant \(C\) depends only on the uniform eighth-moment bound.
\end{lemma}

\begin{proof}
Conditionally on \(\mathcal G\), the matrices \(M,N\) and the vector \(v\)
are deterministic and independent of \(z\). We first record a consequence
of Lemma~\ref{lem_b26_bai_silverstein}.

Let \(A\) be any \(\mathcal G\)-measurable Hermitian nonnegative definite
\(p\times p\) matrix. Conditionally on \(\mathcal G\),
Lemma~\ref{lem_b26_bai_silverstein} gives
\[
\begin{split}
        \E\left[
        \left|
        z^\top Az-\tr(A)
        \right|^4
        \,\middle|\,
        \mathcal G
        \right]
        &\leq
        C\left[
        \|A\|_F^4
        +
        \tr\{(AA^*)^2\}
        \right].
\end{split}
\]
Since \(A=A^*\), we have
\[
        \tr\{(AA^*)^2\}
        =
        \tr(A^4)
        \leq
        \{\tr(A^2)\}^2
        =
        \|A\|_F^4.
\]
It follows that
\begin{equation}
\label{eq:conditional-centered-qf}
        \E\left[
        \left|
        z^\top Az-\tr(A)
        \right|^4
        \,\middle|\,
        \mathcal G
        \right]
        \leq
        C\|A\|_F^4.
\end{equation}
Moreover, because \(A\) is nonnegative definite,
\[
        \|A\|_F^2
        =
        \tr(A^2)
        \leq
        \{\tr(A)\}^2.
\]
Using $|x+y|^4\leq 8|x|^4+8|y|^4$, $x,y\in \mathbb C$, together with \eqref{eq:conditional-centered-qf}, we obtain
\begin{align}\label{eq:conditional-positive-qf}
        \E\left[
        (z^\top Az)^4
        \,\middle|\,
        \mathcal G
        \right]
        &\leq
        8
        \E\left[
        \left|
        z^\top Az-\tr(A)
        \right|^4
        \,\middle|\,
        \mathcal G
        \right]
        +
        8\{\tr(A)\}^4
        \leq
        C\{\tr(A)\}^4.
\end{align}

We apply \eqref{eq:conditional-positive-qf} first with $A=M^*M$.
Since $\|Mz\|_2^2=z^\top M^*Mz$ and $\tr(M^*M)=\|M\|_F^2$,
we obtain
\[
\begin{split}
        \E\left[
        \|Mz\|_2^8
        \,\middle|\,
        \mathcal G
        \right]
        &=
        \E\left[
        (z^\top M^*Mz)^4
        \,\middle|\,
        \mathcal G
        \right]
        \leq
        C\|M\|_F^8.
\end{split}
\]
Furthermore, $\|M\|_F^2\leq p\|M\|_2^2$,
and therefore
\[
        \|M\|_F^8
        \leq
        p^4\|M\|_2^8.
\]
This proves \eqref{eq:conditional-Mz-eighth}.

Next, apply \eqref{eq:conditional-positive-qf} with $A=vv^*$.
Since $|z^\top v|^2=z^\top vv^*z$ and $\tr(vv^*)=\|v\|_2^2$,
we obtain
\[
\begin{split}
        \E\left[
        |z^\top v|^8
        \,\middle|\,
        \mathcal G
        \right]
        &=
        \E\left[
        (z^\top vv^*z)^4
        \,\middle|\,
        \mathcal G
        \right]
        \leq
        C\|v\|_2^8.
\end{split}
\]
This proves \eqref{eq:conditional-linear-eighth}.

Finally, apply \eqref{eq:conditional-positive-qf} with $A=NN^*$.
Since $\|z^\top N\|_2^2=z^\top NN^*z$
and $\tr(NN^*)=\|N\|_F^2$,
we obtain
\[
\begin{split}
        \E\left[
        \|z^\top N\|_2^8
        \,\middle|\,
        \mathcal G
        \right]
        &=
        \E\left[
        (z^\top NN^*z)^4
        \,\middle|\,
        \mathcal G
        \right]
        \leq
        C\|N\|_F^8.
\end{split}
\]
This proves \eqref{eq:conditional-row-eighth}.

It remains to prove the two mixed bounds. By the conditional
Cauchy--Schwarz inequality,
\begin{align*}
\E\left[
        \|Mz\|_2^4|z^\top v|^4
        \,\middle|\,
        \mathcal G
        \right]
&\leq
        \left\{
        \E\left[
        \|Mz\|_2^8
        \,\middle|\,
        \mathcal G
        \right]
        \right\}^{1/2}
        \left\{
        \E\left[
        |z^\top v|^8
        \,\middle|\,
        \mathcal G
        \right]
        \right\}^{1/2}
        \nonumber\\
&\leq
        C\|M\|_F^4\|v\|_2^4
        \leq
        Cp^2\|M\|_2^4\|v\|_2^4.        
\end{align*}
This proves \eqref{eq:conditional-mixed-vector}. Similarly,
\begin{align*}
\E\left[
        \|Mz\|_2^4\|z^\top N\|_2^4
        \,\middle|\,
        \mathcal G
        \right]
        &\leq
        \left\{
        \E\left[
        \|Mz\|_2^8
        \,\middle|\,
        \mathcal G
        \right]
        \right\}^{1/2}
        \left\{
        \E\left[
        \|z^\top N\|_2^8
        \,\middle|\,
        \mathcal G
        \right]
        \right\}^{1/2}\\
        &\leq
        C\|M\|_F^4\|N\|_F^4
        \leq
        Cp^2\|M\|_2^4\|N\|_F^4.
\end{align*}
This proves \eqref{eq:conditional-mixed-matrix}.
\end{proof}

\begin{lemma}
\label{lem:resolvent-denominator-bounds}
Let \(A\) be a real symmetric, nonnegative definite \(p\times p\)
matrix, let \(x\in\mathbb R^p\), let \(c\geq0\), and let $\zeta\in\mathbb C\setminus\mathbb R_{\geq0}$.
Set $d_\zeta\coloneqq\operatorname{dist}\bigl(\zeta,\mathbb R_{\geq0}\bigr)=\inf_{t\geq0}|\zeta-t|$.
Then
\begin{equation}
\label{eq:general-denominator-bound}
        \left|
        1+c\,x^\top(A-\zeta I_p)^{-1}x
        \right|^{-1}
        \leq
        \frac{|\zeta|}{d_\zeta}.
\end{equation}
Moreover, if \(B\) is a real symmetric, nonnegative definite \(p\times p\) matrix, then 
\begin{equation}\label{eq:general-trace-denominator-bound} 
\left| 1+c\,\tr\left\{ B(A-\zeta I_p)^{-1} \right\} \right|^{-1} \leq \frac{|\zeta|}{d_\zeta}. \end{equation} 
Furthermore,
\begin{equation}
\label{eq:general-resolvent-bound}
        \left\|
        (A-\zeta I_p)^{-1}
        \right\|_2
        \leq
        \frac1{d_\zeta}.
\end{equation}

Let, in addition, \(\nu\) be a finite nonnegative measure supported on
\(\mathbb R_{\geq0}\), and define
\[
        v_\nu(\zeta)
        \coloneqq
        \int_0^\infty
        \frac1{t-\zeta}\,d\nu(t).
\]
If \(\Sigma\) is a real symmetric, nonnegative definite
\(p\times p\) matrix, then
\begin{equation}
\label{eq:deterministic-resolvent-bound}
        \left\|
        \left\{
        -\zeta v_\nu(\zeta)\Sigma-\zeta I_p
        \right\}^{-1}
        \right\|_2
        \leq
        \frac1{d_\zeta}.
\end{equation}
\end{lemma}

\begin{proof}
Let $A=\sum_{\ell=1}^p\lambda_\ell u_\ell u_\ell^\top$ be an orthonormal spectral decomposition of \(A\), where
\(\lambda_\ell\geq0\). We first prove \eqref{eq:general-denominator-bound}. Note that,
\begin{align*}
        c\,x^\top(A-\zeta I_p)^{-1}x
        &=
        \sum_{\ell=1}^p
        \frac{c\,\langle x,u_\ell\rangle^2}
        {\lambda_\ell-\zeta}.
\end{align*}
Suppose that
\(\operatorname{Re}(\zeta)<0\). Since \(\lambda_\ell\geq0\),
\[
        \Re
        \left(
        \frac1{\lambda_\ell-\zeta}
        \right)
        =
        \frac{\lambda_\ell-\Re(\zeta)}
        {|\lambda_\ell-\zeta|^2}
        \geq0.
\]
Consequently,
\[
        \Re
        \left\{
        1+c\,x^\top(A-\zeta I_p)^{-1}x
        \right\}
        =
        1+
        \sum_{\ell=1}^p
        c\,\langle x,u_\ell\rangle^2
        \Re
        \left(
        \frac1{\lambda_\ell-\zeta}
        \right)
        \geq1.
\]
It follows that
\begin{align}\label{eq:shm-res:u:bound}
        \frac{1}{\left|
        1+c\,x^\top(A-\zeta I_p)^{-1}x
        \right|}\leq \frac{1}{ \Re
        \left\{
        1+c\,x^\top(A-\zeta I_p)^{-1}x
        \right\}}
        \leq1.
\end{align}

Since \(\Re(\zeta)<0\), we have
\begin{align*}
        |\zeta-t|^2
        &=
        |\Re(\zeta)-t+i\Im(\zeta)|^2
        =
        \{\Re(\zeta)-t\}^2+\Im(\zeta)^2
        \\
        &=
        \Re(\zeta)^2+\Im(\zeta)^2+t^2-2t\Re(\zeta)
        \\
        &\geq
        \Re(\zeta)^2+\Im(\zeta)^2
        =
        |\zeta|^2.
\end{align*}
Therefore,
\[
        |\zeta-t|\geq|\zeta|
        \qquad\text{for every }t\geq0.
\]
On the other hand, equality holds for \(t=0\). Hence, if $\Re(\zeta)<0$,
\[
        d_\zeta
        =
        \inf_{t\geq0}|\zeta-t|
        =
        |\zeta|.
\]
Using \eqref{eq:shm-res:u:bound} and the preceding display, we conclude that

\[
        \left|
        1+c\,x^\top(A-\zeta I_p)^{-1}x
        \right|^{-1}
        \leq 1=\frac{|\zeta|}{d_\zeta},
\]

Suppose now that \(\Re(\zeta)\geq0\). Since \(\zeta\notin\mathbb R_{\geq0}\) by definition, we necessarily have
\(\Im(\zeta)\neq0\). For every \(t\geq0\), we have
\[
|\zeta-t|^2=|\Re(\zeta)-t+i\Im(\zeta)|^2=(\Re(\zeta)-t)^2+\Im(\zeta)^2\geq \Im(\zeta)^2.
\]
Therefore,
\[
        |\zeta-t|\geq
        |\Im(\zeta)|
        \qquad\text{for every }t\geq0.
\]
On the other hand, since \(\Re(\zeta)\geq0\), we may choose \(t=\Re(\zeta)\), and then
$|\zeta-\Re(\zeta)|=|i\Im(\zeta)|=|\Im(\zeta)|$.
Hence,
\begin{align}\label{eq:dist:rekl0}
        d_\zeta
        =
        \inf_{t\geq0}|\zeta-t|
        =
        |\Im(\zeta)|.
\end{align}

Next, note that for every \(\lambda_\ell\geq0\),
\begin{align*}
        \frac{\zeta}{\lambda_\ell-\zeta}
        &=
        \frac{[\Re(\zeta)+i\Im(\zeta)]\{(\lambda_\ell-\Re(\zeta))+i\Im(\zeta)\}}
        {|\lambda_\ell-\zeta|^2}
        \\
        &=
        \frac{
        \Re(\zeta)\{\lambda_\ell-\Re(\zeta)\}-\Im(\zeta)^2
        +
        i\{\Re(\zeta)\Im(\zeta)+\Im(\zeta)[\lambda_\ell-\Re(\zeta)]\}
        }
        {|\lambda_\ell-\zeta|^2}
        \\
        &=
        \frac{
        \Re(\zeta)\{\lambda_\ell-\Re(\zeta)\}-\Im(\zeta)^2
        +
        i\lambda_\ell \Im(\zeta)
        }
        {|\lambda_\ell-\zeta|^2}.
\end{align*}
Hence, it follows that
\begin{align*}
\Im\left(\frac{\zeta}{\lambda_\ell-\zeta}\right)
=\frac{\lambda_\ell\Im(\zeta)}{|\lambda_\ell-\zeta|^2}.
\end{align*}
Thus,
\begin{align*}
\Im
\left[
\zeta
\left\{
1+c\,x^\top(A-\zeta I_p)^{-1}x
\right\}
\right]
&=
\Im(\zeta)
+
\sum_{\ell=1}^p
c\,\langle x,u_\ell\rangle^2
\Im
\left(
\frac{\zeta}{\lambda_\ell-\zeta}
\right)\\
&=
\Im(\zeta)
\left\{
1+
\sum_{\ell=1}^p
\frac{
c\,\lambda_\ell\langle x,u_\ell\rangle^2
}{
|\lambda_\ell-\zeta|^2
}
\right\}.
\end{align*}
Since every term in the sum is nonnegative and using \eqref{eq:dist:rekl0}, we conclude that
\begin{align*}
|\zeta|
\left|
1+c\,x^\top(A-\zeta I_p)^{-1}x
\right|
&=
\left|
\zeta
\left\{
1+c\,x^\top(A-\zeta I_p)^{-1}x
\right\}
\right|
\\
&\geq
\left|
\Im
\left[
\zeta
\left\{
1+c\,x^\top(A-\zeta I_p)^{-1}x
\right\}
\right]
\right|
\\
&\geq
|\Im(\zeta)|
=
d_\zeta.
\end{align*}
This proves \eqref{eq:general-denominator-bound}.

We next prove \eqref{eq:general-trace-denominator-bound}. Note that
\[
        \tr\left\{
        B(A-\zeta I_p)^{-1}
        \right\}
        =
        \sum_{\ell=1}^p
        \frac{
        u_\ell^\top Bu_\ell
        }{
        \lambda_\ell-\zeta
        }.
\]
Since \(B\) is nonnegative definite, $u_\ell^\top Bu_\ell\geq0$, for $\ell\in \{1,\dots,p\}$.

Suppose first that \(\Re(\zeta)<0\). Since
\[
        \Re\left(
        \frac1{\lambda_\ell-\zeta}
        \right)=\frac{\lambda_\ell-\Re(\zeta)}{|\lambda_\ell-\zeta|^2}
        \geq0,
\]
we obtain
\[
\Re\left[
        1+c\,\tr\left\{
        B(A-\zeta I_p)^{-1}
        \right\}
\right]
=
1+
c\sum_{\ell=1}^p
u_\ell^\top Bu_\ell\,
\Re\left(
\frac1{\lambda_\ell-\zeta}
\right)
\geq1.
\]
Consequently, if \(\Re(\zeta)<0\), we have
\[
        \left|
        1+c\,\tr\left\{
        B(A-\zeta I_p)^{-1}
        \right\}
        \right|^{-1}
        \leq1.
\]
As shown above, \(d_\zeta=|\zeta|\) when \(\Re(\zeta)<0\). Hence,
\[
        \left|
        1+c\,\tr\left\{
        B(A-\zeta I_p)^{-1}
        \right\}
        \right|^{-1}
        \leq
        \frac{|\zeta|}{d_\zeta}.
\]

Suppose now that \(\Re(\zeta)\geq0\). Then
\(\Im(\zeta)\neq0\), and $d_\zeta=|\Im(\zeta)|$.
Moreover,
\begin{align*}
\Im\left[
\zeta
\left\{
1+c\,\tr\left(
B(A-\zeta I_p)^{-1}
\right)
\right\}
\right]
&=
\Im(\zeta)
+
c\sum_{\ell=1}^p
u_\ell^\top Bu_\ell\,
\Im\left(
\frac{\zeta}{\lambda_\ell-\zeta}
\right)
\\
&=
\Im(\zeta)
\left\{
1+
c\sum_{\ell=1}^p
\frac{
\lambda_\ell\, u_\ell^\top Bu_\ell
}{
|\lambda_\ell-\zeta|^2
}
\right\}.
\end{align*}
Since $\lambda_\ell\, u_\ell^\top Bu_\ell\geq0$, it follows that
\begin{align*}
|\zeta|
\left|
1+c\,\tr\left\{
B(A-\zeta I_p)^{-1}
\right\}
\right|
&\geq
\left|
\Im\left[
\zeta
\left\{
1+c\,\tr\left(
B(A-\zeta I_p)^{-1}
\right)
\right\}
\right]
\right|
\\
&\geq
|\Im(\zeta)|
=
d_\zeta.
\end{align*}
This proves \eqref{eq:general-trace-denominator-bound}.

To prove \eqref{eq:general-resolvent-bound}, the spectral decomposition
of \(A\) gives
\[
        \left\|
        (A-\zeta I_p)^{-1}
        \right\|_2
        =
        \max_{1\leq\ell\leq p}
        \frac1{|\lambda_\ell-\zeta|}
        \leq
        \frac1{
        \operatorname{dist}
        (\zeta,\mathbb R_{\geq0})
        }
        =
        \frac1{d_\zeta}.
\]

It remains to prove \eqref{eq:deterministic-resolvent-bound}. Let $\Sigma=\sum_{\ell=1}^pt_\ell u_\ell u_\ell^\top$, $t_\ell\geq0$,
be an orthonormal spectral decomposition of \(\Sigma\). Then
\[
        -\zeta v_\nu(\zeta)\Sigma-\zeta I_p
        =
        \sum_{\ell=1}^p
        \left\{
        -\zeta\bigl(1+t_\ell v_\nu(\zeta)\bigr)
        \right\}
        u_\ell u_\ell^\top.
\]
We show that
\begin{equation}
\label{eq:scalar-deterministic-resolvent-bound}
        \left|
        \zeta\bigl(1+t_\ell v_\nu(\zeta)\bigr)
        \right|
        \geq
        d_\zeta,\qquad\text{for every }\ell\in \{1,\dots,p\}
\end{equation}
Suppose first that \(\Re(\zeta)<0\). Then
\[      \Re\bigl\{v_\nu(\zeta)\bigr\}
        =
        \int_0^\infty
        \frac{
        t-\Re(\zeta)
        }{
        |t-\zeta|^2
        }
        \,d\nu(t)
        \geq0,
\]
and $\Re\bigl(1+t_\ell v_\nu(\zeta)\bigr)\geq1$.
Therefore,
\begin{align}\label{eq:bound:real}
        \left|
        \zeta\bigl\{1+t_\ell v_\nu(\zeta)\bigr\}
        \right|
        =
        |\zeta|
        \left|
        1+t_\ell v_\nu(\zeta)
        \right|
        \geq  |\zeta|\,\Re\bigl(1+t_\ell v_\nu(\zeta)\bigr)\geq
        |\zeta|
        =
        d_\zeta.
\end{align}
Suppose now that \(\Re(\zeta)\geq0\). Then, by \eqref{eq:dist:rekl0},
$d_\zeta=|\Im(\zeta)|$.
Moreover,
\begin{align*}
\Im
\left[
\zeta\bigl(1+t_\ell v_\nu(\zeta)\bigr)
\right]
&=
\Im(\zeta)
+
t_\ell
\int_0^\infty
\Im
\left(
\frac{\zeta}{t-\zeta}
\right)
\,d\nu(t)
\\
&=
\Im(\zeta)
\left\{
1+
t_\ell
\int_0^\infty
\frac{t}{|t-\zeta|^2}\,d\nu(t)
\right\}.
\end{align*}
It follows that
\begin{align}\label{eq:bound:real1}
        \left|
        \zeta\bigl(1+t_\ell v_\nu(\zeta)\bigr)
        \right|
        \geq
        \left|
        \Im
        \left[
        \zeta\bigl(1+t_\ell v_\nu(\zeta)\bigr)
        \right]
        \right|
        \geq
        |\Im(\zeta)|
        =
        d_\zeta.
\end{align}
Combining \eqref{eq:bound:real} and \eqref{eq:bound:real1}, this proves \eqref{eq:scalar-deterministic-resolvent-bound}. Consequently,
\[
\left\|
\left\{
-\zeta v_\nu(\zeta)\Sigma-\zeta I_p
\right\}^{-1}
\right\|_2
=
\max_{1\leq\ell\leq p}
\frac1{
\left|
\zeta\bigl(1+t_\ell v_\nu(\zeta)\bigr)
\right|
}\leq
\frac1{d_\zeta},
\]
which proves \eqref{eq:deterministic-resolvent-bound}.
\end{proof}

\begin{lemma}
\label{lem:finite-resolvent-difference-expansions}
Let $x_1,\ldots,x_n\in\mathbb R^p$ be arbitrary
vectors, and let $R_{n,-D}(\zeta)$ be defined as in \eqref{eq:D_del_res_cov}.
Fix $\zeta\in\mathbb C\setminus\mathbb R_{\geq0}$ and write
$R_{-D}=R_{n,-D}(\zeta)$, suppressing $n$ and $\zeta$ in the notation. Let
\(I\subseteq\{1,\ldots,n\}\), let \(H\subseteq I\) satisfy $|H|\in\{1,2,3\}$,
and set $D\coloneqq I\setminus H$.
For every deletion set \(G\subseteq\{1,\ldots,n\}\) and every
\(m\notin G\), define
\[
        \beta_{G,m}
        :=
        \left\{
        1+\frac1n
        x_m^\top R_{-(G\cup\{m\})}x_m
        \right\}^{-1}.
\]
Moreover, define
\[
        \mathcal R_H
        :=
        \sum_{L\subseteq H}
        (-1)^{|H|-|L|}
        R_{-(I\setminus L)}.
\]

If \(H=\{q\}\), then
\begin{equation}
\label{eq:proof-RH-one}
        \mathcal R_H
        =
        -\frac1n
        \beta_{D,q}
        R_{-I}x_qx_q^\top R_{-I}.
\end{equation}

If \(H=\{q,r\}\), then
\begin{equation}
\label{eq:proof-RH-two}
\begin{split}
        \mathcal R_H
        &=
        \frac1{n^2}
        \beta_{D,q}\beta_{D\cup\{q\},r}
        R_{-(D\cup\{q\})}
        x_qx_q^\top
        R_{-I}x_rx_r^\top R_{-I}
        \\
        &\quad+
        \frac1{n^2}
        \beta_{D,r}\beta_{D\cup\{r\},q}
        R_{-(D\cup\{r\})}
        x_rx_r^\top
        R_{-I}x_qx_q^\top R_{-I}.
\end{split}
\end{equation}

Finally, if \(H=\{q,r,s\}\), then
\begin{equation}
\label{eq:proof-RH-three}
\begin{split}
        \mathcal R_H
        &=
        -\frac1{n^3}
        \beta_{D,q}
        \beta_{D\cup\{q\},r}
        \beta_{D\cup\{q,r\},s}
        R_{-(D\cup\{q\})}x_qx_q^\top
        R_{-(D\cup\{q,r\})}x_rx_r^\top
        R_{-I}x_sx_s^\top R_{-I}
        \\
        &\quad-
        \frac1{n^3}
        \beta_{D,q}
        \beta_{D\cup\{q\},s}
        \beta_{D\cup\{q,s\},r}
        R_{-(D\cup\{q\})}x_qx_q^\top
        R_{-(D\cup\{q,s\})}x_sx_s^\top
        R_{-I}x_rx_r^\top R_{-I}
        \\
        &\quad-
        \frac1{n^3}
        \beta_{D,r}
        \beta_{D\cup\{r\},q}
        \beta_{D\cup\{q,r\},s}
        R_{-(D\cup\{r\})}x_rx_r^\top
        R_{-(D\cup\{q,r\})}x_qx_q^\top
        R_{-I}x_sx_s^\top R_{-I}
        \\
        &\quad-
        \frac1{n^3}
        \beta_{D,r}
        \beta_{D\cup\{r\},s}
        \beta_{D\cup\{r,s\},q}
        R_{-(D\cup\{r\})}x_rx_r^\top
        R_{-(D\cup\{r,s\})}x_sx_s^\top
        R_{-I}x_qx_q^\top R_{-I}
        \\
        &\quad-
        \frac1{n^3}
        \beta_{D,s}
        \beta_{D\cup\{s\},q}
        \beta_{D\cup\{q,s\},r}
        R_{-(D\cup\{s\})}x_sx_s^\top
        R_{-(D\cup\{q,s\})}x_qx_q^\top
        R_{-I}x_rx_r^\top R_{-I}
        \\
        &\quad-
        \frac1{n^3}
        \beta_{D,s}
        \beta_{D\cup\{s\},r}
        \beta_{D\cup\{r,s\},q}
        R_{-(D\cup\{s\})}x_sx_s^\top
        R_{-(D\cup\{r,s\})}x_rx_r^\top
        R_{-I}x_qx_q^\top R_{-I}.
\end{split}
\end{equation}
\end{lemma}

\begin{proof}
For every deletion set \(G\subseteq\{1,\ldots,n\}\) and every
\(m\notin G\), the Sherman--Morrison formula gives
\begin{equation}
\label{eq:finite-difference-SM}
        R_{-G}
        =
        R_{-(G\cup\{m\})}
        -
        \frac1n
        \beta_{G,m}
        R_{-(G\cup\{m\})}
        x_mx_m^\top
        R_{-(G\cup\{m\})}.
\end{equation}

\medskip
\noindent
\textbf{Case \(\lvert H\rvert=1\).}

Let \(H=\{q\}\). Since \(D=I\setminus\{q\}\), we have
\[
        D\cup\{q\}=I
\]
and
\[
        \mathcal R_H
        =
        R_{-D}-R_{-I}.
\]
Applying \eqref{eq:finite-difference-SM} with \(G=D\) and \(m=q\)
therefore gives
\[
        \mathcal R_H
        =
        -\frac1n
        \beta_{D,q}
        R_{-I}x_qx_q^\top R_{-I},
\]
which proves \eqref{eq:proof-RH-one}.

\medskip
\noindent
\textbf{Case \(\lvert H\rvert=2\).}

Let \(H=\{q,r\}\). Then
\[
        I=D\cup\{q,r\},
\]
and the definition of \(\mathcal R_H\) gives
\begin{equation}
\label{eq:proof-RH-two-definition}
\begin{split}
        \mathcal R_H
        &=
        R_{-D}
        -
        R_{-(D\cup\{q\})}
        -
        R_{-(D\cup\{r\})}
        +
        R_{-I}.
\end{split}
\end{equation}
We first rewrite this finite difference as
\begin{equation}
\label{eq:proof-RH-two-first}
\begin{split}
        \mathcal R_H
        &=
        -\frac1n
        R_{-D}x_qx_q^\top
        \left\{
        R_{-(D\cup\{q\})}-R_{-I}
        \right\}
        \\
        &\quad-
        \frac1n
        R_{-D}x_rx_r^\top
        \left\{
        R_{-(D\cup\{r\})}-R_{-I}
        \right\}.
\end{split}
\end{equation}
Indeed, expanding the right-hand side of
\eqref{eq:proof-RH-two-first} gives
\[
\begin{split}
&-\frac1n
R_{-D}x_qx_q^\top R_{-(D\cup\{q\})}
-
\frac1n
R_{-D}x_rx_r^\top R_{-(D\cup\{r\})}
\\
&\quad+
\frac1n
R_{-D}
\left(
x_qx_q^\top+x_rx_r^\top
\right)
R_{-I}.
\end{split}
\]
By the resolvent identity,
\[
        -\frac1n
        R_{-D}x_qx_q^\top R_{-(D\cup\{q\})}
        =
        R_{-D}-R_{-(D\cup\{q\})},
\]
\[
        -\frac1n
        R_{-D}x_rx_r^\top R_{-(D\cup\{r\})}
        =
        R_{-D}-R_{-(D\cup\{r\})},
\]
and, since \(I=D\cup\{q,r\}\),
\[
        \frac1n
        R_{-D}
        \left(
        x_qx_q^\top+x_rx_r^\top
        \right)
        R_{-I}
        =
        R_{-I}-R_{-D}.
\]
Adding these three identities gives
\eqref{eq:proof-RH-two-definition}, and hence proves
\eqref{eq:proof-RH-two-first}.

Next, applying \eqref{eq:finite-difference-SM} once more gives
\begin{align}
\label{eq:SM2}
\begin{split}
        R_{-(D\cup\{q\})}-R_{-I}
        &=
        -\frac1n
        \beta_{D\cup\{q\},r}
        R_{-I}x_rx_r^\top R_{-I},
        \\
        R_{-(D\cup\{r\})}-R_{-I}
        &=
        -\frac1n
        \beta_{D\cup\{r\},q}
        R_{-I}x_qx_q^\top R_{-I}.
\end{split}
\end{align}
Moreover, multiplying \eqref{eq:finite-difference-SM}, with \(G=D\)
and \(m=q\), on the right by \(x_qx_q^\top\), we obtain
\begin{equation}
\label{eq:SM3-multiplied}
        R_{-D}x_qx_q^\top
        =
        \beta_{D,q}
        R_{-(D\cup\{q\})}x_qx_q^\top.
\end{equation}
Indeed,
\[
\begin{split}
R_{-D}x_qx_q^\top
&=
\left[
R_{-(D\cup\{q\})}
-
\frac1n
\beta_{D,q}
R_{-(D\cup\{q\})}
x_qx_q^\top
R_{-(D\cup\{q\})}
\right]
x_qx_q^\top
\\
&=
\left[
1-
\frac{
n^{-1}x_q^\top R_{-(D\cup\{q\})}x_q
}{
1+n^{-1}x_q^\top R_{-(D\cup\{q\})}x_q
}
\right]
R_{-(D\cup\{q\})}x_qx_q^\top
\\
&=
\beta_{D,q}
R_{-(D\cup\{q\})}x_qx_q^\top.
\end{split}
\]
The same identity holds with \(q\) replaced by \(r\). Substituting
\eqref{eq:SM2} and \eqref{eq:SM3-multiplied} into
\eqref{eq:proof-RH-two-first} proves
\eqref{eq:proof-RH-two}.

\medskip
\noindent
\textbf{Case \(\lvert H\rvert=3\).}

Let \(H=\{q,r,s\}\). Then
\[
        I=D\cup\{q,r,s\},
\]
and
\begin{equation}
\label{eq:proof-RH-three-definition}
\begin{split}
        \mathcal R_H
        &=
        R_{-D}
        -
        R_{-(D\cup\{q\})}
        -
        R_{-(D\cup\{r\})}
        -
        R_{-(D\cup\{s\})}
        \\
        &\quad+
        R_{-(D\cup\{q,r\})}
        +
        R_{-(D\cup\{q,s\})}
        +
        R_{-(D\cup\{r,s\})}
        -
        R_{-I}.
\end{split}
\end{equation}
We first rewrite this finite difference as
\begin{equation}
\label{eq:proof-RH-three-first}
\begin{split}
        \mathcal R_H
        &=
        -\frac1n
        R_{-D}x_qx_q^\top
        \left\{
        R_{-(D\cup\{q\})}
        -
        R_{-(D\cup\{q,r\})}
        -
        R_{-(D\cup\{q,s\})}
        +
        R_{-I}
        \right\}
        \\
        &\quad-
        \frac1n
        R_{-D}x_rx_r^\top
        \left\{
        R_{-(D\cup\{r\})}
        -
        R_{-(D\cup\{q,r\})}
        -
        R_{-(D\cup\{r,s\})}
        +
        R_{-I}
        \right\}
        \\
        &\quad-
        \frac1n
        R_{-D}x_sx_s^\top
        \left\{
        R_{-(D\cup\{s\})}
        -
        R_{-(D\cup\{q,s\})}
        -
        R_{-(D\cup\{r,s\})}
        +
        R_{-I}
        \right\}.
\end{split}
\end{equation}
To verify \eqref{eq:proof-RH-three-first}, expand its right-hand side:
\[
\begin{split}
-\frac1nR_{-D}
\Big[
&x_qx_q^\top R_{-(D\cup\{q\})}
+
x_rx_r^\top R_{-(D\cup\{r\})}
+
x_sx_s^\top R_{-(D\cup\{s\})}
\\
&-
(x_qx_q^\top+x_rx_r^\top)
R_{-(D\cup\{q,r\})}
\\
&-
(x_qx_q^\top+x_sx_s^\top)
R_{-(D\cup\{q,s\})}
\\
&-
(x_rx_r^\top+x_sx_s^\top)
R_{-(D\cup\{r,s\})}
\\
&+
(x_qx_q^\top+x_rx_r^\top+x_sx_s^\top)R_{-I}
\Big].
\end{split}
\]
By the resolvent identity,
\[
        -\frac1n
        R_{-D}x_qx_q^\top R_{-(D\cup\{q\})}
        =
        R_{-D}-R_{-(D\cup\{q\})},
\]
and the analogous identities hold for \(r\) and \(s\). Moreover,
\[
        \frac1n
        R_{-D}
        (x_qx_q^\top+x_rx_r^\top)
        R_{-(D\cup\{q,r\})}
        =
        R_{-(D\cup\{q,r\})}-R_{-D},
\]
and the analogous identities hold for the pairs
\(\{q,s\}\) and \(\{r,s\}\). Finally,
\[
        -\frac1n
        R_{-D}
        (x_qx_q^\top+x_rx_r^\top+x_sx_s^\top)R_{-I}
        =
        R_{-D}-R_{-I}.
\]
Adding these seven identities gives
\eqref{eq:proof-RH-three-definition}, and therefore proves
\eqref{eq:proof-RH-three-first}.

The first expression in braces on the right-hand side of
\eqref{eq:proof-RH-three-first} is the second-order resolvent
difference in \eqref{eq:proof-RH-two-definition}, with
\(D\cup\{q\}\) in place of \(D\) and \(r,s\) as the two additional
indices. Hence, by \eqref{eq:proof-RH-two},
\[
\begin{split}
&
R_{-(D\cup\{q\})}
-
R_{-(D\cup\{q,r\})}
-
R_{-(D\cup\{q,s\})}
+
R_{-I}
\\
&\quad=
\frac1{n^2}
\beta_{D\cup\{q\},r}
\beta_{D\cup\{q,r\},s}
R_{-(D\cup\{q,r\})}
x_rx_r^\top
R_{-I}x_sx_s^\top R_{-I}
\\
&\qquad+
\frac1{n^2}
\beta_{D\cup\{q\},s}
\beta_{D\cup\{q,s\},r}
R_{-(D\cup\{q,s\})}
x_sx_s^\top
R_{-I}x_rx_r^\top R_{-I}.
\end{split}
\]
The same argument applies to the other two expressions in braces,
with the roles of \(q\), \(r\), and \(s\) permuted accordingly.

Moreover, by \eqref{eq:SM3-multiplied},
\[
        R_{-D}x_mx_m^\top
        =
        \beta_{D,m}
        R_{-(D\cup\{m\})}x_mx_m^\top,
        \qquad
        m\in\{q,r,s\}.
\]
Substituting the preceding identities into
\eqref{eq:proof-RH-three-first} gives
\eqref{eq:proof-RH-three}. This completes the proof.
\end{proof}

\newpage

\stopcontents[appendix]

\end{document}